\documentclass[reqno,  12pt]{amsart}
\usepackage{extarrows}
\usepackage{ragged2e}
\usepackage{amssymb}
\usepackage{amsfonts}
\usepackage{amsmath}
\usepackage{mathrsfs}
\usepackage{diagbox}
\usepackage{color,  soul}
\usepackage{bbm}
\usepackage{amsthm, graphicx,  empheq}
\usepackage{framed}
\usepackage[backref,  colorlinks,  linkcolor=red,  anchorcolor=green,  citecolor=blue]{hyperref}

\usepackage{cancel}
\usepackage{arydshln}

\usepackage{tikz}
\usetikzlibrary{arrows,backgrounds,snakes,shapes}

\usepackage{booktabs}
\usepackage{makecell}
\usepackage[linesnumbered,ruled]{algorithm2e}

\usepackage{lineno}  

\newtheorem{definition}{Definition}[section]
\newtheorem{remark}{Remark}[section]
\newtheorem{proposition}{Proposition}[section] 
\newtheorem{example}{Example}[section] 
\numberwithin{equation}{section}

\makeatletter
\@namedef{subjclassname@2020}{\textup{2020} Mathematics Subject Classification}
\makeatother

\begin{document}


\title[FLUX VECTOR SPLITTING RKDG METHOD]
{Runge-Kutta Discontinuous Galerkin Method Based on Flux Vector Splitting with Constrained Optimization-based TVB(D)-minmod Limiter for Solving Hyperbolic Conservation Laws}
\author{Zhengrong Xie}

\address[Z. Xie]{School of Mathematical Sciences,  
East China Normal University, Shanghai 200241, China}\email{\tt xzr\_nature@163.com; 52265500018@stu.ecnu.edu.cn}

\keywords{Flux Vector Splitting; RKDG; Hyperbolic Conservation Law; TVB(D)-minmod Limiter; Lagrange Multiplier Method; Local Characteristic Decomposition.}

\subjclass[2020]{65M60,  65M99, 35L65}

\date{\today}

\begin{abstract}
  The flux vector splitting (FVS) method, employed for solving hyperbolic conservation systems, 
  is frequently observed in finite difference and finite volume schemes, 
  where it serves to reconstruct the numerical fluxes at the interfaces of adjacent cells. 
  This method now has been incorporated into the discontinuous Galerkin (DG) framework for reconstructing the numerical fluxes required for 
  the spatial semi-discrete formulation, setting it apart from the conventional DG approaches that typically utilize the Lax-Friedrichs flux scheme 
  or classical Riemann solvers such as HLL and HLLC. 
  The FVS method inherently does not introduce any error. 
  Consequently, the control equations of hyperbolic conservation systems are initially reformulated into a flux-split form. 
  Subsequently, a variational approach is applied to this flux-split form, from which a DG spatial semi-discrete scheme based on FVS is derived. 
  When solving hyperbolic conservation laws, DG schemes are prone to numerical pseudo-oscillations, 
  necessitating the incorporation of limiters for correction. Initially, the smoothness measurement function $IS$ from the WENO limiter is integrated into the TVB(D)-minmod limiter, 
  constructing an optimization problem based on the smoothness factor constraint, thereby realizing a TVB(D)-minmod limiter applicable to arbitrary high-order polynomial approximation. 
  Subsequently, drawing on the ``reconstructed polynomial and the original high-order scheme's $L^2$-error constraint'' 
  from the literature \cite{ref1}, combined with our smoothness factor constraint, a bi-objective optimization problem is formulated to 
  enable the TVB(D)-minmod limiter to balance oscillation suppression and high precision. 
  The aforementioned constrained optimization-type TVB(D)-minmod limiter has been extended to two-dimensional scenarios.  
  When solving hyperbolic conservation systems, limiters are typically required to be used in conjunction with local characteristic decomposition. To transform polynomials from the physical space to the characteristic space, an interpolation-based characteristic transformation scheme has been proposed, and its equivalence with the original moment characteristic transformation has been demonstrated in one-dimensional scenarios. 
  For variable-coefficient and nonlinear equations, the transformation matrix utilized in characteristic transformation must be determined through local freezing. To achieve this, we employ the arithmetic mean or Roe average of the integral means on the common interface of adjacent cells for local freezing, rather than the conventional approach of using the arithmetic mean or Roe average of the cell integral means of each adjacent cell for this purpose. 
  Finally, the concept of ``flux vector splitting based on Jacobian eigenvalue decomposition'' has been applied to the conservative linear scalar transport equations and the nonlinear Burgers' equation. This approach has led to the rederivation of the classical Lax-Friedrichs flux scheme and the provision of a Steger-Warming flux scheme for scalar equations. 
\end{abstract}

\allowbreak
\allowdisplaybreaks

\maketitle

\tableofcontents 

\section{Introduction}\label{sec-Introduction}
The Hyperbolic Conservation Law (HCL) is a homogeneous hyperbolic system of 1st-order quasi-linear divergence form equations, which can be uniformly represented as follows: 
\begin{equation}\label{HCL}
\partial_t U + \nabla \cdot \mathcal{F}(U) = 0 \quad \text{in} \quad \Omega \subset \mathbb{R}^{d},
\end{equation}
here,  $d$ denotes the spatial dimension, and {\color{red}{$U$ can represent either a scalar or a vector}}.  \\
In the one-dimensional case, 
$$\nabla\cdot \mathcal{F}(U)=\partial_x(\mathcal{F}(U));$$
For higher-dimensional scenarios, 
$$
\begin{aligned}
\mathcal{F}(U) &= (F_1(U), F_2(U), \ldots, F_d(U)), \\
\nabla\cdot\mathcal{F}(U) &= \sum_{i=1}^{i=d}\partial_{x_i} F_i(U).
\end{aligned}
$$
Hyperbolic conservation laws are utilized to describe physical phenomena where the increment of extensive quantities within a system is balanced by the net flux passing through the system's boundaries \cite{ref2}. 
Common examples of hyperbolic conservation laws include scalar transport equations, Euler equations, and source-free shallow water wave equations. 
\par
Numerical solutions to equation (\refeq{HCL}) commonly employ schemes such as the finite difference method (FD), 
the finite volume method (FV), and the discontinuous Galerkin method (DG). These methods share a common characteristic in that they are all ``local'' approximation schemes: 
first, the FD method, which is derived from a Taylor expansion around the grid framework points, inherently belongs to a local approximation technique; second, 
the FV method can be regarded as a special case of the DG method when using $P^0$-polynomial approximation, 
and in the DG method, the basis functions are entirely discontinuous, i.e., they are piecewise defined according to the grid cells, 
allowing for different orders or types of basis functions on different grid cells. 
Accordingly, the DG method independently and in parallel carries out variation on each grid cell, ultimately yielding a series of weak equivalent integral equations that are as numerous as the grid cells. 
To distinguish from the global variational process of traditional finite element methods, we refer to this as local variation. 
\par
The FD scheme, when reconstructing the flux value $\mathcal{F}_{i+1/2}$ at the midpoint of the grid, 
faces the issue of choosing which side of the cell interface to use for interpolation based on the grid framework points. 
In both FV and DG schemes, each cell generates boundary integral terms during the variational process. 
However, due to the local nature of the variational form, the flux function $\mathcal{F}$ is multivalued at the common interface of adjacent cells. 
Numerical flux techniques are employed in local approximation schemes to ensure the uniqueness of the flux at the common interfaces of neighboring cells. 
That is, numerical flux techniques map a multivalued function to a single-valued function, thereby restoring the mathematical structure and physical connection of the field functions within different cells.
Let the mesh be partitioned as $\mathcal{T}_h$. Consider any cell $\Omega_k \in \mathcal{T}_h$, with its boundary denoted by $\partial\Omega_k$. Let $\Gamma_k$ denote any part of the boundary $\partial\Omega_k$, such that $\Gamma_k \subset \partial\Omega_k$. Denote $C(\Gamma_k)$ as the set of all cells sharing the boundary segment $\Gamma_k$, and suppose there are $q$ cells sharing $\Gamma_k$, i.e., $|C(\Gamma_k)| = q$. The numerical solution $U_h$ is a multivalued function on the cell interface $\Gamma_k$, with $q$ states $\left({\partial U}_1, {\partial U}_2, \cdots, {\partial U}_c, \cdots, {\partial U}_q\right)$, where ${\partial U}_c = \left. U_h^{(c)} \right|_{\Gamma_k}$, $c \in C(\Gamma_k)$. The numerical flux $\widehat{\mathcal{F}}$ is a single-valued function, such that
\begin{equation}
\left\{
\begin{aligned}
& \left. U_h \right|_{\Gamma_k} = \left({\partial U}_1, {\partial U}_2, \cdots, {\partial U}_c, \cdots, {\partial U}_q\right), \quad c \in C(\Gamma_k), \\
& \widehat{\mathcal{F}}: \left. U_h \right|_{\Gamma_k} \mapsto \widehat{\mathcal{F}}_{\Gamma_k}, \quad \forall \Gamma_k \in \partial\Omega_k.
\end{aligned}
\right.
\end{equation}
\par
Upwind numerical flux formats are one of the common numerical flux techniques, 
and this method of constructing numerical fluxes originates from the characteristic line theory of hyperbolic equations: 
along the reverse direction of the characteristic line (which we refer to as the upwind direction), 
the state quantities located in the dependency domain of the hyperbolic equation are selected to participate in the calculation of the numerical flux. 
The upwind direction is easy to determine for scalar hyperbolic equations: since $\frac{\partial f(U)}{\partial x} = \frac{\partial f}{\partial U} \frac{\partial U}{\partial x}$, 
the upwind direction can be determined based on the sign of $\frac{\partial f}{\partial U}$. 
Considering variable coefficients and nonlinear cases, the upwind direction itself is also a local concept, 
meaning that the trend of characteristic lines, and thus the upwind direction, is different in different computational regions. 
However, for systems, since $\frac{\partial f}{\partial U}$ will be a matrix, which we refer to as the Jacobian matrix, 
the aforementioned approach of ``judging the wind direction based on the sign of $\frac{\partial f}{\partial U}$'' cannot be directly generalized. 
Scalar hyperbolic conservation laws involve only a single characteristic wave, 
but the Jacobian matrix of a hyperbolic conservation system has more than one eigenvalue, corresponding to multiple characteristic waves. 
Thus, the propagation direction of each characteristic wave can be entirely inconsistent. 
As long as the Jacobian of the system's flux function is indefinite over the entire flow field (the eigenvalues of the flux function's Jacobian at some point in the flow field are not all of the same sign), 
then from a global flow field perspective, each conservation variable $U^{(i)}$ corresponds to a flux $\mathcal{F}_{(i)}$ 
that is composed of two scalar sub-fluxes with opposite propagation directions $\left\{\mathcal{F}_{(i)}^{+}, \mathcal{F}_{(i)}^{-}\right\}_{i=1}^{i=m}$. Here, $i$ is the system component index, 
and $m$ is the dimension of the system's state. Therefore, applying the ``upwind philosophy'' to hyperbolic conservation systems can be achieved by splitting the flux $\mathcal{F}$ into two sub-fluxes with opposite propagation directions, 
$\mathcal{F}^{+}$ and $\mathcal{F}^{-}$, and then constructing corresponding upwind schemes for each sub-flux. 
\par
Flux splitting techniques, commonly seen in FD schemes and FV schemes, include two categories: flux vector splitting (FVS) and flux difference splitting (FDS) \cite{ref3,ref4}. 
Flux Vector Splitting (FVS) schemes can be further divided into two categories: those based on Jacobian eigenvalue splitting and those based on Mach number splitting. 
The former requires the flux function to satisfy homogeneity, i.e., $\mathcal{F} = \frac{\partial \mathcal{F}}{\partial U} \cdot U$. 
The latter is designed for some practical fluid control equations, such as the compressible Euler equations and shallow water wave equations. Let $A = \frac{\partial \mathcal{F}}{\partial U}$, where $A$ is referred to as the Jacobian matrix of the flux function $\mathcal{F}$. 
In the Jacobian-FVS scheme, the Jacobian matrix $A$ is diagonalized by a similarity transformation, 
and then its eigenvalues are split to obtain $A^{+}$ and $A^{-}$, which satisfy positive definiteness and negative definiteness respectively. 
This allows the flux $\mathcal{F}$ to be split into positive flux $F^{+}$ and negative flux $F^{-}$ \cite{ref5}, 
and subsequently, numerical fluxes are constructed for each according to their respective upwind directions. 
Some hyperbolic conservation systems do not satisfy $\mathcal{F} = \frac{\partial \mathcal{F}}{\partial U} \cdot U$, such as the shallow water wave equations. 
We forcibly modify the wave speed to $a^* = \sqrt{\frac{1}{2} g h}$ to satisfy the homogeneity requirement of the flux function, 
and then perform flux splitting based on the modified Jacobian eigenvalues \cite{ref6}. 
Regardless of whether the flux function of a hyperbolic conservation system satisfies the homogeneity requirement, 
the Mach-FVS scheme can be employed. In the Mach-FVS scheme, the characteristic wave propagation direction is no longer determined by the eigenvalues of the Jacobian matrix. 
Instead, the local Mach number $ M_a = \frac{u}{a} $ (where $ a $ is the local speed of sound) is used to determine the characteristic propagation direction at that point: when $ M_a \geq 1 $, 
the characteristic propagates in the positive direction; when $ M_a \leq -1 $, it propagates in the negative direction; 
and when $ |M_a| < 1 $, both positive and negative characteristic waves coexist. 
Similar to the Jacobi-FVS splitting of the Jacobian eigenvalues, a natural idea is to achieve flux vector splitting by splitting the Mach number $ M_a $. 
Depending on whether the pressure term is split separately, the Mach-FVS can be further divided into the van Leer splitting method \cite{ref3,ref4} and the AUSM splitting method \cite{ref7}. 
The former includes the pressure term in the flux function for splitting, while the latter separates the pressure term from the flux function for individual splitting. 
\par
Reed and Hill \cite{ref8} proposed the first Discontinuous Galerkin method in 1973, for solving the neutron transport equation. 
Since then, the DG method has developed rapidly, with the emergence of local discontinuous Galerkin methods (LDG) \cite{ref9}, 
direct discontinuous Galerkin methods (DDG) \cite{ref10}, and ultra-weak discontinuous Galerkin methods (UWDG) \cite{ref11} for solving higher-order equations. 
In addition, semi-Lagrangian discontinuous Galerkin methods (SLDG) \cite{ref12} and Euler-Lagrange discontinuous Galerkin (ELDG) \cite{ref13} have been developed based on the characteristic line theory of hyperbolic equations from a Lagrangian perspective. 
These methods exhibit high-order structure-preserving and unconditional stability (large time-step stability) characteristics when solving transport equations. 
The DG method discussed in this paper is the Runge-Kutta Discontinuous Galerkin (RKDG) method. RKDG \cite{ref14,ref15,ref16,ref17,ref18} is a commonly used DG scheme for solving hyperbolic conservation systems such as the compressible Euler equations. 
It employs explicit and nonlinearly stable high-order Runge-Kutta methods for temporal discretization and uses DG methods for spatial discretization. The former ensures the nonlinear stability of the method regardless of accuracy, while the latter integrates the concepts of numerical fluxes and slope limiters from high-resolution FD and FV schemes. 
The resulting RKDG method is stable, high-order accurate, and highly parallelizable \cite{ref19}, capable of easily handling complex geometries and boundary conditions. The RKDG method has the advantages of local conservation, arbitrary triangular meshing, good parallel efficiency, h-p adaptive capability \cite{ref20}, and certain superconvergence properties \cite{ref21,ref22,ref23,ref24,ref25}.
\par
Even if the initial conditions are smooth, hyperbolic conservation law may develop discontinuities during their evolution. 
When dealing with discontinuous problems, high-order schemes such as the RKDG method can exhibit numerical pseudo-oscillations (non-physical oscillations) near discontinuities, 
which can lead to instability in the numerical scheme. 
To suppress numerical oscillations, a common practice is to introduce limiters. 
There are currently two types of common limiters: 
TVD/TVB-type limiters \cite{ref1,ref26} and ENO/WENO-type limiters \cite{ref27,ref28,ref29,ref30,ref31,ref32}. 
\par
TVD (Total Variation Diminishing) limiters initially decompose a high-order scheme into a first-order upwind scheme and a linear sum of a correction term (the first-order scheme has sufficient dissipation to automatically satisfy the non-oscillation property). 
Subsequently, according to Harten's lemma \cite{ref33}, the correction term is restricted. 
The requirement is that in smooth regions, the limiter does not affect the correction term, 
while near discontinuities, the limiter sets the correction term to zero, 
thereby automatically degrading to a stable, non-oscillatory first-order upwind scheme. 
TVD limiters tend to lose accuracy near extremum points, leading to the development of TVB (Total Variation Bounded) limiters, 
which only require that the average total variation be bounded. 
WENO (Weighted Essentially Non-Oscillatory) limiters utilize several stencils composed of the troubled cell and its neighboring cells. 
Based on the cell integral mean, several template polynomials are obtained through interpolation. 
Subsequently, a template smoothness measurement function is introduced to non-linearly weight and combine the various template polynomials to obtain the final reconstructed polynomial on the troubled cell. 
The smoother the stencil, the greater the weight of the corresponding template polynomial, thereby achieving a smooth filtering of the pseudo-oscillations near discontinuities \cite{ref27,ref34}. 
\par
From the above, it is known that limiters should act near discontinuities; using limiters in smooth regions can easily lead to loss of accuracy and increase unnecessary computational time costs. Therefore, a discontinuity indicator can be introduced to identify troubled cells, i.e., those that contain discontinuities. Accordingly, the steps for using the limiter are:
(1) The discontinuity indicator identifies troubled cells;
(2) The limiter is applied to troubled cells.

It should be noted that some limiters, such as the TVB(D)-minmod limiter, have built-in discontinuity indication functions, 
while WENO limiters require additional configuration of a discontinuity indicator \cite{ref35}. 
An obvious idea is to remove the slope correction part of the TVB(D)-minmod limiter and keep the remaining part as a discontinuity indicator for WENO. 
In addition, a commonly used discontinuity indicator is the KXRCF discontinuity indicator \cite{ref36}. 
This indicator is based on the conclusion that ``DG algorithms have strong superconvergence on the outflow boundaries of each cell in smooth areas''. 
If a cell contains a discontinuity, then its strong superconvergence on the outflow boundary will be destroyed, 
which in turn can determine that the cell is a troubled cell. 
\par
When applying limiters to a system of equations, there are typically two implementation methods: 
(a) reconstruct each component individually; (b) in conjunction with local characteristic decomposition \cite{ref37}.
The former is easy to implement in a program; reconstructing directly by component is essentially assuming that the components are already decoupled, but in the original system of equations, the unknown components are usually coupled, hence the actual effect may not be good.
The latter leverages the Jacobian of the flux function to transform the system of equations into the characteristic space to achieve true decoupling, and then reconstructs each component individually.
\par
The main work of this paper is as follows: For the first time, the FVS method is introduced into the DG scheme. 
By utilizing $\mathcal{F} = \mathcal{F}^{+} + \mathcal{F}^{-}$, 
an equivalent flux-split form of the original differential equation is obtained:
$ U_t + \nabla \cdot \mathcal{F}^{+} + \nabla \cdot \mathcal{F}^{-} = 0 $. 
A variation is then performed based on the flux-split form of the original equation, 
leading to the derivation of a numerical flux format based on FVS within the DG framework. 
Specifically, for the two-dimensional hyperbolic conservation system, 
FVS is implemented using normal fluxes in the outward normal direction on the interfaces of cells. 
In terms of limiters, we have introduced the smoothness measurement function $IS$ from the WENO limiter into the TVB(D)-minmod limiter. 
An optimization problem based on the smoothness factor constraint has been constructed, 
realizing a TVB(D)-minmod limiter applicable to arbitrary high-order polynomial approximation, which we refer to as the IS-TVB(D)-minmod limiter. 
Thanks to the constraint of the smoothness factor IS, the IS-TVB(D)-minmod limiter possesses strong capabilities in suppressing numerical oscillations, 
but correspondingly, it is prone to losing high-order accuracy. 
To address this issue, we have adopted the approach from the literature \cite{ref1}, 
introducing the $L^2$-error between the reconstructed polynomial and the original high-order scheme into the objective function. 
We have constructed a ``smoothness factor IS-$L^2$-error constraint dual-objective optimization'' problem 
to further improve the TVB(D)-minmod limiter, which we refer to as IS-$L^2$-TVB(D)-minmod limiter.  
Furthermore, when reconstructing for the hyperbolic systems, it is often necessary to accompany it with local characteristic decomposition,   
and the transformation matrix required for characteristic decomposition is determined by local freezing techniques. 
To this end, we propose to use the integral mean on the cell interfaces instead of the traditional approach of using the cell integral mean to achieve local freezing. 
After obtaining the characteristic transformation matrix, it is necessary to transform the DG polynomials from the physical space to the characteristic space. 
The traditional approach is to construct a column vector of modal coefficients for the same degree terms of each system component's DG polynomial, 
and then perform characteristic transformation through matrix-vector multiplication, referred to as ``moment characteristic transformation''. 
In contrast, we have proposed an interpolation-based polynomial characteristic transformation method to achieve local characteristic decomposition, 
and in the one-dimensional case, we have proven the equivalence of the interpolation-based characteristic transformation with the moment characteristic transformation. 
\par
The structure of this paper is organized as follows:
In Section \ref{sec-Overview of RKDG and Roe Average}, a review of the Runge-Kutta Discontinuous Galerkin (RKDG) method and Roe average is presented.
Subsequently, Section \ref{sec-FVS} introduces two implementations of Flux Vector Splitting (FVS) for one-dimensional hyperbolic conservation systems, exemplified by Jacobian eigenvalue-based splitting and Mach number-based splitting. 
Section \ref{sec-1D-FVSDG} details the flux vector splitting formulation for the RKDG scheme applied to one-dimensional hyperbolic conservation systems, 
while the construction process for the two-dimensional case is elaborated in Section \ref{sec-2D-FVSDG}.
Section \ref{sec-Constrained TVB(D)-minmod Limiter} focuses on the improvements made to the TVB(D)-minmod limiter. 
Section \ref{sec-characteristic-reconstruction} is dedicated to showcasing our unique perspectives on local characteristic decomposition.
The effectiveness of the aforementioned work is validated through carefully designed numerical experiments in Section \ref{sec-NumericalTests}.
Finally, Section \ref{sec-Conclusion} provides a summary of the work presented in this paper. 
\par
{\color{blue}{
  Note that a detailed explanation of the symbols and operational rules used in the definition of the FVS-DG weak solution process 
  in Section \ref{sec-1D-FVSDG} and Section \ref{sec-2D-FVSDG} must be referred to Appendix \ref{appendix-notations}.
  Furthermore, Appendix \ref{appendix-SWE-FVS} provides the FVS-DG formulation for the shallow water wave equations.
  Lastly, Appendix \ref{appendix-FVS-scalar} re-examines the construction process and working principles of the classical Lax-Friedrichs flux scheme from the perspective of Jacobian-FVS.
  Inspired by the Steger-Warming splitting method for systems, Appendix \ref{appendix-FVS-scalar} presents the Steger-Warming flux scheme for scalar transport equations and the Burgers' equation, and proves that this flux scheme satisfies compatibility, Lipschitz continuity, and monotonicity, thereby ensuring that the weak solutions obtained for scalar transport equations and the Burgers' equation using this flux scheme meet the cell entropy inequality and $L^2$-stability. 
}}

\section{Overview of RKDG and Roe Average}\label{sec-Overview of RKDG and Roe Average}
The RKDG method serves as the main framework for the algorithm presented in this paper, 
while the Roe average is essential for the local freezing used in the characteristic decomposition 
of hyperbolic conservation systems. A brief review of both techniques is provided below. 
\subsection{Runge-Kutta Discontinuous Galerkin Method}\label{subSec-RKDG}
\hfill\par
Let the computational domain be $\Omega \subset \mathbb{R}^{\mathrm{d}}$, where $d = 1, 2, \cdots$. 
The mesh partitioning is denoted as $\mathcal{T}_h = \Omega_1 \cup \Omega_2 \cup \cdots \cup \Omega_N$, 
with the property that $\Omega_i \cap \Omega_j = \varnothing$, for all $i, j \in \left\{1, 2, \cdots, N \right\}$.
\par
The definition of the local orthogonal basis functions on an arbitrary domain $\Omega_i$ is as follows:
\begin{definition}[Normative Orthogonal System]
  Let $\left\{\phi_{0}^{(\mathcal{D})}, \phi_{1}^{(\mathcal{D})}, \cdots, \phi_{\ell}^{(\mathcal{D})}\right\}$ (where $\ell$ is finite or $\ell \in \infty$) be a set of functions defined on a bounded domain $\mathcal{D}$, with $\phi_{i}^{(\mathcal{D})} \in L^2(\mathcal{D})$, for all $i = 0, 1, 2, \cdots, \ell$.
  We refer to $\left\{\phi_{i}^{(\mathcal{D})}\right\}_{i=0}^{i=\ell}$ as an normative orthonormal system on the bounded domain $\mathcal{D}$, if 
\begin{equation}\label{orthogonality}
  \begin{aligned}
  \left<\phi_j^{(\mathcal{D})}, \phi^{(\mathcal{D})}_k\right>_{L^2(\mathcal{D})}=\int_{\mathcal{D}}\phi_j^{(\mathcal{D})}\phi^{(\mathcal{D})}_k \mathrm{~d}\mathbf{X}=\delta_{jk}=
  \left\{
  \begin{array}{ll}
  0, & j \neq k \\
  1, & j=k
  \end{array}
  \right.
  \end{aligned}
\end{equation}
where,  $\mathbf{X}\in \mathbb{R}^{\mathrm{d}}, d=1, 2, \ldots.$
\end{definition}
\par
Let $P^{K}(\Omega_i)$ be the space of polynomials of degree $K$ defined on $\Omega_i$. It can be spanned by any normative orthonormal system $\left\{\phi_m^{(i)}\right\}_{m=0}^{m=K}$ on it. The set $\left\{\phi_m^{(i)}\right\}_{m=0}^{m=K}$ satisfies:
\begin{itemize}
\item $\left\{\phi_m^{(i)}\right\}_{m=0}^{m=K}$ is an normative orthonormal system in $P^{K}(\Omega_i)$,
\item $\forall q \in P^{K}(\Omega_i),\ \exists!\left\{\alpha_m\right\}_{m=0}^{m=K} \subset \mathbb{R},\ s.t.\ q = \sum_{m=0}^{m=K} \alpha_m \phi_m^{(i)}$.
\end{itemize}
On this basis, the space $\mathbb{V}^{K}_h(\Omega)$, where the DG weak solution resides, can be defined: 
\begin{definition}[DG Weak Solution Space]
\begin{equation}
  \mathbb{V}^{K}_h(\Omega)=\left\{ v \in L^2(\Omega):v\left|\right._{\Omega_i} \in P^{K}(\Omega_i), \ \forall \Omega_i \in \mathcal{T}_h \right\}
  \nonumber
\end{equation}
here, $\mathcal{T}_h\ $is a partitioning of $\Omega \subset \mathbb{R}^\mathrm{d}$. 
\end{definition}
Thus, for Equation (\refeq{HCL}), the definition of the DG weak solution is as follows: 
\begin{definition}[DG Weak Solution]
In the DG weak solution space, a unique polynomial function $U_h \in \mathbb{V}^K_h$ is determined such that for all $V \in \mathbb{V}^K_h$ and for all $\Omega_i \in \mathcal{T}_h$, the following is satisfied:
\begin{equation}\label{definition-DG-space}
\int_{\Omega_i} V \partial_t (U_h) \, \mathrm{d} \mathbf{X} = \int_{\Omega_i} \mathcal{F}(U_h) \cdot (\nabla V) \, \mathrm{d} \mathbf{X} - \int_{\partial\Omega_i} \widehat{\mathcal{F}_{\nu}} V \, \mathrm{d} \sigma.
\end{equation}
Here, $\widehat{\mathcal{F}_{\nu}}$ is the numerical normal flux, which is a single-valued function on $\partial\Omega_i$ and satisfies
$\widehat{\mathcal{F}_{\nu}} \approx \mathbf{n} \cdot \mathcal{F}(U_h)$, 
where $\mathbf{n}$ is the unit outward normal vector on $\partial\Omega_i$.
\end{definition}
\begin{remark}
The construction of the numerical flux $\widehat{\mathcal{F}_{\nu}}$ to approximate $\mathbf{n} \cdot \mathcal{F}(U_h)$ is key to the spatial discretization in the DG algorithm. 
The literature \cite{ref38,ref39} introduces a variety of typical numerical flux formats, 
while our work involves the introduction of flux splitting methods from FD/FV schemes, especially the flux vector splitting method, 
to construct $\widehat{\mathcal{F}_{\nu}}$ in the DG scheme.
\end{remark}
\par
Let $U_h^{(i)} := U_h|_{\Omega_i} \in P^K(\Omega_i)$, then we have
\begin{equation}
U_h^{(i)}(\mathbf{X}, t) = \sum_{\ell=0}^{\ell=K} \alpha_{\ell}^{(i)}(t) \phi_{\ell}^{(i)}(\mathbf{X}),
\end{equation}
The $\alpha_{\ell}^{(i)}$ is referred to as the $\ell$-th moment or the $\ell$-th modal coefficient of $U_h^{(i)}$, satisfying
\begin{equation}
  \alpha_{\ell}^{(i)} = \left\langle U^{(i)}_h, \phi_{\ell}^{(i)} \right\rangle_{L^2(\Omega_i)}.
\end{equation}
Note that the set $\{\phi_{\ell}^{(i)}\}_{\ell=0}^{\ell=K}$ is normative orthonormal.
\par
On $\Omega_i$, take the test function $V = \phi_r^{(i)}$. Note that the basis function $\phi_{\ell}^{(i)}$ depends only on $\mathbf{X}$, that is, it does not evolve with time. Thus, we have
\begin{align}
\int_{\Omega_i} \phi_r^{(i)} \partial_t (U_h) \, \mathrm{d} \mathbf{X} = \int_{\Omega_i} \phi_r^{(i)} \partial_t \left( U_h^{(i)} \right) \, \mathrm{d} \mathbf{X} = \sum_{\ell=0}^{\ell=K} \partial_t \left( \alpha^{(i)}_{\ell} \right) \cdot \left< \phi_r^{(i)}, \phi_{\ell}^{(i)} \right>_{L^2(\Omega_i)}.
\end{align}
Utilizing the orthogonality \eqref{orthogonality}, it immediately follows that
\begin{equation}\label{evolutionary-moment}
  \int_{\Omega_i} \phi_r^{(i)} \partial_t (U_h) \, \mathrm{d} \mathbf{X} = \partial_t \left( \alpha_r^{(i)} \right).
\end{equation}
From equations (\refeq{definition-DG-space}) and (\refeq{evolutionary-moment}), the following semi-discrete spatial format can be derived:
\begin{equation}
\frac{\mathrm{d}}{\mathrm{d} t} \left( \alpha^{(i)}_r \right) = DG\left( U_h^{(i)} \right),
\end{equation}
where $i = 1, 2, \cdots, N; \ r = 0, 1, \cdots, K$.
\par
The following presents the 3rd-order total variation diminishing Runge-Kutta temporal discretization scheme (TVD-RK3)  
and the 4-stage 4th-order non-strong stability preserving Runge-Kutta time discretization format (RK4): 
Let $\Delta t > 0$, and partition the time domain as $t_n = n \cdot \Delta t$, for $0 \leq n \leq M = T / \Delta t$. Let $L_h$ be some spatial discretization operator, such as the previously mentioned discontinuous Galerkin spatial discretization scheme $DG(\cdot)$. Then, 
\par \noindent
$\bullet\ $TVD-RK3
\begin{equation}\label{TVD-RK3}
\left.
\begin{aligned}
u^{(1)}&=u^n+\Delta t \cdot L_h\left(u^n ; t_n\right), \\
u^{(2)}&=\frac{3}{4} u^n+\frac{1}{4} u^{(1)}+\frac{1}{4} \Delta t \cdot L_h(u^{(1)} ; t_n+\Delta t), \\
u^{(3)}&=\frac{1}{3} u^n+\frac{2}{3} u^{(2)}+\frac{2}{3} \Delta t \cdot L_h(u^{(2)} ; t_n+\frac{1}{2} \Delta t), \\
u^{n+1}&=u^{(3)}.
\end{aligned}
\right\}
\end{equation}
$\bullet\ $RK4
\begin{equation}\label{RK4}
\left.
\begin{aligned}
u^{(1)}&=u^n+\frac{1}{2} \Delta t \cdot L_h\left(u^n ; t_n\right), \\
u^{(2)}&=u^n+\frac{1}{2} \Delta t \cdot L_h(u^{(1)} ; t_n+\frac{1}{2} \Delta t), \\
u^{(3)}&=u^n+\Delta t \cdot L_h(u^{(2)} ; t_n+\frac{1}{2} \Delta t), \\
u^{(4)}&=\frac{1}{3}(-u^n+u^{(1)}+2 u^{(2)}+u^{(3)})+\frac{1}{6} \Delta t \cdot L_h(u^{(3)} ; t_n+\Delta t), \\
u^{n+1}&=u^{(4)}.
\end{aligned}
\right\}
\end{equation}
As for the 10-stage 4th-order strong stability preserving Runge-Kutta method ($\operatorname{SSPRK}(10,4)$), please refer to Appendix \ref{appendix-SSPRK(10,4)}.
\par
By combining the RK-type temporal discretization schemes such as (\refeq{TVD-RK3}), (\refeq{RK4}), or (\refeq{SSPRK(10,4)}) 
with the DG spatial discretization format (\refeq{definition-DG-space}), one obtains the fully discrete RKDG scheme, which allows the evolution of $U_h^n$ to $U_h^{n+1}$. This is succinctly denoted as
\begin{equation}
  U_h^{n+1} = RKDG\left(U_h^{n}\right).
\end{equation}
\par
To suppress numerical pseudo-oscillations, it is necessary to introduce a limiter ${\varLambda\varPi}_h$ to restrict $U_h^n$ before the temporal evolution from $t_n$ to $t_{n+1}$, yielding $\widetilde{U}_h^{n}$, that is,
\begin{equation}
  \widetilde{U}_h^{n} = {\varLambda\varPi}_h\left(U_h^n\right).
\end{equation}
Using $\widetilde{U}_h^{n}$ in the temporal evolution gives
\begin{equation}
  U_h^{n+1} = RKDG\left(\widetilde{U}_h^{n}\right).
\end{equation}
\par
Lastly, for convenience, we take the explicit Euler temporal discretization scheme as an example to provide a complete description of the spatial discretization, temporal discretization, and limiter correction processes in the DG algorithm:
\par \noindent
Determine a unique polynomial function $U_h^{n+1} \in \mathbb{V}^K_h$ in the DG weak solution space such that for all $V \in \mathbb{V}^K_h$ and for all $\Omega_i \in \mathcal{T}_h^d$, the following is satisfied:
\begin{equation}\label{2.12}
\int_{\Omega_i} V \frac{U_h^{n+1} - \widetilde{U}_h^{n}}{\Delta t} \, \mathrm{d} \mathbf{X} =
\int_{\Omega_i} \mathcal{F}\left(\widetilde{U}_h^n\right) \cdot (\nabla V) \, \mathrm{d} \mathbf{X}
- \int_{\partial\Omega_i} \widehat{\mathcal{F}_{\nu}}^{n} V \, \mathrm{d} \sigma,
\end{equation}
where $\widetilde{U}_h^{n} = {\varLambda\varPi}_h\left(U_h^n\right)$, and $\widehat{\mathcal{F}_{\nu}}^n$ is the normal numerical flux at time $t_n$.

\subsection{Roe Average}\label{subSec-Roe}
\hfill\par
For control equations of a system, the limiters introduced typically need to be used in conjunction with local characteristic decomposition. 
For variable-coefficient equations or nonlinear equations, local freezing techniques are required to determine the characteristic transformation matrix before performing local characteristic decomposition. 
A simple approach is naturally to take the arithmetic average of the left and right states (in 1D) or the internal and external states of the cell (in 2D), 
but a more effective method, or one that is more in line with the physical process, is the Roe average. 
The approach of Roe averaging is as follows: by utilizing the constant states of the left and right functions, 
$\mathbf{U}^L$ and $\mathbf{U}^{R}$, a reasonable constant matrix ${Q}\left(\boldsymbol{U}^L, \boldsymbol{U}^{R}\right)$ 
is constructed to approximate the original $A(\mathbf{U})$. 
This transformation simplifies the complex nonlinear problem into a linear one. Roe \cite{ref40} achieved this construction 
by ensuring that $Q\left(\mathbf{U}_{i+\frac{1}{2}}^L, \mathbf{U}_{i+\frac{1}{2}}^R\right)$ satisfies
\begin{equation}\label{4.1}
  \mathbf{F}\left(\mathbf{U}_{i+\frac{1}{2}}^R\right) - \mathbf{F}\left(\mathbf{U}_{i+\frac{1}{2}}^L\right) = Q\left(\mathbf{U}_{i+\frac{1}{2}}^L, \mathbf{U}_{i+\frac{1}{2}}^R\right) \left(\mathbf{U}_{i+\frac{1}{2}}^R - \mathbf{U}_{i+\frac{1}{2}}^L\right),
\end{equation}
thereby completing the construction. (Here, $\mathbf{F(U)}$ is flux function in hyperbolic systems.)
In the Roe scheme, the matrix $Q\left(\mathbf{U}_{i+\frac{1}{2}}^L, \mathbf{U}_{i+\frac{1}{2}}^R\right)$ 
is a linear approximation of the Jacobian matrix $A$ at a certain state $\mathbf{U}^{*}_{i+\frac{1}{2}}$. 
The state $\mathbf{U}^{*}_{i+\frac{1}{2}}$ is referred to as the Roe average of $\mathbf{U}^{L}_{i+\frac{1}{2}}$ 
and $\mathbf{U}^{R}_{i+\frac{1}{2}}$.
\par
In the subsequent sections, when employing characteristic reconstruction in the FVS-RKDG scheme for the hyperbolic systems, 
Roe averaging is used to achieve local freezing. Therefore, in this section, 
we provide the Roe averaging formulations for the compressible Euler equations in 1D and 2D/3D as examples.  
\par
The Roe averaging scheme for the 1D-Euler equations is given by:
\begin{equation}
\left.
\begin{aligned}
& \sqrt{\overline{\rho}^{Roe}} = \left(\sqrt{\rho_L} + \sqrt{\rho_R}\right) / 2, \\
& \sqrt{\overline{\rho}^{Roe}} \cdot \overline{u}^{Roe} = \left(\sqrt{\rho_L} u_L + \sqrt{\rho_R} u_R\right) / 2, \\
& \sqrt{\overline{\rho}^{Roe}} \cdot \overline{\widetilde{H}}^{Roe} = \left(\sqrt{\rho_L} \widetilde{H}_L + \sqrt{\rho_R} \widetilde{H}_R\right) / 2.
\end{aligned}
\right\}
\end{equation}
Note: $ \widetilde{H} $ is the specific enthalpy.
\par
The pressure $ \overline{P}^{Roe} $, speed of sound $ \overline{a}^{Roe} $, and total energy $ \overline{E}^{Roe} $ are computed indirectly from $ \overline{\rho}^{Roe} $, $ \overline{u}^{Roe} $, and $ \overline{\widetilde{H}}^{Roe} $ as follows:
\begin{equation}
\left.
\begin{aligned}
\overline{P}^{Roe} &= \frac{\gamma - 1}{\gamma} \left( \overline{\rho}^{Roe} \overline{\widetilde{H}}^{Roe} - \frac{1}{2} (\overline{\rho}^{Roe})^2 \right), \\
(\overline{a}^{Roe})^2 &= (\gamma - 1) \left( \overline{\widetilde{H}}^{Roe} - \frac{(\overline{u}^{Roe})^2}{2} \right), \\
\overline{E}^{Roe} &= \overline{\rho}^{Roe} \overline{\widetilde{H}}^{Roe} - \overline{P}^{Roe}.
\end{aligned}
\right\}
\end{equation}
Note: 
$ \overline{P}^{Roe} $, $ \overline{a}^{Roe} $, and $ \overline{E}^{Roe} $ are not obtained by directly averaging $ P $, $ a $, and $ E $ with density but are computed indirectly.
\par
The Roe matrix $ A^{Roe} $ is constructed as:
\begin{equation}
\left.
\begin{aligned}
& \overline{U}^{Roe}=\left(\overline{\rho}^{Roe},\ (\overline{\rho}^{Roe}\cdot\overline{u}^{Roe}),\ \overline{E}^{Roe}\right)^{\mathrm{T}}, \\
& A^{Roe} = A\left(\overline{\mathbf{U}}^{Roe}\right) = R\left(\overline{\mathbf{U}}^{Roe}\right) \cdot \Lambda\left(\overline{\mathbf{U}}^{Roe}\right) \cdot L\left(\overline{\mathbf{U}}^{Roe}\right) = \overline{R}^{Roe} \cdot \overline{\Lambda}^{Roe} \cdot \overline{L}^{Roe}, \\
& \left|A^{Roe}\right| = \overline{R}^{Roe} \cdot \left|\overline{\Lambda}^{Roe}\right| \cdot \overline{L}^{Roe}, \quad \left|\overline{\Lambda}^{Roe}\right| = \operatorname{diag}\left(\left|\overline{\lambda_1}^{Roe}\right|, \left|\overline{\lambda_2}^{Roe}\right|, \left|\overline{\lambda_3}^{Roe}\right|\right).
\end{aligned}
\right\}
\end{equation}
\par
For the 2D/3D case, the following additional formulas are included:
\begin{equation}
\left.
\begin{aligned}
& \overline{v}^{Roe} = \frac{\left(\sqrt{\rho_L} v_L + \sqrt{\rho_R} v_R\right)}{\sqrt{\rho_L} + \sqrt{\rho_R}}, \\
& \overline{w}^{Roe} = \frac{\left(\sqrt{\rho_L} w_L + \sqrt{\rho_R} w_R\right)}{\sqrt{\rho_L} + \sqrt{\rho_R}}, \\
& \overline{U}^{Roe}=\left(\overline{\rho}^{Roe},\ (\overline{\rho}^{Roe}\cdot\overline{u}^{Roe}),\ (\overline{\rho}^{Roe}\cdot\overline{v}^{Roe}),\ (\overline{\rho}^{Roe}\cdot\overline{w}^{Roe}),\ \overline{E}^{Roe}\right)^{\mathrm{T}}.
\end{aligned}
\right\}
\end{equation}
\begin{remark}
  In the 2D case, the Roe averaging scheme for the compressible Euler equations presented in this paper is based on the normal flux, 
  that is, the normal velocity $q_n = u \cdot n_x + v \cdot n_y$ is introduced. 
  For details, please refer to Section \ref{subSec-2D-DG-Roe} below. 
\end{remark}
\begin{remark}
  The specific form of the Roe average may vary for different fluid motion models. 
  For instance, in the case of the shallow water wave equations, the Roe averaging scheme can be referred to in \cite{ref41}. 
\end{remark}

\section{Preliminaries: Flux Vector Splitting for Hyperbolic Conservative System in One-dimension}\label{sec-FVS}
This section introduces the flux vector splitting method using the one-dimensional hyperbolic conservation system as an example. 
In general, homogeneous systems can be split based on the eigenvalues of the Jacobian matrix of the flux function, 
with the goal of obtaining both a positive definite matrix and a negative definite matrix whose sum equals the original Jacobian matrix. 
This allows the flux function to be divided into positive and negative flux components, 
where the positive flux corresponds to characteristic waves propagating in the positive direction of the coordinate system, 
and the negative flux corresponds to characteristic waves propagating in the negative direction. 
Consequently, the upwind direction for both can be clearly defined, and a stable upwind scheme can be employed for computation. 
For practical fluid equations such as the compressible Euler equations or shallow water wave equations, the Mach number $M_a$ 
(or the Froude number $F_r$ for shallow water waves) can be introduced to rewrite the flux function as a linear function of $M_a$, 
thereby achieving flux function splitting through the splitting of $M_a$.  

\subsection{Based on Jacobian eigenvalue Splitting}\label{subSec-FVS-Jacobian}
The purpose of splitting the flux function through Jacobian eigenvalue splitting is to construct an upwind-type numerical flux scheme. This method requires that the flux function itself be a homogeneous function. The following provides the definition of a homogeneous function:
\begin{definition}[Homogeneous Function]
A function $f$ is homogeneous if it satisfies $f(\alpha u) = \alpha f(u)$, for all $u \in D(f)$, where $D(f)$ is the domain of $f$, and $\alpha$ is any constant such that for a given $u$, $\alpha u \in D(f)$.
\end{definition}
It is easy to see that homogeneous functions satisfy
$ f'(\alpha u) \cdot u = f(u)$. 
\par
By taking $\alpha = 1$, we obtain
\begin{align}
f(u) = f'(u) \cdot u.
\end{align}
\par
Consider the following one-dimensional homogeneous system
$
\partial_t \mathbf{U} + \partial_x f(\mathbf{U}) = 0,
$
where
$
\mathbf{U} = (u_1, u_2, u_3, \cdots, u_m)^{\mathrm{T}},
$
and
$
f(\mathbf{U}) = (f_1(\mathbf{U}), f_2(\mathbf{U}), f_3(\mathbf{U}), \cdots, f_m(\mathbf{U}))^{\mathrm{T}},
$
with the property that
$
f(\mathbf{U}) = \frac{\partial f}{\partial \mathbf{U}} \cdot \mathbf{U}.
$
Let $ A = \frac{\partial f}{\partial \mathbf{U}} $, then $ f = A \mathbf{U} $.
\par
The general process of Jacobian-FVS is as follows: 
\begin{itemize}
  \item Step1.\ First, the flux function needs to be rewritten as $f(\mathbf{U})=A(\mathbf{U}) \cdot \mathbf{U}$. 
  {\color{blue}{
  \begin{remark}
  if $f(\mathbf{U})$ is a homogeneous function (e.g., Euler equations), then take $A(\mathbf{U})=\frac{\partial f}{\partial \mathbf{U}}$. 
  \end{remark}}}
  {\color{blue}{
  \begin{remark}
    The flux function $f$ of the shallow water wave equations is not a homogeneous function and requires special preprocessing. For example, the wave speed $a = \sqrt{gh}$ may be forcibly modified to $a^{*} = \sqrt{\frac{gh}{2}}$. For more details, see Appendix \ref{appendix-SWE-FVS}.  
  \end{remark}}}
  \item  Step2.\ Diagonalize $A$ by similarity transformation.  
    $$
    \begin{aligned}
    &T^{-1} A T=\Lambda, \\
    &A=T \Lambda T^{-1}=R \Lambda L. 
    \end{aligned}
    $$	
    Note: $R L=I_d(\text{identity matrix}),\ T=R,\ T^{-1}=L$. $A, R, L, \Lambda$\ are all matrix-valued functions of $\mathbf{U}$. 
  \item  Step3.\ Split $\Lambda$.
    $$
    \Lambda=\Lambda^{+}+\Lambda^{-}, 
    $$	
    satisfying
    $$\lambda\left(\Lambda^{+}\right) \geq 0, \ \lambda\left(\Lambda^{-}\right) \leq 0,$$
    that is, $\Lambda^+$ is positive definite, $\Lambda^-$ is negative definite. 
  \item  Step4.\ Split $A$. 
    $$
    \begin{aligned}
      A&=R\left(\Lambda^{+}+\Lambda^{-}\right) L=R \Lambda^{+} L+R \Lambda^{-} L=A^{+}+A^{-}, \\
      A^{+}&=R \Lambda^{+} L,\ A^{-}=R \Lambda^{-} L, 
    \end{aligned}
    $$ 
    satisfying$$\lambda\left(\mathrm{A}^{+}\right) \geq 0,\ \lambda\left(A^{-}\right) \leq 0,$$
    that is, $A^+$ is positive definite, $A^-$is negative definite.
  \item  Step5.\ Split the flux vector $f$. 
    $$f=f^+ + f^- \quad ,  \quad f^+=A^+ \mathbf{U} \quad ,  \quad f^-=A^- \mathbf{U}.$$
\end{itemize}
There are two common schemes to implement Step 3, namely the Steger-Warming splitting scheme and the Lax-Friedrichs splitting scheme. 
They are introduced separately below.  
\begin{remark}
  The flux function for the 1D-Euler equations, given by $ f = (\rho u, \rho u^2 + P, (E + u)P)^{\mathrm{T}} $, 
  and the normal flux function for the 2D-Euler equations, given by $ F_n = n_x \cdot F + n_y \cdot G $, 
  are homogeneous functions of the 1-dimensional conservative vector $ \mathbf{U} = (\rho, \rho u, E)^{\mathrm{T}} $ 
  and the 2-dimensional conservative vector $ \mathbf{U} = (\rho, \rho u, \rho v, E)^{\mathrm{T}} $, respectively. 
\end{remark}

\subsubsection{Steger-Warming Splitting}
The Steger-Warming splitting scheme (S-W Splitting) is as follows:
\\
\begin{align}
\lambda_k^{+} = \frac{\lambda_k + |\lambda_k|}{2} \quad , \quad \lambda_k^{-} = \frac{\lambda_k - |\lambda_k|}{2} \quad , \quad k = 1, 2, \cdots, m.
\end{align}
Then,
\begin{align}
\Lambda^{+} = \frac{1}{2} (\Lambda + |\Lambda|) \quad , \quad \Lambda^{-} = \frac{1}{2} (\Lambda - |\Lambda|).
\end{align}
Hence,
\begin{align}
A^{+} = R \frac{(\Lambda + |\Lambda|)}{2} L = \frac{1}{2} (A + |A|) \quad , \quad A^{-} = R \frac{(\Lambda - |\Lambda|)}{2} L = \frac{1}{2} (A - |A|).
\end{align}
Here, $ |A| = R |\Lambda| L $, where $ |\Lambda| $ denotes taking the absolute value of the elements in $ \Lambda $.
\\ 
If smoothness is considered, it can be improved as:
\begin{align}
\lambda_k^{+} = \frac{\lambda_k + \sqrt{|\lambda_k|^2 + \delta^2}}{2} \quad , \quad \lambda_k^{-} = \frac{\lambda_k - \sqrt{|\lambda_k|^2 + \delta^2}}{2},
\end{align}
where $ \delta $ is a small positive constant, for example, $ \delta $ can be taken as $ 10^{-8} $. 



\subsubsection{Lax-Friedrichs Splitting}
The Lax-Friedrichs splitting scheme (L-F Splitting) is as follows: 
\begin{align}
\lambda_k^{+}=\frac{\lambda_k+M}{2},\ \lambda_k^{-}=\frac{\lambda_k-M}{2},\ M>0.
\end{align}
Local L-F splitting: 
\begin{align}
M_j=\max _k\left\{\left|\lambda_k\left(x_j\right)\right|\right\}.
\end{align}
Global L-F splitting: 
\begin{align}
M=\max _j \max _k\left\{\lambda_k\left(x_j\right)\right\}.
\end{align}
Therefore, 
\begin{align}
\Lambda^{+}=\frac{1}{2}(\Lambda+M I_d) \quad ,  \quad \Lambda^{-}=\frac{1}{2}(\Lambda-M I_d).
\end{align}
Note that $R L=I_d$,  hence
\begin{align}
A^{+}=R \frac{(\Lambda+M I_d)}{2} L=\frac{1}{2}(A+M I_d)
\quad ,  \quad 
A^{-}=R \frac{(\Lambda-M I_d)}{2} L=\frac{1}{2}(A-M I_d). 
\end{align}
\begin{remark}
  From this, it can be seen that the L-F flux vector splitting is a ``simplified'' version of the S-W flux vector splitting, as L-F replaces the individual moduli of all eigenvalues with a common upper bound of their moduli, using this upper bound uniformly in the splitting process. 
  Conversely, S-W provides more precise control over artificial viscosity than L-F (S-W has less dissipation than L-F), because S-W involves the splitting with the modulus of each eigenvalue itself, rather than a simplified treatment using a single upper bound value. 
\end{remark}
\begin{remark}
  The S-W flux vector splitting and local L-F flux vector splitting possess adaptivity, meaning that at each cell boundary, the corresponding $|\lambda|_{i+1/2}, M_{i+1/2}$ are taken.
  In contrast, the global Lax-Friedrichs flux vector splitting scheme abandons this adaptivity, but it is more convenient to implement and introduces greater artificial viscosity, which is beneficial for suppressing numerical oscillations (at the cost of reduced accuracy). 
\end{remark}
\begin{remark}
  Noticing the similarity between the Lax-Friedrichs splitting method and the common Lax-Friedrichs flux format in DG schemes, 
  we may also attempt to develop the Steger-Warming splitting method into a flux format for scalar equations. 
  Details can be found in Appendix \ref{appendix-FVS-scalar}. 
\end{remark}

\subsection{Based on Mach Number Splitting}\label{subSec-FVS-Mach}
For practical fluid control equations, such as the compressible Euler equations and shallow water wave equations, a dimensionless parameter, 
the Mach number, can be defined as the ratio of the fluid velocity to the local speed of sound, denoted as $M_a = \frac{u}{a}$. 
The Mach number determines the characteristic propagation direction.
This section illustrates the concept using the 1D-Euler equations as an example, 
and for the Mach-FVS related to the shallow water wave equations, please refer to Appendix \ref{appendix-SWE-FVS}. 
\par
The eigenvalues of the Jacobian for the one-dimensional (1D) Euler equations are given as follows:  
\begin{align}
\lambda_1=u \quad \lambda_2=u-a \quad \lambda_3=u+a.
\end{align}
Note that the speed of sound $a=\sqrt{\frac{\gamma P}{\rho}} > 0$, where $\gamma$ is the specific heat ratio.  
\\
$\bullet\ $if $M_a>1$,  then $u>a>0$,  \\
hence $u>0,  u-a>0,  u+a>0$. \\
This implies that all three characteristic waves propagate in the positive direction.  \\
Consequently, there is no negative flux in the flux vector $\mathbf{F}$,  that is $\mathbf{F}=\mathbf{F}^{+},  \mathbf{F}^{-}=\mathbf{0}$. 
\\
$\bullet\ $if $M_a<-1$,  then $u<-a<0$, \\
hence $u<0,  u-a<0,  u+a<0$. \\
This implies that all three characteristic waves propagate in the negative direction.  \\
Consequently, there is no positive flux in the flux vector $\mathbf{F}$,  that is $\mathbf{F}=\mathbf{F}^{-},  \mathbf{F}^{+}=\mathbf{0}$. 
\\
$\bullet\ $If $|M_a| < 1$, the situation can be discussed in the following two cases:
\\
Case 1: $u > 0$, $u - a < 0$, $u + a > 0$ (both positive and negative fluxes exist);
\\
Case 2: $u < 0$, $u - a < 0$, $u + a > 0$ (both positive and negative fluxes exist).
\begin{remark}
  When $|M_a| > 1$, that is, $M_a > 1$ or $M_a < -1$, it is referred to as supersonic; when $|M_a| < 1$, it is referred to as subsonic.
\end{remark}

\subsubsection{van Leer Splitting}
The split fluxes $f^{+}$ and $f^{-}$ of the Steger-Warming scheme are not continuously differentiable near the sonic point, 
which can lead to errors in calculations. 
In contrast, the van Leer split fluxes are continuously differentiable near the sonic point \cite{ref3}.
\par
Rewrite the flux function $\mathbf{F}$ of the 1D-Euler equations in the van Leer form: 
\begin{align}
\mathbf{F} & =\left[\begin{array}{c}
\rho u \\
\rho u^2+P \\
(E+P) u
\end{array}\right]=\left[\begin{array}{c}
\rho u \\
\rho u^2+P \\
\left(\frac{1}{2} \rho u^2+\frac{P}{\gamma-1}+P\right) u
\end{array}\right]=\left[\begin{array}{c}
\rho u \\
\rho u^2+P \\
\left(\frac{1}{2} \rho u^2+\frac{\gamma P}{\gamma-1}\right) u
\end{array}\right] \nonumber \\
& =\left[\begin{array}{c}
\rho u \\
\rho u^2+\frac{\rho a^2}{\gamma} \\
\left(\frac{1}{2} \rho u^2+\frac{\gamma}{\gamma-1} \frac{\rho a^2}{\gamma}\right) u
\end{array}\right]=\left[\begin{array}{c}
\rho a M_a \\
\rho a M_a u+\frac{\rho a^2}{\gamma} \\
\left(\frac{1}{2} M_a u+\frac{a}{\gamma-1}\right) \rho a u
\end{array}\right] 
=\rho a\left[\begin{array}{l}
M_a \\
u M_a+a / \gamma \\
u^2 M_a / 2+a u /(\gamma-1)
\end{array}\right]. 
\end{align}
Then achieve the splitting of $\mathbf{F}$ by splitting $M_a$: 
\begin{itemize}
\item When $M_a \geq 1$ (characteristics propagate to the right, flux is positive),  
\begin{align}
\mathbf{F}=\mathbf{F}^{+}\left(\mathbf{F}^{-}=0\right);
\end{align}

\item When $M_a \leq-1$ (characteristics propagate to the left, flux is negative),
\begin{align}
\mathbf{F}=\mathbf{F}^{-}\left(\mathbf{F}^{+}=0\right);
\end{align}

\item {\color{red}{$\left|M_a\right|<1$, that is, $-1<M_a<1$, }}
\begin{align}
\mathbf{F}=\mathbf{F}^{+} + \mathbf{F}^{-} \quad ,  \quad
\mathbf{F}^{ \pm}=\left[\begin{array}{c}
1 \\
{[(\gamma-1) u \pm 2 a] / \gamma} \\
{[(\gamma-1) u \pm 2 a]^2 /\left[2\left(\gamma^2-1\right)\right]}
\end{array}\right] \cdot \rho a M_{a}^{\pm},
\end{align}	
here,  
{\color{red}{
\begin{align}
M_a^{+}=\frac{\left(M_a+1\right)^2}{4},  \quad M_a^{-}=-\frac{\left(M_a-1\right)^2}{4},
\end{align}
satisfying
\begin{align}
  M_a=M_a^{+}+M_a^{-}.
\end{align}
}}
\end{itemize}

\subsubsection{Liou-Stenffen Splitting (AUSM-type Methods)}
The fundamental idea of the AUSM (Advection Upstream Splitting Method) scheme is to consider the advection waves 
(related to the characteristic velocity $u$, linear) and the acoustic waves 
(related to the characteristic velocities $u+a$ and $u-a$, nonlinear) as physically distinct processes. 
Therefore, the inviscid flux is split into advection and pressure flux terms for separate treatment. 
In terms of scheme construction, the AUSM scheme is an evolution and improvement of the van Leer scheme. 
However, it distinguishes between linear and nonlinear characteristics, splits the pressure term and the advection flux, 
obtains the ``upwind'' property by judging the specially constructed interface Mach number, and accurately captures shocks through the careful design of a unified sound speed in the adjacent cells at the interface. 
From the analysis of its dissipation term, this is a composite scheme of FVS and FDS.
The AUSM scheme theoretically distinguishes between the linear field (related to the characteristic velocity $u$) 
and the nonlinear field (related to the characteristic velocities $u \pm a$) in the flow's advection characteristics and splits the pressure term and the advection flux separately. 
\par
The flux function $\mathbf{F}$ of the one-dimensional (1D) Euler equations, in the Liou-Stenffen form, is given as follows:  
\begin{align}
\mathbf{F}=\left[\begin{array}{l}
\rho u \\
\rho u^2+P \\
u(E+P)
\end{array}\right]=\left[\begin{array}{l}
\rho u \\
\rho u^2+P \\
\rho u \widetilde{H}
\end{array}\right]=\left[\begin{array}{l}
{\color{red}{\rho a M_a}} \\
{\color{red}{\rho a M_a u}}+{\color{blue}{P}} \\
{\color{red}{\rho a M_a \widetilde{H}}}
\end{array}\right]={\color{red}{\underbrace{\rho a M a\left[\begin{array}{l}
1 \\
u \\
\widetilde{H}
\end{array}\right]}_{\mathbf{F}^a}}}+{\color{blue}{\underbrace{\left[\begin{array}{l}
0 \\
P \\
0
\end{array}\right]}_{\mathbf{F}^P}}},
\end{align}
where $\mathbf{F}^a$ is referred to as the advection flux and $\mathbf{F}^P$ as the pressure flux, 
satisfying $\mathbf{F} = \mathbf{F}^a + \mathbf{F}^P$. 
\par
The Mach number $ M_a $ is split to thereby split the advection flux:
$$
M_a = M_a^{+} + M_a^{-},
$$
where
\begin{subequations}
  \begin{align}
    M_a^{+} &= \begin{cases}
    M_a & , \quad M_a > 1 \\
    \left( M_a + 1 \right)^2 / 4 & , \quad |M_a| \leq 1 \\
    0 & , \quad M_a < -1
    \end{cases} \\
    M_a^{-} &= \begin{cases}
    0 & , \quad M_a > 1 \\
    -\left( M_a - 1 \right)^2 / 4 & , \quad |M_a| \leq 1 \\
    M_a & , \quad M_a < -1
    \end{cases}
  \end{align}
\end{subequations}
Thus,
\begin{align}
\mathbf{F}^{a}_{\pm} = \rho a M_{a}^{\pm} \left[ \begin{array}{c}
1 \\
u \\
\widetilde{H}
\end{array} \right].
\end{align}
\par
Splitting the pressure term thereby leads to the splitting of the pressure flux:
\begin{align}
P = P^{+} + P^{-},
\end{align}
where
\begin{subequations}
\begin{align}
&P^{+} = \left\{
\begin{array}{ccc}
P & , & M_a > 1 \\
\frac{P(1 + M_a)}{2} & , & |M_a| \leq 1 \\
0 & , & M_a < -1
\end{array}
\right.
\\
&P^{-} = \left\{
\begin{array}{ccc}
0 & , & M_a > 1 \\
\frac{P(1 - M_a)}{2} & , & |M_a| \leq 1 \\
P & , & M_a < -1
\end{array}
\right.
\end{align}
\end{subequations}
Thus, there is
\begin{align}
\mathbf{F}^{P}_{\pm} = P^{\pm} \left[ \begin{array}{l}
0 \\
1 \\
0
\end{array} \right].
\end{align}
\par
In summary,
\begin{align}
\mathbf{F} = \mathbf{F}^{+} + \mathbf{F}^{-} \quad , \quad
\mathbf{F}^{\pm} = \mathbf{F}_{\pm}^a + \mathbf{F}_{\pm}^P.
\end{align}
\begin{remark}
  The van Leer splitting method and the Liou-Stenffen splitting method (a class of AUSM methods) share a consistent approach to the splitting of the Mach number $M_a$. 
\end{remark}

\section{DG Based on Flux Vector Splitting in One Dimension (1D-FVS-DG)}\label{sec-1D-FVSDG}
This section will present the flux vector splitting scheme for one-dimensional hyperbolic conservation systems within the DG framework.
Flux vector splitting itself does not introduce any error; therefore, one can first perform a functional split on the flux term in the original control equation and then define the weak solution through variation.
Based on the physical interpretation of the split sub-fluxes and considering the characteristic line theory and the stability of the discrete scheme, the corresponding upwind-type numerical flux format can be derived. 
\par
Consider the one-dimensional hyperbolic conservation system
\begin{align}\label{5.1}
\partial_t(\mathbf{U}) + \partial_x(\mathbf{F}(\mathbf{U}))=0.
\end{align}
\par
We perform the following splitting on the flux vector $\mathbf{F}$:
\begin{align}
\mathbf{F}(\mathbf{U}) = \mathbf{F}^{+}(\mathbf{U}) + \mathbf{F}^{-}(\mathbf{U}),
\end{align}
where $\mathbf{F}^+$ denotes the positive flux, which is the flux carried by characteristic waves propagating in the positive direction of the one-dimensional coordinate system; 
$\mathbf{F}^-$ denotes the negative flux, which is the flux carried by characteristic waves propagating in the negative direction of the one-dimensional coordinate system. 
The construction methods for $\mathbf{F}^+$ and $\mathbf{F}^-$ are discussed in Section \ref{sec-FVS}. 
\par
Thus, we have
\begin{align}
\partial_x(\mathbf{F}(\mathbf{U}))=\partial_x(\mathbf{F}^{+}(\mathbf{U})) + \partial_x(\mathbf{F}^{-}(\mathbf{U})),
\end{align}
\begin{align}\label{5.4}
\partial_t(\mathbf{U}) + \partial_x(\mathbf{F}^{+}(\mathbf{U})) + \partial_x(\mathbf{F}^{-}(\mathbf{U}))=0.
\end{align}
\par
Multiplying both sides of equation (\refeq{5.4}) by an arbitrary test function $\mathbf{V} \in \mathbb{V}^K_h$ and integrating by parts over the cell $I_i$ yields
\begin{align}
& \left<\partial_t(\mathbf{U}) , \mathbf{V}\right>_{L^2(I_i)} + \left<\partial_x(\mathbf{F}^{+}(\mathbf{U})) , \mathbf{V}\right>_{L^2(I_i)} + \left<\partial_x(\mathbf{F}^{-}(\mathbf{U})) , \mathbf{V}\right>_{L^2(I_i)} \nonumber \\
& = \left<\partial_t(\mathbf{U}) , \mathbf{V}\right>_{L^2(I_i)} - \left<\mathbf{F}^{+}(\mathbf{U}) , \partial_x(\mathbf{V})\right>_{L^2(I_i)} - \left<\mathbf{F}^{-}(\mathbf{U}) , \partial_x(\mathbf{V})\right>_{L^2(I_i)} 
+\left(\mathbf{F}^{+}\odot \mathbf{V}\right)^{{i+1/2}}_{{i-1/2}}
+\left(\mathbf{F}^{-}\odot \mathbf{V}\right)^{{i+1/2}}_{{i-1/2}} \nonumber \\
& = \left<\partial_t(\mathbf{U}) , \mathbf{V}\right>_{L^2(I_i)} - \left<\mathbf{F}^{+}(\mathbf{U})+\mathbf{F}^{-}(\mathbf{U}) , \partial_x(\mathbf{V})\right>_{L^2(I_i)}
+\left(\mathbf{F}^{+}\odot \mathbf{V}\right)^{{i+1/2}}_{{i-1/2}}
+\left(\mathbf{F}^{-}\odot \mathbf{V}\right)^{{i+1/2}}_{{i-1/2}} \nonumber \\
& = \left<\partial_t(\mathbf{U}) , \mathbf{V}\right>_{L^2(I_i)} - \left<\mathbf{F}(\mathbf{U}) , \partial_x(\mathbf{V})\right>_{L^2(I_i)}
+\left(\mathbf{F}^{+}\odot \mathbf{V}\right)^{{i+1/2}}_{{i-1/2}}
+\left(\mathbf{F}^{-}\odot \mathbf{V}\right)^{{i+1/2}}_{{i-1/2}}= 0. 
\end{align}
\hspace*{\fill} \\
Thus, the definition of the DG weak solution $\mathbf{U}_h$ based on flux vector splitting can be given: 
\begin{definition}[Weak Solution $\mathbf{U}_h$ for 1D-FVS-DG]
  If $\mathbf{U}_h \in \mathbb{V}^K_h$ and for all $\mathbf{V}_h \in \mathbb{V}_h^K$, it satisfies
\begin{align}
  &\left<\partial_t(\mathbf{U}_h) , \mathbf{V}_h\right>_{L^2(I_i)} - \left<\mathbf{F}(\mathbf{U}_h) , \partial_x(\mathbf{V}_h)\right>_{L^2(I_i)} \nonumber \\
  &+\left[\left(\widehat{\mathbf{F}^{+}}_{i+1/2} \odot \left(\mathbf{V}_h^L\right)_{i+1/2}\right)
  -\left(\widehat{\mathbf{F}^{+}}_{i-1/2} \odot \left(\mathbf{V}_h^R\right)_{i-1/2}\right)\right] \nonumber \\
  & +\left[\left(\widehat{\mathbf{F}^{-}}_{i+1/2} \odot \left(\mathbf{V}_h^L\right)_{i+1/2}\right)
  -\left(\widehat{\mathbf{F}^{-}}_{i-1/2} \odot \left(\mathbf{V}_h^R\right)_{i-1/2}\right)\right]
  = 0,
\end{align}
where $\widehat{\mathbf{F}^{\pm}}$ are the numerical fluxes, then $\mathbf{U}_h$ is referred to as the weak solution 
of the original equation \eqref{5.1} in the sense of FVS-DG. 
\end{definition}
Let
\begin{align}
  \widehat{\mathbf{F}}_{i+1/2} := \widehat{\mathbf{F}^{+}}_{i+1/2}+\widehat{\mathbf{F}^{-}}_{i+1/2},
\end{align}
then we have
\begin{align}\label{5.8}
  &\left<\partial_t(\mathbf{U}_h) , \mathbf{V}_h\right>_{L^2(I_i)} - \left<\mathbf{F}(\mathbf{U}_h) , \partial_x(\mathbf{V}_h)\right>_{L^2(I_i)}
  +\left(\widehat{\mathbf{F}}_{i+1/2}\odot\left(\mathbf{V}_h^L\right)_{i+1/2}\right)
  -\left(\widehat{\mathbf{F}}_{i-1/2}\odot\left(\mathbf{V}_h^R\right)_{i-1/2}\right)
  = 0.
\end{align}
In actual calculations, after obtaining $\widehat{\mathbf{F}^{+}}$ and $\widehat{\mathbf{F}^{-}}$ at the same interface, 
the two are superimposed scalarly to obtain $\widehat{\mathbf{F}}$, which is then substituted into equation (\refeq{5.8}) for computation. 
\par
Based on the characteristic line theory of hyperbolic equations and the stability of the discrete scheme, an upwind numerical flux can be employed, which takes the interface state on the reverse side of the characteristic propagation direction as the numerical state.
This numerical state can be split into two numerical sub-fluxes with opposite directions, and the sub-flux that is consistent with the characteristic propagation direction under consideration is taken as the final numerical flux for this interface in that characteristic propagation direction.
From the definitions of $\mathbf{F}^+$ and $\mathbf{F}^-$, the following can be immediately obtained: 
\begin{subequations}
\begin{align}
& \widehat{\mathbf{F}^+}_{i+1/2}=\mathbf{F}^{L, +}_{i+1/2}, \quad \widehat{\mathbf{F}^+}_{i-1/2}=\mathbf{F}^{L, +}_{i-1/2},\\
& \widehat{\mathbf{F}^-}_{i+1/2}=\mathbf{F}^{R, -}_{i+1/2}, \quad \widehat{\mathbf{F}^-}_{i-1/2}=\mathbf{F}^{R, -}_{i-1/2}.
\end{align}
\end{subequations}
\begin{remark}[DG Based on Jacobian Matrix Eigenvalues Splitting]
  \par
  In the Steger-Warming and Lax-Friedrichs splitting schemes, we have $\mathbf{F}^{+} = A^{+}\mathbf{U}$ and $\mathbf{F}^{-} = A^{-}\mathbf{U}$, so we obtain
\begin{subequations}
\begin{align}
    & \widehat{\mathbf{F}^+}_{i+1/2} = \mathbf{A}^{L, +}(\mathbf{U}^{L}_{i+1/2})\mathbf{U}^L_{i+1/2}, \quad \widehat{\mathbf{F}^+}_{i-1/2} = \mathbf{A}^{L, +}(\mathbf{U}^{L}_{i-1/2})\mathbf{U}^L_{i-1/2},\\
    & \widehat{\mathbf{F}^-}_{i+1/2} = \mathbf{A}^{R, -}(\mathbf{U}^{R}_{i+1/2})\mathbf{U}^R_{i+1/2}, \quad \widehat{\mathbf{F}^-}_{i-1/2} = \mathbf{A}^{R, -}(\mathbf{U}^{R}_{i-1/2})\mathbf{U}^R_{i-1/2}.
\end{align}
\end{subequations}
\end{remark}

\section{DG Based on Flux Vector Splitting in Two-dimension (2D-FVS-DG)}\label{sec-2D-FVSDG}
Consider the two-dimensional hyperbolic conservation system
\begin{align}\label{6.1}
\partial_t(\mathbf{U}) + \partial_x(\mathbf{F}(\mathbf{U})) + \partial_y(\mathbf{G}(\mathbf{U}))=0.
\end{align}
Let
\begin{align}
\mathbb{H}=(\mathbf{F}, \mathbf{G}),
\end{align}
the original equation can be written as
\begin{align}
\partial_t(\mathbf{U}) + \nabla \cdot \mathbb{H}=0.
\end{align}
Now, introduce the normal flux of the two-dimensional hyperbolic conservation system as
\begin{align}
\boldsymbol{\mathcal{F}}_{\nu} :=\mathbf{n}\cdot\mathbb{H}
\end{align}
that is
\begin{align}
\boldsymbol{\mathcal{F}}_{\nu} = n_x\cdot\mathbf{F} + n_y\cdot\mathbf{G}.
\end{align}
\begin{remark}
  For high-dimensional problems, the DG scheme cannot discretize space dimension by dimension like FD or FV schemes, 
  especially on triangular meshes, where the concept of ``dimension-by-dimension discretization'' is invalid.
  Therefore, we introduce the outward normal flux on the cell interfaces, 
  transforming the multi-dimensional problem into several one-dimensional problems along the outward normal directions of the mesh interfaces 
  (the number of normal one-dimensional problems depends on the number of numerical integration points arranged on the interfaces). 
\end{remark}

\subsection{DG Based on Flux Vector Splitting in Two-dimension}\label{subSec-2D-DG-FVS}
In the two-dimensional hyperbolic conservation system, there are flux functions $\mathbf{F}$ and $\mathbf{G}$ in the $x$ and 
$y$ directions, respectively. By using the outward normal vector on the cell interfaces, 
$\mathbf{F}$ and $\mathbf{G}$ are combined into a single flux function, namely the normal flux $\boldsymbol{\mathcal{F}_{\nu}}$. 
On the basis of the normal flux function $\boldsymbol{\mathcal{F}_{\nu}}$, the Jacobian matrix eigenvalue splitting or Mach number splitting is achieved. 
\par
Consider the two-dimensional hyperbolic conservation system
$$
\partial_t(\mathbf{U}) + \partial_x(\mathbf{F}(\mathbf{U})) + \partial_y(\mathbf{G}(\mathbf{U})) = 0,
$$
following the one-dimensional case, split the flux vectors in each dimension ($\mathbf{F}$ is the flux vector in the $x$ direction, $\mathbf{G}$ is the flux vector in the $y$ direction):
\begin{align}
\mathbf{F} = \mathbf{F}^{+} + \mathbf{F}^{-} \quad, \quad \mathbf{G} = \mathbf{G}^{+} + \mathbf{G}^{-}.
\end{align}
$\mathbf{F}^+$ represents the flux component in the $x$ dimension along the $+x$ direction; $\mathbf{F}^-$ represents the flux component in the $x$ dimension along the $-x$ direction; 
$\mathbf{G}^+$ represents the flux component in the $y$ dimension along the $+y$ direction; $\mathbf{G}^-$ represents the flux component in the $y$ dimension along the $-y$ direction. 
The construction methods for $\mathbf{F}^{\pm}, \mathbf{G}^{\pm}$ are discussed in Section \ref{sec-FVS}. 
\par
Let
\begin{align}
\mathbb{H}^{+}:=(\mathbf{F}^{+} ,  \mathbf{G}^{+} )^{\mathrm{T}}
\quad, \quad
\mathbb{H}^{-}:=(\mathbf{F}^{-} ,  \mathbf{G}^{-})^{\mathrm{T}}, 
\end{align}
hence
\begin{align}
\mathbb{H}=\mathbb{H}^{+} + \mathbb{H}^{-},
\end{align}
then
\begin{align}\label{6.9}
\partial_t(\mathbf{U})+\nabla \cdot \mathbb{H}^{+} + \nabla \cdot \mathbb{H}^{-} = 0.
\end{align}
\par
Multiply both sides of equation \eqref{6.9} by an arbitrary test function $\mathbf{V} \in \mathbb{V}^{K}_{h}$, 
and integrate by parts over $\Omega_i$ to obtain
{\footnotesize
\begin{align}
& \left< \partial_t(\mathbf{U}), \mathbf{V}\right>_{L^2(\Omega_i)} + \left<\nabla \cdot \mathbb{H}^{+}\left(\mathbf{U}\right), \mathbf{V}\right>_{L^2(\Omega_i)} + \left<\nabla \cdot \mathbb{H}^{-}\left(\mathbf{U}\right), \mathbf{V}\right>_{L^2(\Omega_i)} \nonumber \\
& = \left< \partial_t(\mathbf{U}), \mathbf{V}\right>_{L^2(\Omega_i)} 
- \left<\mathbb{H}^{+}\left(\mathbf{U}\right):\nabla\mathbf{V}\right>_{L^2(\Omega_i)}
- \left<\mathbb{H}^{-}\left(\mathbf{U}\right):\nabla\mathbf{V}\right>_{L^2(\Omega_i)} 
+ \int_{\Omega_i}\nabla \cdot (\mathbb{H}^{+}\left(\mathbf{U}\right) \odot \mathbf{V})d\mathbf{X}
+ \int_{\Omega_i}\nabla \cdot (\mathbb{H}^{-}\left(\mathbf{U}\right) \odot \mathbf{V})d\mathbf{X} \nonumber \\
& = \left< \partial_t(\mathbf{U}), \mathbf{V}\right>_{L^2(\Omega_i)} 
- \left<\mathbb{H}^{+}\left(\mathbf{U}\right):\nabla\mathbf{V}\right>_{L^2(\Omega_i)}
- \left<\mathbb{H}^{-}\left(\mathbf{U}\right):\nabla\mathbf{V}\right>_{L^2(\Omega_i)}
+ \int_{\partial\Omega_i} (\mathbb{H}^{+}\left(\mathbf{U}\right) \odot \mathbf{V}) \cdot \mathbf{n} dl
+ \int_{\partial\Omega_i} (\mathbb{H}^{-}\left(\mathbf{U}\right) \odot \mathbf{V}) \cdot \mathbf{n} dl \nonumber \\
& = \left< \partial_t(\mathbf{U}), \mathbf{V}\right>_{L^2(\Omega_i)} 
- \left<\mathbb{H}^{+}\left(\mathbf{U}\right):\nabla\mathbf{V}\right>_{L^2(\Omega_i)}
- \left<\mathbb{H}^{-}\left(\mathbf{U}\right):\nabla\mathbf{V}\right>_{L^2(\Omega_i)}
+ \int_{\partial\Omega_i}(\mathbb{H}^{+}\left(\mathbf{U}\right) \cdot \mathbf{n}) \odot \mathbf{V}dl
+ \int_{\partial\Omega_i}(\mathbb{H}^{-}\left(\mathbf{U}\right) \cdot \mathbf{n}) \odot \mathbf{V}dl \nonumber \\
& = \left< \partial_t(\mathbf{U}), \mathbf{V}\right>_{L^2(\Omega_i)} 
- \left<\mathbb{H}^{+}+\mathbb{H}^{-}\left(\mathbf{U}\right):\nabla\mathbf{V}\right>_{L^2(\Omega_i)}
+ \int_{\partial\Omega_i}(\mathbb{H}^{+}\left(\mathbf{U}\right) \cdot \mathbf{n}) \odot \mathbf{V}dl
+ \int_{\partial\Omega_i}(\mathbb{H}^{-}\left(\mathbf{U}\right) \cdot \mathbf{n}) \odot \mathbf{V}dl \nonumber \\
& = \left< \partial_t(\mathbf{U}), \mathbf{V}\right>_{L^2(\Omega_i)} 
- \left<\mathbb{H}\left(\mathbf{U}\right):\nabla\mathbf{V}\right>_{L^2(\Omega_i)}
+ \int_{\partial\Omega_i}(\mathbb{H}^{+}\left(\mathbf{U}\right) \cdot \mathbf{n}) \odot \mathbf{V}dl
+ \int_{\partial\Omega_i}(\mathbb{H}^{-}\left(\mathbf{U}\right) \cdot \mathbf{n}) \odot \mathbf{V}dl
=0.
\end{align}
}
\par
Let
\begin{align}
\boldsymbol{\mathcal{F}}_{\nu}^{+} := \mathbb{H}^{+} \cdot \mathbf{n} \quad, \quad \boldsymbol{\mathcal{F}}_{\nu}^{-} := \mathbb{H}^{-} \cdot \mathbf{n}.
\end{align}
{\color{red}{Note that the outward normal vector has been incorporated into the numerical fluxes.}}
\\
Obviously, 
\begin{align}
\boldsymbol{\mathcal{F}}_{\nu}^{+} = n_x \mathbf{F}^+ + n_y \mathbf{G}^+ \quad, \quad \boldsymbol{\mathcal{F}}_{\nu}^{-} = n_x \mathbf{F}^- + n_y \mathbf{G}^-,
\end{align}
and it holds that
\begin{align}
\boldsymbol{\mathcal{F}_{\nu}}=\boldsymbol{\mathcal{F}}_{\nu}^{+} + \boldsymbol{\mathcal{F}}_{\nu}^{-}.
\end{align}
\begin{remark}
  Given the physical significance of ${\mathbf{F}}^{\pm}, {\mathbf{G}}^{\pm}$, 
  it is understood that $\boldsymbol{\mathcal{F}}_{\nu}^{+}$ represents the flux flowing from the cell interior to the cell exterior 
  in the direction of the unit outward normal vector $\mathbf{n}$. 
  Conversely, $\boldsymbol{\mathcal{F}}_{\nu}^{-}$ represents the flux flowing from the cell exterior to the cell interior 
  in the direction of the unit outward normal vector $\mathbf{n}$.  
\end{remark}
Now, the definition of the DG weak solution $\mathbf{U}_h$ based on flux vector splitting is given:  
\begin{definition}[The weak solution of the 2D-FVS-DG]
  If $\mathbf{U}_h \in \mathbb{V}^K_h$ and for all $\mathbf{V}_h \in \mathbb{V}^K_h$, it satisfies
  \begin{equation}
    \begin{aligned}
    & \left< \partial_t(\mathbf{U}_h), \mathbf{V}\right>_{L^2(\Omega_i)} 
    - \left<\mathbb{H}\left(\mathbf{U}_h\right):\nabla\mathbf{V}_h\right>_{L^2(\Omega_i)}
    + \int_{\partial\Omega_i}\widehat{\boldsymbol{\mathcal{F}}_{\nu}^{+}} \cdot \mathbf{V}_h dl
    + \int_{\partial\Omega_i}\widehat{\boldsymbol{\mathcal{F}}_{\nu}^{-}} \cdot \mathbf{V}_h dl=0,
    \end{aligned}
  \end{equation}
  where $\widehat{\boldsymbol{\mathcal{F}}^{\pm}_{\nu}}$ are the normal numerical fluxes, 
  then $\mathbf{U}_h$ is referred to as the weak solution of the original equation \eqref{6.1} in the sense of FVS-DG. 
\end{definition}
Let
\begin{align}
\widehat{\boldsymbol{\mathcal{F}}_{\nu}}:=\widehat{\boldsymbol{\mathcal{F}}_{\nu}^{+}}+\widehat{\boldsymbol{\mathcal{F}}_{\nu}^{-}},
\end{align}
\\
In actual calculations, the following form is used:
\begin{equation}
  \begin{aligned}
  & \left< \partial_t(\mathbf{U}_h), \mathbf{V}\right>_{L^2(\Omega_i)} 
  - \left<\mathbb{H}\left(\mathbf{U}_h\right):\nabla\mathbf{V}_h\right>_{L^2(\Omega_i)}
  + \int_{\partial\Omega_i}\widehat{\boldsymbol{\mathcal{F}}_{\nu}} \cdot \mathbf{V}_h \, \mathrm{dl}
  = 0.
  \end{aligned}
\end{equation}
Similar to the one-dimensional problem, and still due to considerations of characteristic lines and stability, based on the physical meaning of $\boldsymbol{\mathcal{F}}_{\nu}^{\pm}$, it can be obtained that
\begin{subequations}
\begin{align}
\widehat{\boldsymbol{\mathcal{F}}_{\nu}^+} = n_x \mathbf{F}^{+}(\mathbf{U}^{int}) + n_y \mathbf{G}^{+}(\mathbf{U}^{int}),
\\
\widehat{\boldsymbol{\mathcal{F}}_{\nu}^-} = n_x \mathbf{F}^{-}(\mathbf{U}^{ext}) + n_y \mathbf{G}^{-}(\mathbf{U}^{ext}).
\end{align}
\end{subequations}
Here, ``int'' denotes the interior of the cell, and ``ext'' denotes the exterior of the cell.

\subsection{Normal Roe Average Employed by DG in Two-dimension}\label{subSec-2D-DG-Roe}
Roe average can be adopted for local-freezing when local characteristic decomposition discussed in Section \ref{sec-characteristic-reconstruction}. 
Similar to previous approach in subSection \ref{subSec-2D-DG-FVS}, {\color{blue}{Roe average along the outward normal direction}} on the cell interfaces is utilized to handle {\color{blue}{multi-dimensional cases}}.  
Here, we provides the specific implementation details of {\color{blue}{the noraml Roe averaging}}.  
\par
Firstly, review that the Roe average along the coordinate vectors: 
$$
\begin{aligned}
& \sqrt{\overline{\rho}^{Roe}}=\left(\sqrt{\rho_{int}}+\sqrt{\rho_{ext}}\right) / 2, \\
& \sqrt{\overline{\rho}^{Roe}} \cdot \overline{u}^{Roe}=\left(\sqrt{\rho_{int}} \cdot u_{int}+\sqrt{\rho_{ext}} \cdot u_{ext}\right) / 2, \\
& \sqrt{\overline{\rho}^{Roe}} \cdot \overline{v}^{Roe}=\left(\sqrt{\rho_{int}} \cdot v_{int}+\sqrt{\rho_{ext}} \cdot v_{ext}\right) / 2, \\
& \sqrt{\overline{\rho}^{Roe}} \cdot \overline{\widetilde{H}}^{Roe}=\left(\sqrt{\rho_{int}} \cdot \widetilde{H}_{int}+\sqrt{\rho_{ext}} \cdot \widetilde{H}_{ext}\right) / 2.
\end{aligned}
$$
Note: $ \widetilde{H} $ is the specific enthalpy.  
$$
\begin{aligned}
\overline{P}^{Roe} &= \frac{\gamma-1}{\gamma}\left(\overline{\rho}^{Roe} \overline{\widetilde{H}}^{Roe}-\frac{1}{2} (\overline{\rho}^{Roe})^2\right),\\
(\overline{a}^{Roe})^2 &= (\gamma-1)\left(\overline{\widetilde{H}}^{Roe}-\frac{1}{2}(\overline{u}^{Roe})^2 -\frac{1}{2}(\overline{v}^{Roe})^2\right),\\
\overline{E}^{Roe} &= \overline{\rho}^{Roe} \overline{\widetilde{H}}^{Roe}-\overline{P}^{Roe}.
\end{aligned}
$$
Note that 
$ \overline{P}^{Roe} $, $ \overline{a}^{Roe} $, and $ \overline{E}^{Roe} $ are not obtained by directly averaging $ P $, $ a $, and $ E $ with density but are computed indirectly.
\par
Then, the normal Jacobian and its eigenstructure are substituted as follows:  
\begin{equation}
\left.
\begin{aligned}
&\overline{q_n}^{Roe}=\overline{u}^{Roe} n_x + \overline{v}^{Roe} n_y, \\
&\overline{\lambda_1}^{Roe}=\overline{\lambda_2}^{Roe}=\overline{q_n}^{Roe}, \quad \overline{\lambda_3}^{Roe}=\overline{q_n}^{Roe}-\overline{a}^{Roe}, \quad \overline{\lambda_4}^{Roe}=\overline{q_n}^{Roe}+\overline{a}^{Roe},\\
&\overline{A_n}^{Roe}=A_n\left(\overline{\mathbf{U}}^{Roe}\right)=R_n\left(\overline{\mathbf{U}}^{Roe}\right) \cdot \Lambda_n\left(\overline{\mathbf{U}}^{Roe}\right) \cdot L_n\left(\overline{\mathbf{U}}^{Roe}\right)
=\overline{R_n}^{Roe} \cdot \overline{\Lambda_n}^{Roe} \cdot \overline{L_n}^{Roe}, \\
&\left|\overline{A_n}^{Roe}\right|=\overline{R_n}^{Roe} \cdot \left|\overline{\Lambda_n}^{Roe}\right| \cdot \overline{L_n}^{Roe},  \quad \left|\overline{\Lambda_n}^{Roe}\right|=\operatorname{diag}\left(\left|\overline{\lambda_1}^{Roe}\right|, \left|\overline{\lambda_2}^{Roe}\right|, \left|\overline{\lambda_3}^{Roe}\right|, \left|\overline{\lambda_4}^{Roe}\right|\right). \\
\end{aligned}
\right\}
\end{equation}
\begin{remark}
The normal Roe flux for the 2D-Euler equations is given by:  
\begin{align}
\hat{\mathbf{F}}_{n}^{R o e}=\frac{1}{2}\left(\mathbf{F}_{n}^{ext}+\mathbf{F}_{n}^{int}\right)-\frac{1}{2}\left|\overline{A_n}\right|\left(\mathbf{U}^{ext}-\mathbf{U}^{int}\right),
\end{align}
where    
\begin{subequations}
\begin{align}
\mathbf{F}_{n}^{ext}=\mathbf{F}_{n}\left(\mathbf{U}^{ext}\right)=\mathbf{F}\left(\mathbf{U}^{ext}\right) \cdot n_x + \mathbf{G}\left(\mathbf{U}^{ext}\right) \cdot n_y, \\
\mathbf{F}_{n}^{int}=\mathbf{F}_{n}\left(\mathbf{U}^{int}\right)=\mathbf{F}\left(\mathbf{U}^{int}\right) \cdot n_x + \mathbf{G}\left(\mathbf{U}^{int}\right) \cdot n_y.
\end{align}
\end{subequations}
Roe flux is a kind of Riemann solvers, which is diffrenet from our FVS-based flux schemes introduced in Section \ref{sec-1D-FVSDG} and subSection \ref{subSec-2D-DG-FVS}.
\end{remark}

\section{A Novel TVB(D)-minmod Limiter for Numerical Pseudo-Oscillation Treatment}\label{sec-Constrained TVB(D)-minmod Limiter}
This section will introduce the original TVB(D)-minmod limiter and the discontinuity indicators firstly 
and then focuses on some improvements we have made to the TVB(D)-minmod limiter.

\subsection{Original TVB(D)-minmod Limiter and Commonly Used Discontinuity Indicators}
In this subsection, we introduce two types of discontinuity indicators used in the subsequent numerical experiments, 
namely the TVB(D)-minmod discontinuity indicator and the KXRCF discontinuity indicator.
The TVB(D)-minmod discontinuity indicator is applicable to 1D/2D scalar equations as well as 1D/2D systems. 
{\color{red}{Note that when we utilize the data provided by the TVB(D)-minmod discontinuity indicator to 
further refine the integral mean on the boundary, the discontinuity indicator evolves into 
a TVB(D)-minmod limiter. As mentioned in the ``1.Introduction'', the TVB(D)-minmd limiter is equipped with 
an inherent discontinuity indication function.}} 
We employ this type of limiter solely on rectangular (arbitrary quadrilateral) meshes;
KXRCF discontinuity indicator is used exclusively on triangular meshes for the control equations of fluid systems (the KXRCF is not very convenient for scalar cases).
It is important to note that when discontinuity indicators are applied to system control equations, 
each component's troubled cells should be inspected individually.

\subsubsection{TVB(D)-minmod Discontinuity Indicator and Limiter}\label{subsubSec-classical-TVB(D)-minmod}
The TVB(D)-minmod limiter is a slope limiter that compares the jump between the average of the troubled cell and the averages of the neighboring cells with the jump between the average of the troubled cell and the average of the boundary integrals of the troubled cell. The smaller of these absolute values is taken as the revised average of the boundary integrals of the troubled cell, thereby reducing the average slope of the approximate solution curve within the troubled cell. 
When used as a discontinuity indicator, the working principle of the TVB(D)-minmod limiter is as follows:
\par
$\bullet$ 1D-TVB(D)-minmod Discontinuity Indicator and Limiter
\par
The deviation of the cell left boundary value from the cell average (left deviation) is given by
$$
\widehat{U}_i := \overline{U}_i - U_{i-\frac{1}{2}}^{+},
$$
hence,
$$
U_{i-\frac{1}{2}}^{+} = \overline{U}_i - \widehat{U}_i.
$$
The deviation of the cell right boundary value from the cell average (right deviation) is
$$
\widetilde{U}_i := U_{i+\frac{1}{2}}^{-} - \overline{U}_i,
$$
hence,
$$
U_{i+\frac{1}{2}}^{-} = \overline{U}_i + \widetilde{U}_i.
$$
The forward difference between neighboring cell averages is
$$
\Delta_{-} = \overline{U}_i - \overline{U}_{i-1}.
$$
The backward difference between neighboring cell averages is
$$
\Delta_{+} = \overline{U}_{i+1} - \overline{U}_i.
$$
The modified left deviation is calculated as
$$
\widehat{U}_i^{mod} = \operatorname{minmod}\left(\widehat{U}_i, \Delta_{-}, \Delta_{+}\right).
$$
The modified left boundary value is then
$$
U_{i-\frac{1}{2}}^{+, mod} = \overline{U}_i - \widehat{U}_i^{mod}.
$$
The modified right deviation is calculated as
$$
\widetilde{U}_i^{mod} = \operatorname{minmod}\left(\widetilde{U}_i, \Delta_{-}, \Delta_{+}\right).
$$
The modified right boundary value is then
$$
U_{i+\frac{1}{2}}^{-, mod} = \overline{U}_i + \widetilde{U}_i^{mod}.
$$
If the modification is not equal to the original, it indicates that the cell is a troubled cell, that is:
$$
\text{if} \quad
\widehat{U}_i^{\text{mod}} \neq \widehat{U}_i
\quad \text{or} \quad
\widetilde{U}_i^{\text{mod}} \neq \widetilde{U}_i
\quad \text{then} \quad
I_i \text{ is a troubled cell}.
$$
\begin{remark}[Correct the end-point values of the troubled interval and implement the TVB(D)-minmod limiter]
According to $U_{i-\frac{1}{2}}^{+, mod}$ and $U_{i+\frac{1}{2}}^{-, mod}$, we can correct the DG polynomial on the $I_i$ cell. 
{\color{red}{Executing this correction step can implement the TVB(D)-minmod limiter 
(if not, then it just work as a discontinuity indicator). }}
\end{remark}
$\bullet$ 2D-TVB(D)-minmod Discontinuity Indicator and Limiter
\\
The 2D rectangular cell integral average is given by
$$
\overline{U}_{ij} = \overline{U}_{\Omega_{ij}} = \frac{1}{|\Omega_{ij}|} \int_{\Omega_{ij}} u_{ij}(x, y) \, dxdy.
$$
The left boundary (L) integral average is
$$
\overline{U}^{L}_{ij} = \overline{U}^{L}_{\partial\Omega_{ij}} = \frac{1}{|\partial\Omega_{ij}^{L}|} \int_{\partial\Omega_{ij}^{L}} u_{ij}(x, y) \, dl.
$$
The right boundary (R) integral average is
$$
\overline{U}^{R}_{ij} = \overline{U}^{R}_{\partial\Omega_{ij}} = \frac{1}{|\partial\Omega_{ij}^{R}|} \int_{\partial\Omega_{ij}^{R}} u_{ij}(x, y) \, dl.
$$
The bottom boundary (B) integral average is
$$
\overline{U}^{B}_{ij} = \overline{U}^{B}_{\partial\Omega_{ij}} = \frac{1}{|\partial\Omega_{ij}^{B}|} \int_{\partial\Omega_{ij}^{B}} u_{ij}(x, y) \, dl.
$$
The top boundary (T) integral average is
$$
\overline{U}^{T}_{ij} = \overline{U}^{T}_{\partial\Omega_{ij}} = \frac{1}{|\partial\Omega_{ij}^{T}|} \int_{\partial\Omega_{ij}^{T}} u_{ij}(x, y) \, dl.
$$
The left difference between adjacent cell averages is
$$
\Delta_L = \overline{U}_{ij} - \overline{U}_{i-1, j}.
$$
The right difference between adjacent cell averages is
$$
\Delta_R = \overline{U}_{i+1, j} - \overline{U}_{i, j}.
$$
The bottom difference between adjacent cell averages is
$$
\Delta_B = \overline{U}_{ij} - \overline{U}_{i, j-1}.
$$
The top difference between adjacent cell averages is
$$
\Delta_T = \overline{U}_{i, j+1} - \overline{U}_{ij}.
$$
The left deviation of the boundary average from the cell average is
$$
\widehat{U}^{L}_{ij} = \overline{U}_{ij} - \overline{U}^{L}_{ij},
$$
hence,
$$
\overline{U}^{L}_{ij} = \overline{U}_{ij} - \widehat{U}^{L}_{ij}.
$$
The right deviation of the boundary average from the cell average is
$$
\widetilde{U}^{R}_{ij} = \overline{U}^{R}_{ij} - \overline{U}_{ij},
$$
hence,
$$
\overline{U}^{R}_{ij} = \overline{U}_{ij} + \widetilde{U}^{R}_{ij}.
$$
The bottom deviation of the boundary average from the cell average is
$$
\widehat{U}^{B}_{ij} = \overline{U}_{ij} - \overline{U}^{B}_{ij},
$$
hence,
$$
\overline{U}^{B}_{ij} = \overline{U}_{ij} - \widehat{U}^{B}_{ij}.
$$
The top deviation of the boundary average from the cell average is
$$
\widetilde{U}^{T}_{ij} = \overline{U}^{T}_{ij} - \overline{U}_{ij},
$$
hence,
$$
\overline{U}^{T}_{ij} = \overline{U}_{ij} + \widetilde{U}^{T}_{ij}.
$$
The modified left deviation is
$$
\widehat{U}^{L, mod}_{ij} = \operatorname{minmod}\left(\widehat{U}^{L}_{ij}, \Delta_L, \Delta_R\right).
$$
The modified right deviation is
$$
\widetilde{U}^{R, mod}_{ij} = \operatorname{minmod}\left(\widetilde{U}^{R}_{ij}, \Delta_L, \Delta_R\right).
$$
The modified bottom deviation is
$$
\widehat{U}^{B, mod}_{ij} = \operatorname{minmod}\left(\widehat{U}^{B}_{ij}, \Delta_B, \Delta_T\right).
$$
The modified top deviation is
$$
\widetilde{U}^{T, mod}_{ij} = \operatorname{minmod}\left(\widetilde{U}^{T}_{ij}, \Delta_B, \Delta_T\right).
$$
The 2D-rectangular mesh TVB(D)-minmod discontinuity indicator is as follows:
$$
\begin{aligned}
& \text{if} \quad
\widehat{U}^{L, mod}_{ij} \neq \widehat{U}^{L}_{ij}
\quad \text{or} \quad
\widetilde{U}^{R, mod}_{ij} \neq \widetilde{U}^{R}_{ij}
\quad \text{or} \quad
\widehat{U}^{B, mod}_{ij} \neq \widehat{U}^{B}_{ij}
\quad \text{or} \quad
\widetilde{U}^{T, mod}_{ij} \neq \widetilde{U}^{T}_{ij}, \\
& \text{then} \   
\Omega_{ij} \text{ is a troubled cell}.
\end{aligned}
$$	
\begin{remark}[Correction of Boundary Integral Averages using the 2D Rectangular Cell TVB(D)-minmod Limiter]
\begin{align*}
&\overline{U}^{L,mod}_{ij}=\overline{U}_{ij}+\widehat{U}^{L,mod}_{ij}, \\
&\overline{U}^{R,mod}_{ij}=\overline{U}_{ij}+\widetilde{U}^{R,mod}_{ij}, \\
&\overline{U}^{B,mod}_{ij}=\overline{U}_{ij}+\widehat{U}^{B,mod}_{ij}, \\
&\overline{U}^{T,mod}_{ij}=\overline{U}_{ij}+\widetilde{U}^{T,mod}_{ij}.
\end{align*}
{\color{red}{Executing this correction step can implement the TVB(D)-minmod limiter 
(if not, then it just work as a discontinuity indicator). }}
\end{remark}

\subsubsection{KXRCF Discontinuity Indicator}\label{subsubSec-KXRCF}
Let the inflow boundary of the solution be denoted as $\partial \Omega^{-}$. The \textit{KXRCF} indicator is defined as follows:
$$
J_{\Omega}=\frac{\left|\int_{\partial \Omega^{-}}\left(\left.u\right|_{\Omega}-\left.u\right|_{\Omega_{n b}}\right) \mathrm{d} S_x\right|}{h^{\frac{k+1}{2}}\left|\partial \Omega^{-}\right| \cdot\left\|\left.u\right|_{\Omega}\right\|_{L^{\infty}}}.
$$
If $J_{\Omega}>1$,  then $\Omega$ is identified as a troubled cell. 
Here, {\color{red}{$\Omega_{\text{nb}}$ refers to the neighboring cell sharing the inflow boundary $\partial \Omega^{-}$}} (not all neighboring cells), 
$K$ is the highest degree of the piecewise polynomial, 
and $h$ is the radius of the cell.  
\par
$\bullet\ $1D-KXRCF Discontinuity Indicator: For the 1D fluid systems, if $v \cdot \mathbf{n} < 0$, it is considered an inflow boundary; 
otherwise, it is an outflow boundary. Here $\mathbf{n}$ is the unit outward normal vector; 
For 1D scalar case, we define $v$, taking its value from inside the cell $I_j$ as $f^{\prime}(u)$ and take $u$ as the indicator variable.  
\par
$\bullet\ $2D-KXRCF Discontinuity Indicator: For the 2D fluid systems, if $(v_x,v_y) \cdot \mathbf{n} < 0$, it is classified as an inflow boundary; 
otherwise, it is an outflow boundary. Here $\mathbf{n}$ is the unit outward normal vector; 
For 2D scalar case, we define $v_1,\ v_2$, taking their values from inside the cell $\Omega_j$ as $f^{\prime}(u),\ g^{\prime}(u)$ respectively 
and still take $u$ as the indicator variable.

\subsection{Constrained Optimization-based TVB(D)-minmod Limiter Compatible with High-Order Polynomial Approximation}\label{subSec-Constrained Optimization-based TVB(D)-minmod}
\begin{figure}[htbp]
  \vspace{0.025cm}
  \begin{center}
    \begin{minipage}{0.3\linewidth}
      \centerline{\includegraphics[width=1.25\linewidth]{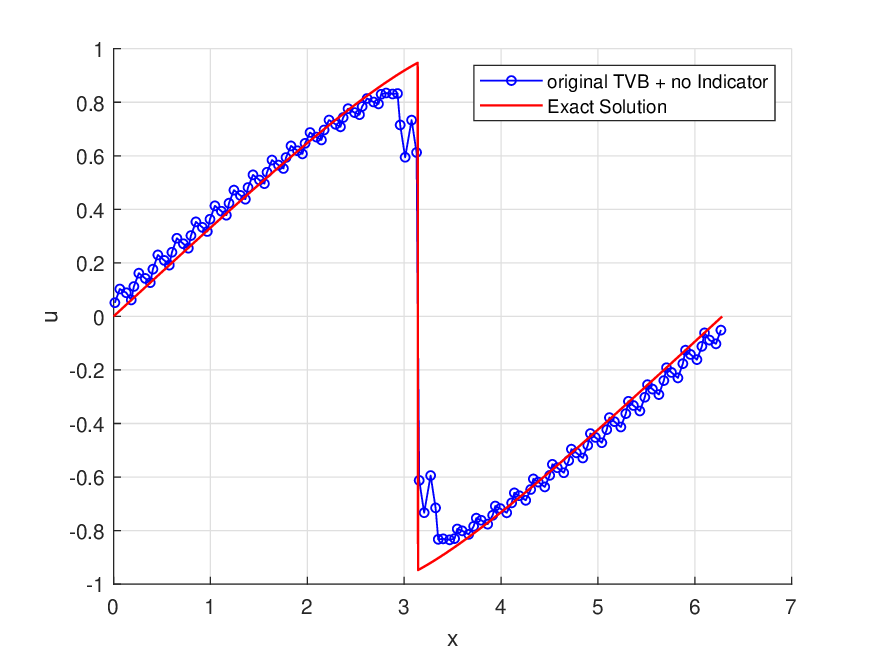}}
    \end{minipage}
    \hfill
    \begin{minipage}{0.3\linewidth}
      \centerline{\includegraphics[width=1.25\linewidth]{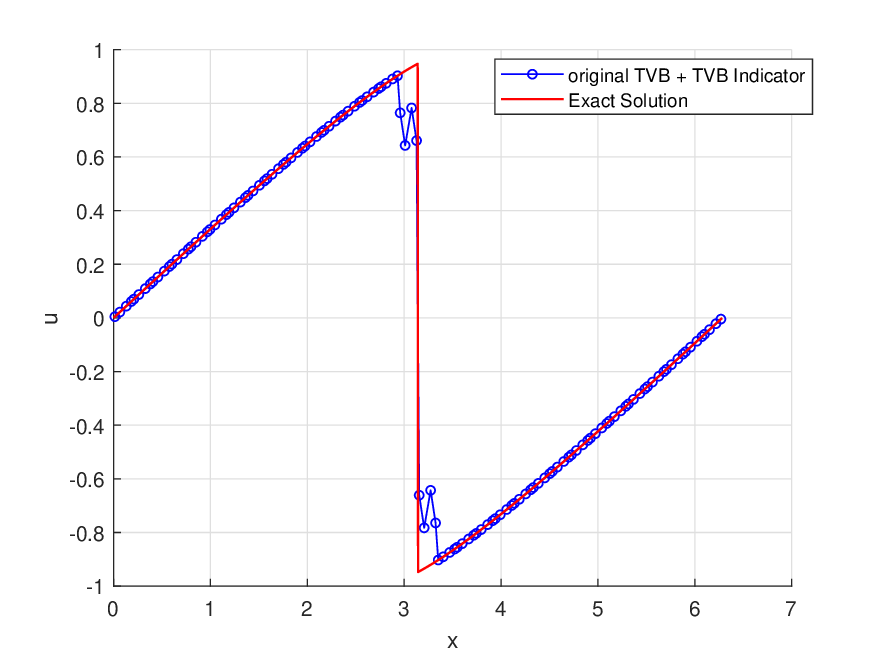}}
    \end{minipage}
    \hfill
    \begin{minipage}{0.3\linewidth}
      \centerline{\includegraphics[width=1.25\linewidth]{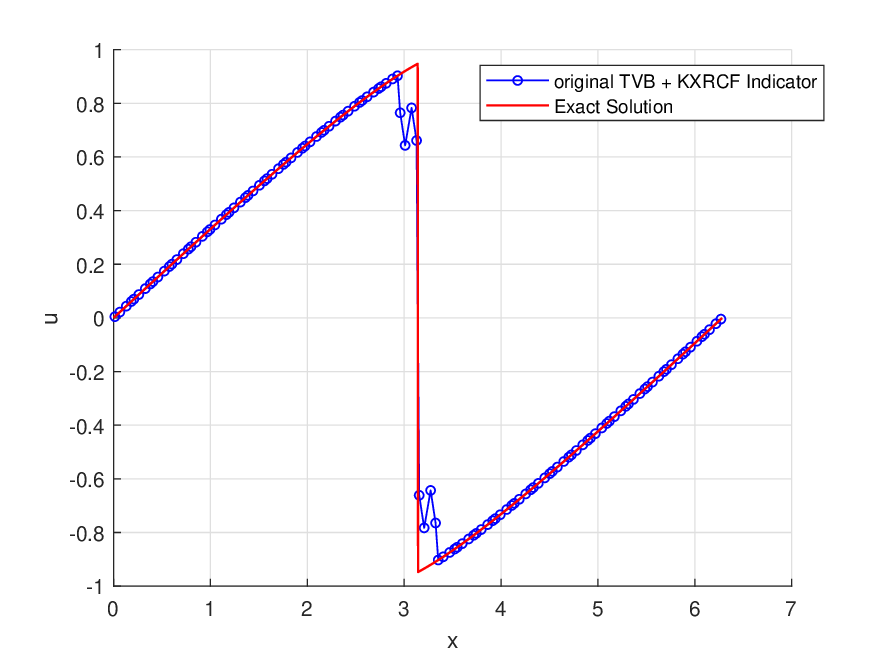}}
    \end{minipage}
    \vfill
    \begin{minipage}{0.3\linewidth}
      \small
      \centerline{\tiny{(a)}}
    \end{minipage}
    \hfill
    \begin{minipage}{0.3\linewidth}
      \small
      \centerline{\tiny{(b)}}
    \end{minipage}
    \hfill
    \begin{minipage}{0.3\linewidth}
      \small
      \centerline{\tiny{(c)}}
    \end{minipage}
  \end{center}
  \caption{\tiny{1D-Burgers' equation $u_t + (\frac{1}{2}u^2)_x = 0$ with initial condition $u_0(x)=\sin(x)$. The simulation is performed up to time $t=2.0$. $P^3$-polynomial approximations and uniquely spaced 32 cells. TVB parameter $M=1$. The $v$ for KXRCF indicator in scalar case is taken as $u_{i\pm1/2}$ from inside the cell $I_i$: 
  (a)\ classical TVB-minmod limiter without discontinuity indicator; (b)\ classical TVB-minmod limiter with TVB discontinuity indicator; 
  (c)\ classical TVB-minmod limiter with KXRCF discontinuity indicator}}
  \label{Fig.classical-TVB}
\end{figure}

In the context of high-order $P^K$-polynomial approximation ($K \geq 3$), the original TVB(D)-minmod limiter results in indeterminate correction equations with non-unique solutions, 
for example, in one-dimensional case, at most 2 correction equations ($U_{i-\frac{1}{2}}^{+, mod},U_{i+\frac{1}{2}}^{-, mod}$) provided by this kind of limiter but more than 3 modal coefficients in high-order polynomial 
while on 2D rectangular meshes, at least 9 modal coefficients ($P^3$) needed to be modified but only 4 conditions ($\widehat{U}^{L,mod}_{ij},\widehat{U}^{R,mod}_{ij},\widehat{U}^{B,mod}_{ij},\widehat{U}^{T,mod}_{ij}$). 
{\color{red}{Direct use of the least squares solution for the above indeterminate correction equations will result in noticeable numerical oscillations when using higher-order polynomial approximation 
(Please refer to Figure \ref{Fig.classical-TVB}). }}
The primary goal of employing the TVB(D)-minmod limiter is to suppress numerical oscillations. 
To better achieve this goal, we select the solution from the set of indeterminate solutions that is ``the smoothest,'' 
serving as the DG solution modified by the TVB(D)-minmod limiter.
Therefore, we introduce the smoothness measurement function $IS$ from WENO reconstruction as the objective function 
to construct an optimization problem. Additionally, while suppressing numerical oscillations, 
it is essential to maintain the solution's high accuracy. 
The approach taken in \cite{ref1} is to introduce an $L^2$-error term to characterize the difference 
between the corrected solution and the original high-order scheme. 
Consequently, a natural idea is to combine the $IS$ and $L^2$-error to form a bi-objective optimization problem 
that balances oscillation suppression and precision.  

\subsubsection{Smoothness Factor $IS$ Constrained TVB(D)-minmod Limiter: One-Dimensional Case}\label{sec-ISTVB-1D}
The one-dimensional smoothness measurement function is defined as:
\begin{align}
I S(Q_i^{K}(x))=\sum_{l=1}^{l=K} \int_{x_{i-1/2}}^{x_{i+1/2}} \Delta x_i^{2 l-1}\left(\frac{\partial^l Q_i^{K}(x)}{\partial x^l}\right)^2 d x.
\end{align}
Let the original DG weak solution on interval $I_i$ without the limiter be denoted as $u_i^{old}$, with modal coefficients given by 
$\mathbf{a}^{i,old}=\left[\alpha_0^{i,old}, \alpha_1^{i,old},  \alpha_2^{i,old},  \alpha_3^{i,old},  \ldots,  \alpha_{K-1}^{i,old},  \alpha_K^{i,old}\right]^{\mathrm{T}}$. 
An optimization problem for the DG weak solution $u_i$ on interval $I_i$ is constructed as follows: 
{\color{red}{
\begin{equation}\label{min(IS(u))1}
\begin{aligned}
& \min IS(u_i) \\
& s.t. \quad 
\left\{\begin{array}{c}
\frac{1}{\left|I_i\right|} \int_{I_i} u_i(x) d x=\overline{U}_i \\
u_i\left(x_{i-1/2}\right)=U_{i-1/2}^{ +, mod } \\
u_i\left(x_{i+1/2}\right)=U_{i+1/2}^{ -, mod }
\end{array}\right.
\end{aligned}
\end{equation}
}}
Let the solution to problem \eqref{min(IS(u))1} be $u_i^{mod}$, with modal coefficients
\begin{align}
  \mathbf{a}^{i,mod}=\left[\alpha_0^{i,mod}, \alpha_1^{i,mod},  \alpha_2^{i,mod},  \alpha_3^{i,mod},  \ldots,  \alpha_{K-1}^{i,mod},  \alpha_K^{i,mod}\right]^{\mathrm{T}}, 
\end{align}
such that
\begin{align}
  u_i^{mod}=\sum_{\ell=0}^{\ell=K}\alpha_{\ell}^{i,mod}\phi_{\ell}^{(I_i)}.
\end{align}
The integral mean remains unchanged, i.e., {\color{red}{$\alpha_0^{i,mod}=\alpha_0^{i,old}$}},  which allows further simplification to
{\color{blue}{
\begin{equation}\label{eq-min(IS(u))2}
\begin{aligned}
& \min IS(u_i)\\
& s.t. \quad 
\left\{\begin{array}{c}
u_i\left(x_{i-1/2}\right)=U_{i-1/2}^{ +, mod } \\
u_i\left(x_{i+1/2}\right)=U_{i+1/2}^{ -, mod }
\end{array}\right.
\end{aligned}
\end{equation}
}}
{\color{red}{For the sake of brevity in notation, we omit the cell indices 
for the modal coefficients and the basis functions. }}
It is known that
$$
\begin{aligned}
u_i^{\prime} & =\alpha_1 \phi_1^{\prime}(x)+\alpha_2 \phi_2^{\prime}(x)+\alpha_3 \phi_3^{\prime}(x)+\cdots+\alpha_K \phi_K^{\prime}(x), \\
u_i^{\prime \prime} & =\alpha_2 \phi_2^{\prime \prime}(x)+\alpha_3 \phi_3^{\prime \prime}(x)+\cdots+\alpha_K \phi_K^{\prime \prime}(x), \\
u_i^{(3)} & =\alpha_3 \phi_3^{(3)}(x)+\cdots+\alpha_K \phi_K^{(3)}(x), \\
\vdots & \\
u_i^{(K-1)} & =\alpha_{K-1} \phi_{K-1}^{(K-1)}(x)+\alpha_K \phi_K^{(K-1)}(x),\\
u_i^{(K)} & =\alpha_K \phi_K^{(K)}(x).
\end{aligned}
$$
Define
\begin{align}
M_{j k}^d := \Delta x^{2 d-1} \int_{I_i} \frac{\partial^d \phi_j}{\partial x^d} \cdot \frac{\partial^d \phi_k}{\partial x^d} d x=\Delta x^{2 d-1}\left\langle\phi_j^{(d)},  \phi_k^{(d)}\right\rangle_{L^2\left(I_i\right)},
\end{align}
then the smoothness factor is given by
\begin{align}
IS_i = \sum_{d=1}^{d=K} \int_{x_{i-1/2}}^{x_{i+1/2}} \Delta x_i^{2 d-1}\left(u_i^{(d)}\right)^2 d x 
= \sum_{d=1}^{d=K} \sum_{j=1}^{K} \sum_{k=1}^{K} M_{jk}^{d} \alpha_j \alpha_k 
= \sum_{d=1}^{d=K} \left(\sum_{j=1}^{K} M^d_{jj} \alpha_j^{{\color{blue}{2}}} +\sum_{j=1}^{K-1}\sum_{k=j+1}^{K} {\color{red}{2}} M_{jk}^{d} \alpha_j \alpha_k \right).
\end{align}
The corresponding relationship between coefficients is as follows: 
{\scriptsize 
\setlength{\arraycolsep}{3.5pt}
\begin{subequations}\label{eq-alpha-M}
\begin{align}
&\left[
\begin{array}{llllll}
\alpha_1 ^{{\color{blue}{2}}} & & & & & \\
\alpha_1 \alpha_2 & \alpha_2 ^{{\color{blue}{2}}} & & & & \\
\alpha_1 \alpha_3 & \alpha_2 \alpha_3 & \alpha_3 ^{{\color{blue}{2}}} & & & \\
\vdots & \vdots & \vdots & \ddots & & \\
\alpha_1 \alpha_{K-1} & \alpha_2 \alpha_{K-1} & \alpha_3 \alpha_{K-1} & \cdots & \alpha_{K-1} ^{{\color{blue}{2}}} & \\
\alpha_1 \alpha_K & \alpha_2 \alpha_K & \alpha_3 \alpha_K & \cdots & \alpha_{K-1} \alpha_K & \alpha_K ^{\color{blue}{2}}
\end{array}
\right],
\\
&\left[
\begin{array}{llllll}
M_{11}^1 & & &  \\
& & & & \\
{\color{red}{2}}M_{12}^1 & M_{22}^1+M_{22}^2 & & &  \\
& & & & \\
{\color{red}{2}}M_{13}^1 & {\color{red}{2}}\left(M_{23}^1+M_{23}^2\right) & M_{33}^1+M_{33}^2+M_{33}^3 & &  \\
& & & & \\
\vdots & \vdots & \vdots & \ddots &  \\
& & & & \\
{\color{red}{2}}M_{1 K-1}^1 & {\color{red}{2}}\left(M_{2 K-1}^1+M_{2 K-1}^2\right) & {\color{red}{2}}\left(M_{3 K-1}^1+M_{3 K-1}^2+M_{3 K-1}^3\right) & \cdots & M_{K-1 K-1}^1+M_{K-1 K-1}^2+\cdots+M_{K-1 K-1}^{K-1}  \\
& & & & \\
{\color{red}{2}}M_{1 K}^1 & {\color{red}{2}}\left(M_{2 K}^1+M_{2 K}^2\right) & {\color{red}{2}}\left(M_{3 K}^1+M_{3 K}^2+M_{3 K}^3\right) & \cdots & {\color{red}{2}}\left(M_{K-1 K}^1+M_{K-1 K}^2+\cdots+M_{K-1 K}^{K-1}\right) & \vartriangle
\end{array}
\right] \label{7.6b}, \\
&\vartriangle = M_{K K}^1+M_{K K}^2+\cdots+M_{K K}^K.
\end{align}
\end{subequations}
}
{\color{red}{From this point on, we restore the cell indices.}} Define
{\scriptsize
\setlength{\arraycolsep}{3.5pt}
\begin{align}
&\mathbb{M}_i:= \nonumber \\
&\left[\begin{array}{llllll}
{\color{blue}{2}}M_{11}^1 & & &  \\
& & & & \\
{\color{red}{2}}M_{12}^1 & {\color{blue}{2}}\left(M_{22}^1+M_{22}^2\right) & & &  \\
& & & & \\
{\color{red}{2}}M_{13}^1 & {\color{red}{2}}\left(M_{23}^1+M_{23}^2\right) & {\color{blue}{2}}\left(M_{33}^1+M_{33}^2+M_{33}^3\right) & &  \\
& & & & \\
\vdots & \vdots & \vdots & \ddots &  \\
& & & & \\
{\color{red}{2}}M_{1 K-1}^1 & {\color{red}{2}}\left(M_{2 K-1}^1+M_{2 K-1}^2\right) & {\color{red}{2}}\left(M_{3 K-1}^1+M_{3 K-1}^2+M_{3 K-1}^3\right) & \cdots & {\color{blue}{2}}\left(M_{K-1 K-1}^1+M_{K-1 K-1}^2+\cdots+M_{K-1 K-1}^{K-1}\right)  \\
& & & & \\
{\color{red}{2}}M_{1 K}^1 & {\color{red}{2}}\left(M_{2 K}^1+M_{2 K}^2\right) & {\color{red}{2}}\left(M_{3 K}^1+M_{3 K}^2+M_{3 K}^3\right) & \cdots & {\color{red}{2}}\left(M_{K-1 K}^1+M_{K-1 K}^2+\cdots+M_{K-1 K}^{K-1}\right) & {\color{blue}{2}}\vartriangle 
\end{array}
\right].
\end{align}
}
{\color{red}{Note that each $M_{jk}^{d}$ of $\mathbb{M}_i$ belongs to the cell $I_i$. }} 
The matrix $\mathbb{M}_i$ is decomposed as follows:
\begin{align}
\mathbb{M}_i=\mathbb{L}_i+\mathbb{D}_i, 
\quad\mathbb{D}_i=\operatorname{diag}(\mathbb{M}_i(1, 1), \mathbb{M}_i(2, 2), \mathbb{M}_i(3, 3), \cdots, \mathbb{M}_i(K, K)).
\end{align}
Let
\begin{align}
\mathbb{U}_i=\mathbb{L}_i^{\mathrm{T}}.
\end{align}
Define the matrix $\mathbb{M}_i^{IS}$ as
\begin{align}\label{MIS}
  \mathbb{M}_i^{IS}:=\mathbb{L}_i+\mathbb{D}_i+\mathbb{U}_i.
\end{align}
Additionally, define 
\begin{align}\label{713-715}
\tilde{\mathbf{a}}^{i,mod}&:=\left[{\color{red}{\alpha_1^{i,mod}}},  \alpha_2^{i,mod},  \alpha_3^{i,mod},  \cdots,  \alpha_{K-1}^{i,mod},  \alpha_K^{i,mod}\right]^{\mathrm{T}},\\
\left[\Phi\right]_i&:=
\left[
\begin{array}{llll}
  \phi_1\left(x_{i-1 / 2}\right) & \phi_2\left(x_{i-1 / 2}\right) & \cdots & \phi_K\left(x_{i-1 / 2}\right) \\
  \phi_1\left(x_{i+1 / 2}\right) & \phi_2\left(x_{i+1 / 2}\right) & \cdots & \phi_K\left(x_{i+1 / 2}\right)
\end{array}
\right]^{\mathrm{T}},
\\
\left[\Lambda^*\right]&:=\left[\lambda_1^*, \lambda^*_2\right]^{\mathrm{T}}.
\end{align}
To solve the constrained optimization problem (\refeq{eq-min(IS(u))2}), the Lagrangian is constructed as follows:  
\begin{align}
L_i(\alpha_1^i, \alpha_2^i, \ldots, \alpha_K^i, \lambda_1, \lambda_2) 
= IS_i
+\lambda_1 \cdot \left(u_i\left(x_{i-1/2}\right)-U_{i-1/2}^{ +, mod } \right)
+ \lambda_2 \cdot \left(u_i\left(x_{i+1/2}\right)-U_{i+1/2}^{ -, mod }\right).
\end{align}
The solution to problem ($\ref{eq-min(IS(u))2}$), denoted as $u_i^{mod}(x;\alpha_0^{i,old},\tilde{\mathbf{a}}^{i,mod})$, 
should satisfy the saddle-point equation, that is:
\renewcommand{\arraystretch}{1.5}
\begin{align}\label{saddle-point equation}
\nabla L_i=\left[\begin{array}{c}
\partial_{\alpha_1} L_i \\
\partial_{\alpha_2} L_i \\
\vdots \\
\partial_{\alpha_K} L_i \\
\hdashline
\partial_{\lambda_1} L_i \\
\partial_{\lambda_2} L_i
\end{array}\right]_{(\tilde{\mathbf{a}}^{mod}_i,\lambda_1^*,\lambda_2^*)}
=
\left[\begin{array}{c}
\mathbb{M}^{IS}_i\tilde{\mathbf{a}}^{mod}_i+\left[\Phi\right]_i\left[\Lambda^*\right]\\
\hdashline
\sum_{k=0}^K \alpha_k^{i,mod} \phi_k\left(x_{i-1/2}\right)-U_{i-1/2}^{+,  \bmod } \\
\sum_{k=0}^K \alpha_k^{i,mod} \phi_k\left(x_{i+1/2}\right)-U_{i+1/2}^{-,  \bmod }
\end{array}\right]
=0.
\end{align}
Taking into account that $\alpha_0^{i,mod} = \alpha_0^{i,old}$ is known, the matrix form of the saddle-point equation is given by:
\begin{align}
  \mathbb{A}^{IS}_i \mathbf{X}^{IS}_i = \mathbf{b}^{IS}_i,
\end{align}
where  
\begin{subequations}
\begin{align}
& \mathbb{A}^{IS}_i=\left[\begin{array}{l:l}
\mathbb{M}^{IS}_i & \left[\Phi\right]_i  \\
\hdashline
\left[\Phi\right]_i^{\mathrm{T}} & \mathbf{0}_{2 \times 2}
\end{array}\right],
\\
& \mathbf{b}^{IS}_i=\left[\begin{array}{c}
\mathbf{0}_{K \times 1} \\
\hdashline
U_{i-1/2}^{+,  mod}-\alpha_0^{i,old} \phi_0^{(I_i)}\left(x_{i-1 / 2}\right) \\
U_{i+1/2}^{-,  mod}-\alpha_0^{i,old} \phi_0^{(I_i)}\left(x_{i+1 / 2}\right)
\end{array}\right],
\\
& \mathbf{X}^{IS}_i
=
\left[\begin{array}{c}
\tilde{\mathbf{a}}^{mod}_i\\
\hdashline
\left[\Lambda^*\right]
\end{array}\right].
\end{align}
\end{subequations}
Thus, the DG modal coefficients modified by the IS-TVB(D)-minmod limiter can be obtained
$\mathbf{a}^{mod}_i=\left[\alpha_0^{i,old},\left[\tilde{\mathbf{a}}^{mod}_i\right]^{\mathrm{T}}\right]^{\mathrm{T}}$.

\subsubsection{Smoothness Factor $IS$ Constrained TVB(D)-minmod Limiter: Two-Dimensional Case}\label{sec-ISTVB-2D}
We consider the DG weak solution on the rectangular element $\Omega_{\mathcal{K}}$: 
\begin{align}
  U_h(x,y,t)=\sum_{k=0}^{K}\alpha_k(t)\phi_k(x,y),\ (x,y)\in\Omega_{\mathcal{K}}.
\end{align}
Let $\mathbf{D}$ be a multi-index of spatial partial derivatives, 
then
\begin{align}
  U_h^{(\mathbf{D})}(x,y,t)=\sum_{k=0}^{K}\alpha_k(t)\phi_k^{(\mathbf{D})}(x,y).
\end{align}
In the two-dimensional case, the computation of the smoothness measurement function $IS$ is as follows:
\begin{align}
  IS\left(Q_K(x,y)\right)=\sum_{1 \leq \left|\mathbf{D}\right| \leq K}\left|\Omega_{\mathcal{K}}\right|^{2\left|\mathbf{D}\right|-1}\int_{\Omega_{\mathcal{K}}}\left(\partial^{\mathbf{D}}Q_K(x,y)\right)^{2}\mathrm{~d}x\mathrm{d}y.
\end{align}
We denote
\begin{align}
  M_{kl}^{\mathbf{D}}:=\left|\Omega_{\mathcal{K}}\right|^{2\left|\mathbf{d}\right|-1}\sum_{\left|\mathbf{d}\right|=\left|D\right|}\langle\phi_k^{(\mathbf{d})},\phi_l^{(\mathbf{d})}\rangle_{L^{2}(\Omega_{\mathcal{K}})},
\end{align}
where $\mathbf{d}$ is a multi-index of spatial partial derivatives.
\\
Thus, we have
\begin{align}
  IS\left(U_h\right) &= \sum_{1 \leq \left|\mathbf{D}\right| \leq K}\left|\Omega_{\mathcal{K}}\right|^{2\left|\mathbf{D}\right|-1}\int_{\Omega_{\mathcal{K}}}\left(\partial^{\mathbf{D}}U_h\right)^{2}\mathrm{~d}x\mathrm{d}y \nonumber \\
  &= \sum_{1 \leq \left|\mathbf{D}\right| \leq K}\left|\Omega_{\mathcal{K}}\right|^{2\left|\mathbf{D}\right|-1}\int_{\Omega_{\mathcal{K}}}\left(\sum_{k={\color{red}{1}}}^{K}\alpha_{k}\phi_k^{(\mathbf{D})}\right)^{2}\mathrm{~d}x\mathrm{d}y \nonumber \\
  &= \sum_{1 \leq \left|\mathbf{D}\right| \leq K}\left|\Omega_{\mathcal{K}}\right|^{2\left|\mathbf{D}\right|-1}\sum_{k={\color{red}{1}}}^{K}\sum_{l={\color{red}{1}}}^{K}\alpha_k\alpha_l\langle\phi_k^{(\mathbf{D})},\phi_l^{(\mathbf{D})}\rangle_{L^{2}(\Omega_{\mathcal{K}})} \nonumber \\
  &= \sum_{1 \leq \left|\mathbf{D}\right| \leq K}\sum_{k={\color{red}{1}}}^{K}\sum_{l={\color{red}{1}}}^{K}\alpha_k\alpha_l M^{\mathbf{D}}_{kl} \nonumber \\
  &= \sum_{1 \leq \left|\mathbf{D}\right| \leq K}\left(\sum_{k={\color{red}{1}}}^{K}\alpha_k^{2} M^{\mathbf{D}}_{kk} + \sum_{k={\color{red}{1}}}^{\color{red}{K-1}}\sum_{l={\color{red}{k+1}}}^{K}2\alpha_k\alpha_l M_{kl}^{\mathbf{D}}\right) \nonumber \\
  &= \sum_{1 \leq \left|\mathbf{D}\right| \leq K}\sum_{k={\color{red}{1}}}^{K}\alpha_k^{2} M^{\mathbf{D}}_{kk} + \sum_{1 \leq \left|\mathbf{D}\right| \leq K}\sum_{k={\color{red}{1}}}^{\color{red}{K-1}}\sum_{l={\color{red}{k+1}}}^{K}2\alpha_k\alpha_l M_{kl}^{\mathbf{D}} \nonumber \\
  &= \sum_{k={\color{red}{1}}}^{K}\alpha_k^{2} \cdot \left(\sum_{1 \leq \left|\mathbf{D}\right| \leq K} M^{\mathbf{D}}_{kk}\right) + \sum_{k={\color{red}{1}}}^{\color{red}{K-1}}\sum_{l={\color{red}{k+1}}}^{K}2\alpha_k\alpha_l \cdot \left(\sum_{1 \leq \left|\mathbf{D}\right| \leq K} M_{kl}^{\mathbf{D}}\right).
\end{align}
Note that 
{\color{red}{when $|\mathbf{d}| > \min(k,l)$, it holds that $\phi_{\min(k,l)}^{(\mathbf{d})}=0$, }}\\
we have
\begin{align}
  \sum_{1 \leq |\mathbf{D}| \leq K}M_{kl}^{\mathbf{D}}=\sum_{1 \leq |\mathbf{D}| \leq \min(k,l)}M_{kl}^{\mathbf{D}},
\end{align}
therefore,
\begin{align}\label{2D-IS}
  IS\left(U_h\right) 
  &= 
  \sum_{k={\color{red}{1}}}^{K}\alpha_k^{2} \cdot \left(\sum_{1 \leq \left|\mathbf{D}\right| \leq k} M^{\mathbf{D}}_{kk}\right) + \sum_{k={\color{red}{1}}}^{\color{red}{K-1}}\sum_{l={\color{red}{k+1}}}^{K}2\alpha_k\alpha_l \cdot \left(\sum_{1 \leq \left|\mathbf{D}\right| \leq {\color{blue}{\min(k,l)}}} M_{kl}^{\mathbf{D}}\right) \nonumber \\
  &=
  \underbrace{\sum_{k={\color{red}{1}}}^{K}\alpha_k^{2} \cdot \left(\sum_{1 \leq \left|\mathbf{D}\right| \leq k} M^{\mathbf{D}}_{kk}\right)}_{\sim \mathbb{D}} + \underbrace{\sum_{k={\color{red}{1}}}^{\color{red}{K-1}}\sum_{l={\color{red}{k+1}}}^{K}2\alpha_k\alpha_l \cdot \left(\sum_{1 \leq \left|\mathbf{D}\right| \leq {\color{blue}{k}}} M_{kl}^{\mathbf{D}}\right)}_{\sim \left(\mathbb{L}+{\color{red}{\mathbb{U}}}\right),\ \mathbb{L}={\color{red}{\mathbb{U}}}^{\mathrm{T}}}.
\end{align}
Note that unlike the 1D case,
\begin{align}\label{2D-Mkl}
  M_{kl}^{\mathbf{D}}=\sum_{|\mathbf{d}|=|\mathbf{D}|}\langle\phi_k^{(\mathbf{d})},\phi_l^{(\mathbf{d})}\rangle_{L^{2}(\Omega_{\mathcal{K}})}.
\end{align}
The remaining derivation is similar to the 1D case.

\subsubsection{TVB(D)-minmod Limiter under Dual Constraints of Smoothness Factor $IS$ and $L^2$-Error}\label{subsubSec-IS-L2-TVB(D)-minmod}
\ {\color{blue}{The TVB(D)-minmod limiter provides three (or two) correction conditions that inherently suppress numerical oscillations, }} 
{\color{blue}{optimizing for the \textit{smoothness factor} still only serves to further suppress numerical oscillations, }} 
{\color{blue}{but \textit{over-suppression} of numerical oscillations is clearly detrimental to accuracy: }} 
the smoothness measurement function is in the form of a \textit{sum of squares},  
with 0 being its optimization ideal point (extreme value). 
Obviously, 0 is quite unfavorable for accuracy 
(large dissipation, with the numerical solution curve tending towards a horizontal line). 
\par
Numerical oscillations are generally more prevalent in high-order schemes.   
{\color{red}{\textit{Considering that the TVB(D)-minmod limiter's three (or two) correction conditions have already played a role in suppressing numerical oscillations}}}, 
we \textit{should select {\color{red}{the solution with higher precision from the weakly oscillatory solution set provided by the original TVB(D)-minmod limiter}}}. 
Therefore, it is appropriate to make the {\color{red}{\textit{corrected numerical solution as close as possible to the oscillating numerical solution of the original high-order scheme under the premise of weak oscillations. }}} 
Hence, one can choose the $L^2$-error to characterize the ``closeness'' of the corrected solution to 
the original high-order scheme's solution. This is the approach taken in the literature \cite{ref1}.
\par
{\color{red}{It is worth noting the following observations:}} \\
$\bullet\ $For the computation of the {\color{red}{Burgers' equation using $P^3$-polynomial approximation,}} the ``smoothness factor constraint'' 
is superior to the ``high-order scheme $L^2$-error constraint.'' 
{\color{red}{(Please refer to Example \ref{example-limiters-1D-Burgers} in subSection \ref{subSec-IS-L2-TVB-scalar}).}}\\
$\bullet\ $In the {\color{red}{1D-Euler Lax shock tube problem}}, the ``smoothness factor constraint'' is superior to the ``high-order scheme $L^2$-error constraint'' 
{\color{red}{(the ``high-order scheme $L^2$-error constraint'' may lead to a ``blow-up'')}}.
\par
The shock tube problem is an example that is relatively extreme and demanding, hence stability (suppression of oscillations) 
is more important than accuracy. {\color{red}{The ``smoothness factor constraint'' introduces additional conditions to suppress oscillations 
beyond those provided by the TVB(D)-minmod limiter (equivalent to additional artificial viscosity),}} thus yielding better results. 
On the other hand, the ``high-order scheme $L^2$-error constraint'' may lead to 
computational divergence due to insufficient suppression of oscillations.
\par
{\color{red}{Therefore, we can further consider the combination of the two and construct and solve the bi-objective optimization problem}} 
according to ``Linear Weighted Sum Method'' from optimization theory. 
\\ \hfill \\
$\blacksquare\ $One-dimensional case 
\par
A common approach to multi-objective optimization is to form a new objective function by taking a linear convex combination of several objective functions, i.e., \\
\begin{equation}
\begin{aligned}
& \min \left\{{\color{blue}{\omega_{IS} \cdot IS(u_i^{mod})}} + {\color{red}{\omega_{L^2} \cdot \left\|u_i^{mod }-u_i^{old}\right\|_{L^2(I_i)}}}\right\}\\
& s.t. \quad 
\left\{\begin{array}{c}
\frac{1}{\left|I_i\right|} \int_{I_i} u_i^{mod}(x) d x=\overline{U}_i \\
u_i^{mod}\left(x_{i-1/2}\right)=U_{i-1/2}^{ +, mod } \\
u_i^{mod}\left(x_{i+1/2}\right)=U_{i+1/2}^{ -, mod }
\end{array}\right.
\end{aligned}
\end{equation}
The integral mean remains unchanged, hence $\alpha_0^{i,mod}=\alpha_0^{i,old}$, which allows further simplification as follows: 
\begin{equation}\label{eq-IS-L2-optimization}
\begin{aligned}
& \min \left\{{\color{blue}{\omega_{IS} \cdot IS(u_i^{mod})}} + {\color{red}{\omega_{L^2} \cdot \left\|u_i^{mod }-u_i^{old}\right\|_{L^2(I_i)}}}\right\}\\
& s.t. \quad 
\left\{\begin{array}{c}
u_i^{mod}\left(x_{i-1/2}\right)=U_{i-1/2}^{ +, mod } \\
u_i^{mod}\left(x_{i+1/2}\right)=U_{i+1/2}^{ -, mod }
\end{array}\right.
\end{aligned}
\end{equation}
We continue to use the notations from the previous Section \ref{sec-ISTVB-1D}: 
\begin{align*}
\tilde{\mathbf{a}}^{mod}_i&:=\left[{\color{red}{\alpha_1^{i,mod}}},  \alpha_2^{i,mod},  \alpha_3^{i,mod},  \cdots,  \alpha_{K-1}^{i,mod},  \alpha_K^{i,mod}\right]^{\mathrm{T}},
\\
\left[\Phi\right]_i&:=
\left[
\begin{array}{llll}
  \phi_1^{(I_i)}\left(x_{i-1 / 2}\right) & \phi_2^{(I_i)}\left(x_{i-1 / 2}\right) & \cdots & \phi_K^{(I_i)}\left(x_{i-1 / 2}\right) \\
  \phi_1^{(I_i)}\left(x_{i+1 / 2}\right) & \phi_2^{(I_i)}\left(x_{i+1 / 2}\right) & \cdots & \phi_K^{(I_i)}\left(x_{i+1 / 2}\right)
\end{array}
\right]^{\mathrm{T}},
\\
\left[\Lambda^*\right]&:=\left[\lambda_1^*, \lambda^*_2\right]^{\mathrm{T}}.
\end{align*}
The solution $u_i^{mod}$ to problem (\refeq{eq-IS-L2-optimization}) satisfies the following saddle-point equation: 
\begin{align}
\mathbb{A}^{IS+L^2}_i \mathbf{X}_i^{IS+L^2} = \mathbf{b}^{IS+L^2}_i,
\end{align}
where  
\begin{subequations}
\begin{align}
&\mathbb{A}^{IS+L^2}_i
=\left[\begin{array}{c:c}
\omega_{IS} \mathbb{M}^{IS}_i + \omega_{L^2} \mathbb{M}^{L^2}_i & \left[\Phi\right]_i  \\
\hdashline
\left[\Phi\right]_i^{\mathrm{T}} & \mathbf{0}_{2 \times 2}
\end{array}\right],
\\
& \mathbf{b}^{IS+L^2}_i=
\left[\begin{array}{c}
\omega_{IS} \cdot \mathbf{0}_{K\times 1} + \omega_{L^2} \cdot {\color{red}{2}}\tilde{\mathbf{a}}^{i,old}\\
\hdashline
U_{i-1/2}^{+,  mod}-\alpha_0^{i,old} \phi_0^{(I_i)}\left(x_{i-1 / 2}\right) \\
U_{i+1/2}^{-,  mod}-\alpha_0^{i,old} \phi_0^{(I_i)}\left(x_{i+1 / 2}\right)
\end{array}\right],
\\
& \mathbf{X}_i^{IS+L^2}=
\left[\begin{array}{c}
\tilde{\mathbf{a}}^{mod}_i\\
\hdashline
\left[\Lambda^*\right]
\end{array}\right].
\end{align}
\end{subequations}
Note: 
\begin{align}\label{eq-a-old}
& \mathbb{M}^{L^2}_i=\operatorname{diag}(2)_{K\times K}, \\
& \tilde{\mathbf{a}}^{i,old}=\left[{\color{red}{\alpha_1^{i,old}}},  \alpha_2^{i,old},  \alpha_3^{i,old},  \ldots,  \alpha_{K-1}^{i,old},  \alpha_K^{i,old}\right]^{\mathrm{T}}.
\end{align}
Finally, the optimal solution to problem (\refeq{eq-IS-L2-optimization}) can be obtained 
$\mathbf{a}^{mod}_i=\left[\alpha_0^{i,old},\left[\tilde{\mathbf{a}}^{mod}_i\right]^{\mathrm{T}}\right]^{\mathrm{T}}$.
\par
\hfill \\
$\blacksquare\ $Two-dimensional case
\par
For the DG weak solution $U_{ij}$ on the two-dimensional rectangular element $\Omega_{ij}$, 
the following bi-objective optimization problem is constructed:
\begin{equation}\label{eq-2D-IS+L2}
\begin{aligned}
& \min \left\{{\color{blue}{\omega_{IS} \cdot IS(U_{ij}^{mod})}} + {\color{red}{\omega_{L^2} \cdot \left\|U_{ij}^{mod }-U_{ij}^{old}\right\|_{L^{2}(\Omega_{ij})}}}\right\}\\
& s.t. \quad 
\left\{
  \begin{aligned}
    & \int_{\partial\Omega_{ij}^{L}}U_{ij}^{mod}(x,y)dl=\overline{U}^{L,mod}_{\partial\Omega_{ij}}=\overline{U}^{L,mod}_{ij} \\
    & \int_{\partial\Omega_{ij}^{R}}U_{ij}^{mod}(x,y)dl=\overline{U}^{R,mod}_{\partial\Omega_{ij}}=\overline{U}^{R,mod}_{ij} \\
    & \int_{\partial\Omega_{ij}^{B}}U_{ij}^{mod}(x,y)dl=\overline{U}^{B,mod}_{\partial\Omega_{ij}}=\overline{U}^{B,mod}_{ij} \\
    & \int_{\partial\Omega_{ij}^{T}}U_{ij}^{mod}(x,y)dl=\overline{U}^{T,mod}_{\partial\Omega_{ij}}=\overline{U}^{T,mod}_{ij}
  \end{aligned}
\right.
\end{aligned}
\end{equation}
For conservation requirement, it still holds that $\alpha_0^{(ij),mod}=\alpha_0^{(ij),old}$. 
\\
Here, we denote
{\renewcommand{\arraystretch}{1.25}
\begin{align}
&  \tilde{\mathbf{a}}^{old}_{(ij)}:=\left[{\color{red}{\alpha_1^{(ij),old}}},\alpha_2^{(ij),old},\ldots,\alpha_{(K+1)(K+2)/2-1}^{(ij),old}\right]^{\mathrm{T}}
,\\
&  \tilde{\mathbf{a}}^{mod}_{(ij)}:=\left[{\color{red}{\alpha_1^{(ij),mod}}},\alpha_2^{(ij),mod},\ldots,\alpha_{(K+1)(K+2)/2-1}^{(ij),mod}\right]^{\mathrm{T}}
,\\
&\left.
\begin{aligned}
& \Gamma^{L,(ij)}_{k}:=\int_{\partial\Omega_{ij}^{L}}\phi^{(ij)}_k(x,y)dl,\quad \Gamma^{R,(ij)}_{k}:=\int_{\partial\Omega_{ij}^{R}}\phi^{(ij)}_k(x,y)dl \\
& \Gamma^{B,(ij)}_{k}:=\int_{\partial\Omega_{ij}^{B}}\phi^{(ij)}_k(x,y)dl,\quad \Gamma^{T,(ij)}_{k}:=\int_{\partial\Omega_{ij}^{T}}\phi^{(ij)}_k(x,y)dl \\
& k={\color{blue}{0}},1,2,3,\cdots,(K+1)(K+2)/2-1
\end{aligned}
\right\}
,\\
&\left[\Gamma\right]_{(ij)}:=\left[\begin{matrix}
  \Gamma^{L,(ij)}_{\color{red}{1}} & \Gamma_2^{L,(ij)} & \cdots & \Gamma_{(K+1)(K+2)/2-1}^{L,(ij)} \\
  \Gamma^{R,(ij)}_{\color{red}{1}} & \Gamma_2^{R,(ij)} & \cdots & \Gamma_{(K+1)(K+2)/2-1}^{R,(ij)} \\
  \Gamma^{B,(ij)}_{\color{red}{1}} & \Gamma_2^{B,(ij)} & \cdots & \Gamma_{(K+1)(K+2)/2-1}^{B,(ij)} \\
  \Gamma^{T,(ij)}_{\color{red}{1}} & \Gamma_2^{T,(ij)} & \cdots & \Gamma_{(K+1)(K+2)/2-1}^{T,(ij)}
  \end{matrix}\right]
,\\
&\left[\Lambda^*\right]:=\left[\lambda_1^*, \lambda^*_2, \lambda^*_3, \lambda^*_4\right]^{\mathrm{T}}. 
\end{align}
}
The saddle-point equation is given by:
\begin{align}
\mathbb{A}^{IS+L^2}_{(ij)} \mathbf{X}^{IS+L^2}_{(ij)} = \mathbf{b}^{IS+L^2}_{(ij)},
\end{align}
where 
\begin{subequations}
\begin{align}
  &\mathbb{A}_{(ij)}^{IS+L^2}
  =
  \left[\begin{array}{c:c}
  \omega_{IS}\mathbb{M}^{IS}_{(ij)} + \omega_{L^2} \mathbb{M}^{L^{2}}_{(ij)} & \left[\Gamma\right]^{\mathrm{T}}_{(ij)}  \\
  \hdashline
  \left[\Gamma\right]_{(ij)} & \mathbf{0}_{4 \times 4}  \\
  \end{array}\right],
  \\
  &\mathbf{b}^{IS+L^2}_{(ij)}=\left[\begin{array}{c}
    \omega_{IS} \cdot \mathbf{0}_{(K+1)(K+2)/2-1\ \times\ 1} + \omega_{L^2} \cdot {\color{red}{2}}\tilde{\mathbf{a}}^{old}_{(ij)} \\
  \hdashline
  \overline{U}_{i j}^{L, mod} -\alpha_0^{(ij),old} \Gamma_0^{L,(ij)} \\
  \overline{U}_{i j}^{R, mod} -\alpha_0^{(ij),old} \Gamma_0^{R,(ij)} \\
  \overline{U}_{i j}^{B, mod} -\alpha_0^{(ij),old} \Gamma_0^{B,(ij)} \\
  \overline{U}_{i j}^{T, mod} -\alpha_0^{(ij),old} \Gamma_0^{T,(ij)}
  \end{array}\right],
  \\
  &\mathbf{X}^{IS+L^2}_{(ij)}=\left[\begin{array}{c}
  \tilde{\mathbf{a}}^{mod}_{(ij)} \\
  \hdashline
  \left[\Lambda^*\right]
  \end{array}\right].
\end{align}
\end{subequations}
Note that in the two-dimensional case,
$\mathbb{M}^{L^2}_{(ij)}=\operatorname{diag}(2)_{(K+1)(K+2)/2-1\ \times\ (K+1)(K+2)/2-1}$, 
while $\mathbb{M}^{IS}_{(ij)}$ must be derived based on the \textit{formulas \eqref{2D-IS} and \eqref{2D-Mkl} in Section \ref{sec-ISTVB-2D}}, 
following the process \textit{\eqref{eq-alpha-M} to \eqref{MIS} outlined in Section \ref{sec-ISTVB-1D}}.  
\\
Then, $\mathbf{a}^{mod}_{(ij)}=\left[\alpha_0^{(ij),old},\left[\tilde{\mathbf{a}}^{mod}_{(ij)}\right]^{\mathrm{T}}\right]^{\mathrm{T}}$ 
is just the optimal solution to problem \eqref{eq-2D-IS+L2}.  
\par
We denote the limiter based on the bi-objective optimization problem concerning the $IS$ function and $L^2$-error 
as the IS-$L^2$-TVB(D)-minmod limiter. It is worth noting that when setting $\omega_{L^2} = 0$ 
and $\omega_{IS} = 1$, the IS-$L^2$-TVB(D)-minmod limiter degenerates into the IS-TVB(D)-minmod limiter discussed 
in the previous sections \ref{sec-ISTVB-1D} and \ref{sec-ISTVB-2D}; On the other hand, if we take $\omega_{L^2} = 1$ and $\omega_{IS} = 0$, 
then the IS-$L^2$-TVB(D)-minmod limiter degenerates into the $L^2$-TVB(D)-minmod limiter originated from literature \cite{ref1}.

\begin{remark}[Hybrid of $P_{IS}^K$ and $P^K_{DG}$]
Recently, the literature \cite{ref42} proposed a hybrid strategy that combines a dissipative 
but low-accuracy limiter (typically a first-order linear polynomial) 
with the high-order polynomials provided by DG methods in a convex combination. 
This strategy effectively balances the goals of ``suppressing numerical oscillations" and ``preserving the high accuracy of the DG solution."
However, {\color{blue}{as mentioned in \cite{ref42}, the classical TVB(D)-minmod limiter ultimately ``replaces high-order schemes with linear polynomials on the troubled cells, thus achieving essentially oscillation-free." }}
This observation {\color{red}{aligns with the point raised in this paper, 
where it is emphasized that ``the classical TVB(D)-minmod limiter is only well-suited for low-order polynomials 
and cannot directly modify all coefficients of higher-order polynomials, 
as the number of correction equations does not match the number of coefficients of the higher-order polynomial." }}
\par
In summary, the literature \cite{ref42} has indeed recognized that {\color{orange}{the classical TVB(D)-minmod limiter 
can only generate low-order polynomials. }}
{\color{blue}{Their approach does not alter the result of the classical TVB(D)-minmod limiter;}} 
instead, it combines it with the original DG high-order polynomial in a convex combination to obtain a higher-order polynomial. 
In contrast, {\color{red}{our approach directly modifies the classical TVB(D)-minmod limiter, }}
transforming it into an optimization problem that can directly generate high-order polynomials. 
\end{remark}

\section{Reconstruction in Characteristic Field}\label{sec-characteristic-reconstruction}
When applying limiters to a system of equations, there are typically two implementation methods: (a) directly reconstructing each component in the original physical space; (b) performing local characteristic decomposition in conjunction to transform into the characteristic space and then reconstructing each component individually.
\par
The former is easier to implement in a program. Directly reconstructing by component is essentially assuming that the components are already decoupled. 
However, in the original system of equations, the unknown components are usually coupled, so the actual effect may not be satisfactory. 
The latter method uses the Jacobian of the flux function to transform the system into the characteristic space, achieving true decoupling, 
and then proceeds to reconstruct each component individually.  
\par
The basic principle of characteristic reconstruction is as follows 
(taking a one-dimensional hyperbolic conservation system as an example):
\\
Consider the equation
$$\partial_t \mathbf{U} + \partial_x (f(\mathbf{U})) = 0,$$
where $A = \frac{\partial f(\mathbf{U})}{\partial \mathbf{U}}$, and assume it is diagonalizable, 
that is, $A = R \Lambda L$ and $RL = I_d$. Then the above equation can be rewritten as
\begin{align}
\partial_t \mathbf{U} + R \Lambda L \partial_x \mathbf{U} = 0.
\end{align}
When $A$ is a constant matrix ( $R, \Lambda, L$ are all constant matrices), 
introduce the characteristic variable $\mathbf{V} = L\mathbf{U}$, and multiply both sides of the equation by $L$ to obtain
\begin{align}
\partial_t \mathbf{V} + \Lambda \partial_x \mathbf{V} = 0.
\end{align}
Since $\Lambda$ is a constant diagonal matrix, 
the above equation is a completely decoupled set of $N$ independent constant-coefficient single-wave equations 
( $\mathbf{U} \in \mathbb{R}^{\mathrm{d}}$ ).
Thus, $\mathbf{V}$ can be reconstructed component by component to obtain $\widetilde{\mathbf{V}}$, 
and then the inverse transformation using $R$ is applied to return to the original physical space, 
that is, $\widetilde{\mathbf{U}} = R\widetilde{\mathbf{V}}$, which completes the characteristic reconstruction. 
\par
Considering the variable coefficients and nonlinear characteristics of the equations (where $A$, $R$, $\Lambda$, 
and $L$ are function matrices), local freezing must be introduced to fix $A$, $R$, $\Lambda$, and $L$. 
Similar to the point values in the FD method and the cell-averaged integrals in the FV schemes, 
the DG scheme also requires determining the ``values at which'' to substitute into the Jacobian of the flux function 
and its characteristic structure, i.e., the ``values at which'' to implement local freezing. 
The approach of this paper can be referred to in Remark\ \ref{remark-OurLocalFreezing}. 
\par
Furthermore, in the DG scheme, each component of $\mathbf{U}$ is a polynomial function; 
therefore, special attention is needed on how to use the eigenstructure of the Jacobian ($L, R$) 
to transform the polynomial into the corresponding characteristic space.  
\par
Specifically, the general procedure for characteristic reconstruction in 2D systems is as follows:
\\
Let the troubled cell be denoted as $\Omega_0$, 
with its boundary composed of $\partial\Omega_0 = \cup_{\ell=1}^{\mathcal{L}}\Gamma_0^{\ell}$, 
and the neighboring cells adjacent to $\Omega_0$ across $\Gamma_0^{\ell}$ are denoted as 
$\Omega_{\ell},\ \ell = 1, 2, \cdots, \mathcal{L}$. 
The outward normal vector on $\Gamma_0^{\ell}$ is $\mathbf{n}^{\ell} = \left(n^{\ell}_x, n_y^{\ell}\right)$, 
and the Jacobian in this direction (the normal Jacobian) is 
$A_n^{\ell} = n^{\ell}_x \cdot \frac{\partial \mathbf{F}}{\partial \mathbf{U}} + n_y^{\ell} \cdot \frac{\partial \mathbf{G}}{\partial \mathbf{U}}$, 
with its eigenstructure being $\mathbb{R}_n^{\ell},\ \mathbb{L}_n^{\ell}$.  
\\
$\bullet\ $step1. Generate $A_n^{\ell}$ and its eigenstructure $\left(\mathbb{L}^{\ell}_n, \mathbb{R}^{\ell}_n\right)$ using the arithmetic mean of the cell integral means of $\Omega_0, \Omega_{\ell}$, i.e., $(\overline{u}_0 + \overline{u}_{\ell}) / 2$. 
\\
$\bullet\ $step2. Perform the characteristic transformation of the template polynomials using $\mathbb{L}^{\ell}_n$ for each edge.  
\\
$\bullet\ $step3. In the characteristic space, perform reconstruction component by component.  
\\
$\bullet\ $setp4. Use the corresponding $\mathbb{R}_n^{\ell}$ to inverse transform back to the physical space, obtaining characteristic reconstruction polynomials $u^{mod, \ell}_0, \ell = 1, 2, \cdots, \mathcal{L}$.  
\\
$\bullet\ $step5. The final result of the characteristic reconstruction is $u^{mod}_0 = \sum_{\ell} w_\ell \cdot u^{mod, \ell}_0$. 
If it is a uniform grid, then $w_{i} = 1 / \mathcal{L}$; 
if it is a non-uniform grid, then $w_{i} = \frac{|\Omega_{i}|}{\sum_{\ell=1}^{\mathcal{L}}|\Omega_{\ell}|}$.  
\begin{remark}\label{remark-OurLocalFreezing}
  It is noted that limiters ultimately serve to correct the interface fluxes, which are calculated on the boundaries of cells. 
  Therefore, unlike the conventional practices mentioned above, in this paper we choose the integral average values of the field functions on the grid edges as candidate states 
  (the cells on either side of an edge can each provide an edge integral average value) and then takes their arithmetic average or Roe average as the local freezing value. 
  \par
  Consider two adjacent cells, $\Omega_0$ and $\Omega_{\ell}$, sharing a common boundary $\Gamma_0^{\ell}$, 
  with DG solutions $u_0$ and $u_{\ell}$, respectively. 
  The integral averages of $u_0$ and $u_{\ell}$ on $\Gamma_0^{\ell}$ are computed component-wise as follows: 
  \begin{align*}
    \overline{\Gamma}_0^{\ell}(i, \Omega_0)=\frac{1}{\left|\Gamma_0^{\ell}\right|}\int_{\Gamma_0^{\ell}}u_0^{(i)}dl, \quad i=1, 2, \cdots, m.
    \\
    \overline{\Gamma}_0^{\ell}(i, \Omega_{\ell})=\frac{1}{\left|\Gamma_0^{\ell}\right|}\int_{\Gamma_0^{\ell}}u_{\ell}^{(i)}dl, \quad i=1, 2, \cdots, m.
  \end{align*}
  Here, $m$ denotes the total number of system components. 
  Then, compute the arithmetic or Roe average of $\overline{\Gamma}_0^{\ell}(i, \Omega_0)$ and $\overline{\Gamma}_0^{\ell}(i, \Omega_{\ell})$ component-wise, 
  and use these averages to generate the normal Jacobian. 
\end{remark}

Next, the one-dimensional hyperbolic conservation system is taken as an example to 
detail a previously used characteristic transformation process for polynomials, referred to as  ``moment characteristic transformation'' or ``modal coefficient transformation'': 
\par \noindent
Consider
$
\partial_t\mathbf{U} + \partial_x \mathbf{F}(\mathbf{U})=0, 
\ \mathbf{U}=(u_1, u_2, u_3, \cdots, u_m)^{\mathrm{T}}
$. 
\\
Let
$
\mathbb{J}=\frac{\partial \mathbf{F}(\mathbf{U})}{\partial \mathbf{U}} =\mathbb{R} \Lambda \mathbb{L}
$.  
Note: $\mathbb{R} \mathbb{L}=\mathbb{I}_d(\text{identity matrix});\ \mathbb{J},  \mathbb{R},  \Lambda,  \mathbb{L}$ are all matrix-valued functions of $\mathbf{U}$. 
\\
The approximate solution for the $i$-th component $u_i$ is given by: 
$$
u_i=\sum_{k=0}^{k=K} \alpha_k^{(i)} \phi_{i,  k},
$$
where $i$ denotes the $i$-th component, $K$ represents the highest degree of the basis functions (polynomials), 
$\phi_{i, k}$ represents the $k$-th basis function (polynomial) of the $i$-th component, 
and $\alpha_k^{(i)}$ is the $k$-th moment (or $k$-th modal coefficient) of the $i$-th component. 
Form a column vector with the basis functions of all orders for the $i$-th system component 
\begin{align}
  \left[\boldsymbol{\Phi}^{(i)}\right]=\left[\phi_{i,0},\phi_{i,1},\phi_{i,2},\ldots,\phi_{i,K}\right]^{\mathrm{T}}.
\end{align}
Now, partition the left eigenvector matrix $\mathbb{L}$ column-wise, such that
\begin{align}
\mathbb{L} &= \left[\begin{array}{ccccc}
l_{11} & l_{12} & l_{13} & \cdots & l_{1 m} \\
l_{21} & l_{22} & l_{23} & \cdots & l_{2 m} \\
l_{31} & l_{32} & l_{33} & \cdots & l_{3 m} \\
\vdots & \vdots & \vdots & \cdots & \vdots \\
l_{m 1} & l_{m 2} & l_{m 3} & \cdots & l_{m m}
\end{array}\right]=\left(\mathbf{L}_1,  \mathbf{L}_2,  \mathbf{L}_3,  \cdots,  \mathbf{L}_m\right).
\end{align}
Form a row vector consisting of the modal coefficients of the $i$-th system component, denoted as ${\left[\mathbf{a}^{(i)}\right]^{\mathrm{T}}}$. 
The modal coefficients of all components constitute the following modal coefficient matrix 
\begin{align}
\mathbb{A} &= \left[\begin{array}{ccccc}
\alpha_0^{(1)} & \alpha_1^{(1)} & \alpha_2^{(1)} & \cdots & \alpha_K^{(1)} \\
\alpha_0^{(2)} & \alpha_1^{(2)} & \alpha_2^{(2)} & \cdots & \alpha_K^{(2)} \\
\alpha_0^{(3)} & \alpha_1^{(3)} & \alpha_2^{(3)} & \cdots & \alpha_K^{(3)} \\
\vdots & \vdots & \vdots & \cdots & \vdots \\
\alpha_0^{(m)} & \alpha_1^{(m)} & \alpha_2^{(m)} & \cdots & \alpha_K^{(m)}
\end{array}\right]
=\left(\begin{array}{c}
{\left[\mathbf{a}^{(1)}\right]^{\mathrm{T}}} \\
{\left[\mathbf{a}^{(2)}\right]^{\mathrm{T}}} \\
{\left[\mathbf{a}^{(3)}\right]^{\mathrm{T}}} \\
\vdots \\
{\left[\mathbf{a}^{(m)}\right]^{\mathrm{T}}}
\end{array}\right).
\end{align}
When the left eigenvector matrix $\mathbb{L}$ acts on the modal coefficient matrix $\mathbb{A}$, 
it produces the matrix of characteristic modal coefficients $\mathbb{B}$, given by
\begin{align}
\mathbb{B}=\mathbb{L} \mathbb{A} & =\left(\mathbf{L}_1,  \mathbf{L}_2,  \mathbf{L}_3,  \cdots,  \mathbf{L}_m\right)\left(\begin{array}{c}
{\left[\mathbf{a}^{(1)}\right]^{\mathrm{T}}} \\
{\left[\mathbf{a}^{(2)}\right]^{\mathrm{T}}} \\
{\left[\mathbf{a}^{(3)}\right]^{\mathrm{T}}} \\
\vdots \\
{\left[\mathbf{a}^{(m)}\right]^{\mathrm{T}}}
\end{array}\right) \nonumber \\
& =\mathbf{L}_1\left[\mathbf{a}^{(1)}\right]^{\mathrm{T}}+\mathbf{L}_2\left[\mathbf{a}^{(2)}\right]^{\mathrm{T}}+\mathbf{L}_3\left[\mathbf{a}^{(3)}\right]^{\mathrm{T}}+\cdots+\mathbf{L}_m\left[\mathbf{a}^{(m)}\right]^{\mathrm{T}}.
\end{align}
The characteristic modal coefficient matrix $\mathbb{B}$ has its $i$-th row denoted as $\mathbb{B}(i, :)$, satisfying
\begin{align}\label{moment-CT-B}
\mathbb{B}(i, :) & =\mathbf{L}_1(i) \cdot\left[\mathbf{a}^{(1)}\right]^{\mathrm{T}}+\mathbf{L}_2(i) \cdot\left[\mathbf{a}^{(2)}\right]^{\mathrm{T}}+\mathbf{L}_3(i) \cdot\left[\mathbf{a}^{(3)}\right]^{\mathrm{T}}+\cdots+\mathbf{L}_m(i) \cdot\left[\mathbf{a}^{(m)}\right]^{\mathrm{T}}.
\end{align}
Let
\begin{subequations}\label{moment-CT-b}
\begin{align}
  & \left[\mathbf{b}^{(i)}\right]^{\mathrm{T}}=\mathbb{B}(i, :), \\
  & \beta_k^{(i)}={\left[\mathbf{b}^{(i)}\right]^{\mathrm{T}}(k)=\mathbb{B}(i, k)}. 
\end{align}
\end{subequations}
Define
\begin{align}
v_i=\left[\mathbf{b}^{(i)}\right]^{\mathrm{T}}\left[\boldsymbol{\Phi}^{(i)}\right]=\sum_{k=0}^{k=K} \beta_k^{(i)} \phi_{i,  k},
\end{align}
then the polynomial $v_i$ is the projection of the polynomial $u_i$ in the original physical space onto the characteristic space.  
\begin{remark}
  The procedure for the inverse characteristic transformation is analogous to that of the characteristic transformation; 
  it only requires substituting $\mathbb{L}$ with $\mathbb{R}$. 
\end{remark}

\subsection{Interpolation-based Characteristic Transformation}
Distinct from the ``moment characteristic transformation'' mentioned earlier, 
this paper proposes an ``interpolation-based characteristic transformation'', 
and illustrates the process using a one-dimensional hyperbolic conservation system as an example. 
\par
Let the set of sampling points be $\mathbf{X} = \left(x_0, x_1, x_2, \cdots, x_K\right)^{\mathrm{T}}$. 
Correspondingly, the sampled values for the $i$-th component are 
$\mathbf{Y}^{(i)} = \left(y_0^{(i)}, y_1^{(i)}, y_2^{(i)}, \cdots, y_K^{(i)}\right)^{\mathrm{T}}$, where $y_k^{(i)} = u_i(x_k)$.
The sampled values of all components together form the sampling point matrix $\mathbb{Y}$: 
\begin{align}
\mathbb{Y}=\left(\begin{array}{c}
{\left[\mathbf{Y}^{(1)}\right]^{\mathrm{T}}} \\
{\left[\mathbf{Y}^{(2)}\right]^{\mathrm{T}}} \\
{\left[\mathbf{Y}^{(3)}\right]^{\mathrm{T}}} \\
\vdots \\
{\left[\mathbf{Y}^{(m)}\right]^{\mathrm{T}}}
\end{array}\right).
\end{align}
Transform the sampling point matrix to the characteristic space
\begin{align}
\widetilde{\mathbb{Y}}=\mathbb{L} \mathbb{Y} & =\left(\mathbf{L}_1,  \mathbf{L}_2,  \mathbf{L}_3,  \cdots,  \mathbf{L}_m\right)\left(\begin{array}{c}
{\left[\mathbf{Y}^{(1)}\right]^{\mathrm{T}}} \\
{\left[\mathbf{Y}^{(2)}\right]^{\mathrm{T}}} \\
{\left[\mathbf{Y}^{(3)}\right]^{\mathrm{T}}} \\
\vdots \\
{\left[\mathbf{Y}^{(m)}\right]^{\mathrm{T}}}
\end{array}\right) \nonumber \\
& =\mathbf{L}_1\left[\mathbf{Y}^{(1)}\right]^{\mathrm{T}}+\mathbf{L}_2\left[\mathbf{Y}^{(2)}\right]^{\mathrm{T}}+\mathbf{L}_3\left[\mathbf{Y}^{(3)}\right]^{\mathrm{T}}+\cdots+\mathbf{L}_m\left[\mathbf{Y}^{(m)}\right]^{\mathrm{T}}.
\end{align}
The $i$-th row of the characteristic sampling point matrix $\widetilde{\mathbb{Y}}$ is denoted as $\widetilde{\mathbb{Y}}(i, :)$, satisfying
\begin{align}
\widetilde{\mathbb{Y}}(i, :)=\mathbf{L}_1(i) \cdot\left[\mathbf{Y}^{(1)}\right]^{\mathrm{T}}+\mathbf{L}_2(i) \cdot\left[\mathbf{Y}^{(2)}\right]^{\mathrm{T}}+\mathbf{L}_3(i) \cdot\left[\mathbf{Y}^{(3)}\right]^{\mathrm{T}}+\cdots+\mathbf{L}_m(i) \cdot\left[\mathbf{Y}^{(m)}\right]^{\mathrm{T}}.
\end{align}
Let
\begin{align}\label{8.13}
\widetilde{\mathbf{Y}}^{(i)}:=[\widetilde{\mathbb{Y}}(i, :)]^{\mathrm{T}}=\mathbf{L}_1(i) \cdot \mathbf{Y}^{(1)}+\mathbf{L}_2(i) \cdot \mathbf{Y}^{(2)}+\mathbf{L}_3(i) \cdot \mathbf{Y}^{(3)}+\cdots+\mathbf{L}_m(i) \cdot \mathbf{Y}^{(m)},
\end{align}
then $\widetilde{\mathbf{Y}}^{(i)}$ is the projection of the original physical space's $\mathbf{Y}^{(i)}$ onto the characteristic space. 
In other words, $\widetilde{\mathbf{Y}}^{(i)}$ represents the values of the characteristic projection $\widetilde{u}_i$ 
of the original physical space's $u_i$ at the same set of sampling points $\mathbf{X} = \left(x_0, x_1, x_2, \cdots, x_K\right)^{\mathrm{T}}$.  
\\
At this point, the polynomial $\widetilde{u}_i$ in the characteristic space that satisfies the interpolation data set $\left(\mathbf{X}, \widetilde{\mathbf{Y}}^{(i)}\right)$ is determined using the method of undetermined coefficients.
\\
Let the modal coefficient set of the characteristic polynomial $\widetilde{u}_i$ be $\widetilde{\mathbf{a}}^{(i)}$, so that
$$
\widetilde{u}_i = \sum_{k=0}^{K} \widetilde{\alpha}_k^{(i)} \phi_{i, k} \quad \text{with} \quad \widetilde{\alpha}_k^{(i)} = \widetilde{\mathbf{a}}^{(i)}(k).
$$
Let
\begin{align}
\mathbb{P}^{(i)}:=\left[\begin{array}{ccccc}
\phi_{i,  0}\left(x_0\right) & \phi_{i,  1}\left(x_0\right) & \phi_{i,  2}\left(x_0\right) & \cdots & \phi_{i,  K}\left(x_0\right) \\
\phi_{i,  0}\left(x_1\right) & \phi_{i,  1}\left(x_1\right) & \phi_{i,  2}\left(x_1\right) & \cdots & \phi_{i,  K}\left(x_1\right) \\
\phi_{i,  0}\left(x_2\right) & \phi_{i,  1}\left(x_2\right) & \phi_{i,  2}\left(x_2\right) & \cdots & \phi_{i,  K}\left(x_2\right) \\
\vdots & \vdots & \vdots & \cdots & \vdots \\
\phi_{i,  0}\left(x_K\right) & \phi_{i,  1}\left(x_K\right) & \phi_{i,  2}\left(x_K\right) & \cdots & \phi_{i,  K}\left(x_K\right)
\end{array}\right].
\end{align}
The characteristic polynomial $\widetilde{u}_i$ should satisfy
\begin{align}\label{8.15}
  \widetilde{u}_i(\mathbf{X})=\widetilde{\mathbf{Y}}^{(i)},
\end{align}
hence
\begin{align}\label{8.16}
\mathbb{P}^{(i)} \widetilde{\mathbf{a}}^{(i)}=\widetilde{\mathbf{Y}}^{(i)}, 
\end{align}
and thus
\begin{align}\label{8.17}
\widetilde{\mathbf{a}}^{(i)}=\left[\mathbb{P}^{(i)}\right]^{-1} \widetilde{\mathbf{Y}}^{(i)}. 
\end{align}
With this, the projection of the polynomial $u_i$ from the original physical space onto the characteristic space, 
$\widetilde{u}_i$, can be obtained.  

\subsection{Equivalence of Interpolation-based Characteristic Transformation and Moment Characteristic Transformation}
\begin{proposition}[Equivalence of Interpolation-based Characteristic Transformation and Moment Characteristic Transformation when Basis Functions are Identical]\label{proposition8.1}
When all components of the system share the same set of basis functions, that is,
$ K_1 = K_2 = K_3 = \cdots = K_m = K $
and
$ \phi_{1, k}(\cdot) \equiv \phi_{2, k}(\cdot) \equiv \phi_{3, k}(\cdot) \equiv \cdots \equiv \phi_{m, k}(\cdot),\ \forall k \in \{1, 2, 3, \cdots, K\}, $
it follows that 
$ \mathbb{P}^{(1)} = \mathbb{P}^{(2)} = \mathbb{P}^{(3)} = \cdots = \mathbb{P}^{(m)} $, 
thus the interpolation-based characteristic transformation is equivalent to the moment characteristic transformation, i.e.,
$ \widetilde{\mathbf{a}}^{(i)} = \mathbf{b}^{(i)},\ \forall i \in \{1, 2, 3, \ldots, m\}. $
\end{proposition}
Taking the one-dimensional hyperbolic conservation system as an example, Proposition \ref{proposition8.1} is proven.
\par
\begin{proof}
By substituting equation (\refeq{8.13}) into equation (\refeq{8.17}), we obtain
\begin{align}\label{8.18}
\left[\widetilde{\mathbf{a}}^{\color{red}{(i)}}\right]^{\mathrm{T}}&=\left(\left[\mathbb{P}^{(i)}\right]^{-1} \widetilde{\mathbf{Y}}^{(i)}\right)^{\mathrm{T}}=\left[\widetilde{\mathbf{Y}}^{(i)}\right]^{\mathrm{T}}\left(\left[\mathbb{P}^{(i)}\right]^{-1}\right)^{\mathrm{T}}=\left[\widetilde{\mathbf{Y}}^{(i)}\right]^{\mathrm{T}}\left(\left[\mathbb{P}^{(i)}\right]^{\mathrm{T}}\right)^{-1} \nonumber \\
& =\left(\mathbf{L}_1(i) \cdot \mathbf{Y}^{(1)}+\mathbf{L}_2(i) \cdot \mathbf{Y}^{(2)}+\cdots+\mathbf{L}_m(i) \cdot \mathbf{Y}^{(m)}\right)^{\mathrm{T}}\left(\left[\mathbb{P}^{(i)}\right]^{-1}\right)^{\mathrm{T}} \nonumber \\
& =\mathbf{L}_1(i) \cdot\left[\mathbf{Y}^{(1)}\right]^{\mathrm{T}}\left(\left[\mathbb{P}^{\color{red}{(i)}}\right]^{-1}\right)^{\mathrm{T}}+\mathbf{L}_2(i) \cdot\left[\mathbf{Y}^{(2)}\right]^{\mathrm{T}}\left(\left[\mathbb{P}^{\color{red}{(i)}}\right]^{-1}\right)^{\mathrm{T}}
+\cdots+\mathbf{L}_m(i) \cdot\left[\mathbf{Y}^{(m)}\right]^{\mathrm{T}}\left(\left[\mathbb{P}^{\color{red}{(i)}}\right]^{-1}\right)^{\mathrm{T}}.
\end{align}
From equations (\refeq{moment-CT-B}) and (\refeq{moment-CT-b}), it is easy to know
\begin{align}\label{8.19}
  \left[\mathbf{b}^{(i)}\right]^{\mathrm{T}} =\mathbf{L}_1(i) \cdot\left[\mathbf{a}^{(1)}\right]^{\mathrm{T}}+\mathbf{L}_2(i) \cdot\left[\mathbf{a}^{(2)}\right]^{\mathrm{T}}+\mathbf{L}_3(i) \cdot\left[\mathbf{a}^{(3)}\right]^{\mathrm{T}}+\cdots+\mathbf{L}_m(i) \cdot\left[\mathbf{a}^{(m)}\right]^{\mathrm{T}}.
\end{align}
Note that in the original physical space, it holds that $u_i(\mathbf{X})=\mathbf{Y}^{(i)}$, similar to equations (\refeq{8.15}) and (\refeq{8.16}), we have
\begin{subequations}\label{8.20}
\begin{align}
\mathbb{P}^{(i)} {\mathbf{a}}^{(i)}&={\mathbf{Y}}^{(i)},\ i=1,2,\ldots,m; \label{8.20a}\\
{\mathbf{a}}^{(i)}&=\left[\mathbb{P}^{(i)}\right]^{-1} {\mathbf{Y}}^{(i)},\ i=1,2,\ldots,m. \label{8.20b}
\end{align}
\end{subequations}
By substituting equation (\refeq{8.20b}) into equation (\refeq{8.19}), we get
\begin{align}
\left[{\mathbf{b}}^{\color{red}{(i)}}\right]^{\mathrm{T}}
=
\mathbf{L}_1(i) \cdot\left[\mathbf{Y}^{(1)}\right]^{\mathrm{T}}\left(\left[\mathbb{P}^{\color{red}{(1)}}\right]^{-1}\right)^{\mathrm{T}}+\mathbf{L}_2(i) \cdot\left[\mathbf{Y}^{(2)}\right]^{\mathrm{T}}\left(\left[\mathbb{P}^{\color{red}{(2)}}\right]^{-1}\right)^{\mathrm{T}}
+\cdots+\mathbf{L}_m(i) \cdot\left[\mathbf{Y}^{(m)}\right]^{\mathrm{T}}\left(\left[\mathbb{P}^{\color{red}{(m)}}\right]^{-1}\right)^{\mathrm{T}}.
\end{align}
Comparing ${\left[\widetilde{\mathbf{a}}^{(i)}\right]^{\mathrm{T}} }$ and $\left[\mathbf{b}^{(i)}\right]^{\mathrm{T}}$, 
and noting the superscripts of the matrix $\mathbb{P}$, \\
it can be seen that \\
if $\phi_{1, k}=\phi_{2, k}=\phi_{3, k}=\cdots=\phi_{m, k}, \forall k \in\{1, 2, 3, \cdots, K\}$, \\ 
then $\mathbb{P}^{(1)}=\mathbb{P}^{(2)}=\mathbb{P}^{(3)}=\cdots=\mathbb{P}^{(m)}$, \\
and thus $\widetilde{\mathbf{a}}^{(i)}=\mathbf{b}^{(i)}, i=1, 2, 3, \cdots, m$. \\
That is, at this point, the ``moment characteristic transformation'' and the ``interpolation-based characteristic transformation'' 
are equivalent. 
\end{proof}

\section{Numerical Results}\label{sec-NumericalTests}
This section first conducts a accuracy test on the FVS-DG method for the compressible Euler equations and the shallow water wave equations, 
including one-dimensional and two-dimensional test cases. 
Subsequently, it demonstrates the advantage of the IS-$L^2$-TVB(D)-minmod limiter in suppressing numerical oscillations 
near discontinuities on the linear scalar transport equation, the nonlinear Burgers' equation, and the non-convex Buckley-Leverett equation. 
Finally, we solve several classical Riemann problems for the compressible Euler equations and the shallow water wave equations 
using the FVS-DG method coupled with the IS-$L^2$-TVB(D)-minmod limiter, 
including the Sod problem, Lax problem, Shu-Osher problem, Blast problem, 
three typical Riemann test cases from the literature \cite{ref43}, 
and the Dam Break problem for the shallow water wave equations \cite{ref38}.
\par
In all numerical experiments of this section, the adiabatic index for the compressible Euler equations is set to $\gamma = 1.4$ (simulating the dynamical process of air), 
and the acceleration due to gravity for the shallow water wave equations is set to $ g = 9.8120\ m/s^2$.
\par
All numerical experiments in this section use a uniform mesh partitioning, 
where the two-dimensional test cases include both rectangular and structured triangular meshes.

\subsection{Accuracy Tests for FVS-DG}
\subsubsection{One-dimensional Test Cases}
\begin{example}
1D-Euler Compressible Equations.
\\
Control Eqs: 
$$
\frac{\partial}{\partial t}\left(\begin{array}{c}
\rho \\
\rho u \\
E
\end{array}\right)+\frac{\partial}{\partial x}\left(\begin{array}{c}
\rho u \\
\rho u^2+p \\
u(E+p)
\end{array}\right)=0;
$$
Computational domain:  
$$
\Omega \times [0,T_{end}] = [0,1] \times [0,1];
$$
$\mathrm{I.C.}\ $
$$
\begin{aligned}
    \rho_0(x) = 1 + 0.2 \cos(\pi x) ,\quad
    u_0(x) = -0.7 ,\quad
    P_0(x) = 1;
\end{aligned}
$$
$\mathrm{B.C.}\ ${\color{blue}{periodic boundary conditions; }}
\\ \hspace*{\fill} \\
True solutions: 
$$
\begin{aligned}
    \rho(x,t) = 1 + 0.2 \cos(\pi (x+0.7t)) ,\quad
    u(x,t) = -0.7 ,\quad
    P(x,t) = 1.
\end{aligned}
$$

\begin{itemize}
  \item Numderical flux format: AUSM; \\
  Temporal discretization format: TVD-RK3; \\
  CFL=0.1. \\
  $L^{\infty},\ L^2,\ L^1$\ numerical errors and convergence orders  for conservative variables $\rho,\ \rho u,\ E$
  with $P^2$- and $P^3$-polynomial approximations 
  are summarized respectively in Table \ref{table-1D-Euler-Accuracy-AUSM-P2} and Table \ref{table-1D-Euler-Accuracy-AUSM-P3}. 
  {\tiny
  \begin{table}[htbp]
    \centering
    \caption{\tiny $P^2$-FVS(AUSM)-DG using equally spaced cells. 1D-Compressible Euler Equations with smooth initial conditions: 
    $\rho_0(x) = 1 + 0.2 \cos(\pi x) ,\ u_0(x) = -0.7 ,\ P_0(x) = 1$. $L^{\infty},\ L^2,\ L^1$\ errors and convergence orders  for 
    conservative variables $\rho,\ \rho u,\ E$. }
    \label{table-1D-Euler-Accuracy-AUSM-P2}
    \begin{tabular}{lllllllll}
    \hline  & $Mesh$ & 10 & 20 & 40 & 80 & 160 \\
    \hline  & $L^{\infty}$-error & 1.1024E-03   & 1.4025E-04   & 1.7693E-05   & 2.2156E-06   & 2.7701E-07 \\
            & $L^{\infty}$-order & —        & 2.9746   & 2.9868   & 2.9974   & 2.9997 \\
    $\rho$  & $L^{2}$-error      & 2.4737E-04   & 3.1521E-05   & 3.9607E-06   & 4.9575E-07   & 6.1989E-08  \\
            & $L^{2}$-order      & —        & 2.9723    & 2.9925    & 2.9981    & 2.9995 \\
            & $L^{1}$-error      & 1.8559E-04   & 2.3295E-05   & 2.9126E-06   & 3.6395E-07   & 4.5479E-08 \\
            & $L^{1}$-order      & —        & 2.9940    & 2.9997    & 3.0005    & 3.0005 \\
    \hline    & $L^{\infty}$-error & 5.3443E-04   & 6.6924E-05   & 8.4195E-06   & 1.0543E-06   & 1.3184E-07 \\
              & $L^{\infty}$-order & —          & 2.9974    & 2.9907    & 2.9975    & 2.9994 \\
    $\rho u$  & $L^{2}$-error      & 2.4737E-04   & 3.1521E-05   & 3.9607E-06   & 4.9575E-07   & 6.1989E-08  \\
              & $L^{2}$-order      & —          & 2.9723    & 2.9925    & 2.9981    & 2.9995 \\
              & $L^{1}$-error      & 1.8559E-04   & 2.3295E-05   & 2.9126E-06   & 3.6395E-07   & 4.5479E-08 \\
              & $L^{1}$-order      & —          & 2.9940    & 2.9997    & 3.0005    & 3.0005 \\
    \hline    & $L^{\infty}$-error & 8.2911E-04   & 1.0494E-04   & 1.3268E-05   & 1.6616E-06   & 2.0772E-07 \\
              & $L^{\infty}$-order & —        & 2.9819    & 2.9836    & 2.9972    & 2.9999 \\
    $   E  $  & $L^{2}$-error      & 2.4737E-04   & 3.1521E-05   & 3.9607E-06   & 4.9575E-07   & 6.1989E-08  \\
              & $L^{2}$-order      & —        & 2.9723    & 2.9925    & 2.9981    & 2.9995 \\
              & $L^{1}$-error      & 1.8559E-04   & 2.3295E-05   & 2.9126E-06   & 3.6395E-07   & 4.5479E-08 \\
              & $L^{1}$-order      & —        & 2.9940    & 2.9997    & 3.0005    & 3.0005 \\
    \hline  
    \end{tabular}
    \end{table}

  \begin{table}[htbp]
    \centering
    \caption{\tiny $P^3$-FVS(AUSM)-DG using equally spaced cells. 1D-Compressible Euler Equations with smooth initial conditions: 
    $\rho_0(x) = 1 + 0.2 \cos(\pi x) ,\ u_0(x) = -0.7 ,\ P_0(x) = 1$. $L^{\infty},\ L^2,\ L^1$\ errors and convergence orders  for 
    conservative variables $\rho,\ \rho u,\ E$. }
    \label{table-1D-Euler-Accuracy-AUSM-P3}
    \begin{tabular}{lllllllll}
    \hline  & $Mesh$ & 10 & 20 & 40 & 80 & 160 \\
    \hline  & $L^{\infty}$-error & 4.1850E-05   & 2.7313E-06   & 1.7220E-07   & 1.0778E-08   & 6.7472E-10 \\
            & $L^{\infty}$-order & —        & 3.9376    & 3.9875    & 3.9978    & 3.9977 \\
    $\rho$  & $L^{2}$-error      & 7.8755E-06   & 4.9119E-07   & 3.0563E-08   & 1.9097E-09   & 1.1935E-10  \\
            & $L^{2}$-order      & —        & 4.0030    & 4.0064    & 4.0004    & 4.0001 \\
            & $L^{1}$-error      & 5.4508E-06   & 3.3906E-07   & 2.0980E-08   & 1.3095E-09   & 8.1857E-11 \\
            & $L^{1}$-order      & —        & 4.0069    & 4.0145    & 4.0019    & 3.9998 \\
    \hline    & $L^{\infty}$-error & 2.3771E-05   & 1.5271E-06   & 9.6370E-08   & 6.0382E-09   & 3.7818E-10 \\
              & $L^{\infty}$-order & —          & 3.9604    & 3.9861    & 3.9964    & 3.9970 \\
    $\rho u$  & $L^{2}$-error      & 7.8755E-06   & 4.9119E-07   & 3.0563E-08   & 1.9097E-09   & 1.1935E-10  \\
              & $L^{2}$-order      & —          & 4.0030    & 4.0064    & 4.0004    & 4.0001 \\
              & $L^{1}$-error      & 5.4508E-06   & 3.3906E-07   & 2.0980E-08   & 1.3095E-09   & 8.1857E-11 \\
              & $L^{1}$-order      & —          & 4.0069    & 4.0145    & 4.0019    & 3.9998 \\
    \hline    & $L^{\infty}$-error & 3.2236E-05   & 2.0668E-06   & 1.2930E-07   & 8.0974E-09   & 5.0454E-10 \\
              & $L^{\infty}$-order & —        & 3.9632    & 3.9986    & 3.9971    & 4.0044 \\
    $   E  $  & $L^{2}$-error      & 7.8755E-06   & 4.9119E-07   & 3.0563E-08   & 1.9097E-09   & 1.1935E-10  \\
              & $L^{2}$-order      & —        & 4.0030    & 4.0064    & 4.0004    & 4.0001 \\
              & $L^{1}$-error      & 5.4508E-06   & 3.3906E-07   & 2.0980E-08   & 1.3095E-09   & 8.1857E-11 \\
              & $L^{1}$-order      & —        & 4.0069    & 4.0145    & 4.0019    & 3.9998 \\
    \hline  
    \end{tabular}
    \end{table}
    }
    \hfill \\
    \item Numderical flux format: Steger-Warming; \\
    Temporal discretization format: TVD-RK3; \\
    CFL=0.1. \\
    $L^{\infty},\ L^2,\ L^1$\ numerical errors and convergence orders  for conservative variables $\rho,\ \rho u,\ E$ 
    with $P^2$- and $P^3$-polynomial approximations   
    are summarized respectively in Table \ref{table-1D-Euler-Accuracy-Steger-Warming-P2} and Table \ref{table-1D-Euler-Accuracy-Steger-Warming-P3}.
    {\tiny
    \begin{table}[htbp]
      \centering
      \caption{\tiny $P^2$-FVS(Steger-Warming)-DG using equally spaced cells. 1D-Compressible Euler Equations with smooth initial conditions: 
    $\rho_0(x) = 1 + 0.2 \cos(\pi x) ,\ u_0(x) = -0.7 ,\ P_0(x) = 1$. $L^{\infty},\ L^2,\ L^1$\ errors and convergence orders  for 
    conservative variables $\rho,\ \rho u,\ E$. }
    \label{table-1D-Euler-Accuracy-Steger-Warming-P2}
      \begin{tabular}{lllllllll}
      \hline  & $Mesh$ & 10 & 20 & 40 & 80 & 160 \\
      \hline  & $L^{\infty}$-error & 1.1728E-03   & 1.5873E-04   & 2.0237E-05   & 2.5390E-06   & 3.1747E-07 \\
              & $L^{\infty}$-order & —        & 2.8853    & 2.9716    & 2.9947    & 2.9996 \\
      $\rho$  & $L^{2}$-error      & 2.7279E-04   & 3.6698E-05   & 4.6856E-06   & 5.8895E-07   & 7.3721E-08  \\
              & $L^{2}$-order      & —        & 2.8940    & 2.9694    & 2.9920    & 2.9980 \\
              & $L^{1}$-error      & 2.0508E-04   & 2.7383E-05   & 3.4908E-06   & 4.3859E-07   & 5.4894E-08 \\
              & $L^{1}$-order      & —        & 2.9048    & 2.9717    & 2.9926    & 2.9981 \\
      \hline    & $L^{\infty}$-error & 5.2307E-04   & 6.8032E-05   & 8.5503E-06   & 1.0695E-06   & 1.3369E-07 \\
                & $L^{\infty}$-order & —          & 2.9427    & 2.9922    & 2.9991    & 3.0000 \\
      $\rho u$  & $L^{2}$-error      & 2.7279E-04   & 3.6698E-05   & 4.6856E-06   & 5.8895E-07   & 7.3721E-08  \\
                & $L^{2}$-order      & —          & 2.8940    & 2.9694    & 2.9920    & 2.9980 \\
                & $L^{1}$-error      & 2.0508E-04   & 2.7383E-05   & 3.4908E-06   & 4.3859E-07   & 5.4894E-08 \\
                & $L^{1}$-order      & —          & 2.9048    & 2.9717    & 2.9926    & 2.9981 \\
      \hline    & $L^{\infty}$-error & 8.8839E-04   & 1.2148E-04   & 1.5573E-05   & 1.9591E-06   & 2.4527E-07 \\
                & $L^{\infty}$-order & —        & 2.8705    & 2.9636    & 2.9908    & 2.9978 \\
      $   E  $  & $L^{2}$-error      & 2.7279E-04   & 3.6698E-05   & 4.6856E-06   & 5.8895E-07   & 7.3721E-08  \\
                & $L^{2}$-order      & —        & 2.8940    & 2.9694    & 2.9920    & 2.9980 \\
                & $L^{1}$-error      & 2.0508E-04   & 2.7383E-05   & 3.4908E-06   & 4.3859E-07   & 5.4894E-08 \\
                & $L^{1}$-order      & —        & 2.9048    & 2.9717    & 2.9926    & 2.9981 \\
      \hline  
      \end{tabular}
      \end{table}

      \begin{table}[htbp]
        \centering
        \caption{\tiny $P^3$-FVS(Steger-Warming)-DG using equally spaced cells. 1D-Compressible Euler Equations with smooth initial conditions: 
        $\rho_0(x) = 1 + 0.2 \cos(\pi x) ,\ u_0(x) = -0.7 ,\ P_0(x) = 1$. $L^{\infty},\ L^2,\ L^1$\ errors and convergence orders  for 
        conservative variables $\rho,\ \rho u,\ E$. }
        \label{table-1D-Euler-Accuracy-Steger-Warming-P3}
        \begin{tabular}{lllllllll}
        \hline  & $Mesh$ & 10 & 20 & 40 & 80 & 160 \\
        \hline  & $L^{\infty}$-error & 3.6110E-05   & 2.3683E-06   & 1.4967E-07   & 9.3550E-09   & 5.8578E-10 \\
                & $L^{\infty}$-order & —        & 3.9304    & 3.9840    & 3.9999    & 3.9973 \\
        $\rho$  & $L^{2}$-error      & 6.6707E-06   & 4.1281E-07   & 2.5553E-08   & 1.5954E-09   & 9.9687E-11  \\
                & $L^{2}$-order      & —        & 4.0143    & 4.0139    & 4.0015    & 4.0004 \\
                & $L^{1}$-error      & 4.6767E-06   & 2.8866E-07   & 1.7743E-08   & 1.1087E-09   & 6.9260E-11 \\
                & $L^{1}$-order      & —        & 4.0180    & 4.0240    & 4.0003    & 4.0008 \\
        \hline    & $L^{\infty}$-error & 2.1895E-05   & 1.4053E-06   & 8.8435E-08   & 5.5424E-09   & 3.4715E-10 \\
                  & $L^{\infty}$-order & —          & 3.9616    & 3.9901    & 3.9960    & 3.9969 \\
        $\rho u$  & $L^{2}$-error      & 6.6707E-06   & 4.1281E-07   & 2.5553E-08   & 1.5954E-09   & 9.9687E-11  \\
                  & $L^{2}$-order      & —          & 4.0143    & 4.0139    & 4.0015    & 4.0004 \\
                  & $L^{1}$-error      & 4.6767E-06   & 2.8866E-07   & 1.7743E-08   & 1.1087E-09   & 6.9260E-11 \\
                  & $L^{1}$-order      & —          & 4.0180    & 4.0240    & 4.0003    & 4.0008 \\
        \hline    & $L^{\infty}$-error & 2.4615E-05   & 1.5544E-06   & 9.6544E-08   & 6.0418E-09   & 3.7597E-10 \\
                  & $L^{\infty}$-order & —        & 3.9851    & 4.0090    & 3.9981    & 4.0063 \\
        $   E  $  & $L^{2}$-error      & 6.6707E-06   & 4.1281E-07   & 2.5553E-08   & 1.5954E-09   & 9.9687E-11  \\
                  & $L^{2}$-order      & —        & 4.0143    & 4.0139    & 4.0015    & 4.0004 \\
                  & $L^{1}$-error      & 4.6767E-06   & 2.8866E-07   & 1.7743E-08   & 1.1087E-09   & 6.9260E-11 \\
                  & $L^{1}$-order      & —        & 4.0180    & 4.0240    & 4.0003    & 4.0008 \\
        \hline  
        \end{tabular}
        \end{table}
      }
\end{itemize}
\end{example}
\begin{example}
  1D-Shallow Water Wave Equations.
  \\
  Control Eqs: 
  $$
  \begin{aligned}
  \partial_t
  \left[\begin{array}{l}
  h \\
  h u
  \end{array}\right]
  +
  \partial_x
  \left[\begin{array}{c}
  h u \\
  h u^2+\frac{1}{2} g h^2
  \end{array}\right]
  =
  \left[\begin{array}{c}
  0 \\
  -g h\left(z_0\right)_x
  \end{array}\right]
  \end{aligned};
  $$
  Bottom elevation: 
  $$
  Z_0 = \sin^2(\pi x);
  $$
  Computational domain:  
  $$
  \Omega \times [0,T_{end}] = [0,1] \times [0,0.075];
  $$
  Note: The equation will develop a discontinuity near $t = 0.1$. To test for accuracy, we set the simulation time to $t = 0.075$, 
  at which point the solution is still smooth. This example does not have an analytical true solution; 
  instead, we use numerical solutions on finer meshes as references to obtain the errors and convergence orders of the algorithm. \\
  $\mathrm{I.C.}\ $
  $$
  h_0(x)=5+e^{\cos (2 \pi x)}, \quad u_0(x)=\frac{\sin (\cos (2 \pi x))}{5+e^{\cos (2 \pi x)}};
  $$
  $\mathrm{B.C.}\ ${\color{blue}{periodic boundary conditions}}; \\
  Numderical flux format: van Leer; \\
  Temporal discretization format: TVD-RK3; \\
  $L^{\infty},\ L^2,\ L^1$\ numerical errors and convergence orders  for conservative variables $h,\ hu$ 
  with $P^2$-polynomial approximation   
  are summarized respectively in Table \ref{table-1D-SWE-Accuracy-vanLeer-P2}.  
  {\tiny
  \begin{table}[htbp]
    \centering
    \caption{\tiny $P^2$-FVS(vanLeer)-DG using equally spaced cells with TVD-RK3 and CFL=$0.01$. 1D-Shallow Water Equations with smooth initial conditions: 
    $h_0(x)=5+e^{\cos (2 \pi x)}, \quad u_0(x)=\frac{\sin (\cos (2 \pi x))}{5+e^{\cos (2 \pi x)}}$. $L^{\infty},\ L^2,\ L^1$\ errors and convergence orders  for 
    conservative variables $h,\ h u$. }
    \label{table-1D-SWE-Accuracy-vanLeer-P2}
    \begin{tabular}{lllllllll}
    \hline  & $Mesh$ & 50 & 100 & 200 & 400 & 800 \\
    \hline  & $L^{\infty}$-error & 2.0111E-02   & 6.3551E-03   & 2.0769E-03   & 2.9931E-04   & 4.3458E-05 \\
            & $L^{\infty}$-order & — & 1.6620    & 1.6135    & 2.7947    & 2.7840 \\
    $ h $   & $L^{2}$-error      & 2.9563E-03   & 5.3791E-04   & 9.2442E-05   & 1.1241E-05   & 1.4936E-06  \\
            & $L^{2}$-order      & — &  2.4584    & 2.5407    & 3.0397    & 2.9119 \\
            & $L^{1}$-error      & 1.0235E-03   & 1.3895E-04   & 2.0637E-05   & 2.6006E-06   & 3.5516E-07 \\
            & $L^{1}$-order      & — &  2.8809    & 2.7513    & 2.9883    & 2.8723 \\
    \hline    & $L^{\infty}$-error & 7.4667E-02   & 2.5399E-02   & 8.3407E-03   & 1.2216E-03   & 1.7540E-04 \\
              & $L^{\infty}$-order & — &  1.5557    & 1.6065    & 2.7714    & 2.8000 \\
    $ h u $   & $L^{2}$-error      & 1.3649E-02   & 2.4733E-03   & 4.0608E-04   & 4.7905E-05   & 7.1252E-06  \\
              & $L^{2}$-order      & — &  2.4643    & 2.6066    & 3.0835    & 2.7492 \\
              & $L^{1}$-error      & 5.6059E-03   & 7.4207E-04   & 1.0431E-04   & 1.3226E-05   & 1.8956E-06 \\
              & $L^{1}$-order      & — &  2.9173    & 2.8307    & 2.9795    & 2.8026 \\
    \hline  
    \end{tabular}
    \end{table}
    }
\end{example}

\subsubsection{Two-dimensional Test Cases}
\begin{example}
2D-Compressible Euler Equations.
\\
Control Eqs: 
$$
\frac{\partial}{\partial t}\left(\begin{array}{c}
\rho \\
\rho u \\
\rho v \\
E
\end{array}\right)+\frac{\partial}{\partial x}\left(\begin{array}{c}
\rho u \\
\rho u^2+p \\
\rho u v \\
u(E+p)
\end{array}\right)+\frac{\partial}{\partial y}\left(\begin{array}{c}
\rho v \\
\rho u v \\
\rho v^2+p \\
v(E+p)
\end{array}\right)=0;
$$
Computational domain:  
$$
\Omega \times [0,T_{end}] = \left\{ [0,2] \times [-1,1] \right\} \times [0,1];
$$
$\mathrm{I.C.}\ $
$$
\begin{aligned}
    \rho_0(x,y) = 1 + 0.2 \cos(\pi x + \pi y) ,\quad
    u_0(x,y) = -0.7 ,\quad
    v_0(x,y) = 0.3 ,\quad
    P_0(x,y) = 1
\end{aligned};
$$
$\mathrm{B.C.}\ ${\color{blue}{periodic boundary conditions; }}
\\ \hspace*{\fill} \\
True solutions:
$$
\left\{
\begin{aligned}
    \rho(x,y,t) &= 1 + 0.2 \cos(\pi (x+0.7t) + \pi (y-0.3t)), \\
    u(x,y,t) &= -0.7, \\
    v(x,y,t) &= 0.3, \\
    P(x,y,t) &= 1.
\end{aligned}
\right.
$$
\begin{itemize}
  \item Numderical flux format: AUSM; \\
  Temporal discretization format: TVD-RK3; \\
  CFL=0.025; \\
  $L^{\infty},\ L^2,\ L^1$\ numerical errors and convergence orders  for conservative variables $\rho,\ \rho u,\ \rho v,\ E$ 
  with $P^2$-polynomial approximation  
  are shown in Table \ref{table-2D-Euler-Accuracy-AUSM-P2}.
  {\tiny
  \begin{table}[htbp]
    \centering
    \caption{\tiny $P^2$-FVS(AUSM)-DG using {\color{blue}{uniform structured triangular meshes}}. 2D-Compressible Euler Equationswith smooth initial conditions: 
    $\rho_0(x,y) = 1 + 0.2 \cos(\pi x + \pi y) ,\ u_0(x,y) = -0.7 ,\ v_0(x,y) = 0.3 ,\ P_0(x,y) = 1$. $L^{\infty},\ L^2,\ L^1$\ errors and convergence orders  for 
    conservative variables $\rho,\ \rho u,\ \rho v,\ E$. }
    \label{table-2D-Euler-Accuracy-AUSM-P2}
    \begin{tabular}{lllllllll}
    \hline  & $Mesh$ & 8$\times$8 & 16$\times$16 & 32$\times$32 & 64$\times$64 \\
    \hline  & $L^{\infty}$-error &  8.9167E-04   & 1.2342E-04   & 1.5657E-05   & 1.9642E-06 \\
            & $L^{\infty}$-order & —        & 2.8530    & 2.9786    & 2.9948 \\
    $\rho$  & $L^{2}$-error      & 7.0277E-04   & 9.0059E-05   & 1.1425E-05   & 1.4334E-06  \\
            & $L^{2}$-order      & —        & 2.9641    & 2.9787    & 2.9947 \\
            & $L^{1}$-error      & 1.1249E-03   & 1.4124E-04   & 1.7661E-05   & 2.2048E-06 \\
            & $L^{1}$-order      & —        & 2.9936    & 2.9995    & 3.0018 \\
    \hline    & $L^{\infty}$-error & 5.5158E-04   & 7.6466E-05  &  9.7164E-06   & 1.2195E-06 \\
              & $L^{\infty}$-order & —          & 2.8507    & 2.9763    & 2.9942 \\
    $\rho u$  & $L^{2}$-error      & 4.5169E-04   & 5.7971E-05   & 7.3577E-06   & 9.2326E-07  \\
              & $L^{2}$-order      & —          & 2.9619    & 2.9780    & 2.9944 \\
              & $L^{1}$-error      & 7.2709E-04   & 9.1365E-05   & 1.1470E-05   & 1.4332E-06 \\
              & $L^{1}$-order      & —          & 2.9924    & 2.9938    & 3.0005 \\
    \hline    & $L^{\infty}$-error & 3.4939E-04   & 4.5836E-05   & 5.8122E-06   & 7.3063E-07 \\
              & $L^{\infty}$-order & —          & 2.9303    & 2.9793    & 2.9919 \\
    $\rho v$  & $L^{2}$-error      & 2.6454E-04   & 3.4709E-05   & 4.4310E-06   & 5.5688E-07  \\
              & $L^{2}$-order      & —          & 2.9301    & 2.9696    & 2.9922 \\
              & $L^{1}$-error      & 4.2513E-04   & 5.4671E-05   & 6.8958E-06   & 8.6307E-07 \\
              & $L^{1}$-order      & —          & 2.9590    & 2.9870    & 2.9982 \\          
    \hline    & $L^{\infty}$-error & 4.1011E-04   & 5.2523E-05   & 6.6514E-06   & 8.3368E-07 \\
              & $L^{\infty}$-order & —        & 2.9650    & 2.9812    & 2.9961 \\
    $   E  $  & $L^{2}$-error      & 3.2336E-04   & 4.2374E-05   & 5.4071E-06   & 6.7950E-07  \\
              & $L^{2}$-order      & —        & 2.9319    & 2.9702    & 2.9923 \\
              & $L^{1}$-error      & 4.9739E-04   & 6.5928E-05   & 8.3778E-06   & 1.0506E-06 \\
              & $L^{1}$-order      & —        & 2.9154    & 2.9762    & 2.9954 \\
    \hline  
    \end{tabular}
    \end{table}
  }
  \hfill \\
  \item Numderical flux format: AUSM; \\
  Temporal discretization format: RK4; \\
  CFL=0.01; \\
  $L^{\infty},\ L^2,\ L^1$\ numerical errors and convergence orders  for conservative variables $\rho,\ \rho u,\ \rho v,\ E$ 
  with $P^3$-polynomial approximation  
  are shown in Table \ref{table-2D-Euler-Accuracy-AUSM-P3}.
  {\tiny
  \begin{table}[htbp]
    \centering
    \caption{\tiny $P^3$-FVS(AUSM)-DG using {\color{blue}{uniform structured rectangular meshes}}. 2D-Compressible Euler Equationswith smooth initial conditions: 
    $\rho_0(x,y) = 1 + 0.2 \cos(\pi x + \pi y) ,\ u_0(x,y) = -0.7 ,\ v_0(x,y) = 0.3 ,\ P_0(x,y) = 1$. $L^{\infty},\ L^2,\ L^1$\ errors and convergence orders  for 
    conservative variables $\rho,\ \rho u,\ \rho v,\ E$. }
    \label{table-2D-Euler-Accuracy-AUSM-P3}
    \begin{tabular}{lllllllll}
    \hline  & $Mesh$ & 10$\times$10 & 20$\times$20 & 40$\times$40 & 80$\times$80 \\
    \hline  & $L^{\infty}$-error & 1.5314E-04   & 9.0975E-06   & 5.6476E-07   & 3.2868E-08 \\
            & $L^{\infty}$-order & —        & 4.0732    & 4.0098    & 4.1029 \\
    $\rho$  & $L^{2}$-error      & 1.2555E-04   & 8.8175E-06   & 5.9757E-07   & 2.8805E-08  \\
            & $L^{2}$-order      & —        & 3.8317    & 3.8832    & 4.3747 \\
            & $L^{1}$-error      & 2.0791E-04   & 1.4274E-05   & 1.0451E-06   & 4.8381E-08 \\
            & $L^{1}$-order      & —        & 3.8645    & 3.7717    & 4.4331 \\
    \hline    & $L^{\infty}$-error & 9.2719E-05   & 5.8423E-06   & 3.6815E-07   & 1.9127E-08 \\
              & $L^{\infty}$-order & —          & 3.9883    & 3.9882    & 4.2666 \\
    $\rho u$  & $L^{2}$-error      & 8.0222E-05   & 5.7501E-06   & 3.9938E-07   & 1.8463E-08  \\
              & $L^{2}$-order      & —          & 3.8024    & 3.8477    & 4.4351 \\
              & $L^{1}$-error      & 1.3091E-04   & 9.3984E-06   & 6.9878E-07   & 3.0859E-08 \\
              & $L^{1}$-order      & —          & 3.8000    & 3.7495    & 4.5010 \\
    \hline    & $L^{\infty}$-error & 5.6795E-05   & 3.1587E-06   & 1.7577E-07   & 1.1372E-08 \\
              & $L^{\infty}$-order & —          & 4.1684    & 4.1675    & 3.9502 \\
    $\rho v$  & $L^{2}$-error      & 4.2407E-05   & 2.7470E-06   & 1.8497E-07   & 9.3140E-09  \\
              & $L^{2}$-order      & —          & 3.9484    & 3.8925    & 4.3118 \\
              & $L^{1}$-error      & 7.0188E-05   & 4.3302E-06   & 3.2201E-07   & 1.5252E-08 \\
              & $L^{1}$-order      & —          & 4.0187    & 3.7493    & 4.4000 \\          
    \hline    & $L^{\infty}$-error & 7.6232E-05   & 4.8337E-06   & 2.8666E-07   & 1.6498E-08 \\
              & $L^{\infty}$-order & —        & 3.9792    & 4.0757    & 4.1190 \\
    $   E  $  & $L^{2}$-error      & 5.7506E-05   & 3.8493E-06   & 2.3815E-07   & 1.2923E-08  \\
              & $L^{2}$-order      & —        & 3.9011    & 4.0146    & 4.2038 \\
              & $L^{1}$-error      & 9.5498E-05   & 6.1273E-06   & 4.1033E-07   & 2.1644E-08 \\
              & $L^{1}$-order      & —        & 3.9621    & 3.9004    & 4.2448 \\
    \hline  
    \end{tabular}
    \end{table}
  }
  \hfill \\
  \item Numderical flux format: Steger-Warming; \\
  Temporal discretization format: TVD-RK3; \\
  CFL=0.05; \\
  $L^{\infty},\ L^2,\ L^1$\ numerical errors and convergence orders  for conservative variables $\rho,\ \rho u,\ \rho v,\ E$ 
  with $P^2$-polynomial approximation  
  are shown in Table \ref{table-2D-Euler-Accuracy-StegerWarming-P2}.
  {\tiny
  \begin{table}[htbp]
    \centering
    \caption{\tiny $P^2$-FVS(Steger-Warming)-DG using {\color{blue}{uniform structured triangular meshes}}. 2D-Compressible Euler Equationswith smooth initial conditions: 
    $\rho_0(x,y) = 1 + 0.2 \cos(\pi x + \pi y) ,\ u_0(x,y) = -0.7 ,\ v_0(x,y) = 0.3 ,\ P_0(x,y) = 1$. $L^{\infty},\ L^2,\ L^1$\ errors and convergence orders  for 
    conservative variables $\rho,\ \rho u,\ \rho v,\ E$. }
    \label{table-2D-Euler-Accuracy-StegerWarming-P2}
    \begin{tabular}{lllllllll}
    \hline  & $Mesh$ & 8$\times$8 & 16$\times$16 & 32$\times$32 & 64$\times$64 \\
    \hline  & $L^{\infty}$-error &  9.3377E-04   & 1.3671E-04   & 1.7590E-05   & 2.1995E-06 \\
            & $L^{\infty}$-order & —        & 2.7719    & 2.9583    & 2.9996 \\
    $\rho$  & $L^{2}$-error      & 8.1292E-04   & 1.1096E-04   & 1.4344E-05   & 1.8091E-06  \\
            & $L^{2}$-order      & —        & 2.8730    & 2.9515    & 2.9872 \\
            & $L^{1}$-error      & 1.3097E-03   & 1.7368E-04   & 2.1929E-05   & 2.7399E-06 \\
            & $L^{1}$-order      & —        & 2.9147    & 2.9856    & 3.0006 \\
    \hline    & $L^{\infty}$-error & 5.9550E-04   & 8.5540E-05   & 1.1053E-05   & 1.3838E-06 \\
              & $L^{\infty}$-order & —          & 2.7994    & 2.9522    & 2.9977 \\
    $\rho u$  & $L^{2}$-error      & 5.1786E-04   & 7.0930E-05   & 9.1898E-06   & 1.1599E-06  \\
              & $L^{2}$-order      & —          & 2.8681    & 2.9483    & 2.9861 \\
              & $L^{1}$-error      & 8.4123E-04   & 1.1217E-04   & 1.4225E-05   & 1.7796E-06 \\
              & $L^{1}$-order      & —          & 2.9069    & 2.9791    & 2.9988 \\
    \hline    & $L^{\infty}$-error & 3.7102E-04   & 5.5142E-05   & 7.1893E-06   & 9.1021E-07 \\
              & $L^{\infty}$-order & —          & 2.7503    & 2.9392    & 2.9816 \\
    $\rho v$  & $L^{2}$-error      & 3.0768E-04   & 4.3852E-05   & 5.7438E-06   & 7.2709E-07  \\
              & $L^{2}$-order      & —          & 2.8107    & 2.9326    & 2.9818 \\
              & $L^{1}$-error      & 4.9562E-04   & 6.8821E-05   & 8.8443E-06   & 1.1098E-06 \\
              & $L^{1}$-order      & —          & 2.8483    & 2.9600    & 2.9945 \\          
    \hline    & $L^{\infty}$-error & 4.2025E-04   & 6.0872E-05   & 7.8618E-06   & 9.9419E-07 \\
              & $L^{\infty}$-order & —        & 2.7874    & 2.9528    & 2.9833 \\
    $   E  $  & $L^{2}$-error      & 3.4478E-04   & 4.9342E-05   & 6.4757E-06   & 8.2032E-07  \\
              & $L^{2}$-order      & —        & 2.8048    & 2.9297    & 2.9808 \\
              & $L^{1}$-error      & 5.4155E-04   & 7.7310E-05   & 1.0117E-05   & 1.2801E-06 \\
              & $L^{1}$-order      & —        & 2.8084    & 2.9339    & 2.9824 \\
    \hline  
    \end{tabular}
    \end{table}
  }
\end{itemize}
\end{example}
From the accuracy test results of the aforementioned one-dimensional and two-dimensional examples, 
it can be seen that both the $P^K$-FVS-DG scheme based on the Jacobian eigenvalue splitting 
and the $P^K$-FVS-DG scheme based on the Mach number splitting can achieve the optimal spatial convergence order of $(K+1)$. In addition,  
The examples for the Euler system and the shallow water wave system demonstrate the universality of the FVS-DG 
for general hyperbolic conservation laws.

\subsection{Performance of the IS-$L^2$-TVB(D)-minmod Limiter for Scalar Conservation Law}\label{subSec-IS-L2-TVB-scalar}
In this subsection, IS-$L^2$-TVB(D)-minmod Limiter is applied to 1D/2D-scalar equations, 
including inviscid Burgers' equation (Example \ref{example-limiters-1D-Burgers} and Example \ref{example-limiters-2D-Burgers}), non-convex scalar Buckley-Leverett problem (Example \ref{example-limiters-Buckley-Leverett}) and linear variable coefficient transport equation 
(swirling deformation flow in Example \ref{example-2D-linear-Transport}).

\subsubsection{One-dimensional Test Cases for Scalar Equations}
To compare IS-$L^2$-TVB(D)-minmod Limiter with other classical limiters, we choose local Lax-Friedrichs (LLF) flux format for all of them. 
The LLF flux format is given as follows: 
\begin{align*}
  & \hat{f}^{\ LF} = \frac{1}{2} \cdot\left(f\left(U^{L}\right)+f\left(U^{R}\right)-\alpha \cdot\left(U^{R}-U^{L}\right)\right), \\
  & {\color{red}\alpha=\max _{u \in I(u^L,u^R)}\left\{\left|f^{\ \prime}\left(u\right)\right|\right\}},\ I(u^L , u^R) = \left(\min \left\{u^{L},u^{R}\right\},\max \left\{u^{L},u^{R}\right\}\right). \\
\end{align*}
\begin{itemize}
  \item If $\ f(u)=c \cdot u\ $ then $\ f^{\ \prime}(u)=c$, $\alpha=|c|$; 
  \item If $\ f(u)=c \cdot \frac{1}{2}u^2\ $ then $\ f^{\ \prime}(u)=c\cdot u$, $\alpha=|c| \cdot \max \left\{\left|u^{L}\right|,\left|u^{R}\right|\right\}$.  
\end{itemize}
The relationship between Local Lax-Friedrichs splitting method previously mentioned and LLF flux format here should be referred to Appendix \ref{appendix-LF-FVS-LF-flux}.
\begin{example}\label{example-limiters-1D-Burgers}
1D inviscid Burgers' problem with smooth initial conditions and a shock during the evolution. 
\\
Control Eq:  
  $
  u_t + (\frac{1}{2}u^2)_x = 0;
  $\\
  Computational domain:   
  $
  \Omega \times [0,T_{end}] = [0,2\pi] \times [0,2.0];
  $\\
  $\mathrm{I.C.}\ $ 
  $
  \begin{aligned}
      u_0(x)=\sin(x);
  \end{aligned}
  $\\
  Note: a discontinuity occurs at $t=1$ and $x=\pi$ under this given initial condition; \\
  $\mathrm{B.C.}\ ${\color{blue}{periodic boundary conditions; }}
  \\
  True solutions: 
  $
  \begin{aligned}
    u(x,t)=u_0(x^*),\ x^* \text{\ satisfies\ } x^*+u_0(x^*)\cdot t=x
  \end{aligned}
  $;\\
  Flux format: local Lax-Friedrichs flux; \\
  Temporal discretization format: TVD-RK3; \\
  CFL=0.1; \\
  $P^3$ and $P^5$-polynomial approximations and uniquely spaced cells are utilized; 
  \\
  Note that in order to fully demonstrate the influence of various limiters on numerical results, 
  {\color{blue}{no discontinuity indicators are used in this Example \ref{example-limiters-1D-Burgers}}}.
  \par
  Numerical results based on different limiters including IS-TVB-minmod limiter, $L^2$-TVB-minmod limiter, SimpleWENO and OEDG \cite{ref44}
  are demonstrated in Figure \ref{Fig.limiters-Burgers-1D}.
  \begin{figure}[htbp]
  \vspace{0.025cm}
  \begin{center}
    \begin{minipage}{1\linewidth}
      \centerline{\includegraphics[width=1\linewidth]{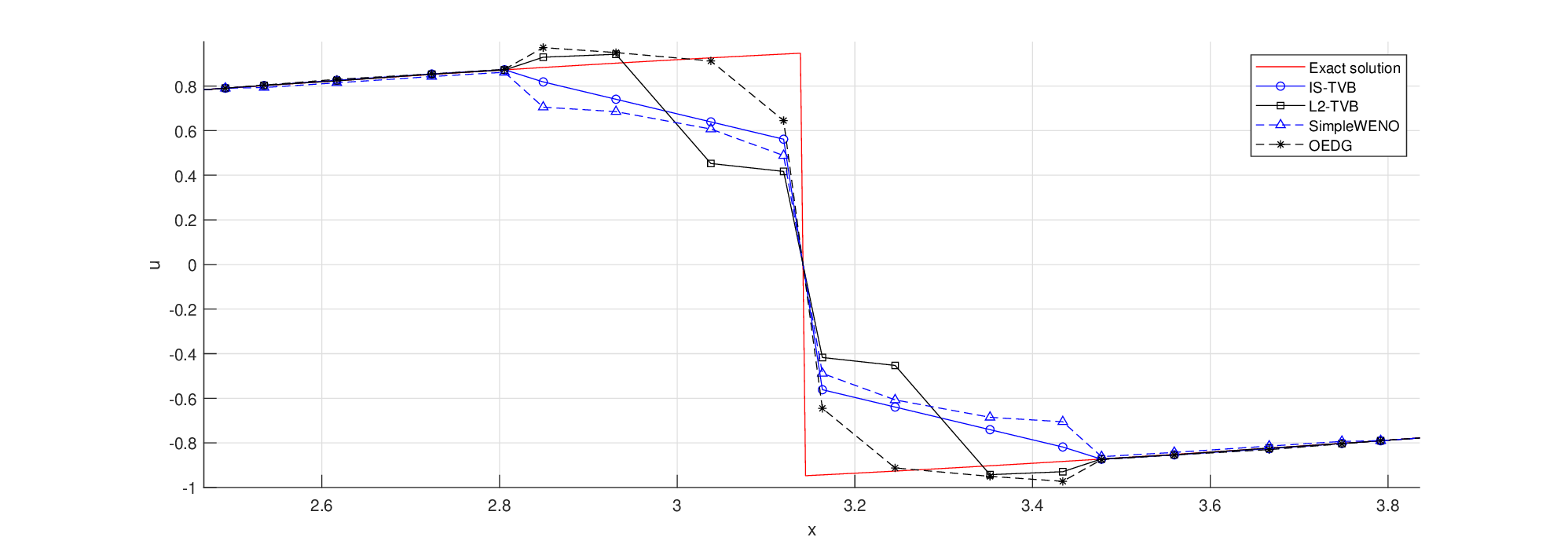}}
    \end{minipage}
    \vfill
    \begin{minipage}{1\linewidth}
      \small
      \centerline{\tiny{(a)\ 20 cells, $P^3$-polynomial approximation with different limiters}}
    \end{minipage}
    \vfill
    \vspace{0.025cm}
    \begin{minipage}{1\linewidth}
      \centerline{\includegraphics[width=1\linewidth]{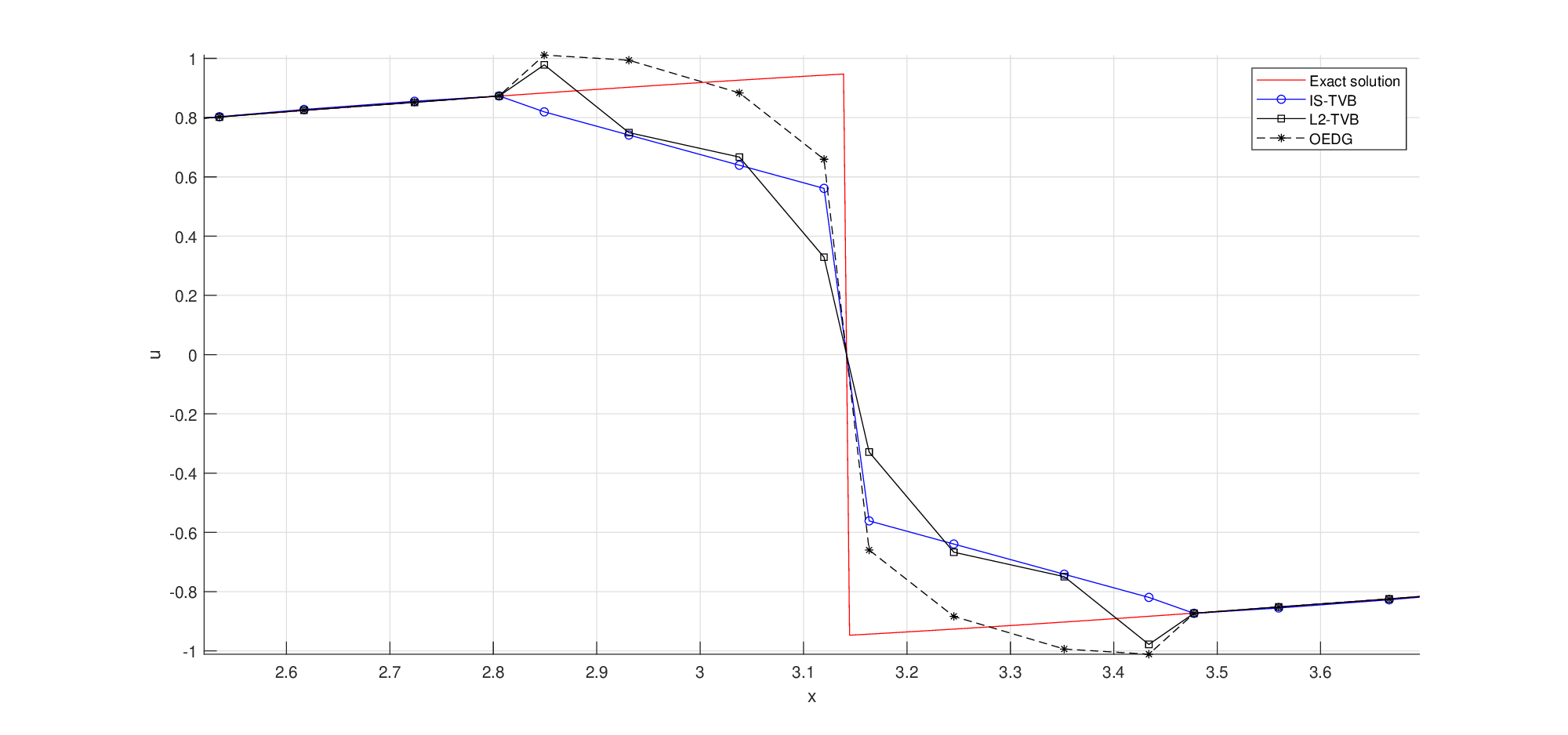}}
    \end{minipage}
    \vfill
    \begin{minipage}{1\linewidth}
      \small
      \centerline{\tiny{(b)\ 20 cells, $P^5$-polynomial approximation with different limiters}}
    \end{minipage}
    \vfill
    \vspace{0.025cm}
    \begin{minipage}{1\linewidth}
      \centerline{\includegraphics[width=1\linewidth]{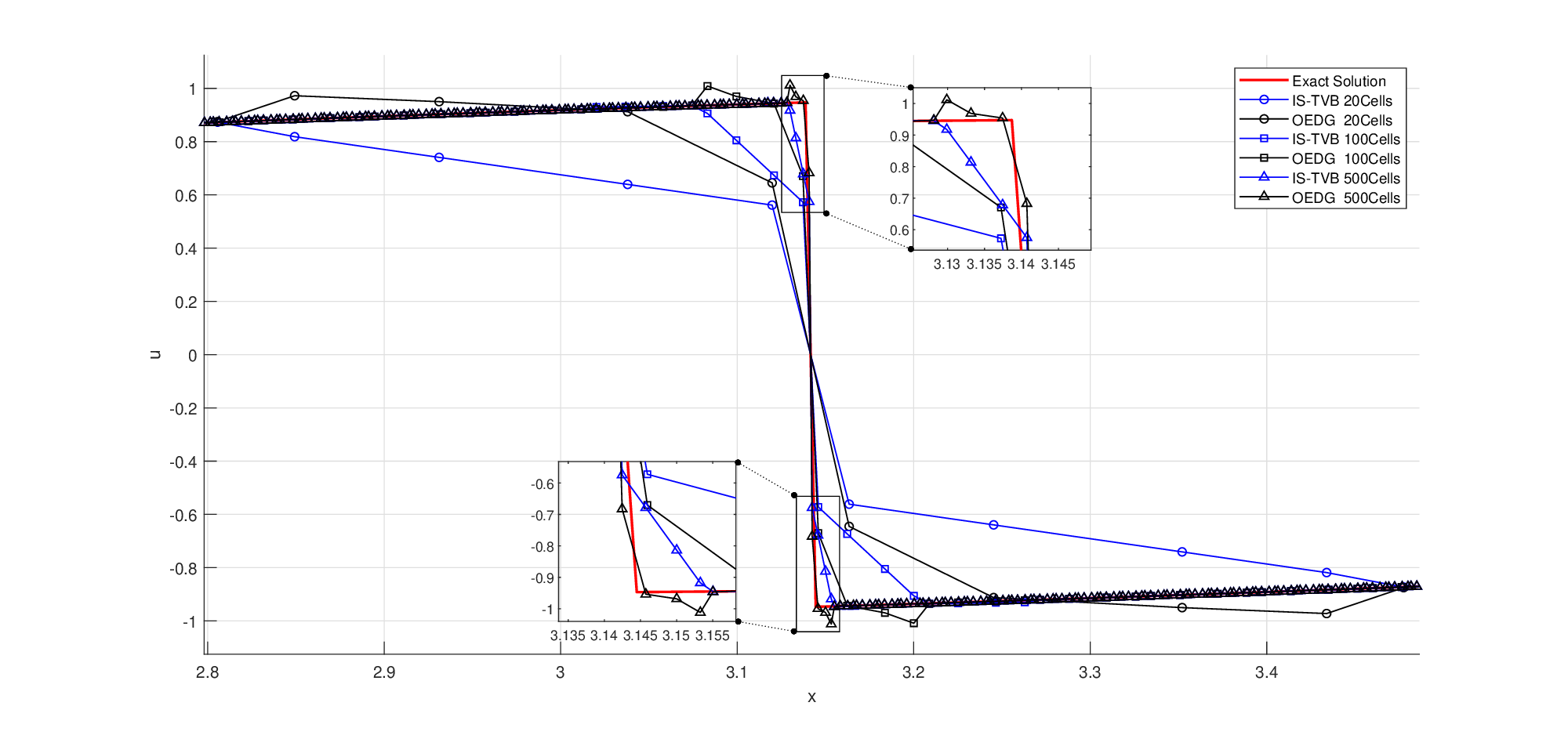}}
    \end{minipage}
    \vfill
    \begin{minipage}{1\linewidth}
      \small
      \centerline{\tiny{(c)\ $P^3$-polynomial approximation with IS-TVB-minmod limiter and OEDG on different meshes}}
    \end{minipage}
  \end{center}
  \caption{\tiny{1D-Burgers' equation $u_t + (\frac{1}{2}u^2)_x = 0$ with initial condition $u_0(x)=\sin(x)$. The simulation is performed up to $t=2.0$. $P^3$ and $P^5$-polynomial approximations and uniquely spaced cells. 
  Numerical results based on different limiters including IS-TVB-minmod limiter, $L^2$-TVB-minmod limiter, SimpleWENO and OEDG are compared with each other. Localized magnification has
  been applied to all sub-figures.}}
  \label{Fig.limiters-Burgers-1D}
\end{figure}
Note: when using $P^5$-polynomial approximation, SimpleWENO doesn't work with 20 cells. (the numerical results are ``NaN'')
\par
According to the Figure \ref{Fig.limiters-Burgers-1D}, we can see that 
when approximating with high-order polynomials, the IS-TVB-minmod limiter exhibits no overshoot near discontinuities, 
whereas both the $L^2$-TVB-minmod limiter and OEDG show overshoot phenomena; 
When approximating with $P^3$-polynomial, SimpleWENO, although not overshoot near discontinuities, has slight oscillations 
and is not as smooth as the IS-TVB-minmod limiter; 
Under the same CFL number (0.1) and the same coarse mesh (20 cells) conditions, when approximating with $P^5$-polynomial, 
SimpleWENO does not work (in fact , even if adding up cells to 500, SimpleWENO still doesn't work), 
while the IS-TVB-minmod limiter remains effective, 
demonstrating that the mesh resolution of the IS-TVB-minmod limiter is superior to that of SimpleWENO 
and IS-TVB-minmod limiter is more suitable for high-order approximation than SimpleWENO. 
Compared with the IS-TVB-minmod limiter, OEDG requires more cells to mitigate overshoot phenomena, in fact, 
OEDG still has overshoot even with 500 cells, while IS-TVB-minmod limiter consistently maintains no overshoot 
during mesh refinement, even on a very coarse mesh (20 cells), there is no overshoot phenomenon with IS-TVB-minmod limiter. 
Therefore, the artificial viscosity and the mesh resolution of the IS-TVB-minmod limiter is also superior to that of OEDG. 
Finally, it needs to be emphasized that as the mesh is refined, the discontinuities captured by IS-TVB-minmod limiter become increasingly clear and sharp, yet there is still no overshoot. 
\end{example}
\begin{example}\label{example-limiters-Buckley-Leverett}
nonlinear non-convex scalar Buckley-Leverett problem. \\
  Control Eq:  
  $
  u_t + (\frac{4u^2}{4u^2+(1-u)^2})_x = 0;
  $\\
  Computational domain:   
  $
  \Omega \times [0,T_{end}] = [-1,1] \times [0,0.4];
  $\\
  $\mathrm{I.C.}\ $ 
  $
  u_0(x)=
  \begin{cases}
      1 \ , \ -\frac{1}{2} \leq x \leq 0,\\
      0 \ , \ \text{elsewhere};
  \end{cases}
  $\\
  $\mathrm{B.C.}\ ${\color{blue}{periodic boundary conditions; }} \\
  Flux format: {\color{blue}{global Lax-Friedrichs flux with fixed $\alpha=2.4$;}} \\
  Temporal discretization format: TVD-RK3; \\
  CFL=0.1; \\
  $P^3$-polynomial approximation and uniquely spaced 80 cells are utilized; \\
  TVB-minmod discontinuity indicator with parameter $M=1$ is utilized for all limiters in this Example \ref{example-limiters-Buckley-Leverett}.
  \par
  Numerical results based on different limiters including IS-TVB-minmod limiter, $L^2$-TVB-minmod limiter, WENO5-JS, SimpleWENO and OEDG 
  are demonstrated in Figure \ref{Fig.limiters-Buckley-Leverett-1D(2)}.
  \begin{figure}[htbp]
  \vspace{0.025cm}
  \begin{center}
    \begin{minipage}{1\linewidth}
      \centerline{\includegraphics[width=1\linewidth]{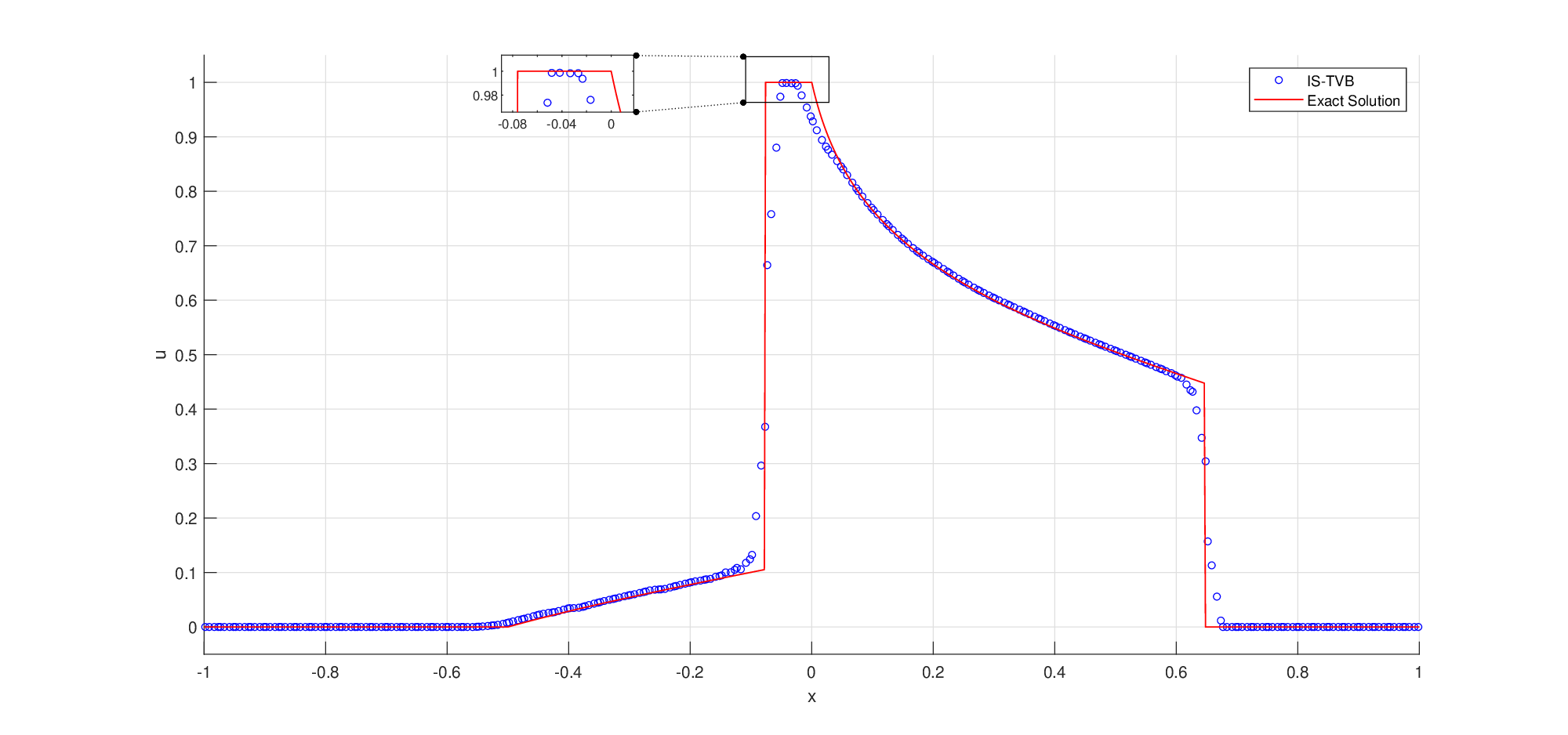}}
    \end{minipage}
    \vfill
    \begin{minipage}{1\linewidth}
      \small
      \centerline{\tiny{(a)\ 80 cells, $P^3$-polynomial approximation with IS-TVB}}
    \end{minipage}
    \vfill
    \begin{minipage}{1\linewidth}
      \centerline{\includegraphics[width=1\linewidth]{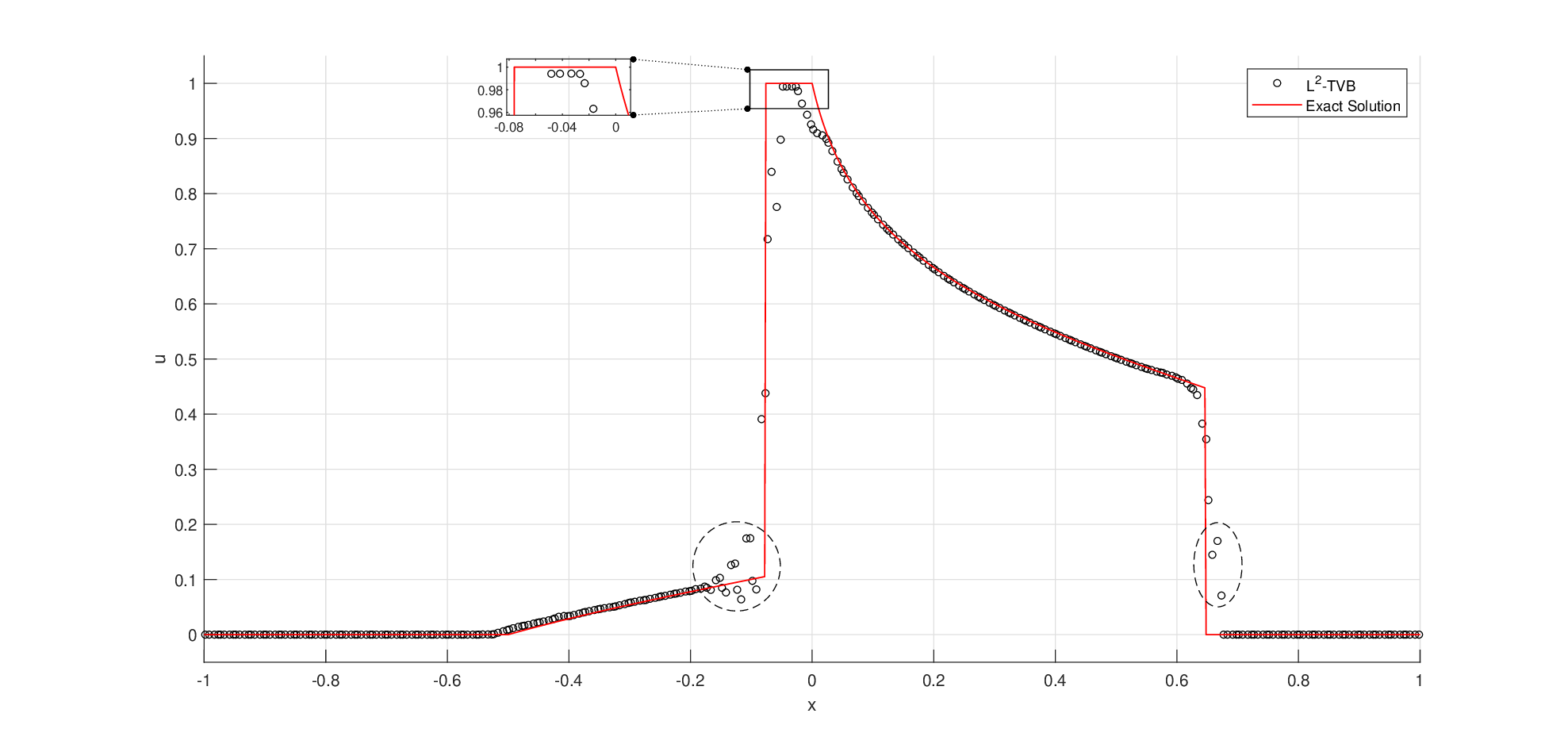}}
    \end{minipage}
    \vfill
    \begin{minipage}{1\linewidth}
      \centerline{\tiny{(b)\ 80 cells, $P^3$-polynomial approximation with $L^2$-TVB}}
    \end{minipage}
    \vfill
    \begin{minipage}{1\linewidth}
      \small
      \centerline{\includegraphics[width=1\linewidth]{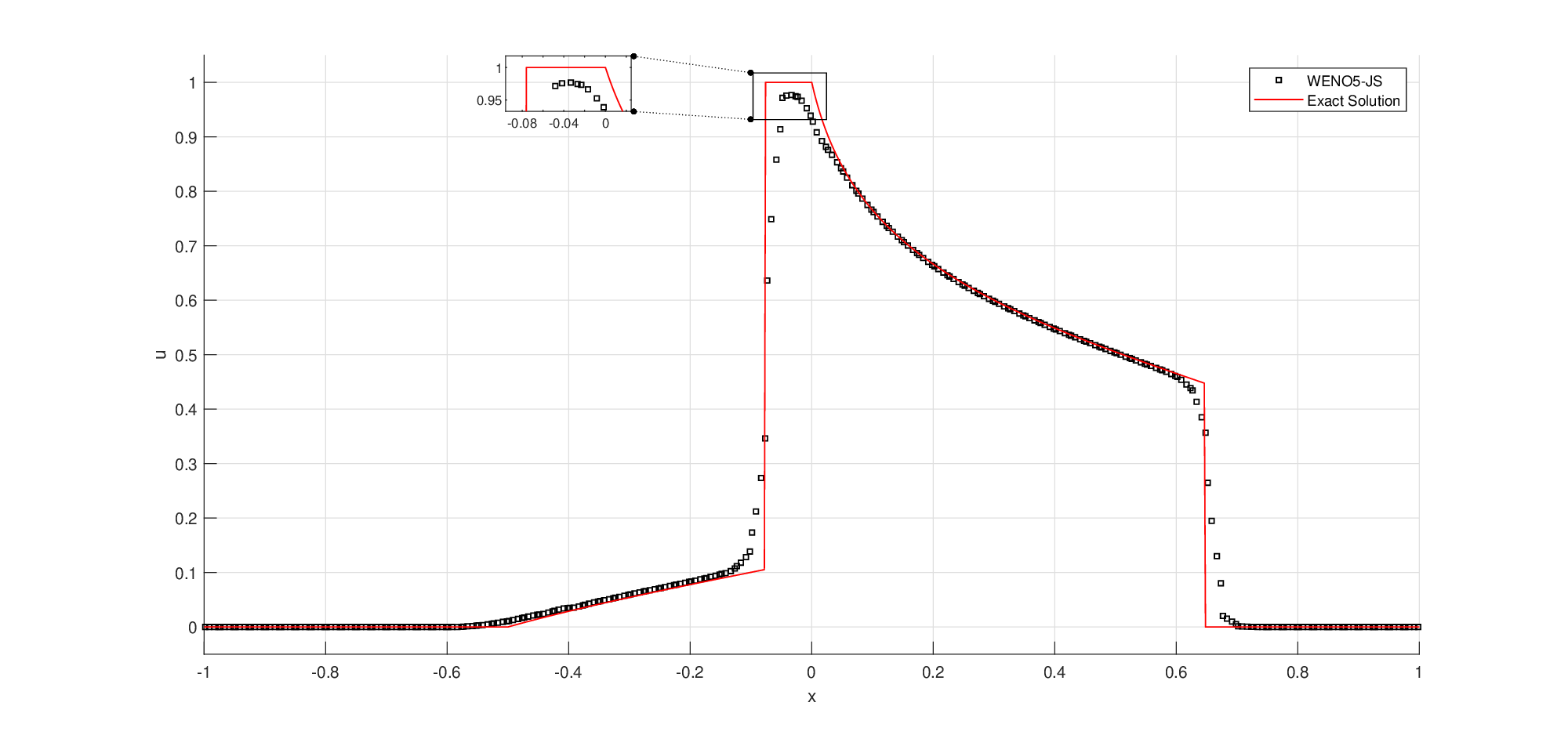}}
    \end{minipage}
    \vfill
    \begin{minipage}{1\linewidth}
      \centerline{\tiny{(c)\ 80 cells, $P^3$-polynomial approximation with WENO5-JS}}
    \end{minipage}
  \end{center}
  \label{Fig.limiters-Buckley-Leverett-1D(1)}
\end{figure}
\begin{figure}[htbp]
  \vspace{0.025cm}
  \begin{center}
    \begin{minipage}{0.49\linewidth}
      \centerline{\includegraphics[width=1\linewidth]{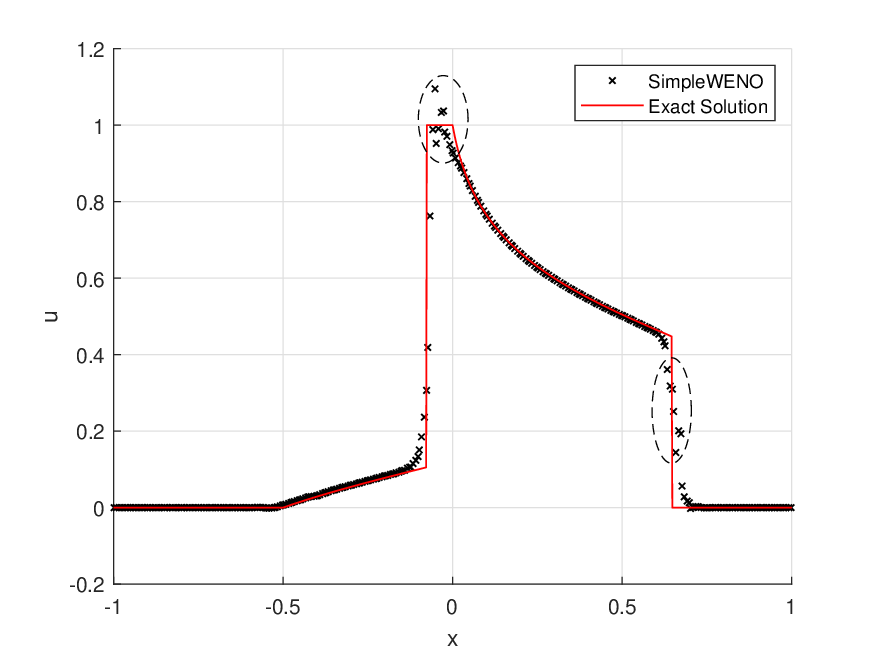}}
    \end{minipage}
    \hfill
    \begin{minipage}{0.49\linewidth}
      \centerline{\includegraphics[width=1\linewidth]{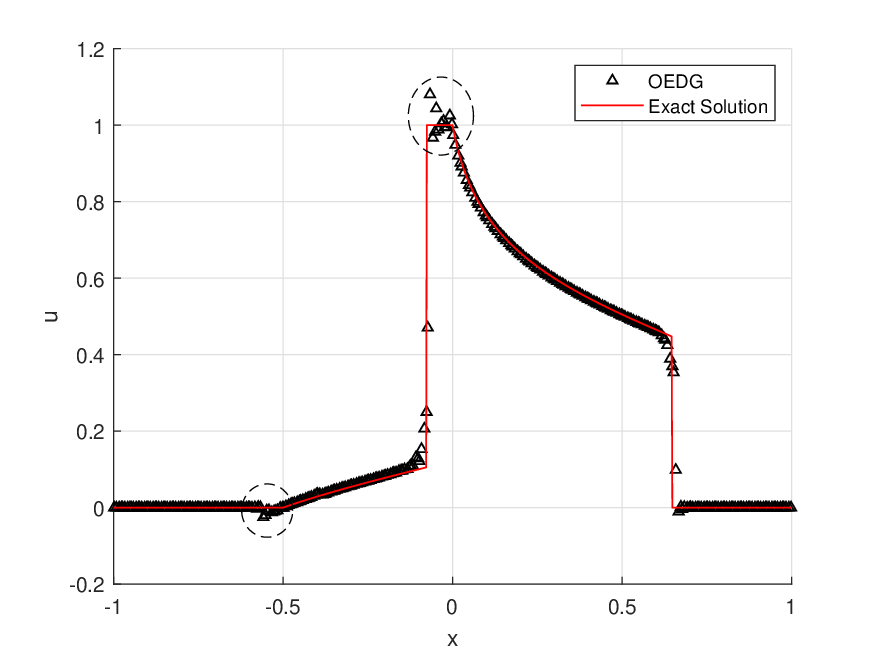}}
    \end{minipage}
    \vfill
    \begin{minipage}{0.49\linewidth}
      \small
      \centerline{\tiny{(d)\ 80 cells, $P^3$-polynomial approximation with SimpleWENO}}
    \end{minipage}
    \hfill
    \begin{minipage}{0.49\linewidth}
      \centerline{\tiny{(e)\ 80 cells, $P^3$-polynomial approximation with OEDG}}
    \end{minipage}
  \end{center}
  \caption{\tiny{1D-Buckley-Leverett problem $u_t + (\frac{4u^2}{4u^2+(1-u)^2})_x = 0$ with initial condition 
  $u=1$ when $-\frac{1}{2} \leq x \leq 0$ and $u=0$ elsewhere. The simulation is performed up to time $t=0.4$. $P^3$-polynomial approximation and uniquely spaced 80 cells. 
  Numerical results based on different limiters including IS-TVB-minmod limiter, $L^2$-TVB-minmod limiter, WENO5-JS, SimpleWENO and OEDG are compared with each other. 
  }}
  \label{Fig.limiters-Buckley-Leverett-1D(2)}
\end{figure}
\par
The IS-TVB exhibits no overshoot behavior, while the SimpleWENO and OEDG exhibit noticeable overshoot phenomena; 
the IS-TVB has no numerical spurious oscillations, 
whereas the $L^2$-TVB shows noticeable numerical oscillations within the range $-0.2 \leq x \leq 0$; 
although the WENO5-JS does not have numerical oscillations, 
its dissipation is significantly larger 
and its precision is inferior to that of the IS-TVB within the range $-0.06 \leq x \leq -0.02$.
\end{example}

\subsubsection{Two-dimensional Test Cases for Scalar Equations}
In this subsection, in addition to LLF flux, we used {\color{red}{``scalar Steger-Warming flux'' inspired by Steger-Warming splitting for systems.}}
More details about ``scalar Steger-Warming flux'', please refer to Appendix \ref{appendix-2D-scalar-SW-flux}. 

\begin{example}\label{example-limiters-2D-Burgers}
2D nonlinear inviscid Burgers' problem. \\
Control Eq:  
$
U_t + (\frac{1}{2}U)_x +  (\frac{1}{2}U)_y = 0;
$
\\ \hfill \\
$\bullet\ $2D-Burgers' problem with smooth initial conditions and a shock during the evolution\\
Computational domain:   
$
\Omega \times [0,T_{end}] = \left\{[0,4] \times [0,4] \right\} \times [0,\frac{1.5}{\pi}];
$\\
$\mathrm{I.C.}\ $ 
$
U_0(x,y)=
\sin\left(\frac{\pi}{2}(x+y)\right)
$\\
$\mathrm{B.C.}\ ${\color{blue}{periodic boundary conditions; }} \\
Flux format: local Lax-Friedrichs flux; \\
Temporal discretization format: TVD-RK3; \\
CFL=0.1; \\
$P^3$-polynomial approximation and uniform rectangular 50 $\times$ 50 cells are utilized; \\
TVB-minmod discontinuity indicator with parameter $M=1$ is utilized in this test case.
\par
Numerical results based on IS-TVB-minmod limiter ($\omega_{IS}=1,\ \omega_{L^2}=0$) are demonstrated in Figure \ref{Fig.limiters-2D-Burgers-1}.
\begin{figure}[htbp]
  \begin{center}
    \begin{minipage}{0.49\linewidth}
      \centerline{\includegraphics[width=1\linewidth]{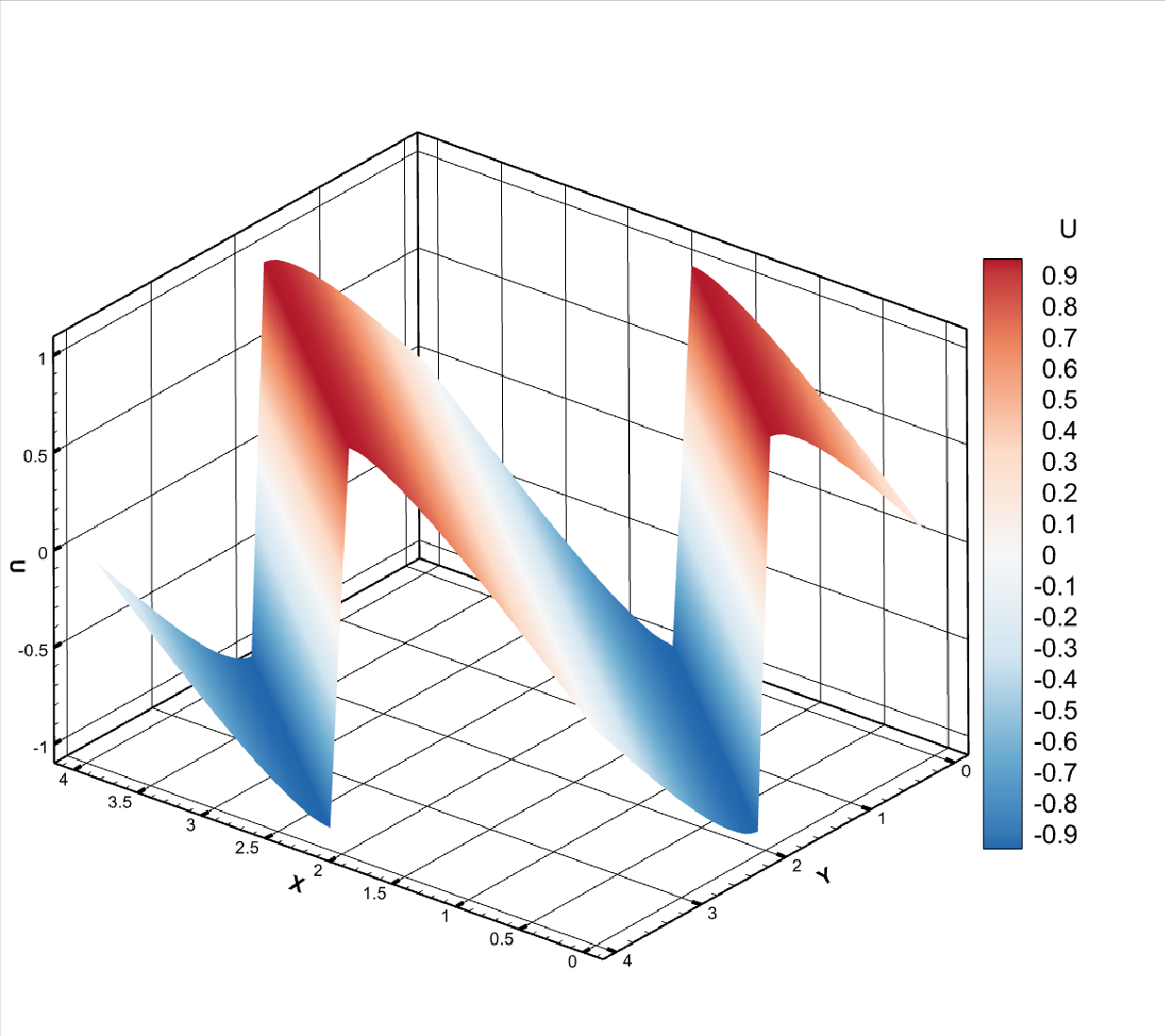}}
    \end{minipage}
    \hfill
    \begin{minipage}{0.49\linewidth}
      \centerline{\includegraphics[width=1\linewidth]{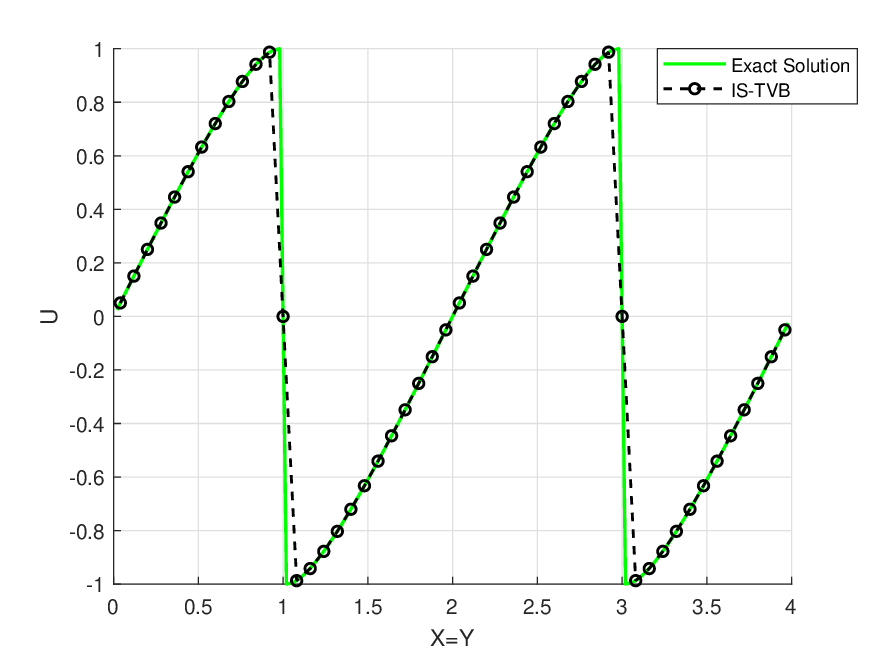}}
    \end{minipage}
    \vfill
    \begin{minipage}{0.49\linewidth}
      \centerline{\tiny{(a)\ 3D view for IS-TVB numerical solution}}
    \end{minipage}
    \hfill
    \begin{minipage}{0.49\linewidth}
      \centerline{\tiny{(b)\  IS-TVB solution vs. exact solution along the cut $y=x$}}
    \end{minipage}
  \end{center}
\caption{\tiny{2D-Burgers' problem $U_t + (\frac{1}{2}U)_x +  (\frac{1}{2}U)_y = 0$ with smooth initial conditions $U_0(x,y)=\sin\left(\frac{\pi}{2}(x+y)\right)$ while a shock during the evolution. The simulation is performed up to time $t=\frac{1.5}{\pi}$ when a shock has already appeared. 
$P^3$-polynomial approximation and uniquely spaced rectangular 50 $\times$ 50 cells. 
Discontinuity is captured based on IS-TVB-minmod limiter ($\omega_{IS}=1,\ \omega_{L^2}=0$). 
}}
\label{Fig.limiters-2D-Burgers-1}
\end{figure}
\par
\hfill \\
$\bullet\ $2D-Burgers' problem with discontinuous initial conditions\\
Computational domain:   
$
\Omega \times [0,T_{end}] = \left\{[0,0.1] \times [0,0.1] \right\} \times [0,0.05];
$\\
$\mathrm{I.C.}\ $ 
$
u_0(x, y)= \begin{cases}0.5, & x<0.05, y<0.05 \\ 0.8, & x>0.05, y<0.05 \\ -1, & x>0.05, y>0.05 \\ -0.2, & x<0.05, y>0.05\end{cases}
$\\
$\mathrm{B.C.}\ ${\color{blue}{free boundary conditions are imposed on all edges of $\Omega$; }} \\
Flux format: Steger-Warming flux (please refer to Appendix \ref{appendix-2D-scalar-SW-flux}); \\
Temporal discretization format: TVD-RK3; \\
CFL=0.1; \\
$P^3$-polynomial approximation and uniform rectangular 50 $\times$ 50 cells are utilized; \\
TVB-minmod discontinuity indicator with parameter $M=1$ is utilized in this test case.
\par
Numerical results based on IS-TVB-minmod limiter ($\omega_{IS}=1,\ \omega_{L^2}=0$) are demonstrated in Figure \ref{Fig.limiters-2D-Burgers-2}.
\begin{figure}[htbp]
\begin{center}
  \begin{minipage}{0.3\linewidth}
    \centerline{\includegraphics[width=1\linewidth]{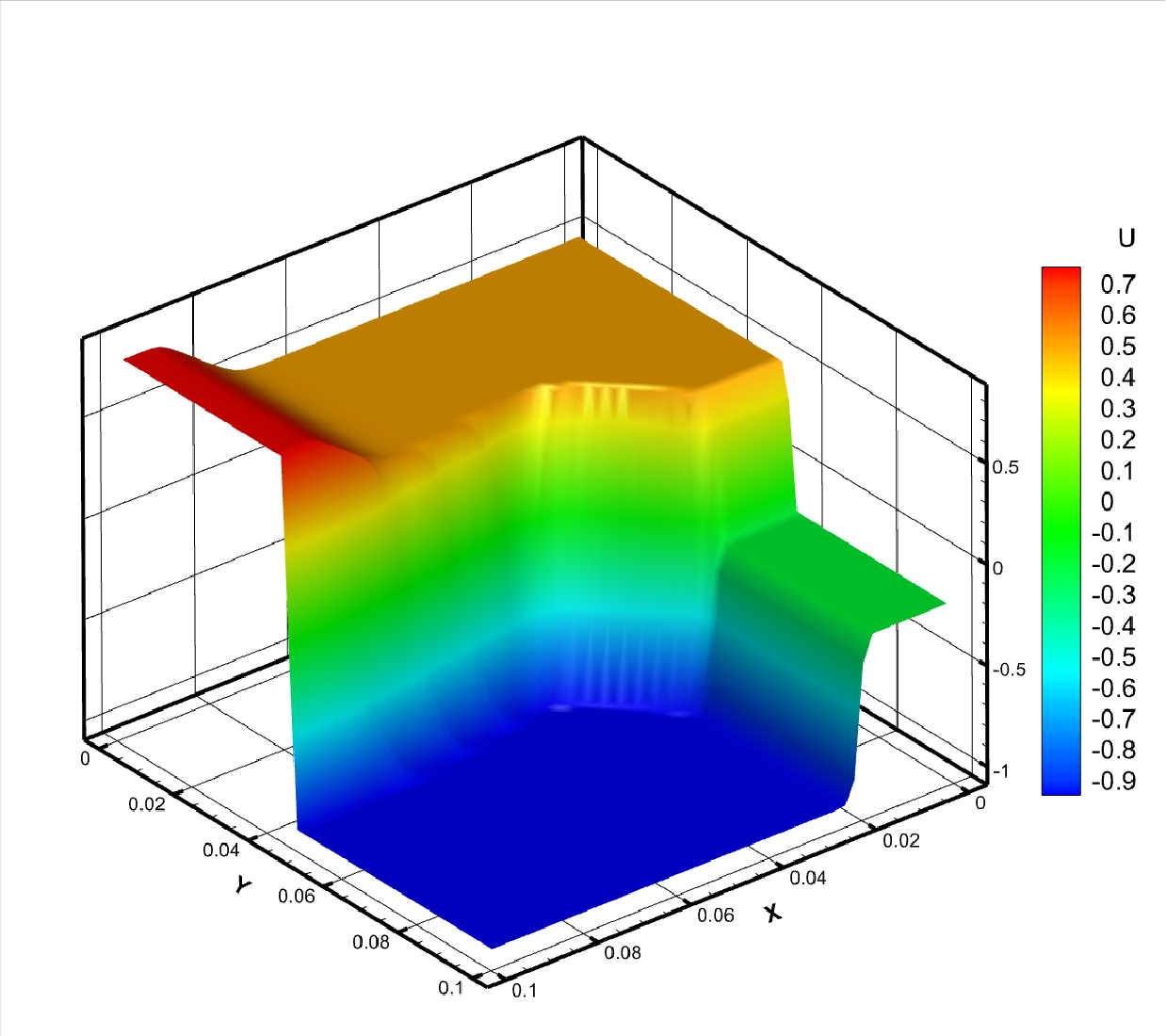}}
  \end{minipage}
  \hfill
  \begin{minipage}{0.3\linewidth}
    \centerline{\includegraphics[width=1\linewidth]{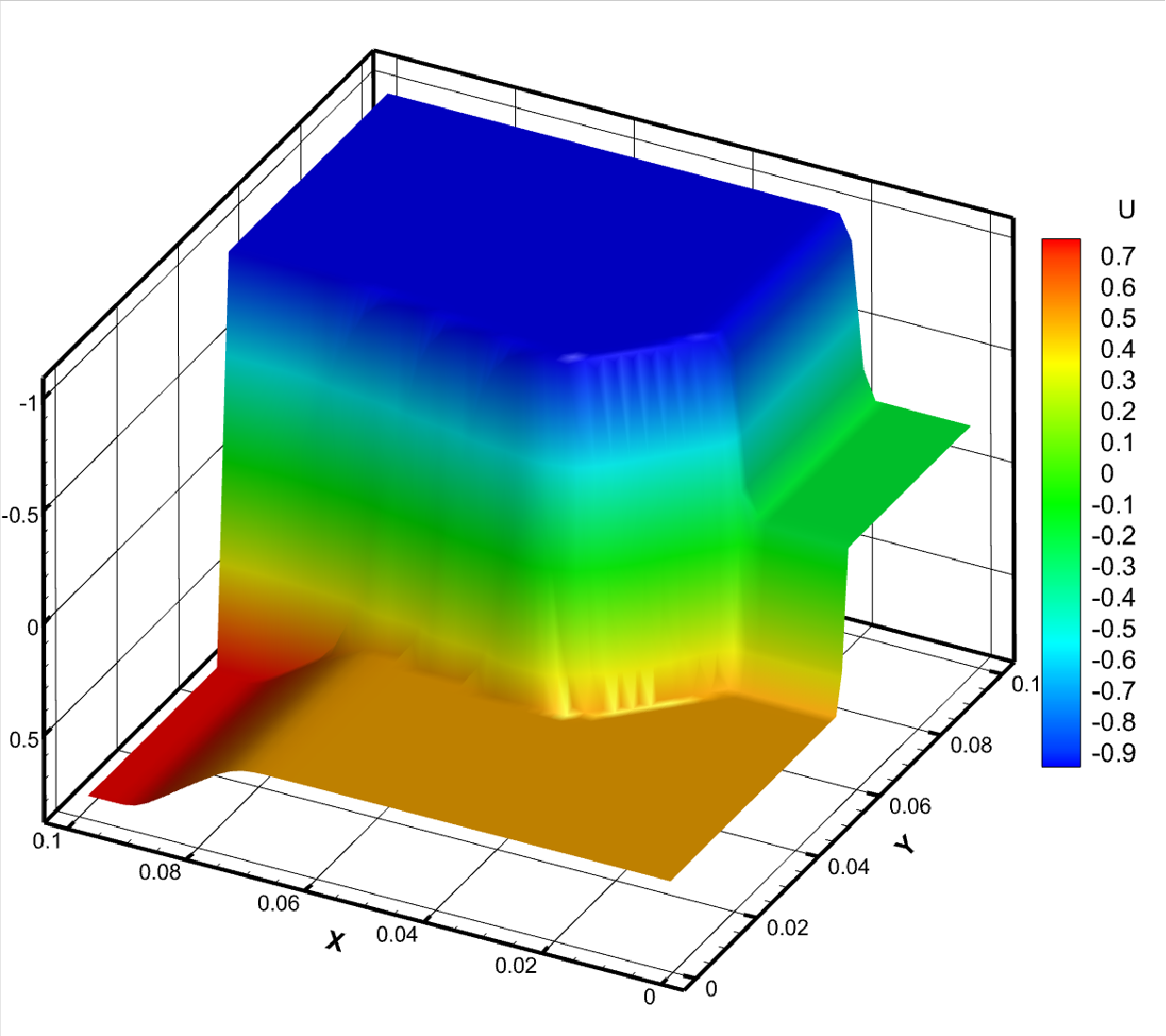}}
  \end{minipage}
  \hfill
  \begin{minipage}{0.3\linewidth}
    \centerline{\includegraphics[width=1\linewidth]{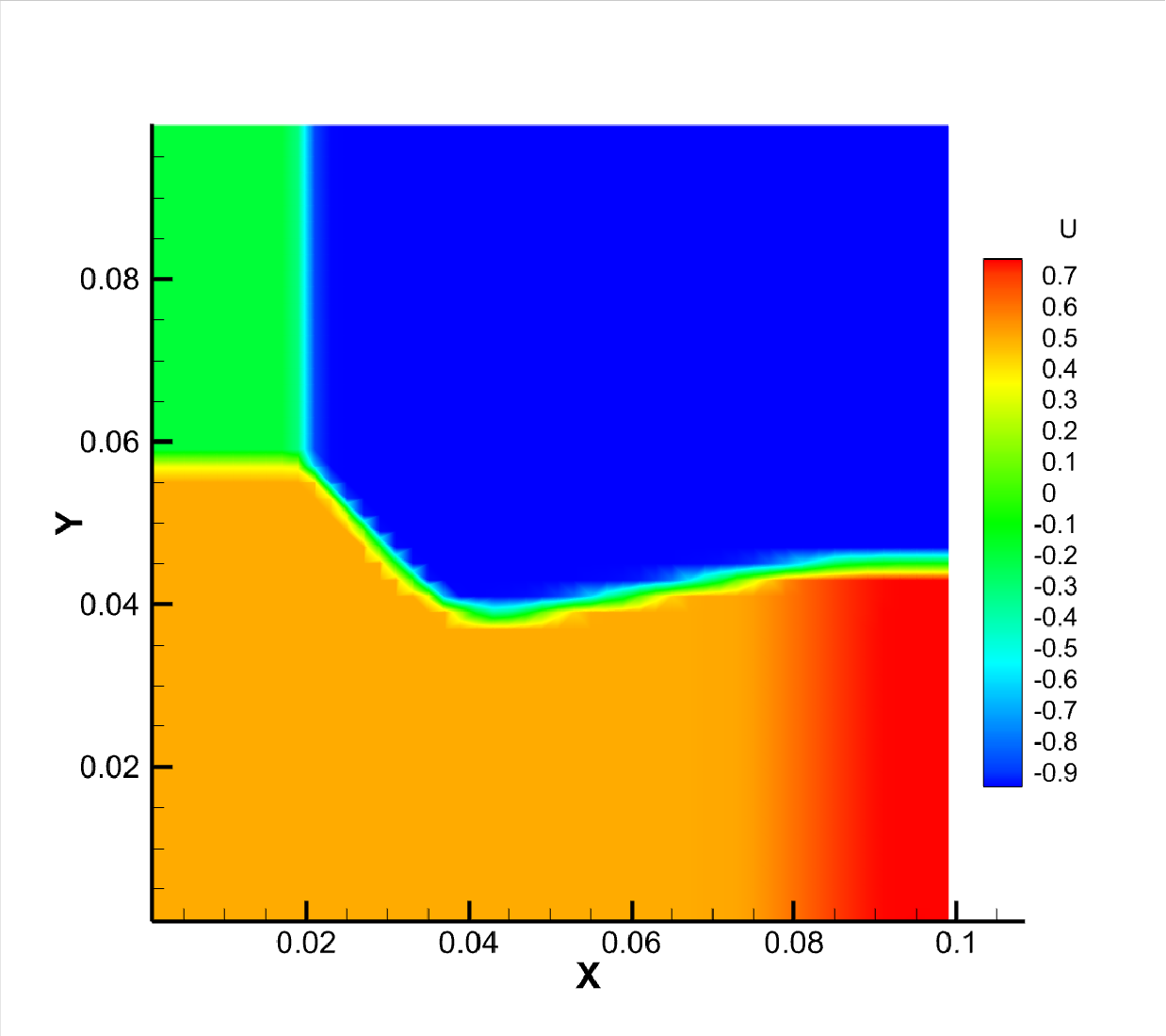}}
  \end{minipage}
  \vfill
  \begin{minipage}{0.3\linewidth}
    \centerline{\tiny{(a)\  3D view}}
  \end{minipage}
  \hfill
  \begin{minipage}{0.3\linewidth}
    \centerline{\tiny{(b)\  3D view (handstand)}}
  \end{minipage}
  \hfill
  \begin{minipage}{0.3\linewidth}
    \centerline{\tiny{(c)\  2D top-down view}}
  \end{minipage}
\end{center}
\caption{\tiny{2D-Burgers' problem with discontinuous initial conditions.  
$U_t + (\frac{1}{2}U)_x +  (\frac{1}{2}U)_y = 0$. 
At $t=0$, $U_0(x,y)=0.5$ when $x<0.05 \text{\ and\ } y<0.05$, $U_0(x,y)=0.8$ when $x>0.05 \text{\ and\ } y<0.05$, 
$U_0(x,y)=-1$ when $x>0.05 \text{\ and\ } y>0.05$, $U_0(x,y)=-0.2$ when $x<0.05 \text{\ and\ } y>0.05$. 
The simulation is performed up to time $t=0.05$. 
$P^3$-polynomial approximation and uniform rectangular 50 $\times$ 50 cells. 
Discontinuity is dealed with IS-TVB-minmod limiter ($\omega_{IS}=1,\ \omega_{L^2}=0$). 
}}
\label{Fig.limiters-2D-Burgers-2}
\end{figure}
\end{example}

\begin{example}\label{example-2D-linear-Transport}
2D linear advection equation. 
\begin{figure}[htbp]
  \begin{center}
    \begin{minipage}{0.49\linewidth}
      \centerline{\includegraphics[width=1\linewidth]{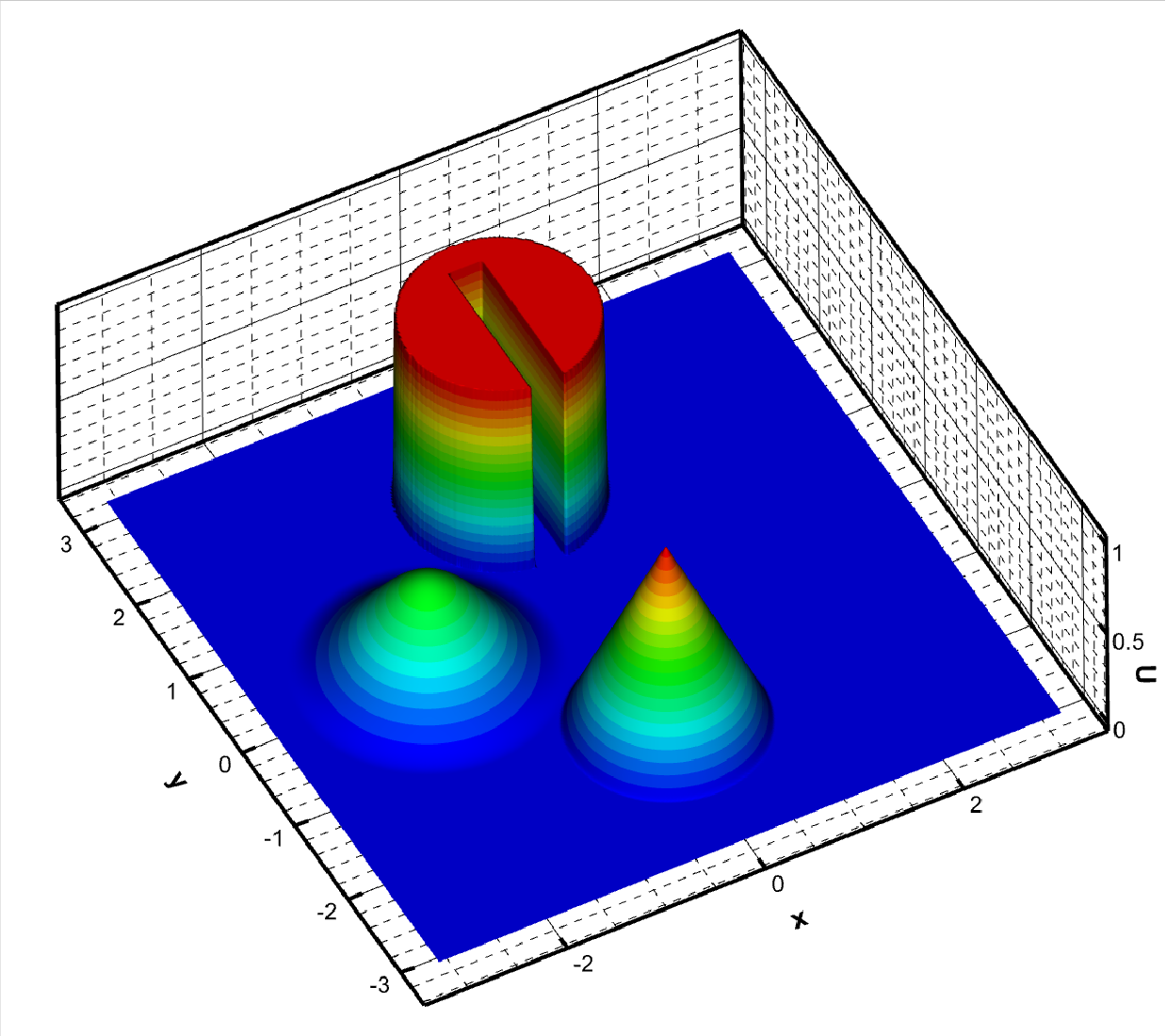}}
    \end{minipage}
    \hfill
    \begin{minipage}{0.49\linewidth}
      \centerline{\includegraphics[width=1\linewidth]{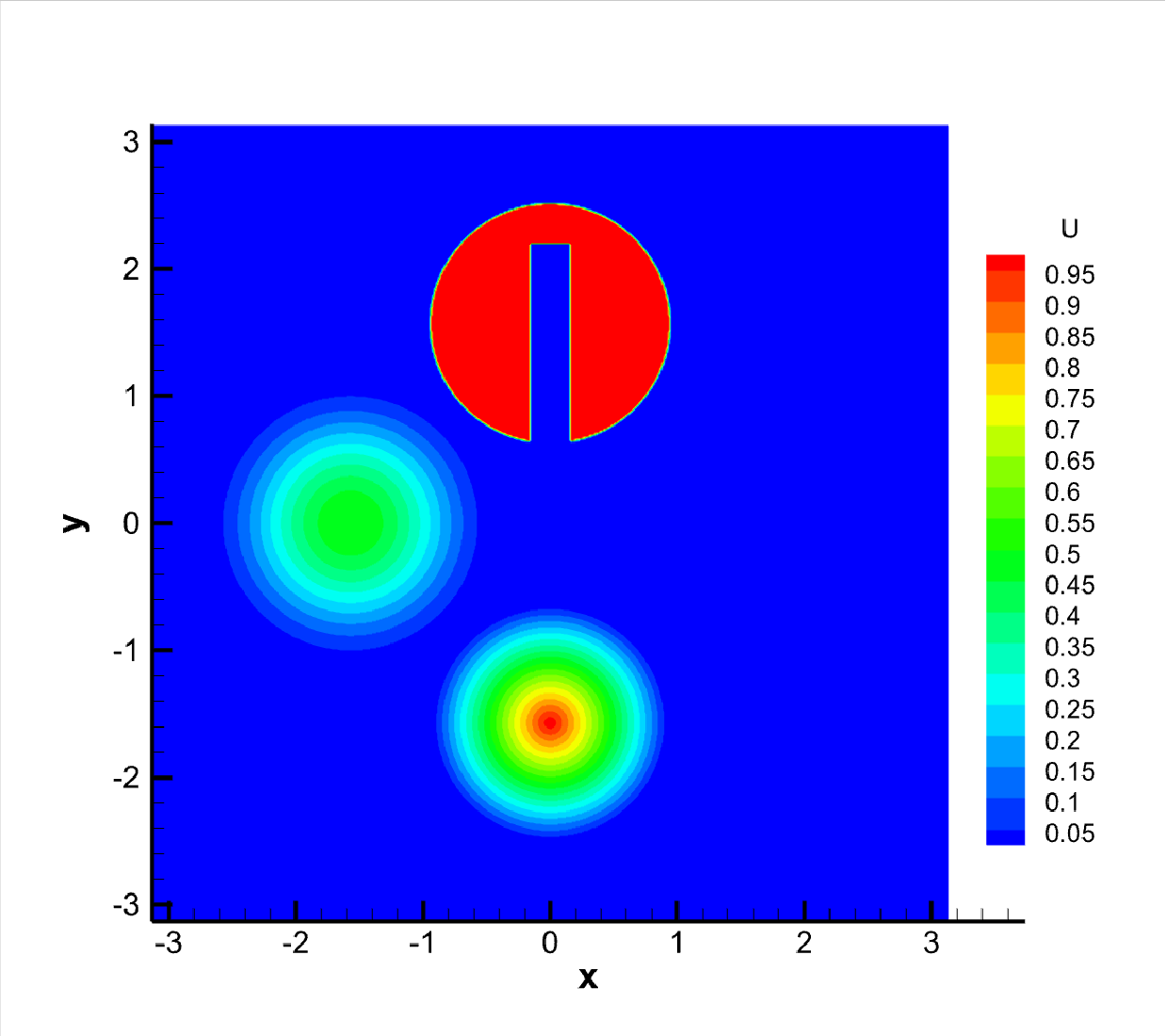}}
    \end{minipage}
    \vfill
    \begin{minipage}{0.49\linewidth}
      \small
      \centerline{\tiny{(a)\ 3D view}}
    \end{minipage}
    \hfill
    \begin{minipage}{0.49\linewidth}
      \small
      \centerline{\tiny{(b)\ 2D top-down view}}
    \end{minipage}
  \end{center}
  \caption{\tiny{initial profile for rigid body rotation and swirling deformation flow}}
  \label{Fig-initial data}
\end{figure}
\\
$\bullet\ $Swirling deformation flow \cite{ref45}.
\\
Control Eqs:
$
U_t-\left(\cos ^2\left(\frac{x}{2}\right) \sin (y) g(t) U\right)_x+\left(\sin (x) \cos ^2\left(\frac{y}{2}\right) g(t) U\right)_y=0, 
$ 
\\
where $g(t)=2 \pi \cos (\pi t / T) ,\quad {\color{red}{T=0.75}}$; 
\\
Computational Domain: 
$
\Omega \times [0,T_{end}] = \left\{ [-\pi,\pi] \times [-\pi,\pi] \right\} \times [0,{\color{red}{0.75}}]; 
$
\\
$\mathrm{I.C.}\ $ 
initial condition plotted in Figure.\ref{Fig-initial data} which
consists of a slotted disk, a cone as well as a smooth hump, similar to the one used in \cite{ref45};
\\
$\mathrm{B.C.}\ ${\color{blue}{periodic boundary conditions; }} \\
Flux format: Steger-Warming flux (please refer to Appendix \ref{appendix-2D-scalar-SW-flux}); \\
Temporal discretization format: TVD-RK3; \\
CFL=0.1; \\
$P^3$-polynomial approximation and uniform rectangular 120 $\times$ 120 cells are utilized; \\
TVB-minmod discontinuity indicator with parameter $M=1$ is utilized in this test case.
\par
Numerical results based on IS-$L^2$-TVB-minmod limiter ($\omega_{IS}=0.75,\ \omega_{L^2}=0.25$) 
are demonstrated in Figure \ref{Fig.limiters-2D-Swirling deformation flow-T/2} and Figure \ref{Fig.limiters-2D-Swirling deformation flow-T}.
\begin{figure}[htbp]
  \begin{center}
    \begin{minipage}{0.49\linewidth}
      \centerline{\includegraphics[width=1\linewidth]{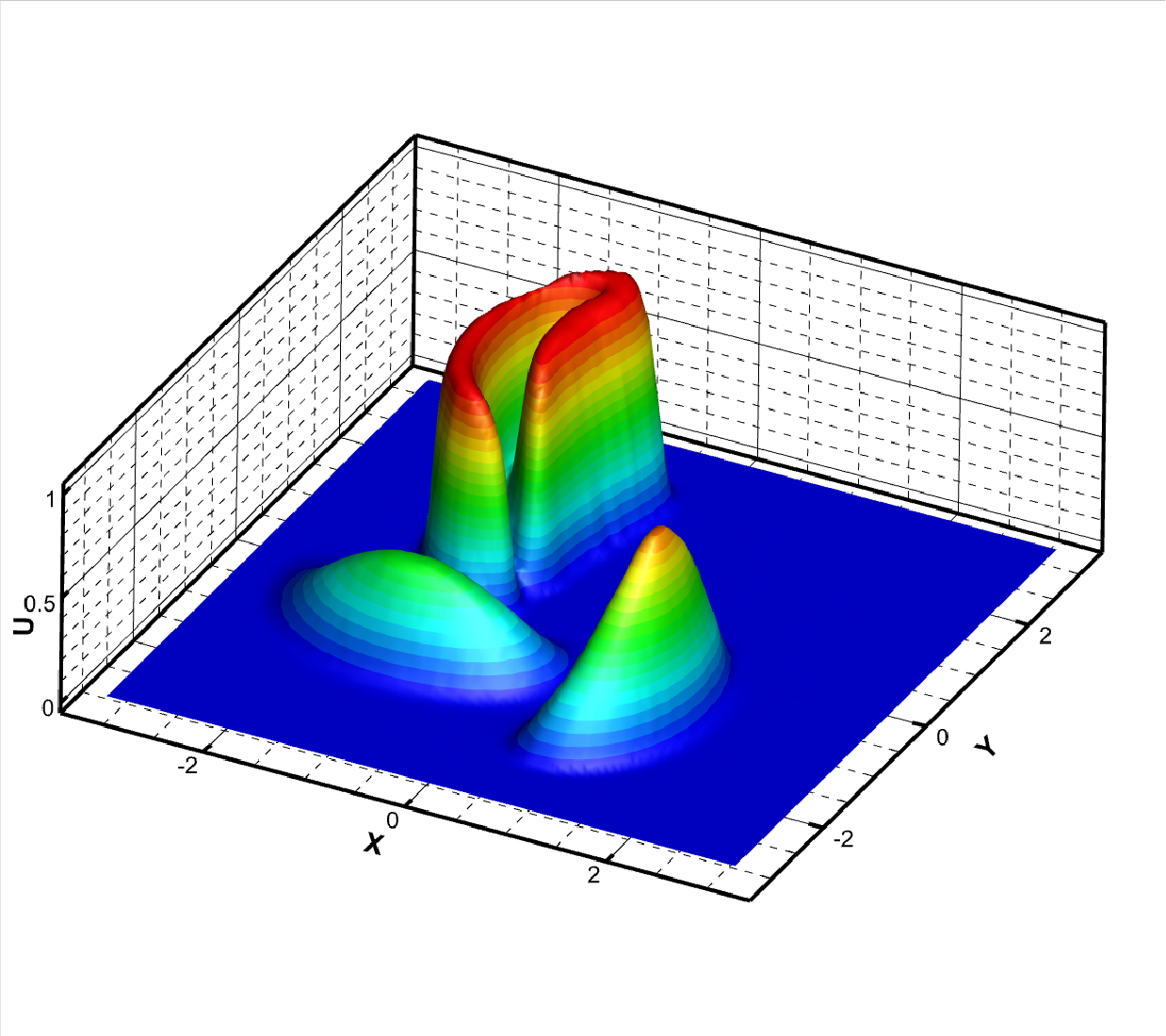}}
    \end{minipage}
    \hfill
    \begin{minipage}{0.49\linewidth}
      \centerline{\includegraphics[width=1\linewidth]{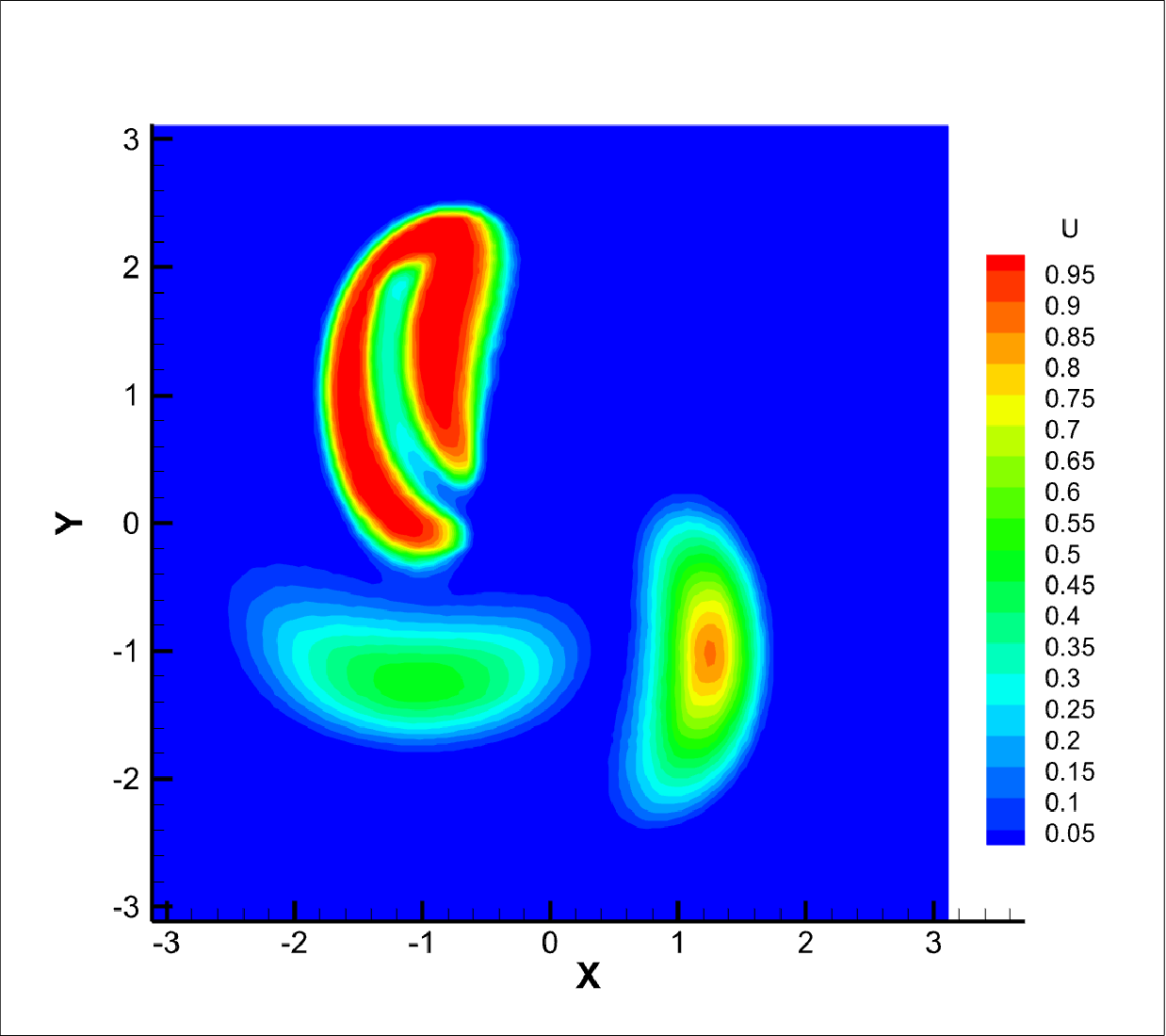}}
    \end{minipage}
    \vfill
    \begin{minipage}{0.49\linewidth}
      \small
      \centerline{\tiny{(a)\ 3D view}}
    \end{minipage}
    \hfill
    \begin{minipage}{0.49\linewidth}
      \centerline{\tiny{(b)\ 2D top-down view}}
    \end{minipage}
  \end{center}
  \caption{\tiny{2D-swirling deformation flow 
  $U_t-\left(\cos ^2\left(\frac{x}{2}\right) \sin (y) g(t) U\right)_x+\left(\sin (x) \cos ^2\left(\frac{y}{2}\right) g(t) U\right)_y=0,\ g(t)=2 \pi \cos (\pi t / T)$  
  with initial data in Fig.\ref{Fig-initial data}:  
  {\color{red}{the period $T=0.75$ and the final integration time 0.375.}} $P^3$-polynomial approximation with IS-$L^2$-TVB-minmod limiter ($\omega_{IS}=0.75,\ \omega_{L^2}=0.25$). 
  Uniform rectangular 120 $\times$ 120 cells. 
  }}
  \label{Fig.limiters-2D-Swirling deformation flow-T/2}
\end{figure} 
\begin{figure}[htbp]
  \begin{center}
    \begin{minipage}{0.49\linewidth}
      \centerline{\includegraphics[width=1\linewidth]{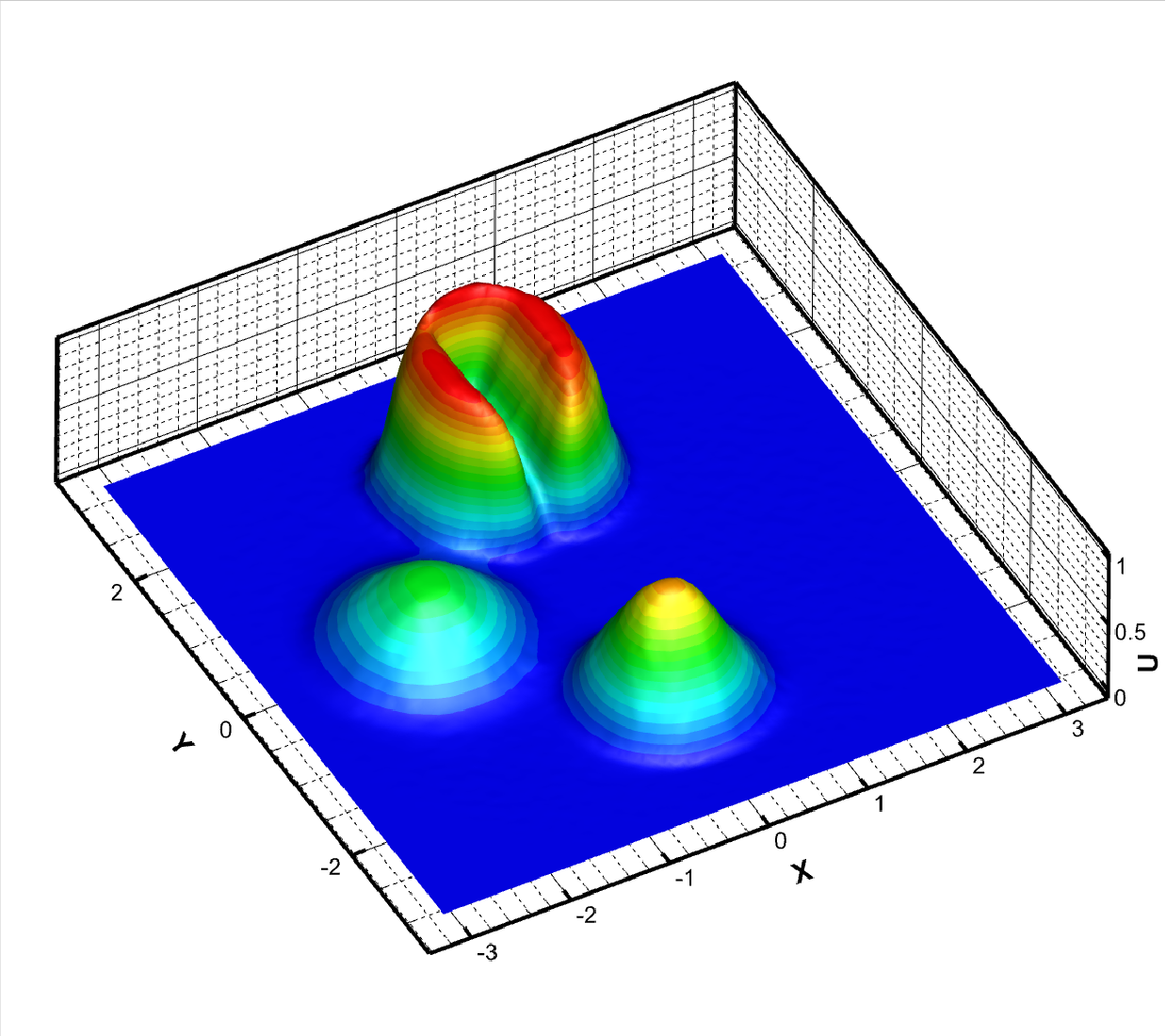}}
    \end{minipage}
    \hfill
    \begin{minipage}{0.49\linewidth}
      \centerline{\includegraphics[width=1\linewidth]{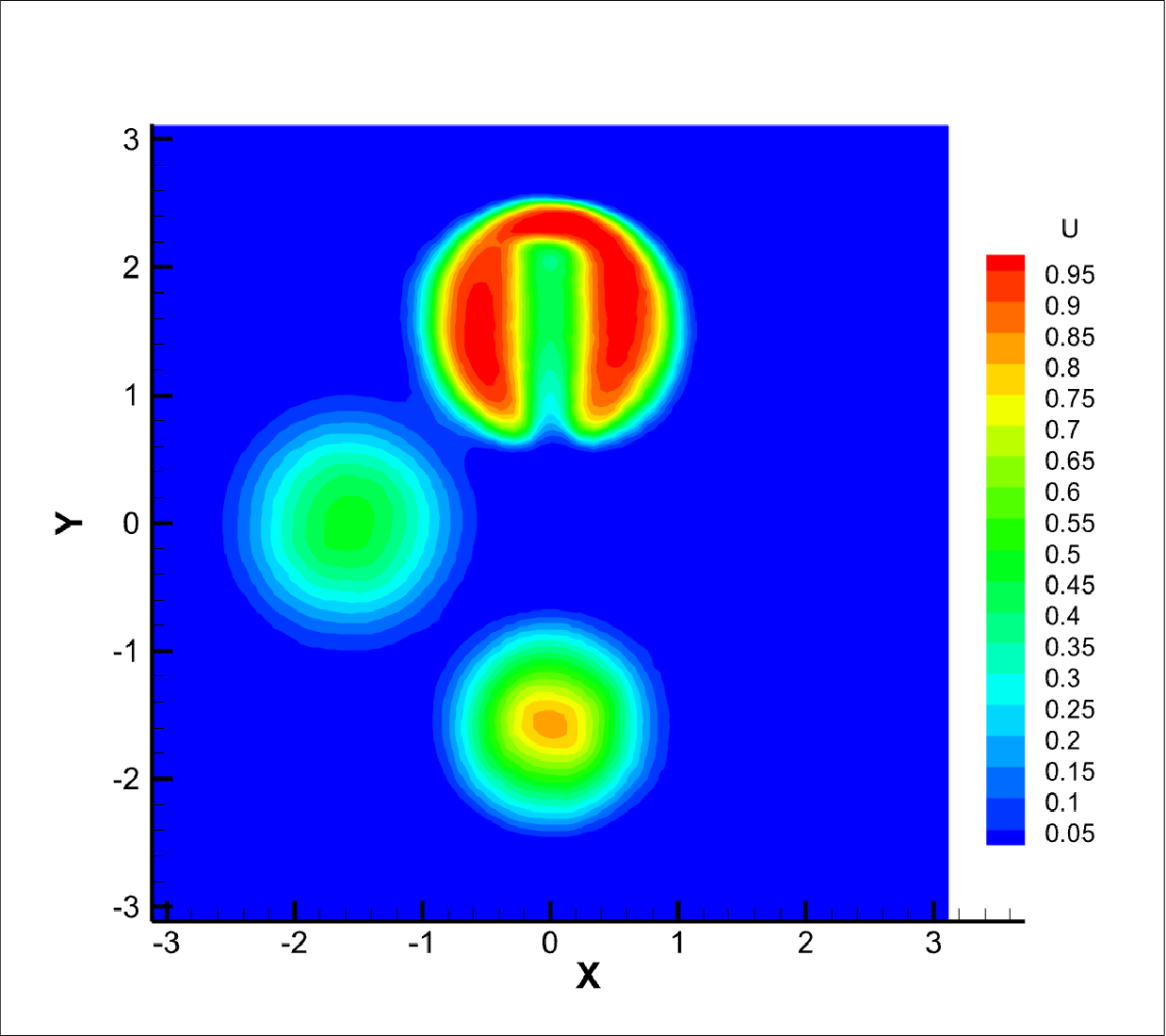}}
    \end{minipage}
    \vfill
    \begin{minipage}{0.49\linewidth}
      \small
      \centerline{\tiny{(a)\ 3D view}}
    \end{minipage}
    \hfill
    \begin{minipage}{0.49\linewidth}
      \centerline{\tiny{(b)\ 2D top-down view}}
    \end{minipage}
  \end{center}
  \caption{\tiny{2D-swirling deformation flow 
  $U_t-\left(\cos ^2\left(\frac{x}{2}\right) \sin (y) g(t) U\right)_x+\left(\sin (x) \cos ^2\left(\frac{y}{2}\right) g(t) U\right)_y=0,\ g(t)=2 \pi \cos (\pi t / T)$ 
  with initial data in Fig.\ref{Fig-initial data}:  
  {\color{red}{the period $T=0.75$ and the final integration time 0.75.}} $P^3$-polynomial approximation with IS-$L^2$-TVB-minmod limiter ($\omega_{IS}=0.75,\ \omega_{L^2}=0.25$). 
  Uniform rectangular 120 $\times$ 120 cells. 
  }}
  \label{Fig.limiters-2D-Swirling deformation flow-T}
\end{figure} 
\end{example}

\subsection{Solving the Riemann Problems for Hyperbolic Conservative Systems Using FVS-DG with IS-$L^2$-TVB(D)-minmod Limiter}
\subsubsection{One-dimensional Test Cases for Hyperbolic Systems}
In this subsection, several classical one-dimensional Riemann problems related to the compressible Euler equations and shallow water wave equations 
are numerically simulated to test the suitability of the IS-$L^2$-TVB(D)-minmod limiter with high-order polynomial approximations ($P^K,\ K \geq 3$)
and its ability to suppress numerical oscillations.   
\par
$\bullet\ $ Riemann Problems for one-dimensional compressible Euler system. 
\begin{example}\label{example-1D-Euler-Sod}
Sod problem with compressible Euler equations. \\
$\mathrm{I.C.}\ $
$$
(\rho, u, p)= 
\begin{cases}(1,0,1), & -1<x<0, \\ 
  (0.125,0,0.1), & 0<x<1;
\end{cases}
$$
{\color{blue}{Free}} boundary conditions are applied on both the left and right boundaries; \\
Computational Domain: 
$$
\Omega \times [0,T_{end}] = [-1,1] \times [0,0.2];
$$
Limiter parameters: IS-TVB-minmod limiter,$\ \omega_{IS}=1.0,\ \omega_{L^2}=0.0$; \\
TroubledCell-Indicator: TVB-minmod discontinuous indicator and take $M=1$ into effect; \\
$P^3$-polynomial approximation was applied; \\
Numderical flux format: Steger-Warming; \\
Temporal discretization format: TVD-RK3; \\
CFL=0.05; \\
Mesh with 400 cells was applied. \\
At $t=0.2$, the density $\rho$, velocity $u$, pressure $P$ and total energy $E$ are ploted in Figure \ref{Fig.1D-Euler-Sod}.
\begin{figure}[htbp]
  \vspace{0.3cm}
  \begin{center}
    \begin{minipage}{0.49\linewidth}
      \centerline{\includegraphics[width=1\linewidth]{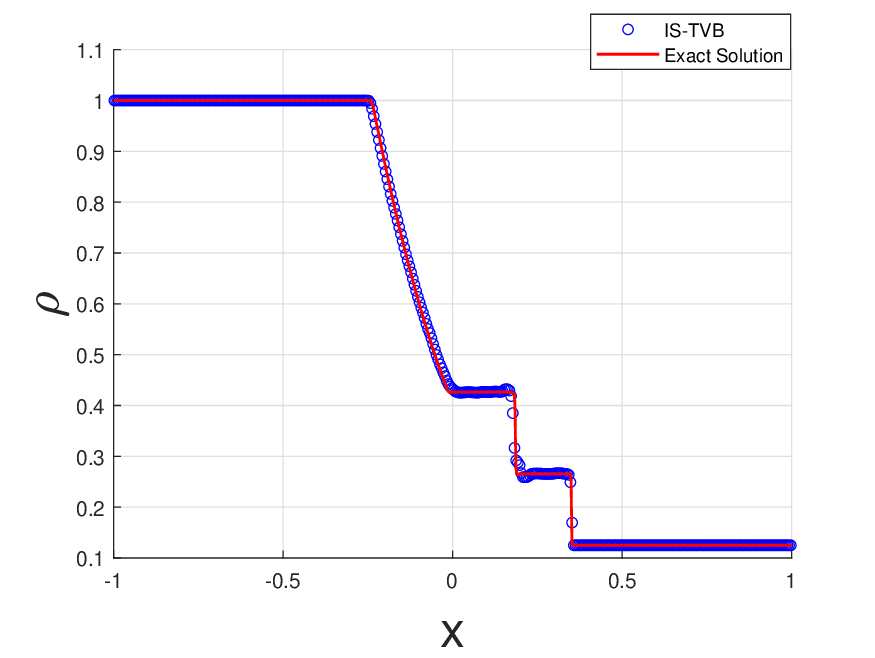}}
    \end{minipage}
    \hfill
    \begin{minipage}{0.49\linewidth}
      \centerline{\includegraphics[width=1\linewidth]{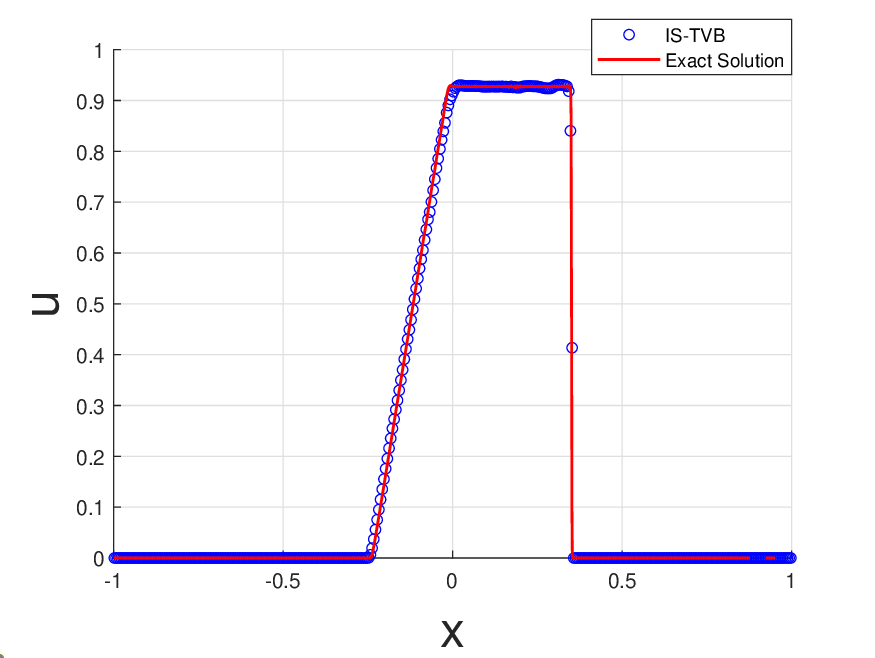}}
    \end{minipage}
    \vfill
    \vspace{0.2cm}
    \begin{minipage}{0.49\linewidth}
      \small
      \centerline{(a)\ Density $\rho$}
    \end{minipage}
    \hfill
    \begin{minipage}{0.49\linewidth}
      \small
      \centerline{(b)\ Velocity $u$}
    \end{minipage}
    \hfill
    \begin{minipage}{0.49\linewidth}
      \centerline{\includegraphics[width=1\linewidth]{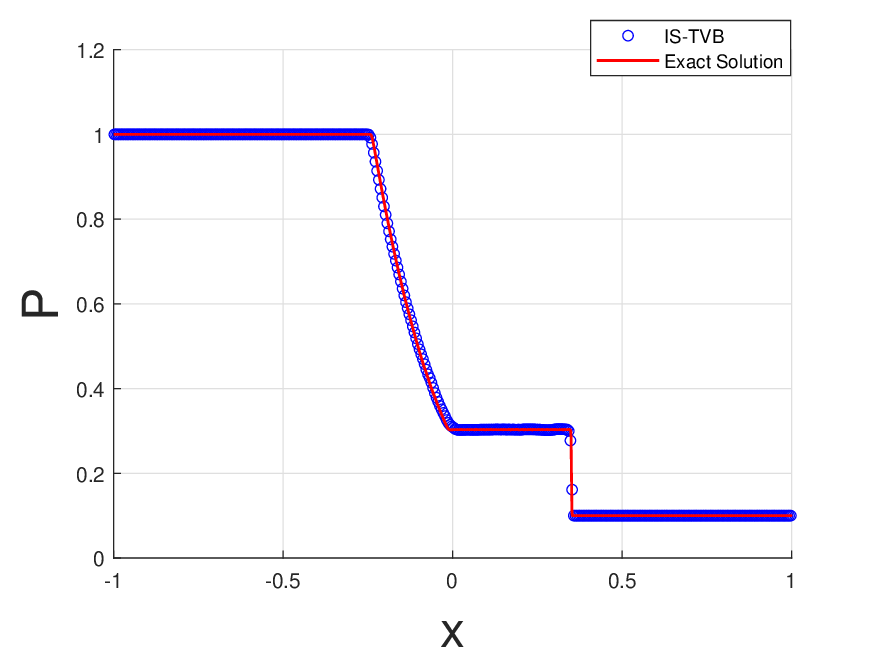}}
    \end{minipage}
    \hfill
    \begin{minipage}{0.49\linewidth}
      \centerline{\includegraphics[width=1\linewidth]{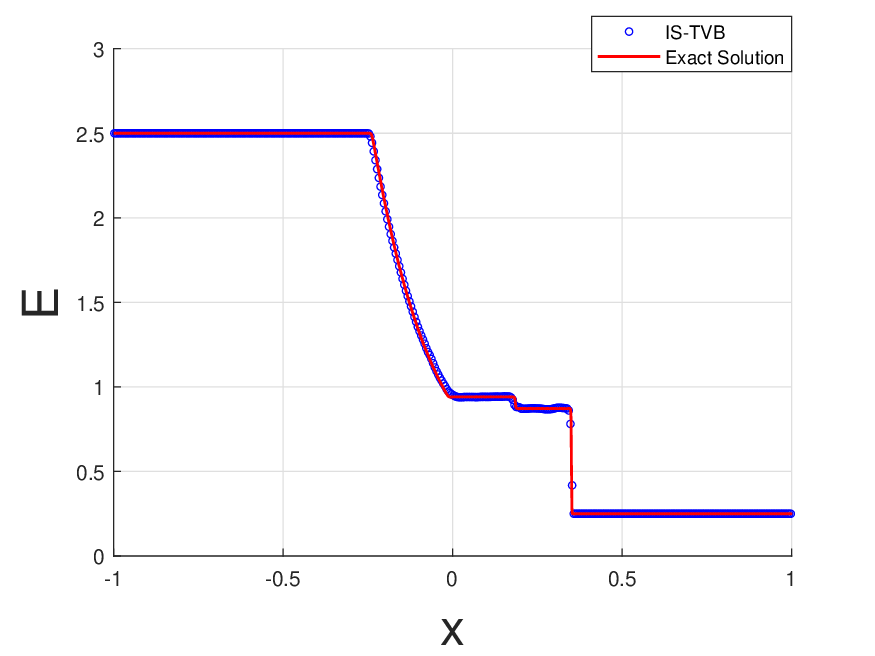}}
    \end{minipage}
    \vfill
    \vspace{0.2cm}
    \begin{minipage}{0.49\linewidth}
      \small
      \centerline{(c)\ Pressure $P$}
    \end{minipage}
    \hfill
    \begin{minipage}{0.49\linewidth}
      \small
      \centerline{(d)\ Total energy $E$}
    \end{minipage}
  \end{center}
  \caption{Sod problem for 1D-Compressible Euler Equations. $P^3$-FVS(Steger-Warming)-DG with IS-TVB-minmod limiter using equally spaced 400 cells.  
  At $t=0.2$, primitive physical variables $\rho,\ u,\ P,\ E$ are demonstrated.}
  \vspace{0.3cm}
  \label{Fig.1D-Euler-Sod}
\end{figure}
\end{example}

\begin{example}\label{example-1D-Euler-Lax}
Lax problem with compressible Euler equations. \\
$\mathrm{I.C.}\ $
$$
(\rho, u, p)= 
\begin{cases}(0.445,0.698,3.528), & -5<x<0 \\ 
  (0.5,0,0.571), & 0<x<5
\end{cases}
$$
{\color{blue}{Free}} boundary conditions are applied on both the left and right boundaries; \\
Computational Domain: 
$$
\Omega \times [0,T_{end}] = [-5,5] \times [0,1.3];
$$
Limiter parameters: IS-TVB-minmod limiter,$\ \omega_{IS}=1.0,\ \omega_{L^2}=0.0$; \\
TroubledCell-Indicator: TVB-minmod discontinuous indicator and take $M=1$ into effect; \\
$P^3$-polynomial approximation was applied; \\
Numderical flux format: van Leer; \\
Temporal discretization format: TVD-RK3; \\
CFL=0.1; \\
Uniform mesh with 2000 cells was applied. \\
At $t=1.3$, the density $\rho$, velocity $u$, pressure $P$ and total energy $E$ are ploted in Figure \ref{Fig.1D-Euler-Lax}.
\begin{figure}[htbp]
  \vspace{0.3cm}
  \begin{center}
    \begin{minipage}{0.49\linewidth}
      \centerline{\includegraphics[width=1\linewidth]{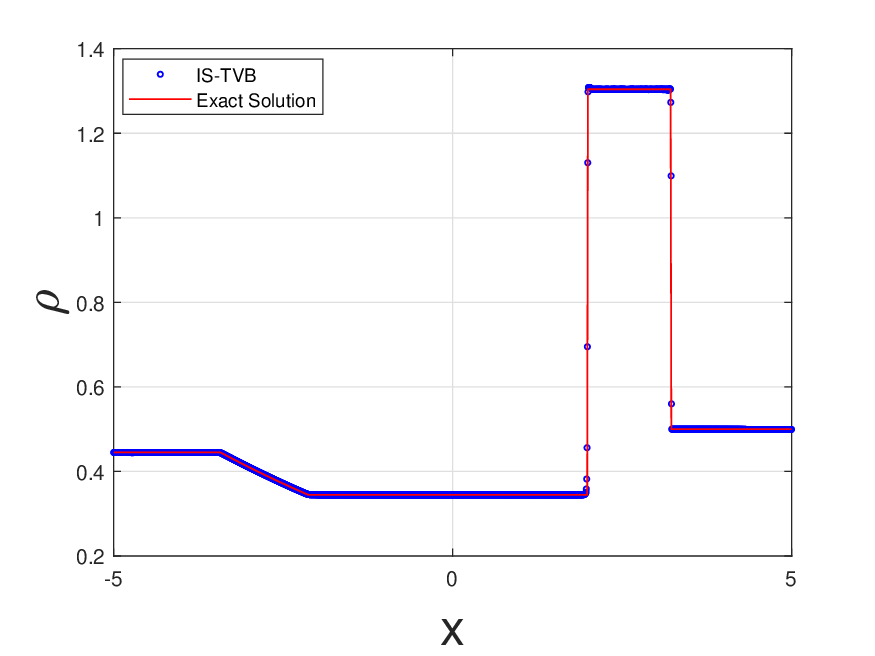}}
    \end{minipage}
    \hfill
    \begin{minipage}{0.49\linewidth}
      \centerline{\includegraphics[width=1\linewidth]{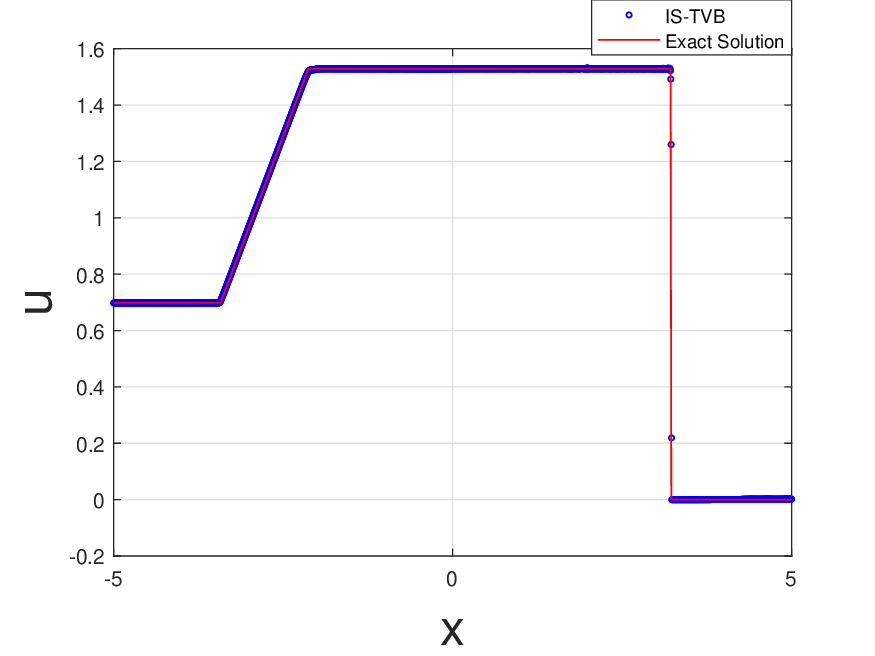}}
    \end{minipage}
    \vfill
    \vspace{0.2cm}
    \begin{minipage}{0.49\linewidth}
      \small
      \centerline{(a)\ Density $\rho$}
    \end{minipage}
    \hfill
    \begin{minipage}{0.49\linewidth}
      \small
      \centerline{(b)\ Velocity $u$}
    \end{minipage}
    \hfill
    \begin{minipage}{0.49\linewidth}
      \centerline{\includegraphics[width=1\linewidth]{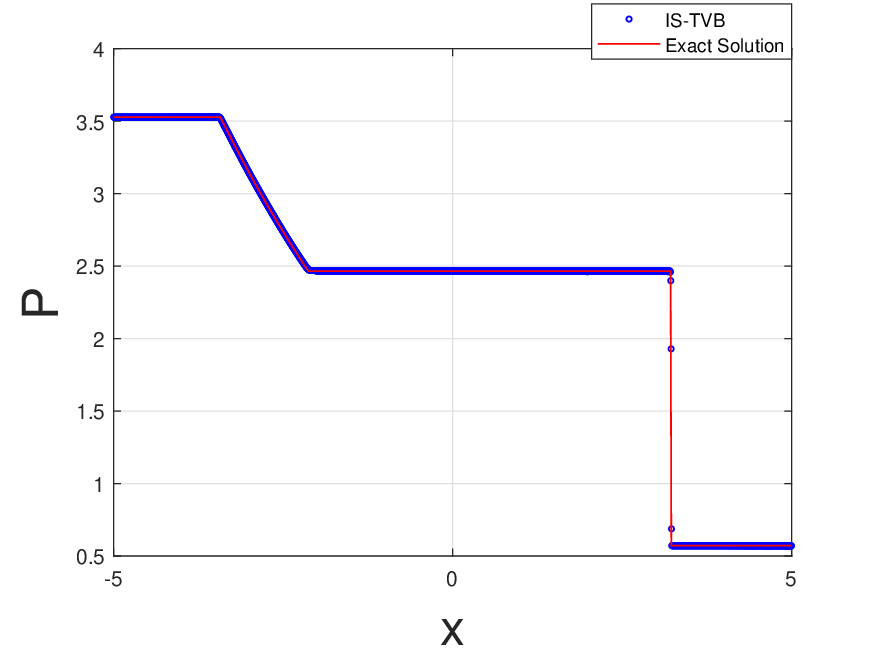}}
    \end{minipage}
    \hfill
    \begin{minipage}{0.49\linewidth}
      \centerline{\includegraphics[width=1\linewidth]{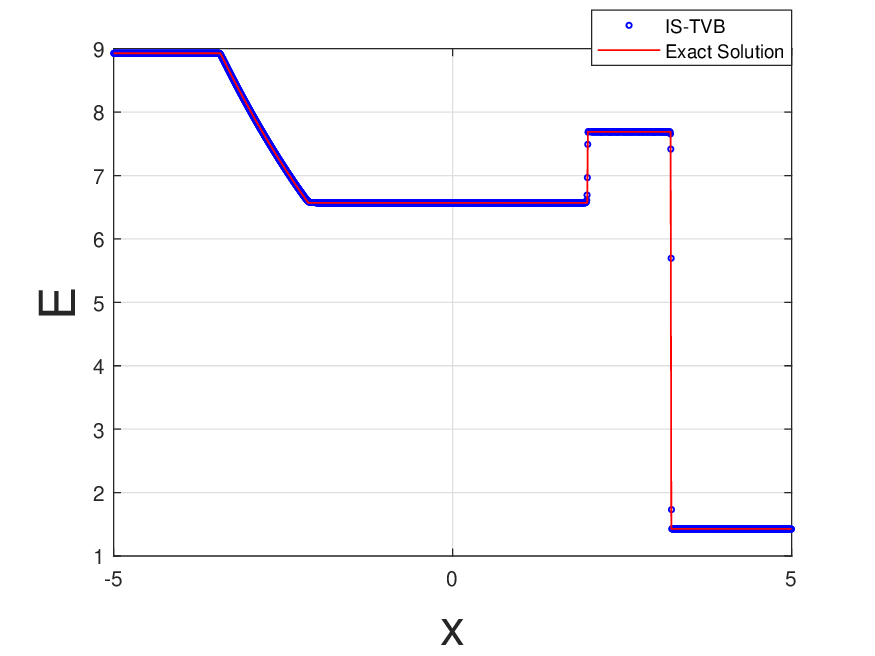}}
    \end{minipage}
    \vfill
    \vspace{0.2cm}
    \begin{minipage}{0.49\linewidth}
      \small
      \centerline{(c)\ Pressure $P$}
    \end{minipage}
    \hfill
    \begin{minipage}{0.49\linewidth}
      \small
      \centerline{(d)\ Total energy $E$}
    \end{minipage}
  \end{center}
  \caption{Lax problem for 1D-Compressible Euler Equations. $P^3$-FVS(van Leer)-DG with IS-TVB-minmod limiter using equally spaced 2000 cells.  
  At $t=1.3$, primitive physical variables $\rho,\ u,\ P,\ E$ are demonstrated.}
  \vspace{0.3cm}
  \label{Fig.1D-Euler-Lax}
\end{figure}
\end{example}

\begin{example}\label{example-1D-Euler-Shu-Osher}
Shu-Osher problem with compressible Euler equations. \\
$\mathrm{I.C.}\ $
$$
(\rho, u, p)= 
\begin{cases}(3.857143,2.629369,10.333333), & -5<x<-4 \\ 
  (1+0.2\sin(5x),0,1), & -4<x<5
\end{cases}
$$
{\color{blue}{Free}} boundary conditions are applied on both the left and right boundaries; \\
Computational Domain: 
$$
\Omega \times [0,T_{end}] = [-5,5] \times [0,1.8];
$$
Limiter parameters: IS-$L^2$-TVB-minmod limiter,$\ \omega_{IS}=0.75,\ \omega_{L^2}=0.25$; \\
TroubledCell-Indicator: TVB-minmod discontinuous indicator and take $M=1$ into effect; \\
$P^5$-polynomial approximation was applied; \\
Numderical flux format: Steger-Warming; \\
Temporal discretization format: TVD-RK3; \\
CFL=0.1; \\
Uniform mesh with 500 cells was applied. \\
At $t=1.8$, the density $\rho$, velocity $u$, pressure $P$ and total energy $E$ are ploted in Figure \ref{Fig.1D-Euler-Shu-Osher}.
\begin{figure}[htbp]
  \vspace{0.3cm}
  \begin{center}
    \begin{minipage}{0.49\linewidth}
      \centerline{\includegraphics[width=1\linewidth]{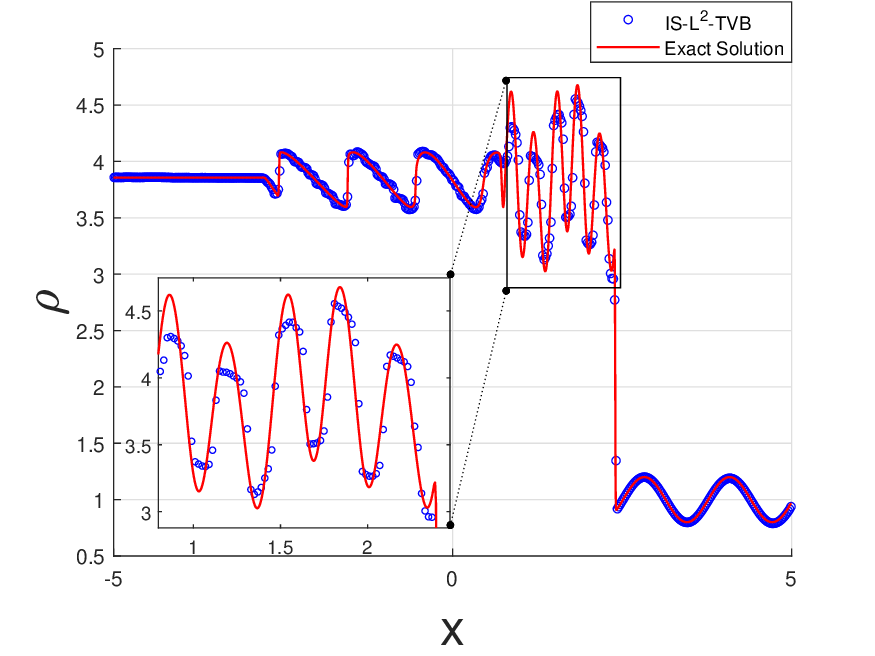}}
    \end{minipage}
    \hfill
    \begin{minipage}{0.49\linewidth}
      \centerline{\includegraphics[width=1\linewidth]{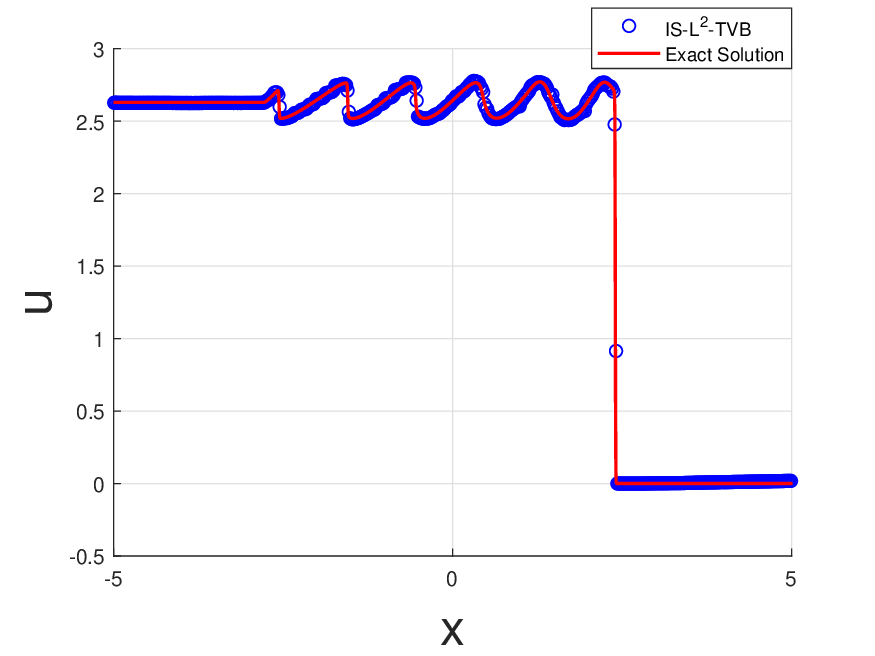}}
    \end{minipage}
    \vfill
    \vspace{0.2cm}
    \begin{minipage}{0.49\linewidth}
      \small
      \centerline{(a)\ Density $\rho$}
    \end{minipage}
    \hfill
    \begin{minipage}{0.49\linewidth}
      \small
      \centerline{(b)\ Velocity $u$}
    \end{minipage}
    \hfill
    \begin{minipage}{0.49\linewidth}
      \centerline{\includegraphics[width=1\linewidth]{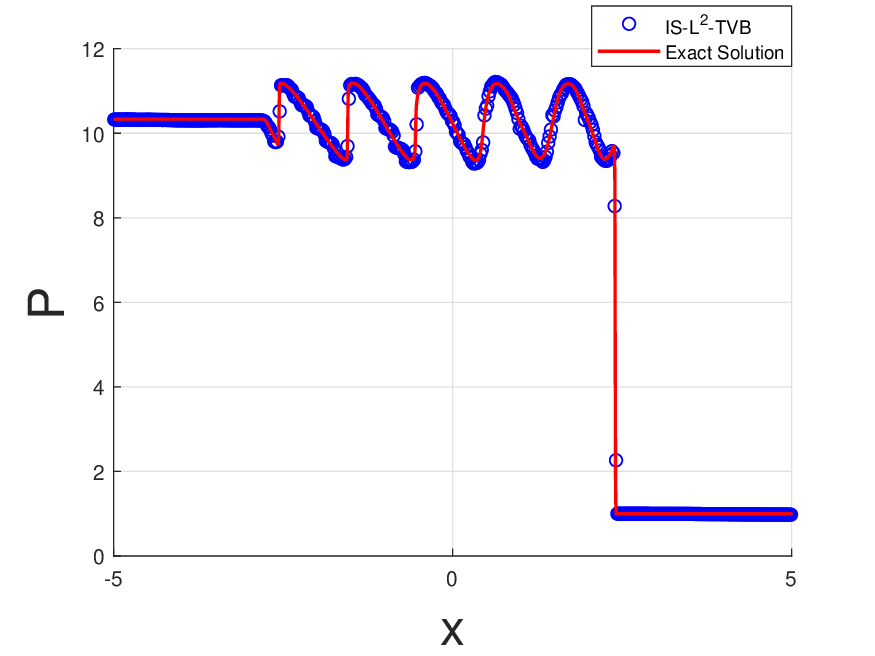}}
    \end{minipage}
    \hfill
    \begin{minipage}{0.49\linewidth}
      \centerline{\includegraphics[width=1\linewidth]{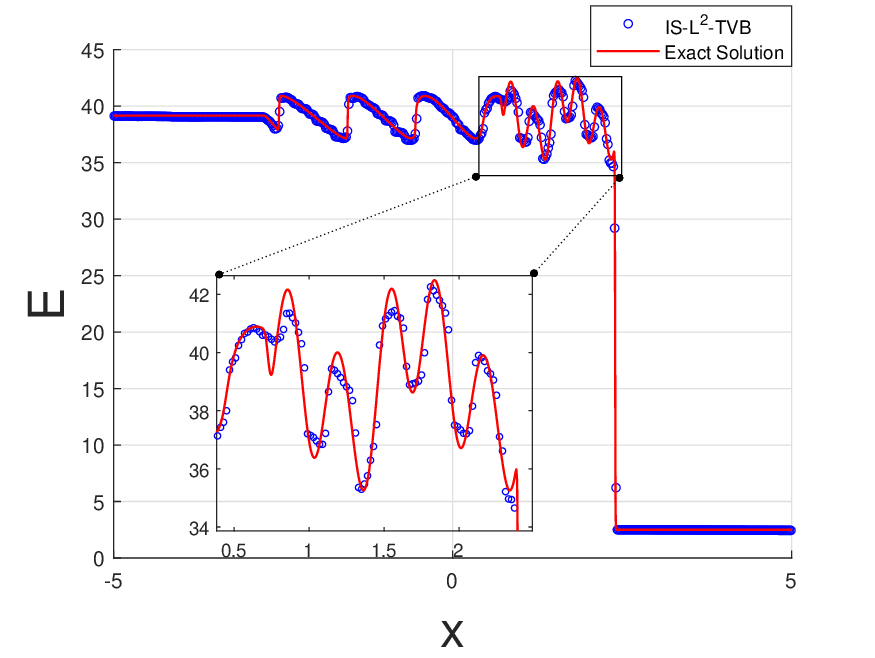}}
    \end{minipage}
    \vfill
    \vspace{0.2cm}
    \begin{minipage}{0.49\linewidth}
      \small
      \centerline{(c)\ Pressure $P$}
    \end{minipage}
    \hfill
    \begin{minipage}{0.49\linewidth}
      \small
      \centerline{(d)\ Total energy $E$}
    \end{minipage}
  \end{center}
  \caption{Shu-Osher problem for 1D-Compressible Euler Equations. $P^5$-FVS(Steger-Warming)-DG with IS-$L^2$-TVB-minmod limiter using equally spaced 500 cells.  
  At $t=1.8$, primitive physical variables $\rho,\ u,\ P,\ E$ are demonstrated.}
  \vspace{0.3cm}
  \label{Fig.1D-Euler-Shu-Osher}
\end{figure}
\end{example}

\begin{example}\label{example-1D-Euler-Blast}
Blast problem with compressible Euler equations. \\
$\mathrm{I.C.}\ $
$$
(\rho, u, p)= 
\begin{cases}\left(1,0,10^3\right), & 0 \leq x<0.1 \\ 
  \left(1,0,10^{-2}\right), & 0.1<x \leq 0.9 \\ 
  \left(1,0,10^2\right), & 0.9<x \leq 1
\end{cases}
$$
{\color{blue}{Reflect}} boundary conditions are applied on both the left and right boundaries; \\
Computational Domain: 
$$
\Omega \times [0,T_{end}] = [0,1] \times [0,0.026];
$$
Limiter parameters: IS-$L^2$-TVB-minmod limiter,$\ \omega_{IS}=0.8,\ \omega_{L^2}=0.2$; \\
TroubledCell-Indicator: TVB-minmod discontinuity indicator and take $M=1$ into effect; \\
$P^2$-polynomial approximation was applied; \\
Numderical flux format: AUSM; \\
Temporal discretization format: TVD-RK3; \\
CFL=0.005; \\
uniform mesh with 800 cells was applied. \\
At $t=0.026$, the density $\rho$ is ploted in Figure \ref{Fig.1D-Euler-Blast}.
\begin{figure}[htbp]
  \begin{center}
      \centerline{\includegraphics[width=0.55\linewidth]{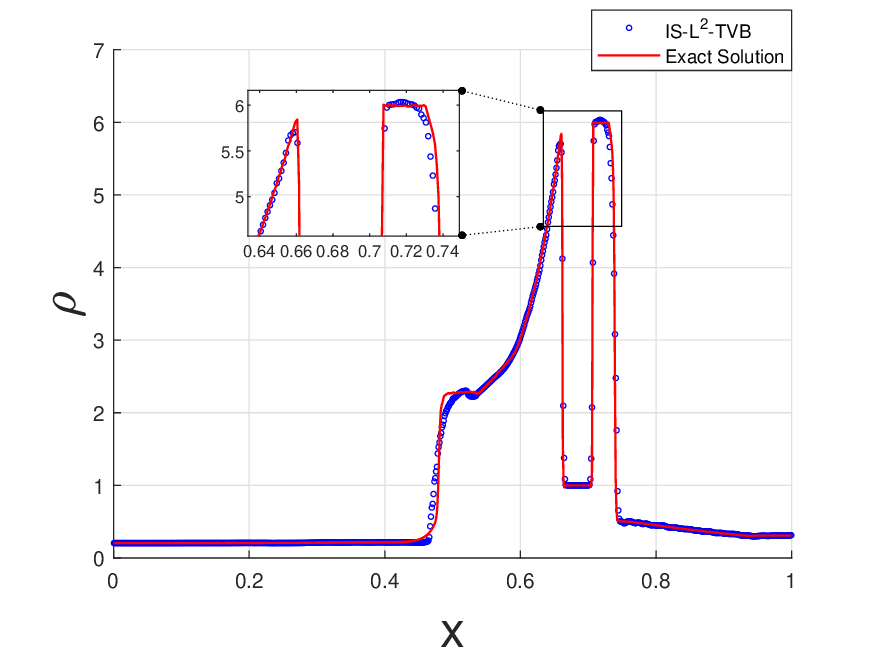}}
  \end{center}
  \caption{Blast problem for 1D-Compressible Euler Equations. $P^2$-FVS(AUSM)-DG with IS-$L^2$-TVB-minmod limiter ($\omega_{IS}=0.8,\ \omega_{L^2}=0.2$) using equally spaced 800 cells.  
  At $t=0.026$, physical variable $\rho$ is demonstrated.}
  \label{Fig.1D-Euler-Blast}
\end{figure}
\end{example}
$\bullet\ $ Riemann Problems for one-dimensional shallow water wave system. 
\begin{example}\label{1D-SWE-DamBreak}
  Dam break on a flat bed with shallow water wave equations. \\
  Bottom elevation: 
  $$
  Z_0\equiv 0;
  $$
  Computational domain:  
  $$
  \Omega \times [0,T_{end}] = [-1,1] \times [0,0.2];
  $$
  Riemann $\mathrm{I.C.}\ $
  $$
  (h, u)= 
  \begin{cases}\left(1,0\right), & -1 < x < 0, \\ 
    \left(0.1,0\right), & 0 < x < 1; 
  \end{cases}
  $$
  $\mathrm{B.C.}\ ${\color{blue}{free}} boundary conditions are applied on both the left and right boundaries; \\
    $P^4$-polynomial approximation was applied; \\
  Limiter parameters: IS-$L^2$-TVB-minmod limiter,$\ \omega_{IS}=0.75,\ \omega_{L^2}=0.25$; \\
  TroubledCell-Indicator: TVD-minmod discontinuous indicator; \\
  Numderical flux format: van Leer; \\
  Temporal discretization format: TVD-RK3; \\
  CFL=0.1; \\
  Uniform mesh with 200 cells was applied. \\
  At $t=0.2$, the surface level $h$ and discharge $hu$ are ploted in Figure \ref{Fig.1D-SWE-DamBreak}.
  \begin{figure}[htbp]
    \vspace{0.3cm}
    \begin{center}
      \begin{minipage}{0.49\linewidth}
        \centerline{\includegraphics[width=1\linewidth]{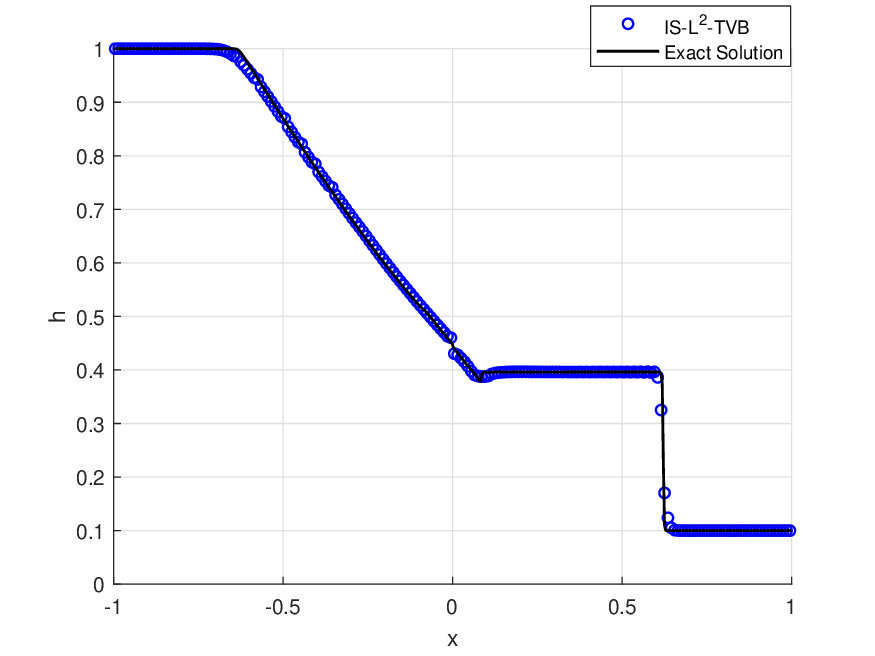}}
      \end{minipage}
      \hfill
      \begin{minipage}{0.49\linewidth}
        \centerline{\includegraphics[width=1\linewidth]{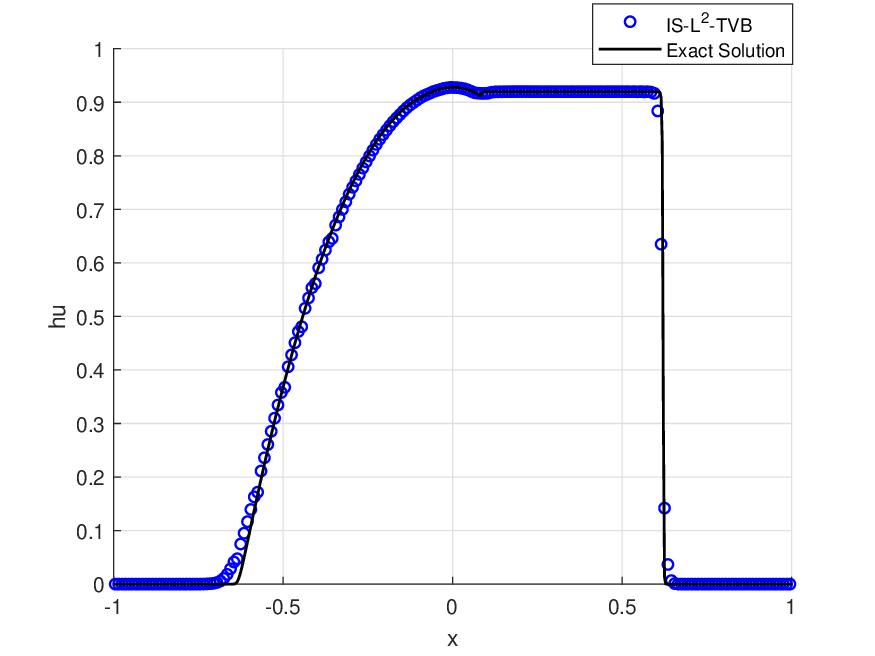}}
      \end{minipage}
      \vfill
      \vspace{0.2cm}
      \begin{minipage}{0.49\linewidth}
        \small
        \centerline{(a)\ DamBreak: surface level $h$ at $t=0.2$}
      \end{minipage}
      \hfill
      \begin{minipage}{0.49\linewidth}
        \small
        \centerline{(b)\ DamBreak: discharge $hu$ at $t=0.2$}
      \end{minipage}
    \end{center}
    \caption{Dam break at a flat bed with shallow water wave equations. $P^4$-FVS(vanLeer)-DG with IS-$L^2$-TVB-minmod limiter using equally spaced 200 cells.  
    At $t=0.2$, conservative variables $h,\ hu$ are demonstrated.}
    \vspace{0.3cm}
    \label{Fig.1D-SWE-DamBreak}
  \end{figure}
\end{example}

\subsubsection{Two-dimensional Test Cases for Hyperbolic Systems}
Three Riemann problems in 2D case for compressible Euler equations are simulated here to demonstrate the performance of IS-$L^2$-TVB(D)-Limiter in conjunction with FVS-DG. 
\begin{example}\label{2D-Riemann-T1}
\hfill \\
Computational Domain: 
$
\Omega \times [0,T_{end}] = \left\{ [0,0.1] \times [0,0.1] \right\} \times [0,0.022]; 
$
\\
$\mathrm{I.C.1}$
$$
(\rho, u, v, p)^T= \begin{cases}(0.5313,0,0,0.4)^T, & x>0.05, y>0.05, \\ (1,0.7276,0,1)^T, & x<0.05, y>0.05, \\ (0.8,0,0,1)^T, & x<0.05, y<0.05, \\ (1,0,0.7276,1)^T, & x>0.05, y<0.05,\end{cases}
$$
$\mathrm{B.C.}\ ${\color{blue}{free boundary conditions are imposed on all edges of $\Omega$; }}
\\
$P^3$-polynomial approximation was applied; \\
Limiter parameters: IS-$L^2$-TVB-minmod limiter,$\ \omega_{IS}=0.8,\ \omega_{L^2}=0.2$; \\
TroubledCell-Indicator: TVB-minmod discontinuity indicator with $M=1$; \\
Numderical Flux Format: Steger-Warming; \\
Temporal Discretization: TVD-RK3; \\
CFL: 0.2; \\
Uniform rectangular mesh with $40 \times 40$ cells. 
\par
The numerical result of density at $t=0.022$ is plotted in Figure \ref{Fig.2D-Riemann-T1}. 
\begin{figure}[htbp]
  \begin{center}
    \begin{minipage}{1\linewidth}
      \centerline{\includegraphics[width=0.55\linewidth]{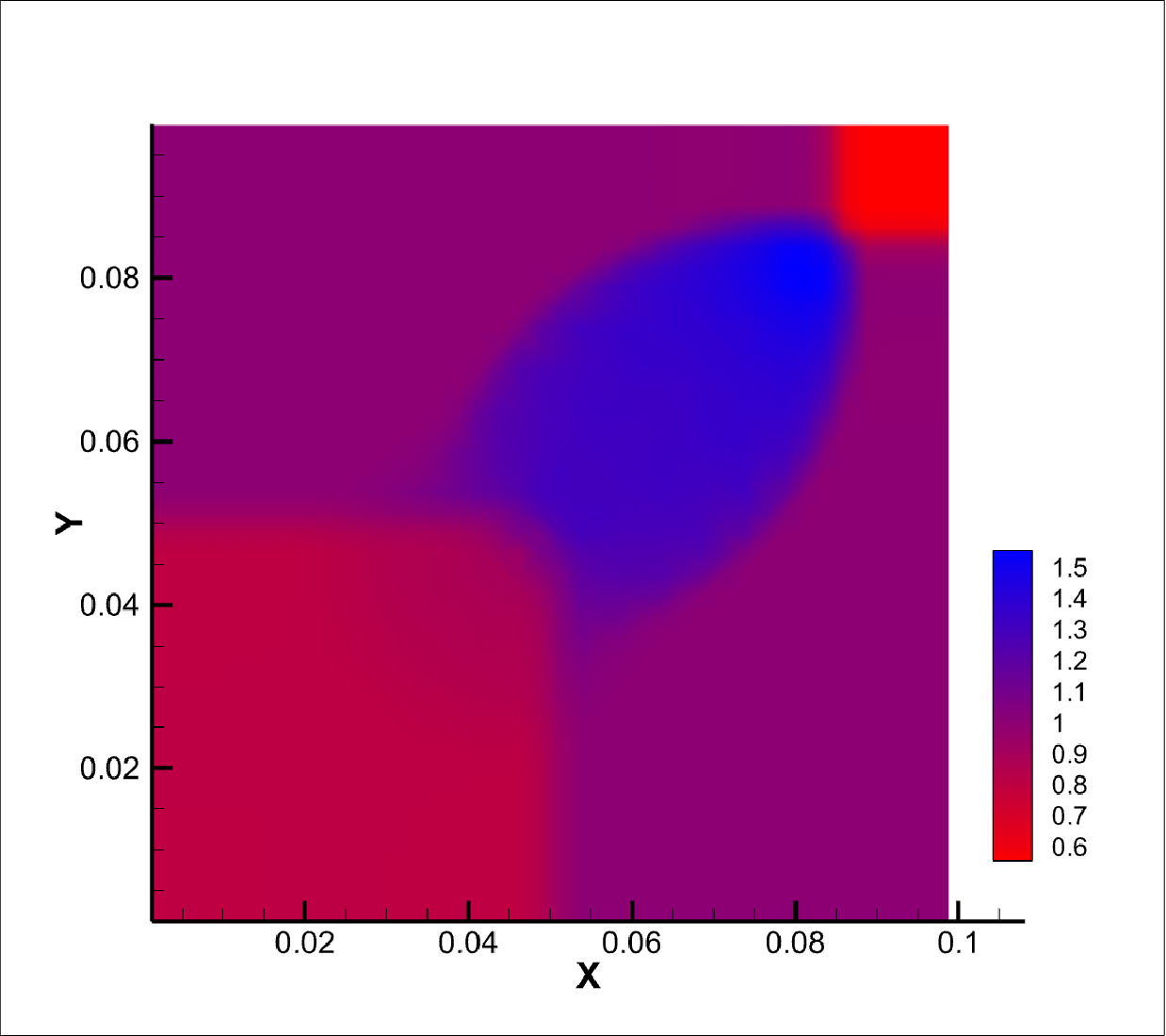}}
    \end{minipage}
  \end{center}
  \caption{\tiny{2D Euler equations for Riemann problem with $\mathrm{I.C.1}$. $P^3$-FVS(Steger-Warming)-DG with IS-$L^2$-TVB-minmod limiter ($\omega_{IS}=0.8,\ \omega_{L^2}=0.2$). 
  Uniform rectangular mesh with $40 \times 40$ cells. The density distribution at $t=0.022$}}
  \label{Fig.2D-Riemann-T1}
\end{figure}
\end{example}
\begin{example}\label{2D-Riemann-T2}
\hfill \\
Computational Domain: 
$
\Omega \times [0,T_{end}] = \left\{ [0,1] \times [0,1] \right\} \times [0,0.25]; 
$
\\
$\mathrm{I.C.2}$
$$
(\rho, u, v, p)^T= \begin{cases}(1.1,0,0,1.1)^T, & x>0.5, y>0.5 \\ (0.5065,0.8939,0,0.35)^T, & x<0.5, y>0.5 \\ (1.1,0.8939,0.8939,1.1)^T, & x<0.5, y<0.5 \\ (0.5065,0,0.8939,0.35)^T, & x>0.5, y<0.5\end{cases}
$$
$\mathrm{B.C.}\ ${\color{blue}{free boundary conditions are imposed on all edges of $\Omega$; }}
\\
$P^2$-polynomial approximation was applied; \\
Limiter parameters: IS-TVB-minmod limiter,$\ \omega_{IS}=1.0,\ \omega_{L^2}=0.0$; \\
TroubledCell-Indicator: TVB-minmod discontinuity indicator with $M=1$; \\
Numderical Flux Format: AUSM; \\
Temporal Discretization: TVD-RK3; \\
CFL: 0.2; \\
Uniform rectangular mesh with $100 \times 100$ cells. 
\par
The numerical result of density at $t=0.25$ is plotted in Figure \ref{Fig.2D-Riemann-T2}. 
\begin{figure}[htbp]
  \begin{center}
    \begin{minipage}{1\linewidth}
      \centerline{\includegraphics[width=0.55\linewidth]{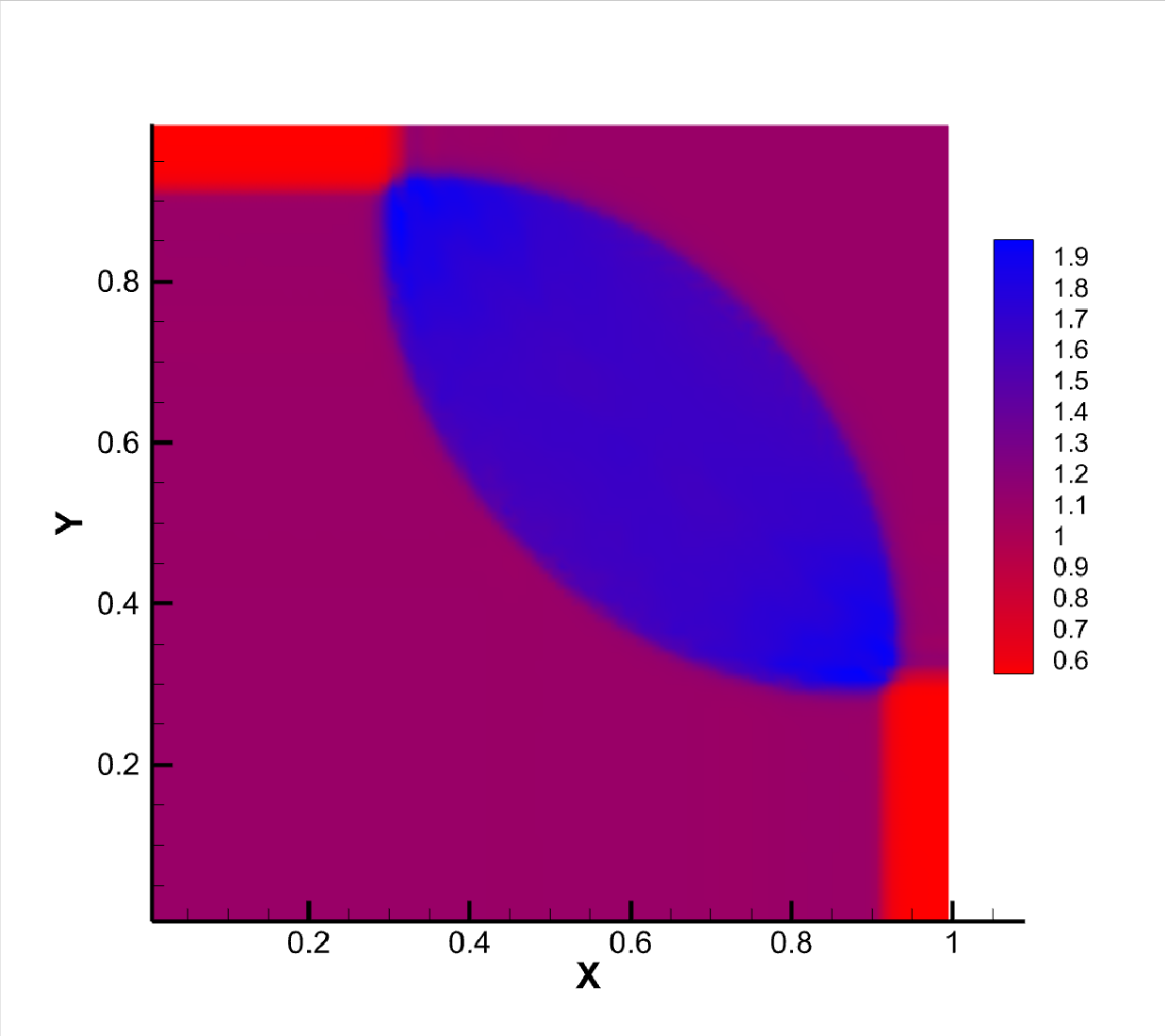}}
    \end{minipage}
  \end{center}
  \caption{\tiny{2D Euler equations for Riemann problem with $\mathrm{I.C.2}$. $P^2$-FVS(AUSM)-DG with IS-TVB-minmod limiter ($\omega_{IS}=1.0,\ \omega_{L^2}=0.0$). 
  Uniform rectangular mesh with $100 \times 100$ cells. The density distribution at $t=0.25$}}
  \label{Fig.2D-Riemann-T2}
\end{figure}
\end{example}
\begin{example}\label{2D-Riemann-T3}
  \hfill \\
  Computational Domain: 
  $
  \Omega \times [0,T_{end}] = \left\{ [0,1] \times [0,1] \right\} \times [0,0.3]; 
  $
  \\
$\mathrm{I.C.3}$
  $$
  (\rho, u, v, p)^T= \begin{cases}(0.5313,0,0,0.4)^T, & x>0.5, y>0.5, \\ (1,0.7276,0,1)^T, & x<0.5, y>0.5, \\ (0.8,0,0,1)^T, & x<0.5, y<0.5, \\ (1,0,0.7276,1)^T, & x>0.5, y<0.5,\end{cases}
  $$
  $\mathrm{B.C.}\ ${\color{blue}{free boundary conditions are imposed on all edges of $\Omega$; }}
  \\
  $P^2$-polynomial approximation was applied; \\
  Limiter parameters: IS-$L^2$-TVB-minmod limiter,$\ \omega_{IS}=0.8,\ \omega_{L^2}=0.2$; \\
  TroubledCell-Indicator: TVB-minmod discontinuity indicator with $M=1$; \\
  Numderical Flux Format: Steger-Warming; \\
  Temporal Discretization: TVD-RK3; \\
  CFL: 0.2; \\
  Uniform rectangular mesh with $100 \times 100$ cells. 
  \par
  The numerical result of density at $t=0.3$ is plotted in Figure \ref{Fig.2D-Riemann-T3}. 
  \begin{figure}[htbp]
    \begin{center}
      \begin{minipage}{1\linewidth}
        \centerline{\includegraphics[width=0.55\linewidth]{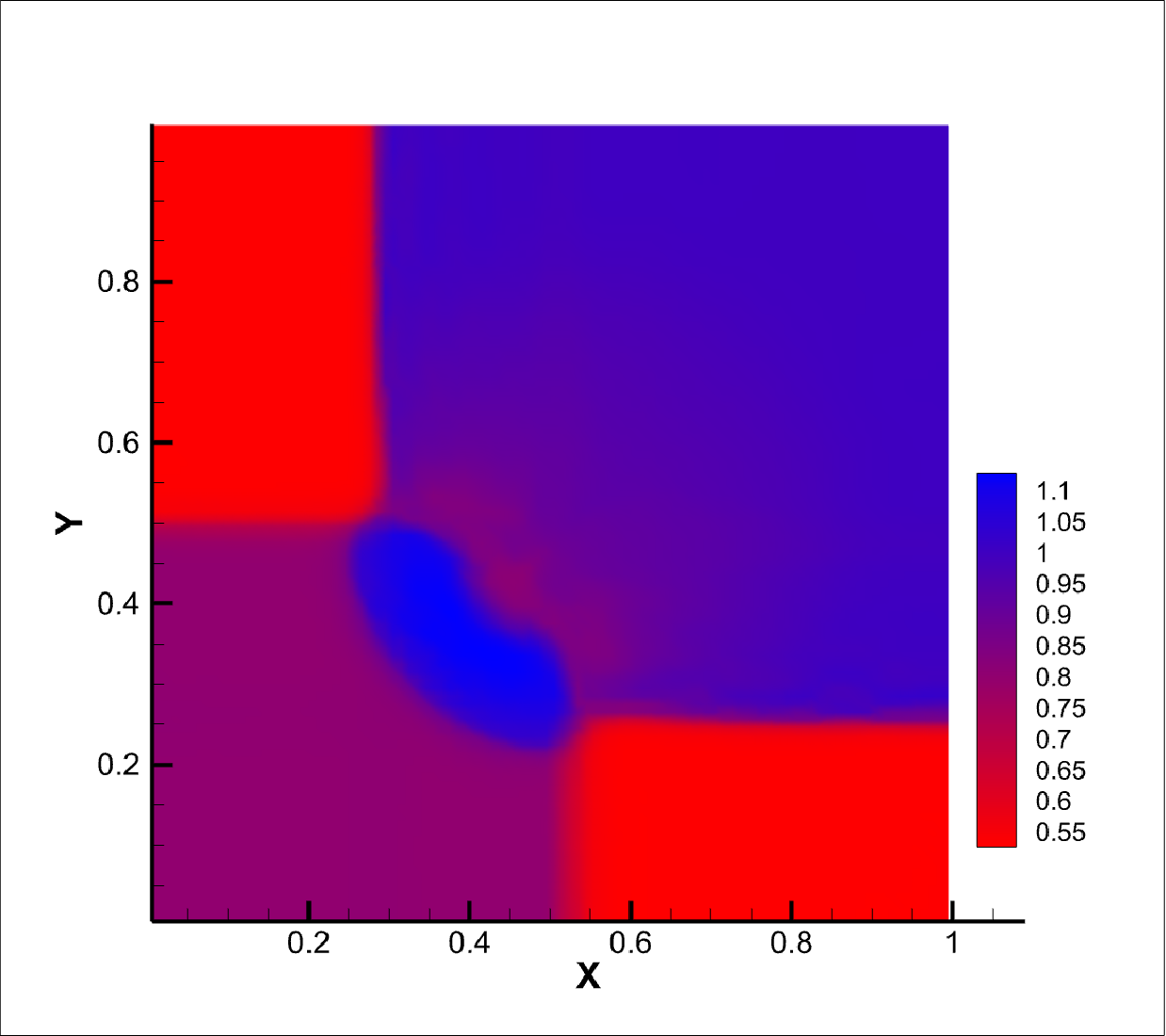}}
      \end{minipage}
    \end{center}
    \caption{\tiny{2D Euler equations for Riemann problem with $\mathrm{I.C.3}$. $P^2$-FVS(Steger-Warming)-DG with IS-$L^2$-TVB-minmod limiter ($\omega_{IS}=0.8,\ \omega_{L^2}=0.2$). 
    Uniform rectangular mesh with $100 \times 100$ cells. The density distribution at $t=0.3$}}
    \vspace{0.3cm}
    \label{Fig.2D-Riemann-T3}
  \end{figure}
  \end{example}

\section{Conclusion}\label{sec-Conclusion}
This paper constructs the numerical flux required for the DG spatial discretization scheme based on the flux vector splitting method, which is implemented in two-dimensional problems by introducing the normal flux on cell interfaces. 
The TVB(D) limiter and the WENO limiter are currently the two main streams of limiters; we utilize the smoothness factor $IS$ from the WENO limiter in the TVB(D)-minmod limiter, 
successfully overcoming the ill-posed problem of the TVB(D)-minmod limiter in the case of high-order polynomial approximation by constructing an optimization problem based on the smoothness factor constraint. 
Furthermore, by drawing on existing practices and introducing the $L^2$-error constraint, the TVB(D)-minmod limiter is able to balance the suppression of oscillations and the maintenance of high precision.  
For local characteristic decomposition, 
we use an arithmetic average or Roe average of the integral mean on the common interface of adjacent cells for local freezing, 
which is different from the usual practice. 
We also use an interpolation-based characteristic transformation to achieve the mutual conversion between the one-dimensional model 
in physical space and characteristic space. Numerical experiments have confirmed the effectiveness of the aforementioned work.

\section*{Acknowledgments}
The authors would like to thank the anonymous referees for their very valuable comments and suggestions.

\appendix
\section{Vector, Matrix, and Tensor Operations in FVS-DG} \label{appendix-notations}
The following notations and operational rules are used for the derivation of the FVS-DG spatial discretization scheme (definition of the FVS-DG weak solution). 
\begin{itemize}
  \item The Hadamard product between two vector-valued functions: $\mathbf{a} \odot \mathbf{b}=(a_1 b_1, a_2 b_2, \ldots, a_m b_m)^{\mathrm{T}}$. 
  \item The quasi-Hadamard product between a matrix-valued function and a vector-valued function is defined as
  $$
  \mathbb{H} \odot \mathbf{V} =
  \left[
  \begin{array}{l}
  \mathbf{H}^{(1)} \cdot v_1 \\
  \mathbf{H}^{(2)} \cdot v_2 \\
  \vdots \\
  \mathbf{H}^{(m)} \cdot v_m \\
  \end{array}
  \right],
  $$
  where $\mathbf{H}^{(i)} = (h_{i1}, h_{i2}, \ldots, h_{id})$, that is, the matrix $\mathbb{H}$ is partitioned row-wise.
  \item
  The inner product between vector-valued functions acts component-wise, i.e.,
  $$
  \langle\mathbf{U}, \mathbf{V}\rangle =
  \begin{bmatrix}
  \langle u_1, v_1 \rangle \\
  \langle u_2, v_2 \rangle \\
  \vdots \\
  \langle u_m, v_m \rangle
  \end{bmatrix}.
  $$ 
  \item The Frobenius inner product between vector functions is defined as
  $\left<\mathbf{U} : \mathbf{V}\right> = \sum_{i=1}^{m} \left<u_i, v_i\right>$. 
  \item Let the matrix-valued function
  $\mathbb{H}
  =
  \begin{bmatrix}
    \mathbf{H}^{(1)} \\
    \mathbf{H}^{(2)} \\
    \vdots \\
    \mathbf{H}^{(m)}
  \end{bmatrix}
  $,  where
  $\mathbf{H}^{(i)}=(h_{i1}, h_{i2}, \ldots, h_{id})$, that is, the matrix $\mathbb{H}$ is partitioned into row blocks. \\
  Let the matrix-valued function
  $\mathbb{V}
  =
  \begin{bmatrix}
    \mathbf{V}^{(1)} \\
    \mathbf{V}^{(2)} \\
    \vdots \\
    \mathbf{V}^{(m)}
  \end{bmatrix}
  $,  where
  $\mathbf{V}^{(i)}=(v_{i1}, v_{i2}, \ldots, v_{id})$, that is, the matrix $\mathbb{V}$ is partitioned into row blocks. \\
  The Frobenius inner product between the matrix-valued functions is defined to act component-wise along the rows, i.e.,
  $
  \left<\mathbb{H}:\mathbb{V}\right>
  =
  \left[
  \begin{array}{l}
    \left<\mathbf{H}^{(1)}:\mathbf{V}^{(1)}\right> \\
    \left<\mathbf{H}^{(2)}:\mathbf{V}^{(2)}\right> \\
    \vdots \\
    \left<\mathbf{H}^{(m)}:\mathbf{V}^{(m)}\right>
  \end{array}
  \right]
  $. 
\end{itemize}
\begin{itemize}
  \item When the divergence operator acts on a matrix-valued function, the matrix-valued function is partitioned column-wise for the computation, such that
  $$
  \nabla \cdot \mathbb{H} = \partial_{x_1} \mathbf{H}_{(1)} + \partial_{x_2} \mathbf{H}_{(2)} + \cdots + \partial_{x_d} \mathbf{H}_{(d)},
  $$
  where $\mathbf{H}_{(j)} = \left(h_{1j}, h_{2j}, \ldots, h_{mj}\right)^{\mathrm{T}}$ represents the $j$-th column of matrix $\mathbb{H}$. 
  \item In this paper, Kronecker product between two vector-valued functions, still denoted by $\otimes$, is redefined as: 
  $$
  \mathbf{a} \otimes \mathbf{b}=\mathbf{a b}^{\mathrm{T}}
  =\left[\begin{array}{l}
  a_1 \\
  a_2 \\
  a_3
  \end{array}\right]
  \left[b_1,  b_2,  b_3\right]
  =\left[\begin{array}{lll}
  a_1 b_1 & a_1 b_2 & a_1 b_3 \\
  a_2 b_1 & a_2 b_2 & a_2 b_3 \\
  a_3 b_1 & a_3 b_2 & a_3 b_3
  \end{array}\right]. 
  $$
  In the dummy index summation notation, the new Kronecker product between two vector-valued functions is expressed as: 
  $
  \mathbf{a} \otimes \mathbf{b} = a_i\mathbf{e}_i \otimes b_j\mathbf{e}_j = a_ib_j\mathbf{e}_i\mathbf{e}_j
  $, 
  where $\mathbf{e}_i$ and $\mathbf{e}_j$ are the standard basis vectors in their respective spaces. 
  \item The gradient of a vector function is a second-order tensor function. Let $\mathbf{V}=(v_1, v_2, \ldots, v_m)^{\mathrm{T}}$, then we have \\
  $$
  \begin{aligned}
  &\nabla \mathbf{V}=
  \begin{bmatrix}
    \partial_{x_1} v_1 & \partial_{x_2} v_1 & \cdots & \partial_{x_d} v_1 \\
    \partial_{x_1} v_2 & \partial_{x_2} v_2 & \cdots & \partial_{x_d} v_2 \\
    \vdots & \vdots & \vdots &\vdots \\
    \partial_{x_1} v_m & \partial_{x_2} v_m & \cdots & \partial_{x_d} v_m
  \end{bmatrix},
  \\
  &\nabla \mathbf{V}=(\partial_{i}v_j)\mathbf{e}_i\mathbf{e}_j.
  \end{aligned}
  $$
\end{itemize}

\section{Flux Vector Splitting Method for Shallow Water System}\label{appendix-SWE-FVS}
\renewcommand{\arraycolsep}{6pt}
\renewcommand{\arraystretch}{1.25}
Detailed introduction to the flux functions (including normal flux), Jacobian matrixs and their eigenstructures of the shallow water wave equations,  
readers are referred to \cite{ref46}.

\subsection{Flux Vector Splitting for Shallow Water Wave Equations in One-dimension}
\subsubsection{Mach Number Splitting Based on True Wave Speed}
\hspace*{\fill} \\
True wave speed: $a=\sqrt{gh}$. \\
True local Mach number: $M_a=u/a$. \\
Pressure: $P=\frac{1}{2}gh^2$; Adiabatic index: $\gamma=2$. 
$$
\begin{aligned}
\mathbf{F} 
=\left[\begin{array}{c}
  h u \\
  h u^2+\frac{1}{2} g h^2
  \end{array}\right]
=h a\left[\begin{array}{l}
M_a \\
u M_a+a / \gamma
\end{array}\right].
\end{aligned}
\leqno{(vanLeer)}
$$
$$
\mathbf{F}=\left[\begin{array}{l}
hu \\
hu^2+\frac{1}{2}gh^2
\end{array}\right]
=
\left[\begin{array}{l}
{\color{red}{h a M_a}} \\
{\color{red}{h a M_a u}}+{\color{blue}{P}}
\end{array}\right]
=
{\color{red}{\underbrace{h a M_a\left[\begin{array}{l}
1 \\
u
\end{array}\right]}_{\mathbf{F}^a}}}
+
{\color{blue}{\underbrace{P \left[\begin{array}{l}
0 \\
1  
\end{array}\right]}_{\mathbf{F}^P}}}.
\leqno{(AUSM)}
$$
The splitting method for $M_a$ and $P$ is the same as for the compressible Euler equations.

\subsubsection{Jacobian Eigenvalue Splitting Based on Modified Wave Speed\cite{ref6}}
\hspace*{\fill} \\
True wave speed: $a=\sqrt{gh}$. 
\\
Modified wave speed: $a^{*}=\sqrt{gh/2}$. 
\\
The Jacobian and its eigenstructure after forced modification are:  
\begin{align*}
&\mathbf{A}^*
=
\left[\begin{array}{ll}
  0 &  1 \\
  {\color{red}{\frac{g h}{2}}}-u^2 & 2 u
  \end{array}\right] ,
\\
&\boldsymbol{\Lambda}^*=
\begin{bmatrix}
u-{\color{red}{\sqrt{g h/2}}} & 0 \\
0 & u+{\color{red}{\sqrt{g h/2}}}
\end{bmatrix},  \\
& \mathbf{R}^*=\left[\begin{array}{cc}
1 & 1 \\
u-{\color{red}{\sqrt{g h/2}}} & u+{\color{red}{\sqrt{g h/2}}}
\end{array}\right],  
\\
& \mathbf{L}^*=\left[\begin{array}{ll}
\frac{u+{\color{red}{\sqrt{g h/2}}}}{2 {\color{red}{\sqrt{g h/2}}}} & -\frac{1}{2 {\color{red}{\sqrt{g h/2}}}} \\
-\frac{u-{\color{red}{\sqrt{g h/2}}}}{2 {\color{red}{\sqrt{g h/2}}}} & \frac{1}{2 {\color{red}{\sqrt{g h/2}}}}
\end{array}\right].
\end{align*}
{\color{red}{After the modification of the wave speed, it holds that $\mathbf{F} = \mathbf{A}^* \mathbf{U}$}}.
\\
Subsequently, the Steger-Warming splitting or Lax-Friedrichs splitting can be applied.

\subsection{Flux Vector Splitting for Shallow Water Wave Equations in Two-dimension}
\subsubsection{Mach Number Splitting Based on True Wave Speed}
\hspace*{\fill} \\
True wave speed: $a=\sqrt{gh}$. \\
True local normal Mach number: $M_a^n=q_n/a$. \\
Pressure: $P=\frac{1}{2}gh^2$; Adiabatic index: $\gamma=2$.
$$
\boldsymbol{\mathcal{F}}_{n}=\left[\begin{array}{l}
h q_n \\
h u q_n + P n_x \\
h v q_n +P n_y
\end{array}\right]
=
\left[\begin{array}{l}
{\color{red}{h a M_a}} \\
{\color{red}{h a M_a u}}+{\color{blue}{P n_x}} \\
{\color{red}{h a M_a v}}+{\color{blue}{P n_y}}
\end{array}\right]
=
{\color{red}{\underbrace{h a M_a\left[\begin{array}{l}
1 \\
u \\
v
\end{array}\right]}_{\mathbf{F}^a}}}
+
{\color{blue}{\underbrace{P \left[\begin{array}{l}
0 \\
n_x \\
n_y \\
\end{array}\right]}_{\mathbf{F}^P}}}.
\leqno{(AUSM)}
$$
The splitting method for $M_a$ and $P$ is the same as for the compressible Euler equations.  

\subsubsection{Jacobian Eigenvalue Splitting Based on Modified Wave Speed\cite{ref6}}
\hspace*{\fill} \\
True wave speed: $a=\sqrt{gh}$. 
\\
Modified wave speed: $a^{*}=\sqrt{gh/2}$. 
\\
The Jacobian and its eigenstructure after forced modification are: 
$$
\begin{aligned}
& \boldsymbol{\Lambda}_n^*=\left[\begin{array}{ccc}
q_n-{\color{red}{\sqrt{g h/2}}} & 0 & 0 \\
0 & q_n & 0 \\
0 & 0 & q_n+{\color{red}{\sqrt{g h/2}}}
\end{array}\right],  
\\
& \mathbf{R}_n^*=\left[\begin{array}{ccc}
1 & 0 & 1 \\
u-{\color{red}{\sqrt{g h/2}}} n_x & \ell_x & u+{\color{red}{\sqrt{g h/2}}} n_x \\
v-{\color{red}{\sqrt{g h/2}}} n_y & \ell_y & v+{\color{red}{\sqrt{g h/2}}} n_y
\end{array}\right], 
\mathbf{L}_n^*=\left[\begin{array}{ccc}
\frac{{\color{red}{\sqrt{g h/2}}}+q_n}{2 {\color{red}{\sqrt{g h/2}}}} & -\frac{n_x}{2 {\color{red}{\sqrt{g h/2}}}} & -\frac{n_y}{2 {\color{red}{\sqrt{g h/2}}}} \\
-q_{\ell} & \ell_x & \ell_y \\
\frac{{\color{red}{\sqrt{g h/2}}}-q_n}{2 {\color{red}{\sqrt{g h/2}}}} & \frac{n_x}{2 {\color{red}{\sqrt{g h/2}}}} & \frac{n_y}{2 {\color{red}{\sqrt{g h/2}}}}
\end{array}\right], 
\\
&\mathbf{A}_n^*=\mathbf{R}_n^*\boldsymbol{\Lambda}_n^*\mathbf{L}_n^*.
\end{aligned}
$$ 
{\color{red}{After the modification of the wave speed, it holds that $\boldsymbol{\mathcal{F}}_n = \mathbf{A}_n^* \mathbf{U}$}}. 
\\
Subsequently, the Steger-Warming splitting or Lax-Friedrichs splitting can be applied.

\section{Development of Numerical Flux Formats for Scalar Equations Based on the System's Jacobian-FVS Method}\label{appendix-FVS-scalar}
\subsection{The Classical Lax-Friedrichs Flux Format from the Jacobian-FVS Perspective}\label{appendix-LF-FVS-LF-flux}
This section attempts to derive the classical Lax-Friedrichs flux format from the L-F splitting method within the Jacobian-FVS.
\par
Consider a one-dimensional hyperbolic conservation system as an example. The L-F splitting process is shown below:
\\
Take
$$
\begin{aligned}
& A^{L, +} = \frac{1}{2}(A^L+M^L I_d)=\frac{1}{2}(A(U^L)+M^L I_d), \\
& A^{R, -} = \frac{1}{2}(A^R-M^R I_d)=\frac{1}{2}(A(U^R)-M^R I_d),
\end{aligned}
$$
where $M^L,M^R$ are positive constants that should ensure
$$
\left.
\begin{aligned}
  & \lambda(A^{L,+})=\lambda(A^L + M^L I_d) \geq 0, \\
  & \lambda(A^{R,-})=\lambda(A^R - M^R I_d) \leq 0,
\end{aligned}
\right\}
$$
hold at every cell interface. 
\\
At $x_{i+1/2}$, there is
$$
\begin{aligned}
\hat{f}_{i+1 / 2} 
& = A^{L, +}_{i+1/2}\mathbf{U}^L_{i+1/2} + A^{R, -}_{i+1/2}\mathbf{U}^R_{i+1/2} \\
& =\frac{1}{2}\left(A\left(\mathbf{U}_{i+1 / 2}^L\right)+M^L_{i+1/2} I\right) \mathbf{U}_{i+1 / 2}^L+\frac{1}{2}\left(A\left(\mathbf{U}_{i+1 / 2}^R\right)-M^R_{i+1/2} I\right) \mathbf{U}_{i+1 / 2}^R \\
& =\frac{1}{2}\left(A\left(\mathbf{U}_{i+1 / 2}^L\right) \mathbf{U}_{i+1 / 2}^L+M^L_{i+1/2} \mathbf{U}_{i+1 / 2}^L\right)+\frac{1}{2}\left(A\left(\mathbf{U}_{i+1 / 2}^R\right) \mathbf{U}_{i+1 / 2}^R-M^R_{i+1/2} \mathbf{U}_{i+1 / 2}^R\right) \\
& =\frac{1}{2}\left(f\left(\mathbf{U}_{i+1 / 2}^L\right)+M^L_{i+1/2} \mathbf{U}_{i+1 / 2}^L\right)+\frac{1}{2}\left(f\left(\mathbf{U}_{i+1 / 2}^R\right)-M^R_{i+1/2} \mathbf{U}_{i+1 / 2}^R\right) \\
& =\frac{1}{2}\left(f\left(\mathbf{U}_{i+1 / 2}^R\right)+f\left(\mathbf{U}_{i+1 / 2}^L\right)-\left(M^R_{i+1/2} \mathbf{U}_{i+1 / 2}^R-M^L_{i+1/2} \mathbf{U}_{i+1 / 2}^L\right)\right)
.\end{aligned}
$$
The classic L-F flux format is given by 
$$
\hat{\hat{f}}_{i+1/2}=\frac{1}{2}\left(f\left(\mathbf{U}_{i+1 / 2}^{R}\right)+f\left(\mathbf{U}_{i+1 / 2}^{L}\right)-\alpha_{i+1/2}\left(\mathbf{U}_{i+1 / 2}^{R}-\mathbf{U}_{i+1 / 2}^{L}\right)\right)
,$$
where $\alpha_{i+1/2}$ is usually taken as  
$$
\alpha_{i+1/2}=\rho\left(\left.\frac{\partial f}{\partial u}\right|_{x_{i+1/2}}\right),
$$
which means $\alpha_{i+1/2}$ is the spectral radius of the Jacobian of the flux function $f$ with respect to the field variable $u$ evaluated at $x_{i+1/2}$.  
\\  
It is noted that if we take
$$
M_{i+1/2}=\max\{M^L_{i+1/2},M^R_{i+1/2}\}, 
$$
then $M_{i+1/2}$ also satisfies
$$
\left.
\begin{aligned}
  & \lambda(A^{L,+}_{i+1/2})=\lambda(A^L_{i+1/2} + M_{i+1/2} I_d) \geq 0, \\
  & \lambda(A^{R,-}_{i+1/2})=\lambda(A^R_{i+1/2} - M_{i+1/2} I_d) \leq 0.
\end{aligned}
\right\}
$$
Thus, $M_{i+1/2}$ can be used to simultaneously replace $M^L_{i+1/2}$ and $M^R_{i+1/2}$ in the numerical flux $\hat{f}$ to obtain
$$
\begin{aligned}
\tilde{\hat{f}}_{i+1 / 2} 
&=\frac{1}{2}\left(f\left(\mathbf{U}_{i+1 / 2}^R\right)+f\left(\mathbf{U}_{i+1 / 2}^L\right)-\left(M_{i+1/2} \mathbf{U}_{i+1 / 2}^R-M_{i+1/2} \mathbf{U}_{i+1 / 2}^L\right)\right) \\
&=\frac{1}{2}\left(f\left(\mathbf{U}_{i+1 / 2}^R\right)+f\left(\mathbf{U}_{i+1 / 2}^L\right)-M_{i+1/2}\left(\mathbf{U}_{i+1 / 2}^R-\mathbf{U}_{i+1 / 2}^L\right)\right).
\end{aligned}
$$
By comparing $\tilde{\hat{f}}$ and $\hat{\hat{f}}$, 
it can be observed that the $\alpha_{i+1/2}$ in the classical Lax-Friedrichs flux formula is a specific choice for the $M_{i+1/2}$ 
in the Lax-Friedrichs splitting method.

\subsection{Constructing Steger-Warming flux for the scalar equation based on Steger-Warming splitting for system}
\hspace*{\fill} \par
From \ref{appendix-LF-FVS-LF-flux}, it is known that the Lax-Friedrichs splitting method can derive the classical Lax-Friedrichs flux format, 
which is not only used for system control equations but also commonly applied to scalar equations. 
Inspired by this, this section attempts to apply the Steger-Warming splitting method to scalar equations in order to 
derive a new numerical flux format suitable for scalar equations.  

\subsubsection{Steger-Warming Flux for One-Dimensional Scalar Equations}
The Jacobian-FVS for a system of equations first requires rewriting the flux function in the form
$F(\mathbf{U}) = A(\mathbf{U}) \cdot \mathbf{U}$.  
Similarly, for the 1D scalar equation, we make a similar rewrite:
\\
Let
\begin{align}
  a=K f^{\prime}(u),
\end{align} 
satisfy
\begin{align}
  f(u)=a u=K f^{\prime}(u) u,
\end{align}
where $K$ is an undetermined constant that varies with different equations (different $f(u)$). 
Generally, $a=a(x,t,u)$. 
\par
It should be emphasized that in scalar equations, $A$, $\Lambda$, and $\lambda$ are ``the same'', all ``equal to $a$''. 
Thus, the scalar form of the Steger-Warming splitting method is as follows: 
\begin{align}
a^{+}=\frac{a+\left|a\right|}{2},\ a^{-}=\frac{a-\left|a\right|}{2}.
\end{align}
If smoothness is considered: 
\begin{align}
a^{+}=\frac{a+\sqrt{\left|a\right|^2+\delta^2}}{2},\ a^{-}=\frac{a-\sqrt{\left|a\right|^2+\delta^2}}{2},
\end{align}
where $\delta$ is a small positive quantity, for example, $\delta$ can be taken as $10^{-8}$. 
\par
Taking the interface $x_{i+1/2}$ as an example, due to the multi-valuedness of $u$ at the interface $(u^L_{i+1/2}, u^R_{i+1/2})$, 
$a$ is also multi-valued at the interface, i.e., $(a^L_{i+1/2}, a^R_{i+1/2})$.
Accordingly, the Steger-Warming flux format for the 1D scalar equation is as follows:  
\begin{align}
\hat{f}_{i+1 / 2}^{S W} 
& =a_{i+1 / 2}^{L, +} u_{i+1 / 2}^L + a_{i+1 / 2}^{R, -} u_{i+1 / 2}^R \nonumber \\
& =\left(\frac{a_{i+1 / 2}^{L}+\left|a_{i+1 / 2}^{L}\right|}{2}\right) u_{i+1 / 2}^{L}+\left(\frac{a_{i+1 / 2}^{R}-\left|a_{i+1 / 2}^{R}\right|}{2}\right) u_{i+1 / 2}^{R} \nonumber \\
& = \frac{1}{2}\left(f\left(u_{i+1 / 2}^{R}\right)+f\left(u_{i+1 / 2}^{L}\right)-\left(\left|a_{i+1 / 2}^{R}\right| u_{i+1 / 2}^{R}-\left|a_{i+1 / 2}^{L}\right| u_{i+1 / 2}^{L}\right)\right)
.\end{align}
\par
The following discussion will focus on the linear advection equation and the nonlinear Burgers' equation in detail: 
\begin{itemize}
	\item Case 1:$\quad u_t + (c(x,t)u)_x=0$. 
  \\
	$$f(u)=c(x,t) u,\ f^{\prime}(u)=c(x,t),$$
	thus we take
	$$
	K=1,\ a=c(x,t).
	$$
  At this point, with $a_{i+1/2}^L=a_{i+1/2}^R=c(x_{i+1/2},\ \cdot)$, 
  the Steger-Warming flux splitting scheme for the linear scalar advection equation is given by:
  \begin{align}
    \hat{f}_{i+1 / 2}^{S W} 
    = \frac{1}{2}\left(c(x_{i+1/2},t)\left(u_{i+1 / 2}^{R} + u_{i+1 / 2}^{L}\right)-\left|c(x_{i+1/2},t)\right|\left( u_{i+1 / 2}^{R} - u_{i+1 / 2}^{L}\right)\right)
  .\end{align}
  Let's examine the time step $t^n$: 
  \begin{itemize}
    \item If $c(x_{i+1/2},t^n) > 0$, abbreviated as $c_{i+1/2}^n > 0$, \\
    then $a_{i+1/2} = c_{i+1/2}^n > 0$, 
    and thus $|a_{i+1/2}| = a_{i+1/2}$. \\
    Consequently, $a^{L,+}_{i+1/2} = a_{i+1/2}$, $a_{i+1/2}^{R,-} = 0$, \\ 
    and the Steger-Warming numerical flux at the interface is given by
    $$
    \hat{f}_{i+1 / 2}^{S W} = a_{i+1 / 2}^{L, +} u_{i+1 / 2}^L = c_{i+1/2}^{n} u_{i+1 / 2}^L.
    $$
    Given that $c_{i+1/2}^n > 0$, the local upwind direction at $x_{i+1/2}$ should be in the direction of $-x$, 
    meaning that the numerical flux should be calculated using the left state of the field function. 
    That is, take $\hat{f}_{i+1/2} = f(u^L_{i+1/2}) = c_{i+1/2}^{n} u_{i+1 / 2}^L$. \\ 
    It can be seen that $\hat{f}_{i+1 / 2}^{S W}$ satisfies the upwind property.  \\

    \item If $c(x_{i+1/2},t^n) < 0$, abbreviated as $c_{i+1/2}^n < 0$, \\
    then $a_{i+1/2} = c_{i+1/2}^n < 0$, and thus $|a_{i+1/2}| = -a_{i+1/2}$. \\
    Consequently, $a^{L,+}_{i+1/2} = 0$, $a_{i+1/2}^{R,-} = a_{i+1/2}$, \\
    and the Steger-Warming numerical flux at the interface is given by
    $$
    \hat{f}_{i+1 / 2}^{S W} = a_{i+1 / 2}^{R, -} u_{i+1 / 2}^R = c_{i+1/2}^{n} u_{i+1 / 2}^R.
    $$
    Given that $c_{i+1/2}^n < 0$, the local upwind direction at $x_{i+1/2}$ should be in the direction of $+x$, 
    meaning that the numerical flux should be calculated using the right state of the field function. 
    That is, take $\hat{f}_{i+1/2} = f(u^R_{i+1/2}) = c_{i+1/2}^{n} u_{i+1 / 2}^R$. \\
    It can be seen that $\hat{f}_{i+1 / 2}^{S W}$ satisfies the upwind property.  
  \end{itemize}
  \begin{remark}
    The fact that $c$ is a function of $x$ and $t$ implies that the upwind direction in different regions of the flow field 
    is not consistent and may change over time. When $c$ is only a function of $t$, 
    or more so when $c = \text{const}$, the above discussion results are still applicable, 
    but the entire flow field will share the same upwind direction (if $c = c(t)$, this upwind direction will change over time). 
  \end{remark}
  \hspace*{\fill}
	\item Case 2:$\quad u_t+\left(\frac{1}{2} u^2\right)_x=0$. 
	$$
  f(u)=\frac{1}{2} u^2,\ f^{\prime}(u)=u, 
	$$
	thus we take
	$$
	{\color{blue}{K=\frac{1}{2}}},\ a=Kf'(u)={\color{blue}{\frac{1}{2}}} u.
	$$
  The Steger-Warming flux for the Burgers' equation is given by 
  \begin{align}
    \hat{f}_{i+1 / 2}^{S W} 
     &= \frac{1}{2}\left(\frac{1}{2}\left|u_{i+1 / 2}^{R}\right|^2+\frac{1}{2}\left|u_{i+1 / 2}^{L}\right|^2-\left(\frac{1}{2}\left|u_{i+1 / 2}^{R}\right| u_{i+1 / 2}^{R}-\frac{1}{2}\left|u_{i+1 / 2}^{L}\right| u_{i+1 / 2}^{L}\right)\right) \nonumber \\
     &= \frac{1}{4}\left(\left|u_{i+1 / 2}^{R}\right|^2+\left|u_{i+1 / 2}^{L}\right|^2-\left(\left|u_{i+1 / 2}^{R}\right| u_{i+1 / 2}^{R}-\left|u_{i+1 / 2}^{L}\right| u_{i+1 / 2}^{L}\right)\right) \nonumber \\
     &= \frac{1}{4}\left(\left|u_{i+1 / 2}^{R}\right|\left(\left|u_{i+1 / 2}^{R}\right|-u_{i+1 / 2}^{R}\right)+\left|u_{i+1 / 2}^{L}\right|\left(\left|u_{i+1 / 2}^{L}\right|+u_{i+1 / 2}^{L}\right)\right)
  .\end{align}
  When applied to the nonlinear Burgers' equation, the Steger-Warming flux satisfies ``compatibility, Lipschitz continuity, and monotonicity.'' 
  This ensures that the spatial semi-discrete weak solution of the Burgers' equation, based on the Steger-Warming flux, 
  satisfies the cell entropy inequality and $L^2$ stability, as detailed below in section \ref{subsubSec-scalar-Steger-Warming-property}. 
\end{itemize}

\subsubsection{Steger-Warming Flux for Two-Dimensional Scalar Equations}\label{appendix-2D-scalar-SW-flux}
Similar to the techniques used in Section \ref{subSec-2D-DG-FVS}, the 2D problem is converted into several 1D problems 
along the outward normal direction by introducing the normal flux. 
\par
Consider the following two-dimensional scalar equation: 
\begin{align*}
  U_t + (f(U))_x +(g(U))_y = 0.
\end{align*}
Let $\mathbf{F}(U)=(f(U),g(U))$, and suppose the unit outward normal vector at a point on the cell boundary is $\mathbf{n}=(n_x,n_y)$, 
The normal flux is then defined as
\begin{align*}
  F_n(U)=\mathbf{n} \cdot \mathbf{F} = n_x f(U) + n_y g(U).
\end{align*}
To use the Jacobian-FVS, it is necessary to construct $a_n(x,t,U)$ such that
\begin{align*}
  F_n(U) = a_n \cdot U.
\end{align*}
Similar to the one-dimensional case, $a_n$ is chosen in the following form:
\begin{align*}
  a_n = K F_n^{\prime}(U),
\end{align*}
where $K$ is a constant to be determined. 
\\
The splitting method for $a_n$ is consistent with the previously described one-dimensional case. 
Consequently, the Steger-Warming flux for the two-dimensional scalar equation at a specific point on the cell interface, 
denoted by $\mathbf{G_k}$, is delineated as follows: 
{\scriptsize 
\begin{align}
\hat{F}_{n}^{S W}(U^{int}(\mathbf{G_k},t),U^{ext}(\mathbf{G_k},t)) 
& =a_{n}^{int, +} \cdot U^{int}(\mathbf{G_k},t) + a_{n}^{ext, -} \cdot U^{ext}(\mathbf{G_k},t) \nonumber \\
& =\left(\frac{a_{n}^{int}+\left|a_{n}^{int}\right|}{2}\right) U^{int}(\mathbf{G_k},t)+\left(\frac{a_{n}^{ext}-\left|a_{n}^{ext}\right|}{2}\right) U^{ext}(\mathbf{G_k},t) \nonumber \\
& = \frac{1}{2}\left(F_n\left(U(\mathbf{G_k},t)^{ext}\right)+F_n\left(U^{int}(\mathbf{G_k},t)\right)-\left(\left|a_{n}^{ext}\right| U^{ext}(\mathbf{G_k},t)-\left|a_{n}^{int}\right| U^{int}(\mathbf{G_k},t)\right)\right)
.\end{align}
}
\begin{itemize}
  \item The two-dimensional conservative linear scalar transport equation \\
  $U_t + (\alpha(x,y,t)U)_x + (\beta(x,y,t)U)_y = 0$. 
  \begin{align*}
    F_n(U) &= n_x \cdot \alpha(x,y,t) U + n_y \cdot \beta(x,y,t) U,\\ 
    F_n^{\prime}(U) &= n_x \cdot \alpha(x,y,t) + n_y \cdot \beta(x,y,t),
  \end{align*}
	hence we take
	$$
	K=1,\ a_n=n_x \cdot \alpha(x,y,t) + n_y \cdot \beta(x,y,t).
	$$
  Therefore, $a^{int}(\mathbf{G_k},\cdot)=a^{ext}(\mathbf{G_k},\cdot)=n_x \cdot \alpha(\mathbf{G_k},\cdot) + n_y \cdot \beta(\mathbf{G_k},\cdot)$, \\
  the Steger-Warming flux for the two-dimensional linear scalar transport equation is given by
  \begin{align}
    &\hat{F}_{n}^{S W}(U^{int}(\mathbf{G_k},t),U^{ext}(\mathbf{G_k},t)) \nonumber \\
    &= \frac{1}{2}\bigg(\left(n_x \cdot \alpha(\mathbf{G_k},t) + n_y \cdot \beta(\mathbf{G_k},t)\right) \cdot \left(U^{ext}(\mathbf{G_k},t)+U^{int}(\mathbf{G_k},t)\right) \nonumber \\
    &-\left|n_x \cdot \alpha(\mathbf{G_k},t) + n_y \cdot \beta(\mathbf{G_k},t)\right|\left(U^{ext}(\mathbf{G_k},t)- U^{int}(\mathbf{G_k},t)\right)\bigg)
  .\end{align}
  \item The two-dimensional nonlinear Burgers' equation\ $U_t + (\frac{1}{2}U^2)_x + (\frac{1}{2}U^2)_y = 0$. 
  \begin{align*}
    F_n(U) &= \left(n_x + n_y\right) \cdot \frac{1}{2}U^2,\\ 
    F_n^{\prime}(U) &= \left(n_x + n_y\right) \cdot U,
  \end{align*}
	hence we take 
	$$
	{\color{blue}{K=\frac{1}{2}}},\ a_n= {\color{blue}{\frac{1}{2}}} \left(n_x + n_y\right) \cdot U.
	$$
  the Steger-Warming flux for the 2D Burgers' equation is given by
  \begin{align}
    &\hat{F}_{n}^{S W}(U^{int}(\mathbf{G_k},t),U^{ext}(\mathbf{G_k},t))  \\
    &= \frac{1}{4}\left(n_x + n_y\right)\bigg(\left(\left(U^{ext}(\mathbf{G_k},t)\right)^2+\left(U^{int}(\mathbf{G_k},t)\right)^2\right) \nonumber \\
    &- \left(\left|U^{ext}(\mathbf{G_k},t)\right| U^{ext}(\mathbf{G_k},t) - \left|U^{int}(\mathbf{G_k},t)\right| U^{int}(\mathbf{G_k},t)\right)\bigg)
.\end{align}
\end{itemize}

\subsubsection{Properties of the Steger-Warming Flux for the Scalar Equation}\label{subsubSec-scalar-Steger-Warming-property}
This section discusses the properties of the Steger-Warming flux 
within the scope of one-dimensional conservative linear scalar transport equations 
and nonlinear Burgers' equations.
\par
The one-dimensional scalar equation Steger-Warming flux is given by 
\begin{align*}
  \hat{f}_{i+1 / 2}^{S W} 
  = \frac{1}{2}\left(f\left(u_{i+1 / 2}^{R}\right)+f\left(u_{i+1 / 2}^{L}\right)-\left(\left|a_{i+1 / 2}^{R}\right| u_{i+1 / 2}^{R}-\left|a_{i+1 / 2}^{L}\right| u_{i+1 / 2}^{L}\right)\right)
.\end{align*}

\begin{itemize}
  \item Conservative linear scalar transport equation $u_t+(c(x,t)u)_x=0$, 
  $$a=c(x,t);$$
  \item Nonlinear Burgers' equation $u_t+(\frac{1}{2}u^2)_x=0$, 
  $$a=\frac{1}{2}u.$$
\end{itemize}
\par
When applied to conservative linear scalar transport equations and nonlinear Burgers' equations, 
the Steger-Warming flux scheme satisfies ``consistency, Lipschitz continuity, and monotonicity'', 
which will be introduced and proven below:
\par
Since all discussions are at $x_{i+1/2}$, the subscript ``$i+1/2$'' is omitted for all variables in the following processes.
\begin{itemize}
  \item Consistency, that is, it satisfies $\hat{f}^{SW}(u,u) = f(u)$. 
  \begin{proof}
    Let
    $$
    u^L=u^{R}=u^{*}.
    $$
    Since $a=a(x,t,u)$, then
    $$
    a^L(\cdot,\ \cdot,\ u^L)=a^R(\cdot,\ \cdot,\ u^R)=a(\cdot,\ \cdot,\ u^*).
    $$
    Let us denote $a^*=a(\cdot,\ \cdot,\ u^*)$, 
    thus we have
    \begin{align*}
      \hat{f}^{S W} (u^*,u^*) &= \frac{1}{2}\left(f\left(u^*\right)+f\left(u^*\right)-\left(\left|a^{*}\right| u^*-\left|a^{*}\right| u^*\right)\right) \\
      &= \frac{1}{2} \cdot 2f\left(u^*\right) \\
      &= f\left(u^*\right).
    \end{align*}
  \end{proof}
  \item When $f(u;x,t)$ is Lipschitz continuous function of $u$, the Steger-Warming numerical flux $\hat{f}^{SW}(u^L,u^R)$ is Lipschitz continuous with respect to both $u^L$ and $u^R$.  
  \begin{proof}
    \hspace*{\fill}\\
      We first establish the Lipschitz continuity of $\hat{f}^{SW}(u^L,u^R)$ with respect to $u^L$. \\
      Let $\Delta$ is a constant that can be either positive or negative, and define 
      \begin{align}a^{L+\Delta}:=a(\cdot,\ \cdot,\ u^L+\Delta),\end{align}
      then
      \begin{align}
        \hat{f}^{SW}(u^L+\Delta,u^R)=
        \frac{1}{2}\left(f\left(u^{R}\right)+f\left(u^{L}+\Delta\right)-\left(\left|a^{R}\right| u^{R}-\left|a^{L+\Delta}\right| \left(u^{L}+\Delta\right)\right)\right),
      \end{align}
      hence
      \begin{align}\label{SWflux-uL-Lipschitz-1}
        &\hat{f}^{SW}(u^L+\Delta,u^R)-\hat{f}^{SW}(u^L,u^R) \nonumber \\
        &=
        \frac{1}{2}\left(f\left(u^{L}+\Delta\right)+\left|a^{L+\Delta}\right| \left(u^{L}+\Delta\right)-f\left(u^{L}\right)-\left|a^{L}\right| u^{L}\right).
      \end{align}
      Perform the following simple decomposition:
      \begin{align}\label{SWflux-uL-Lipschitz-2}
        &\hat{f}^{SW}(u^L+\Delta,u^R)-\hat{f}^{SW}(u^L,u^R) \nonumber \\
        &=\frac{1}{2}\left(f\left(u^{L}+\Delta\right)-f\left(u^{L}\right) + 
        \left|a^{L+\Delta}\right| \left(u^{L}+\Delta\right)
        {\color{red}{-\left|a^{L}\right| \left(u^{L}+\Delta\right) + \left|a^{L}\right| \left(u^{L}+\Delta\right)}} 
        - \left|a^{L}\right| u^{L} \right) \nonumber \\ 
        &=\frac{1}{2}\left(\underbrace{\left[f\left(u^{L}+\Delta\right)-f\left(u^{L}\right)\right]}_{T_1} + 
        \underbrace{\left(\left|a^{L+\Delta}\right|- \left|a^{L}\right|\right) \cdot \left(u^{L}+\Delta\right)}_{T_2}
        +\underbrace{\left|a^{L}\right| \cdot \Delta}_{T_3} \right). 
      \end{align}
      Since $f(u;x,t)$ is Lipschitz continuous with respect to $u$, then
     $$
      \exists L_1>0, s.t.\ |T_1| \leq L_1 |\Delta|.
     $$
      $f$ being Lipschitz continuous implies that its derivative $f'$ is bounded. Consequently, $a$ is bounded as well (given $a = K f'(u)$, where $K$ is a constant), and thus we have
     $$
     \exists L_3>0, s.t.\ |T_3| \leq L_3 |\Delta|.
     $$
      Both $a^{L+\Delta}$ and $a^L$ are bounded, which implies that the difference $|a^{L+\Delta}| - |a^L|$ is also bounded. Noting that $u^L$ is a finite state, we deduce that
     $$
     \exists L_2>0, s.t.\ |T_2| \leq L_2 |\Delta|.
     $$
      Hence, there exists a constant $M > 0$ (e.g.\ take $M=\frac{3}{2} \cdot \max\{L_1,L_2,L_3\}$) such that 
      \begin{align}
        \left|\hat{f}^{SW}(u^L+\Delta,u^R)-\hat{f}^{SW}(u^L,u^R)\right|
        \leq
        \frac{1}{2}(|T_1|+|T_2|+|T_3|)
        \leq
        M\left|\Delta\right|,
      \end{align}
      which implies that $\hat{f}^{SW}(u^L,u^R)$ is Lipschitz continuous with respect to $u^L$ 
      when $f(u;x,t)$ is Lipschitz continuous function of $u$.  
      \par \noindent
      Similarly, it can be shown that $\hat{f}^{SW}(u^L,u^R)$ is Lipschitz continuous with respect to $u^R$ under the same conditions.  
  \end{proof}
  \item Monotonicity: The Steger-Warming numerical flux $\hat{f}^{SW}(u^L,u^R)$ is non-decreasing with respect to the first variable $u^L$ 
  and non-increasing with respect to the second variable $u^R$. 
  This property is succinctly denoted as $\hat{f}^{SW}(\uparrow,\downarrow)$.  
  \begin{proof}
    \hspace*{\fill} \\
    We first prove that $\hat{f}^{SW}(u^L,u^R)$ is non-decreasing with respect to the first variable $u^L$.
    \\
    Following the proof of Lipschitz continuity, we have
    \begin{align*}
      &\hat{f}^{SW}(u^L+\Delta,u^R)-\hat{f}^{SW}(u^L,u^R) \nonumber \\
      &=
      \frac{1}{2}\left(f\left(u^{L}+\Delta\right)+\left|a^{L+\Delta}\right| \left(u^{L}+\Delta\right)-f\left(u^{L}\right)-\left|a^{L}\right| u^{L}\right).
    \end{align*}
    However, here we set $\Delta>0$. 
    \begin{itemize}
      \item Conservative linear scalar transport equation: $u_t+(c(x,t)u)_x=0$, $f(u)=c(x,t)u$, $a=c(x,t)$. \\
      Since $a$ is independent of $u$, it follows that $a^{L+\Delta} = a^L$. Consequently,
      \begin{align*}
        &\hat{f}^{SW}(u^L+\Delta,u^R)-\hat{f}^{SW}(u^L,u^R) \nonumber \\
        &=
        \frac{1}{2}\left(c\cdot\left(u^{L}+\Delta\right)+\left|c\right| \left(u^{L}+\Delta\right)-c\cdot\left(u^{L}\right)-\left|c\right| u^{L}\right) \\
        &= \frac{1}{2}\left(c+|c|\right)\Delta
        \geq 0.
      \end{align*}
      This implies that $\hat{f}^{SW}(u^L,u^R)$ is non-decreasing with respect to the first variable $u^L$ 
      when applied to the conservative linear scalar transport equation.
      \item Nonlinear Burgers' equation: $u_t+(\frac{1}{2}u^2)_x=0$, $f(u)=\frac{1}{2}u^2$, $a=\frac{1}{2}u$. 
      \begin{align*}
        &\hat{f}^{SW}(u^L+\Delta,u^R)-\hat{f}^{SW}(u^L,u^R) \nonumber \\
        &=\frac{1}{2}\left(\frac{1}{2}\left(u^{L}+\Delta\right)^2 + \frac{1}{2}\left|u^{L}+\Delta\right|\left(u^{L}+\Delta\right)-\frac{1}{2}\left(u^{L}\right)^2 - \frac{1}{2}\left|u^{L}\right|u^{L}\right) \nonumber \\
        &=\frac{1}{4}\left(\left|u^L+\Delta\right|\left(\left|u^L+\Delta\right|+u^L+\Delta\right) - \left|u^L\right|\left(\left|u^L\right|+u^L\right)\right).
      \end{align*}
      We introduce an auxiliary function $g(u) = u^2 + |u|u$. Then, we have
      \begin{align*}
        \hat{f}^{SW}(u^L+\Delta,u^R)-\hat{f}^{SW}(u^L,u^R)
        =\frac{1}{4}\left(g(u^L+\Delta) - g(u^L)\right). 
      \end{align*}
      It is clear that $g(u)$ is non-decreasing with respect to $u$. \\
      Given that $\Delta > 0$ and $u^L + \Delta > u^L$, \\
      it follows that 
      $\hat{f}^{SW}(u^L+\Delta,u^R) \geq \hat{f}^{SW}(u^L,u^R)$, 
      \\
      which demonstrates the non-decreasing nature of $\hat{f}^{SW}(u^L,u^R)$ with respect to the first variable $u^L$.
    \end{itemize}
    \par \noindent
    By similar reasoning, it can be proven that the Steger-Warming numerical flux $\hat{f}^{SW}(u^L,u^R)$, 
    when applied to either the conservative linear scalar transport equation or the nonlinear Burgers' equation, 
    is non-increasing with respect to the second variable $u^R$. 
  \end{proof}
\end{itemize}
\par
When a numerical flux scheme satisfies the aforementioned properties of ``consistency, Lipschitz continuity, and monotonicity,'' 
it can be proven that the DG (Discontinuous Galerkin) weak solution based on this flux scheme satisfies the cell entropy inequality 
and $L^2$ stability. This proof process is standardized, and readers are referred to the literature \cite{ref47,ref48} for details.

\subsubsection{Numerical Experiments}
We conducted numerical experiments to compare the Steger-Warming flux scheme discussed in this appendix 
with the classical Lax-Friedrichs flux scheme within the context of conservative linear scalar transport equations 
and nonlinear Burgers' equation. 
\\
$\bullet$ Accuracy test
\begin{example}
Linear transport equations.
\par \noindent
\begin{itemize}
  \item Control Eqs: 
  $
  u_t + (\sin(\omega t)u)_x = 0,\ \omega=\pi;
  $\\
  Computational domain:  
  $
  \Omega \times [0,T_{end}] = [0,2\pi] \times [0,20];
  $\\
  $\mathrm{I.C.}\ $
  $
  \begin{aligned}
      u_0(x)=\sin(x);
  \end{aligned}
  $\\
  $\mathrm{B.C.}\ ${\color{blue}{periodic boundary conditions; }}\\
  True solutions:
  $
  \begin{aligned}
    u(x,t)=u_0(x+\frac{1}{\omega}(\cos(\omega t) - 1))
  \end{aligned}
  $;\\
  Note: The smallest positive period of the solution $u$ is $T^*=\frac{2\pi}{\omega}=2$, hence, it implies simulating up to 10 periods that $T_{end}$ is taken as $20$. \\
  Temporal discretization format: TVD-RK3; \\
  CFL=0.1. \\
  $L^2,\ L^1$\ numerical errors and convergence orders
  with $P^2$-polynomial approximation 
  are summarized respectively in Table \ref{table-1D-scalar-Transport-sin(t)-Accuracy-SW-P2} 
  and Figure \ref{Fig.SW-LF-sin(t)}.  
  {\tiny
\begin{table}[htbp]
  \centering
  \caption{\tiny  $P^2$-Steger-Warming-scalar-DG and $P^2$-Lax-Friedrichs-scalar-DG both using equally spaced cells. $u_t+(\sin(t)u)_x=0$ with smooth initial conditions: 
  $u_0(x) = \sin(x)$. $L^2,\ L^1$\ errors for $u$. }
  \label{table-1D-scalar-Transport-sin(t)-Accuracy-SW-P2}
  \begin{tabular}{l|llllllll}  
    \hline  & $Mesh$ & 20 & 40 & 80 & 160 & 320 \\ 
    \hline 
    S-W & $L^{2}$-error      & 5.0283E-04     & 6.3601E-05     & 0.7973E-05     & 0.9975E-06     & 1.2471E-07 \\
    L-F & $L^{2}$-error      & 6.1941E-04     & 8.8553E-05     & 1.2804E-05     & 1.8863E-06     & 2.8289E-07  \\  
    \hline 
    S-W & $L^{1}$-error      & 1.1090E-03     & 1.4070E-04     & 1.7664E-05     & 2.2117E-06     & 2.7663E-07 \\ 
    L-F & $L^{1}$-error      & 1.3487E-03     & 1.9093E-04     & 2.7256E-05     & 3.9575E-06     & 5.8503E-07 \\
    \hline 
    \end{tabular}
    \end{table}
  }
  \begin{figure}[htbp]
    \vspace{0.3cm}
    \begin{center}
      \begin{minipage}{0.49\linewidth}
        \centerline{\includegraphics[width=1\linewidth]{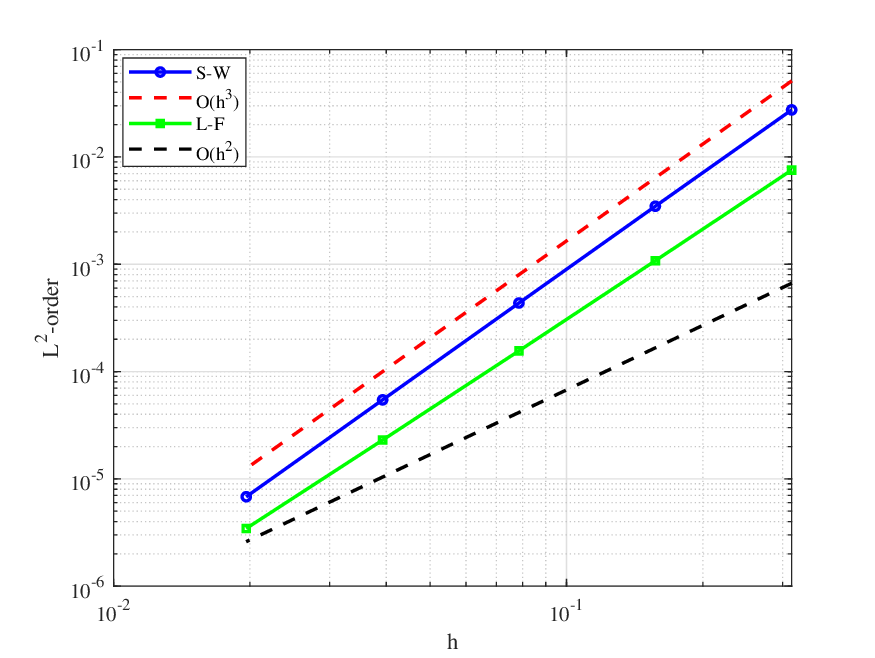}}
      \end{minipage}
      \hfill
      \begin{minipage}{0.49\linewidth}
        \centerline{\includegraphics[width=1\linewidth]{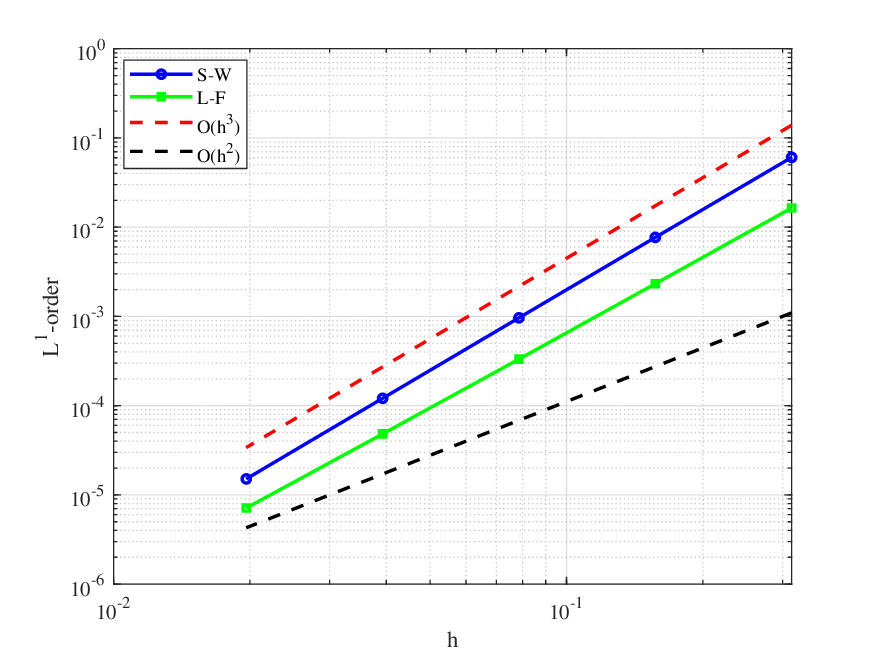}}
      \end{minipage}
      \vfill
      \vspace{0.2cm}
      \begin{minipage}{0.49\linewidth}
        \small
        \centerline{(a)\ $L^2$-convergence order}
      \end{minipage}
      \hfill
      \begin{minipage}{0.49\linewidth}
        \small
        \centerline{(b)\ $L^1$-convergence order}
      \end{minipage}
    \end{center}
    \caption{$P^2$-Steger-Warming-scalar-DG and $P^2$-Lax-Friedrichs-scalar-DG both using equally spaced cells. $u_t+(\sin(t)u)_x=0$ with smooth initial conditions: 
  $u_0(x) = \sin(x)$. $L^2,\ L^1$\ convergence orders for $u$.}
    \vspace{0.3cm}
    \label{Fig.SW-LF-sin(t)}
  \end{figure}
  \hfill \par \noindent
  \item Control Eqs:  
  $
  u_t + (\sin(x)u)_x = 0;
  $\\
  Computational domain:   
  $
  \Omega \times [0,T_{end}] = [0,2\pi] \times [0,1];
  $\\
  $\mathrm{I.C.}\ $ 
  $
  \begin{aligned}
      u_0(x)=1;
  \end{aligned}
  $\\
  $\mathrm{B.C.}\ ${\color{blue}{periodic boundary conditions; }}
  \\
  True solutions: 
  $
  \begin{aligned}
    u(x,t)=\frac{\sin(2\arctan(e^{-t}\tan(\frac{x}{2})))}{\sin(x)}
  \end{aligned}
  $;\\
  Temporal discretization format: RK4; \\
  CFL=0.05. \\
  $L^2,\ L^1$\ numerical errors 
  with $P^5$-polynomial approximation 
  are summarized respectively in Table \ref{table-1D-scalar-Transport-sin(x)-Accuracy-SW-P5}
  and Figure \ref{Fig.SW-LF-sin(x)} .
  {\tiny
  \begin{table}[htbp]
      \centering
      \caption{\tiny  $P^5$-Steger-Warming-scalar-DG and $P^5$-Lax-Friedrichs-scalar-DG both using equally spaced cells. $u_t+(\sin(x)u)_x=0$ with smooth initial conditions: 
      $u_0(x) = 1$. $L^2,\ L^1$\ errors and convergence orders  for  $u$. }
      \label{table-1D-scalar-Transport-sin(x)-Accuracy-SW-P5}
    \begin{tabular}{l|llllllll}  
      \hline  & $Mesh$ & 20 & 40 & 80 & 160 & 320 \\ 
      \hline 
      S-W & $L^{2}$-error      & 1.1214E-06     & 3.1954E-08     & 5.3956E-10     & 9.4229E-12     & 5.7994E-13 \\  
      L-F & $L^{2}$-error      & 1.2660E-06     & 3.3854E-08     & 4.7961E-10     & 7.8348E-12     & 7.5041E-13 \\  
      \hline 
      S-W & $L^{1}$-error      & 1.1410E-06     & 3.0824E-08     & 4.8678E-10     & 1.0050E-11     & 0.8858E-12 \\ 
      L-F & $L^{1}$-error      & 1.3414E-06     & 3.2669E-08     & 4.6412E-10     & 0.9267E-11     & 1.0540E-12 \\ 
      \hline 
      \end{tabular}
      \end{table}
    }
    \begin{figure}[htbp]
      \vspace{0.3cm}
      \begin{center}
        \begin{minipage}{0.49\linewidth}
          \centerline{\includegraphics[width=1\linewidth]{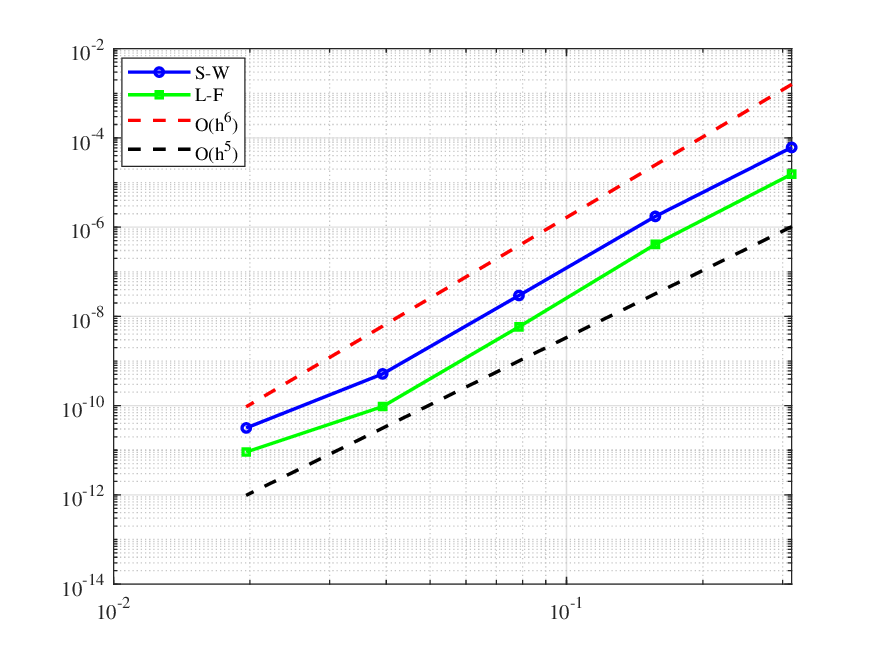}}
        \end{minipage}
        \hfill
        \begin{minipage}{0.49\linewidth}
          \centerline{\includegraphics[width=1\linewidth]{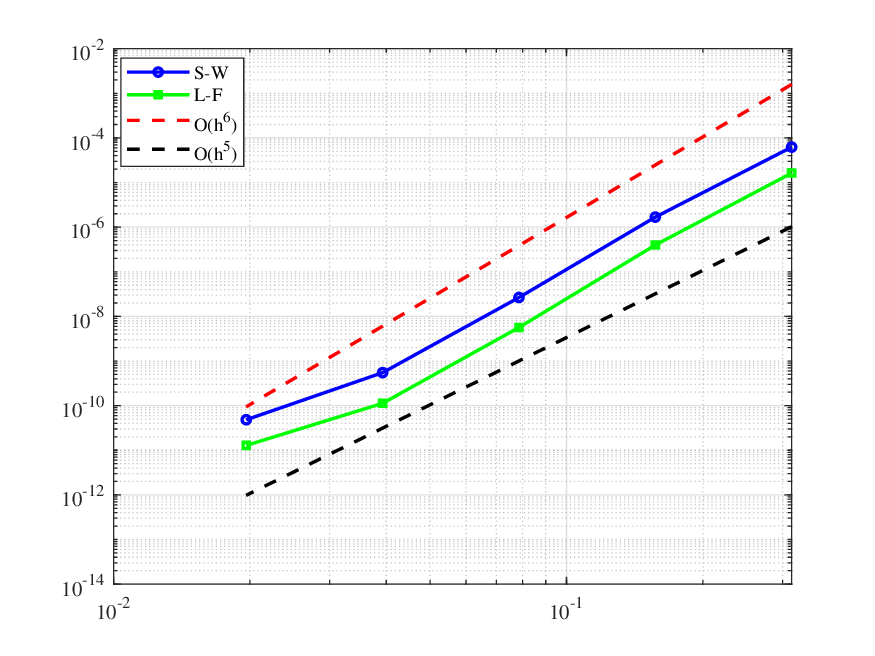}}
        \end{minipage}
        \vfill
        \vspace{0.2cm}
        \begin{minipage}{0.49\linewidth}
          \small
          \centerline{(a)\ $L^2$-convergence order}
        \end{minipage}
        \hfill
        \begin{minipage}{0.49\linewidth}
          \small
          \centerline{(b)\ $L^1$-convergence order}
        \end{minipage}
      \end{center}
      \caption{$P^5$-Steger-Warming-scalar-DG and $P^5$-Lax-Friedrichs-scalar-DG both using equally spaced cells. $u_t+(\sin(x)u)_x=0$ with smooth initial conditions: 
    $u_0(x) = 1$. $L^2,\ L^1$\ convergence orders for $u$.}
      \vspace{0.3cm}
      \label{Fig.SW-LF-sin(x)}
    \end{figure}
\end{itemize}
\end{example}
\begin{example}\label{example-Burgers}
Nonlinear Burgers' equation. 
\\
$\bullet$\ 1D-Control Eqs:  
  $
  u_t + (\frac{1}{2}u^2)_x = 0;
  $\\
  Computational domain:   
  $
  \Omega \times [0,T_{end}] = [0,2\pi] \times [0,0.6];
  $\\
  Note: the exact solution remains smooth during $t\in [0,0.6]$; 
  \\
  $\mathrm{I.C.}\ $ 
  $
  \begin{aligned}
      u_0(x)=\sin(x);
  \end{aligned}
  $\\
  $\mathrm{B.C.}\ ${\color{blue}{periodic boundary conditions; }}
  \\
  True solutions: 
  $
  \begin{aligned}
    u(x,t)=u_0(x^*),\ x^* \text{\ satisfies\ } x^*+u_0(x^*)\cdot t=x
  \end{aligned}
  $;\\
  Temporal discretization format: RK4; \\
  CFL=0.05. \\
  $L^2,\ L^1$\ numerical errors and convergence orders
  with $P^3$-polynomial approximation 
  are summarized respectively in Table \ref{table-1D-scalar-Burgers-Accuracy-SW-P5}
  and Figure \ref{Fig.SW-LF-Burgers-1D}.
{\tiny

    \begin{table}[htbp]
      \centering
      \caption{\tiny $P^5$-Steger-Warming-scalar-DG and $P^5$-Lax-Friedrichs-scalar-DG both using equally spaced cells. $u_t+(\frac{1}{2}u^2)_x=0$ with smooth initial conditions: 
      $u_0(x) = \sin(x)$. $L^2,\ L^1$\ errors and convergence orders  for  $u$. }
      \label{table-1D-scalar-Burgers-Accuracy-SW-P5}
      \begin{tabular}{l|llllllll}
      \hline  & $Mesh$ & 20 & 40 & 80 & 160 & 320 \\ 
      \hline 
      S-W & $L^{2}$-error      & 5.342170E-06     & 1.62220980E-07     & 2.285972076E-09     & 3.9838927E-11     & 8.2960E-13  \\ 
      L-F & $L^{2}$-error      & 5.342175E-06     & 1.62220986E-07     & 2.285972077E-09     & 3.9838921E-11     & 8.2980E-13  \\ 
      \hline 
      S-W & $L^{1}$-error      & 4.25502E-06     & 1.0785780E-07     & 1.4120418E-09     & 2.2865E-11     & 1.2853E-12 \\ 
      L-F & $L^{1}$-error      & 4.25507E-06     & 1.0785784E-07     & 1.4120415E-09     & 2.2864E-11     & 1.2861E-12 \\ 
      \hline 
      \end{tabular}
      \end{table}
    }
    \begin{figure}[htbp]
      \vspace{0.3cm}
      \begin{center}
        \begin{minipage}{0.49\linewidth}
          \centerline{\includegraphics[width=1\linewidth]{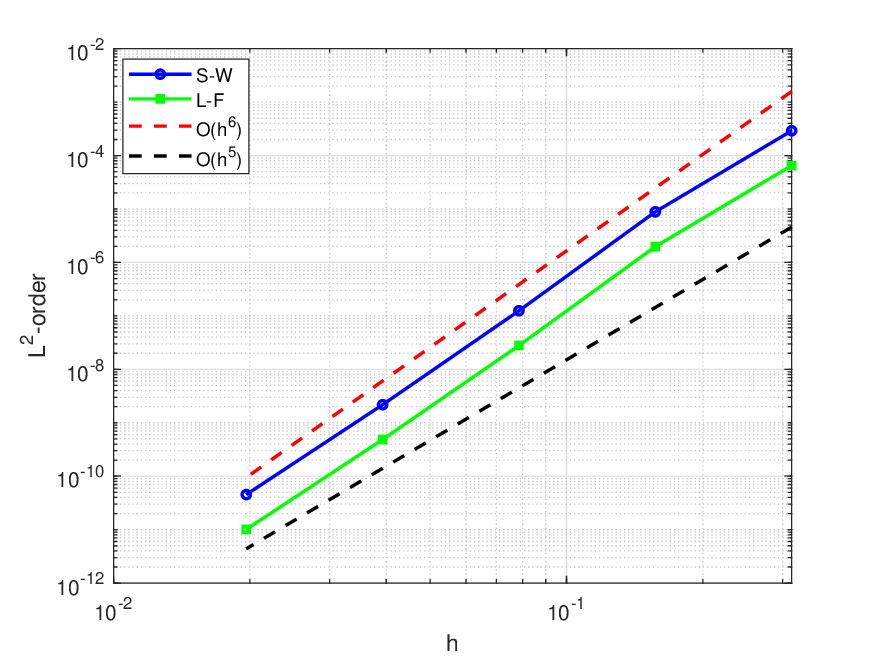}}
        \end{minipage}
        \hfill
        \begin{minipage}{0.49\linewidth}
          \centerline{\includegraphics[width=1\linewidth]{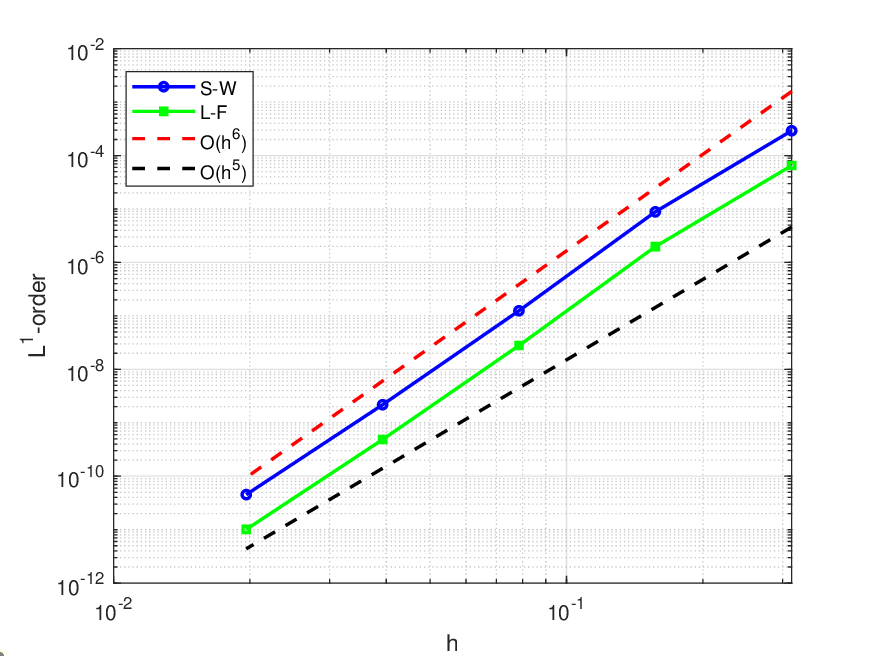}}
        \end{minipage}
        \vfill
        \vspace{0.2cm}
        \begin{minipage}{0.49\linewidth}
          \small
          \centerline{(a)\ $L^2$-convergence order}
        \end{minipage}
        \hfill
        \begin{minipage}{0.49\linewidth}
          \small
          \centerline{(b)\ $L^1$-convergence order}
        \end{minipage}
      \end{center}
      \caption{$P^5$-Steger-Warming-scalar-DG and $P^5$-Lax-Friedrichs-scalar-DG both using equally spaced cells. $u_t+(\frac{1}{2}u^2)_x=0$ with smooth initial conditions: 
    $u_0(x) = \sin(x)$. $L^2,\ L^1$\ convergence orders for $u$.}
      \vspace{0.3cm}
      \label{Fig.SW-LF-Burgers-1D}
    \end{figure}
    \hfill \\
  $\bullet$\ 2D-Control Eqs:  
  $
  U_t + (\frac{1}{2}U^2)_x + (\frac{1}{2}U^2)_y = 0;
  $\\
  Computational domain:   
  $
  \Omega \times [0,T_{end}] = \{[0,4] \times [0,4]\} \times [0,\frac{0.5}{\pi}];
  $\\
  Note: the exact solution remains smooth during $t\in [0,\frac{0.5}{\pi}]$; 
  \\
  $\mathrm{I.C.}\ $ 
  $
  \begin{aligned}
      U_0(x)=\sin(\frac{\pi}{2}x+\frac{\pi}{2}y);
  \end{aligned}
  $\\
  $\mathrm{B.C.}\ ${\color{blue}{periodic boundary conditions; }}
  \\
  By performing the variable substitution $\xi = x + y$ on $U_0(x, y)$, we obtain $\tilde{u}_0(\xi)=\sin(\frac{\pi}{2}\xi)$, 
  then the true solution is given as: \\
  $
  \begin{aligned}
    U(x,y,t)=\tilde{u}_0(\xi^*),\ \xi^* \text{\ satisfies\ } 
    \xi= \xi^* + 2\tilde{u}_0(\xi^*)\cdot t
  \end{aligned}
  $, where $\xi=x+y$; \\
  Temporal discretization format: TVD-RK3; \\
  CFL=0.05. \\
  $L^{\infty},\ L^2,\ L^1$\ numerical errors and convergence orders
  with $P^2$-polynomial approximation based on  scalar Steger-Warming flux
  are summarized in Table \ref{table-2D-Burgers-Accuracy-SW-P2}. 
  {\tiny
  \begin{table}[htbp]
    \centering
    \caption{\tiny $P^2$-Steger-Warming-scalar-DG using uniform {\color{blue}{rectangular}} meshes. 2D-Burgers' Equation with smooth initial condition:  
    $U_0(x)=\sin(\frac{\pi}{2}x+\frac{\pi}{2}y)$. $L^{\infty},\ L^2,\ L^1$\ errors and convergence orders  for $U$. }
    \label{table-2D-Burgers-Accuracy-SW-P2}
    \begin{tabular}{lllllllll}
    \hline  & $Mesh$ & 15$\times$15 & 30$\times$30 & 60$\times$60 & 120$\times$120 & 240$\times$240 \\
    \hline  & $L^{\infty}$-error & 6.5224E-02   & 1.1657E-02   & 1.5466E-03   & 1.8680E-04   & 3.0712E-05  \\
            & $L^{\infty}$-order & —        & 2.4843    & 2.9140    & 3.0495    & 2.6046  \\
      S-W   & $L^{2}$-error      & 5.1671E-02   & 6.5497E-03   & 8.5550E-04   & 1.1426E-04   & 1.5193E-05   \\
            & $L^{2}$-order      & —        & 2.9798    & 2.9366    & 2.9044    & 2.9109  \\
            & $L^{1}$-error      & 1.0390E-01   & 1.1668E-02   & 1.6231E-03   & 2.1695E-04   & 2.7951E-05 \\
            & $L^{1}$-order      & —        & 3.1546    & 2.8456    & 2.9034    & 2.9564 \\
    \hline  
    \end{tabular}
    \end{table}
    }
\end{example}
The one-dimensional and two-dimensional numerical results indicate that the $P^K$-DG method, based on the Steger-Warming flux scheme, 
is capable of achieving the optimal convergence order of $K+1$. Furthermore,
the accuracy test results in one-dimensional case indicate that for scalar equations, the dissipation of the Steger-Warming flux is slightly lower than 
that of the classical Lax-Friedrichs flux. This is consistent with the statement that ``in system control equations, 
the dissipation of the Steger-Warming splitting is less than that of the Lax-Friedrichs splitting''.
\begin{remark}
  Due to the excessive dissipation of the Lax-Friedrichs flux, it paradoxically results in a slightly smaller $L^{\infty}$-error 
  compared to the Steger-Warming flux. This is primarily because the $L^{\infty}$-error is more sensitive to peaks than $L^2$-error and $L^1$-error; 
  hence, a scheme with higher dissipation may often achieve a smaller $L^{\infty}$-error. 
  The degree of dissipation of numerical fluxes becomes more pronounced in high-order schemes, 
  which is why this paper selects the $P^5$-DG method for testing, 
  in order to more clearly reflect the dissipative characteristics of the L-F and S-W fluxes respectively. 
  Furthermore, the dissipative properties of the flux scheme become evident in long-time numerical simulations. 
  At such times, using lower-order polynomial approximation, such as $P^2$-polynomial approximation, can distinguish the degree of dissipation between 
  the S-W flux and the L-F flux, as we have done in the test case 
  ``$u_t + (sin(\omega t)u)_x = 0,\ t\in[0,20]$''.
\end{remark}
\par \noindent
\hfill \\
$\bullet$ Numerical oscillation test
\begin{example}
We continue to utilize the one-dimensional test case presented in Example \ref{example-Burgers}, 
with the simulation duration extended to {\color{red}{$T_{\text{end}} = 2.6$}}, 
at which point discontinuities have developed. No additional limiters are incorporated to amend the numerical solutions, 
allowing for observing the inherent performance differences of the Steger-Warming and Lax-Friedrichs numerical fluxes near discontinuities. 
Contrastive tests were performed under four distinct conditions: coarse and fine meshes, 
as well as low- and high-order polynomial approximations ($P^1$ and $P^3$). The results are demonstrated in Figure \ref{Fig.SW-LF-Burgers-discontinuity-1D}. 
\begin{figure}[htbp]
  \vspace{0.3cm}
  \begin{center}
    \begin{minipage}{0.49\linewidth}
      \centerline{\includegraphics[width=1\linewidth]{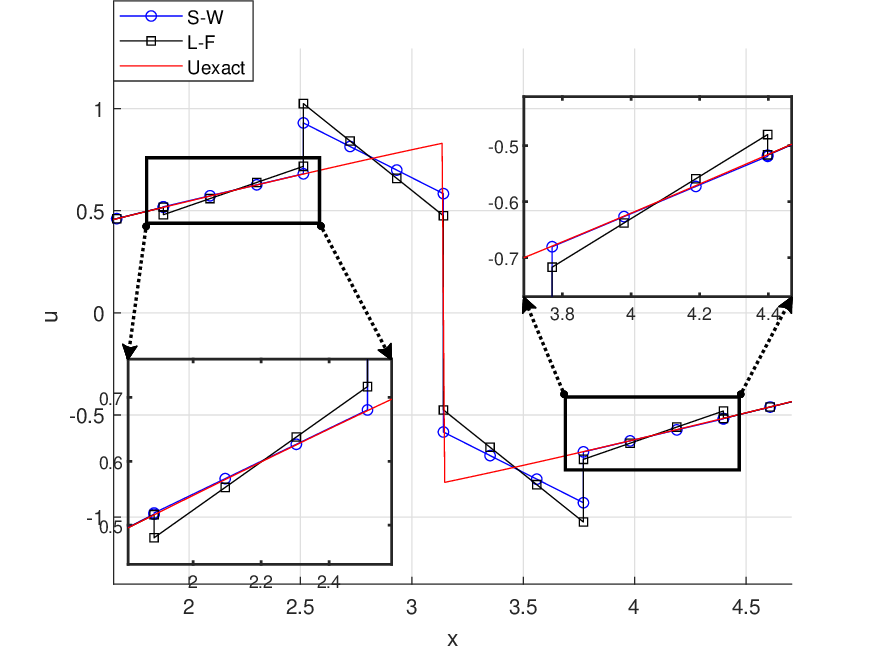}}
    \end{minipage}
    \hfill
    \begin{minipage}{0.49\linewidth}
      \centerline{\includegraphics[width=1\linewidth]{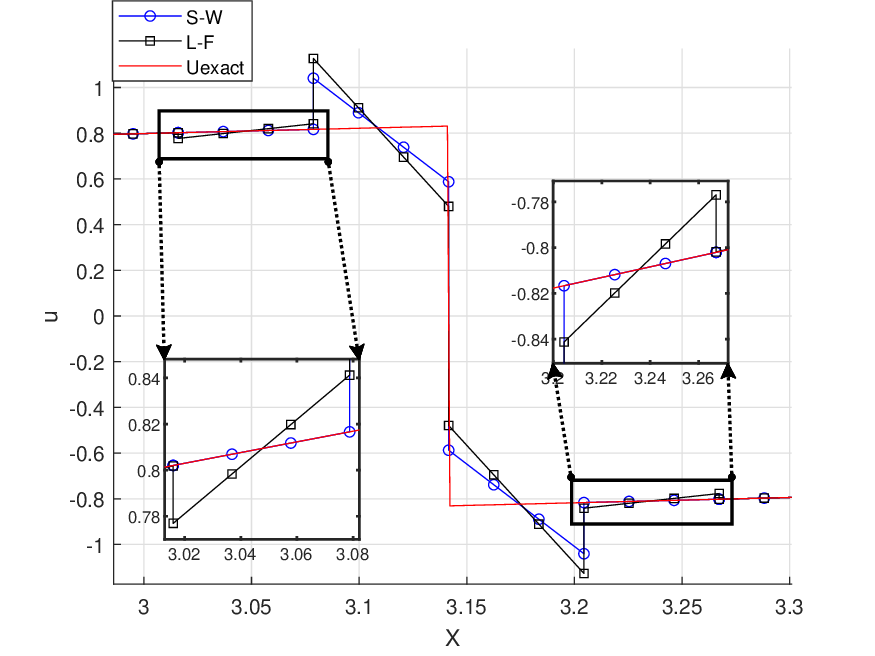}}
    \end{minipage}
    \vfill
    \vspace{0.2cm}
    \begin{minipage}{0.49\linewidth}
      \small
      \centerline{(a)\ $P^1$, Mesh=10}
    \end{minipage}
    \hfill
    \begin{minipage}{0.49\linewidth}
      \small
      \centerline{(b)\ $P^1$, Mesh=100}
    \end{minipage}
    \begin{minipage}{0.49\linewidth}
      \small
      \centerline{\includegraphics[width=1\linewidth]{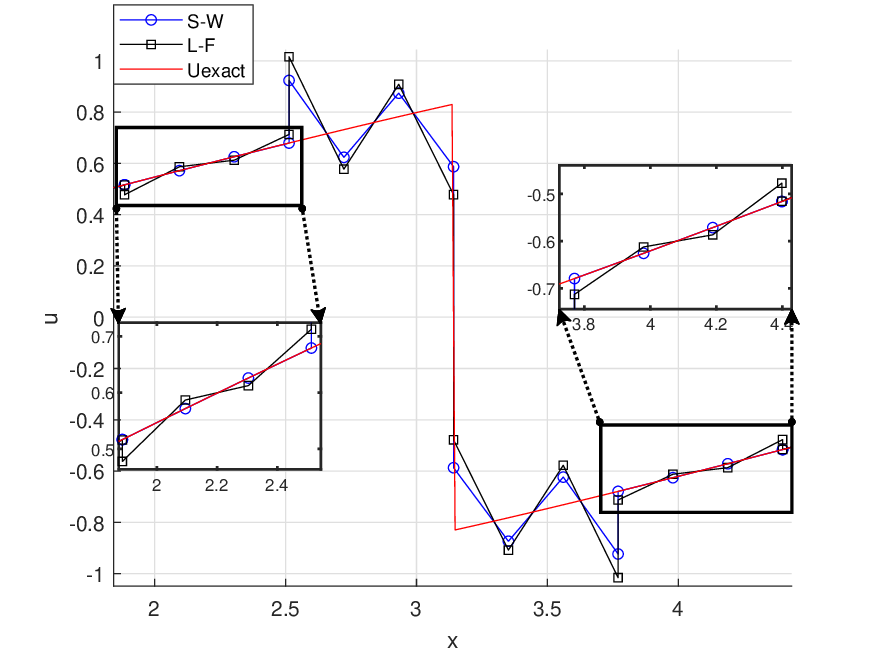}}
    \end{minipage}
    \hfill
    \begin{minipage}{0.49\linewidth}
      \small
      \centerline{\includegraphics[width=1\linewidth]{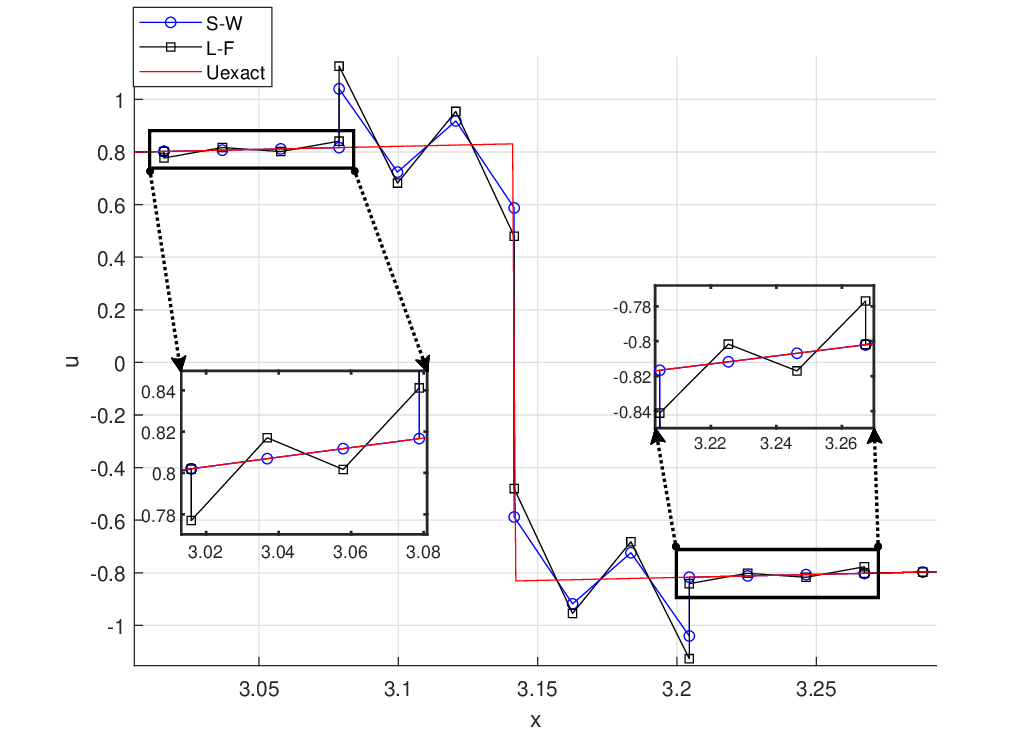}}
    \end{minipage}
    \begin{minipage}{0.49\linewidth}
      \small
      \centerline{(c)\ $P^3$, Mesh=10}
    \end{minipage}
    \hfill
    \begin{minipage}{0.49\linewidth}
      \small
      \centerline{(d)\ $P^3$, Mesh=100}
    \end{minipage}
  \end{center}
  \caption{Steger-Warming-scalar-DG and Lax-Friedrichs-scalar-DG both using equally spaced cells. $u_t+(\frac{1}{2}u^2)_x=0$ with 
  $u_0(x) = \sin(x)$ and $T_{end}=2.6$ (Discontinuity has been developed). {\color{red}{No limiters}}.  Localized magnification has been applied to all sub-figures.}
  \vspace{0.3cm}
  \label{Fig.SW-LF-Burgers-discontinuity-1D}
\end{figure}
\par
Without the introduction of any limiters, near discontinuities, the number of oscillations 
and the degree of overshoot of the L-F flux are both greater than those of the S-W flux. 
Additionally, by observing the jump in the DG approximation solution $u_h$ at the cell interfaces, 
it is found that $[\![u_h]\!]^{L-F}_{x_{i+1/2}} > [\![u_h]\!]^{S-W}_{x_{i+1/2}}$. 
In summary, the S-W flux scheme inherently has a better capability for handling discontinuities 
than the L-F flux, and the DG approximation based on the S-W flux is more ``smooth'' 
(with smaller jumps at interfaces).  
\end{example}

\section{10-Stage 4th-Order Strong Stability Preserving Runge-Kutta Method}\label{appendix-SSPRK(10,4)}
SSPRK(10,4) is a commonly used strong stability preserving Runge-Kutta method, whose specific form can be found in many literatures. 
However, few literatures directly provide the time parameter in the spatial discrete operator $L_h(\cdot)$ under the SSPRK(10,4) framework, 
which is necessary for time-dependent Dirichlet boundary condition and KXRCF discontinuity indicator. 
For the convenience of future readers, we hereby specifically provide the moment parameters of the spatial discrete operator $L_h(\cdot)$ 
at each sub-time step of SSPRK(10,4):
\begin{equation}\label{SSPRK(10,4)}
\left.
\begin{aligned}
 u^{(1)}&=u^n+\frac{1}{6} \Delta t L_h(u^n;t_n), \\
 u^{(2)}&=u^{(1)}+\frac{1}{6} \Delta t L_h(u^{(1)};t_n+\frac{1}{6}\Delta t), \\
 u^{(3)}&=u^{(2)}+\frac{1}{6} \Delta t L_h(u^{(2)};t_n+\frac{1}{3}\Delta t), \\
 u^{(4)}&=u^{(3)}+\frac{1}{6} \Delta t L_h(u^{(3)};t_n+\frac{1}{2}\Delta t), \\
 u^{(5)}&=\frac{3}{5} u^n+\frac{2}{5} u^{(4)}+\frac{1}{15} \Delta t L_h(u^{(4)};t_n+\frac{2}{3}\Delta t), \\
 u^{(6)}&=u^{(5)}+\frac{1}{6} \Delta t L_h(u^{(5)};t_n+\frac{1}{3}\Delta t), \\
 u^{(7)}&=u^{(6)}+\frac{1}{6} \Delta t L_h(u^{(6)};t_n+\frac{1}{2}\Delta t), \\
 u^{(8)}&=u^{(7)}+\frac{1}{6} \Delta t L_h(u^{(7)};t_n+\frac{2}{3}\Delta t), \\
 u^{(9)}&=u^{(8)}+\frac{1}{6} \Delta t L_h(u^{(8)};t_n+\frac{5}{6}\Delta t), \\
 u^{n+1}&=\frac{1}{25} u^n+\frac{9}{25} u^{(4)}+\frac{3}{5} u^{(9)}+\frac{3}{50} \Delta t L_h(u^{(4)};t_n+\frac{2}{3}\Delta t)+\frac{1}{10} \Delta t L_h (u^{(9)};t_n+\Delta t).
\end{aligned}
\right\}
\end{equation}
More SSPRK schemes, please refer to \cite{ref49}.

\end{document}